\documentclass[12pt,a4paper]{report}
\let\openright=\clearpage

\usepackage[utf8]{inputenc}

\usepackage{bibunits}

\usepackage{colortbl} 

\usepackage[utf8]{inputenc}
\usepackage{amsmath,amsfonts,amssymb,amsthm}
	\numberwithin{equation}{section} 
\usepackage{enumerate}
\usepackage{graphicx}

\newtheorem{theorem}{\protect\theoremname}[section]
\newtheorem{definition}[theorem]{\protect\definitionname}
\newtheorem{lemma}[theorem]{\protect\lemmaname}
\newtheorem{proposition}[theorem]{\protect\propositionname}
\newtheorem*{proposition*}{\protect\propositionname}
\newtheorem{example}[theorem]{\protect\examplename}
\newtheorem{corollary}[theorem]{\protect\corollaryname}
\newtheorem{remark}[theorem]{\protect\remarkname}
\newtheorem{problem}[theorem]{\protect\problemname}
\newtheorem{notation}[theorem]{\protect\notationname}
\newtheorem{claim}[theorem]{\protect\claimname}
\newtheorem{conjecture}[theorem]{\protect\conjecturename}

\providecommand{\corollaryname}{Corollary}
\providecommand{\claimname}{Claim}
\providecommand{\definitionname}{Definition}
\providecommand{\lemmaname}{Lemma}
\providecommand{\notationname}{Notation}
\providecommand{\remarkname}{Remark}
\providecommand{\problemname}{Problem}
\providecommand{\propositionname}{Proposition}
\providecommand{\examplename}{Example}
\providecommand{\theoremname}{Theorem}
\providecommand{\conjecturename}{Conjecture}

\newcommand{\N}{\mathbb{N}}

\newcommand{\Q}{\mathbb{Q}}

\newcommand{\slantfrac}[2]{\,^#1\!/_#2}
\newcommand{\mc}{\mathcal}

\newcommand{\fsd}{F_{\sigma\delta}}
\newcommand{\bd}{^{\star\star}}
\newcommand{\seq}{\omega^{<\omega}}
\newcommand{\Fa}{\mc F_\alpha}
\newcommand{\ext}{\hat{\ }}

\newcommand{\baire}{\omega^\omega}

\newcommand{\Amg}[2]{\textnormal{Amg}\left(#1,#2\right)}
\newcommand{\Compl}[2]{\textnormal{Compl}\left(#1,#2\right)}
\newcommand{\rank}{r_{\textnormal{iie}}}
\newcommand{\Tr}{\textnormal{Tr}}
\newcommand{\D}{D_{\textnormal{iie}}}
\newcommand{\cltr}[1]{\textnormal{cl}_{\Tr}\left(#1\right)}
\newcommand{\ims}[2]{\textnormal{ims}_{#1}\left(#2\right)}

\newcounter{vkNoteCounter}

\newcounter{okNoteCounter}

\newcommand{\errata}[1]{{\color{purple} {#1}[Err]}}
\newcommand{\err}[1]{{\color{purple} {#1}}}
\renewcommand{\errata}[1]{{#1}}
\renewcommand{\err}[1]{{#1}}

\makeatletter
\def\@makechapterhead#1{
  {\parindent \z@ \raggedright \normalfont
   \Huge\bfseries \thechapter. #1
   \par\nobreak
   \vskip 20\p@
}}
\def\@makeschapterhead#1{
  {\parindent \z@ \raggedright \normalfont
   \Huge\bfseries #1
   \par\nobreak
   \vskip 20\p@
}}
\makeatother

\def\chapwithtoc#1{
\chapter*{#1}
\addcontentsline{toc}{chapter}{#1}
}

\usepackage{etoolbox}
\patchcmd{\thebibliography}{\chapter*}{\section*}{}{}

\renewcommand{\thechapter}{\Roman{chapter}}

\begin{document}

\lefthyphenmin=2
\righthyphenmin=2


\pagestyle{empty}
\begin{center}

\large

Charles University

\medskip

Faculty of Mathematics and Physics

\vfill

{\bf\Large DOCTORAL THESIS}

\vfill

\centerline{\mbox{\includegraphics[width=60mm]{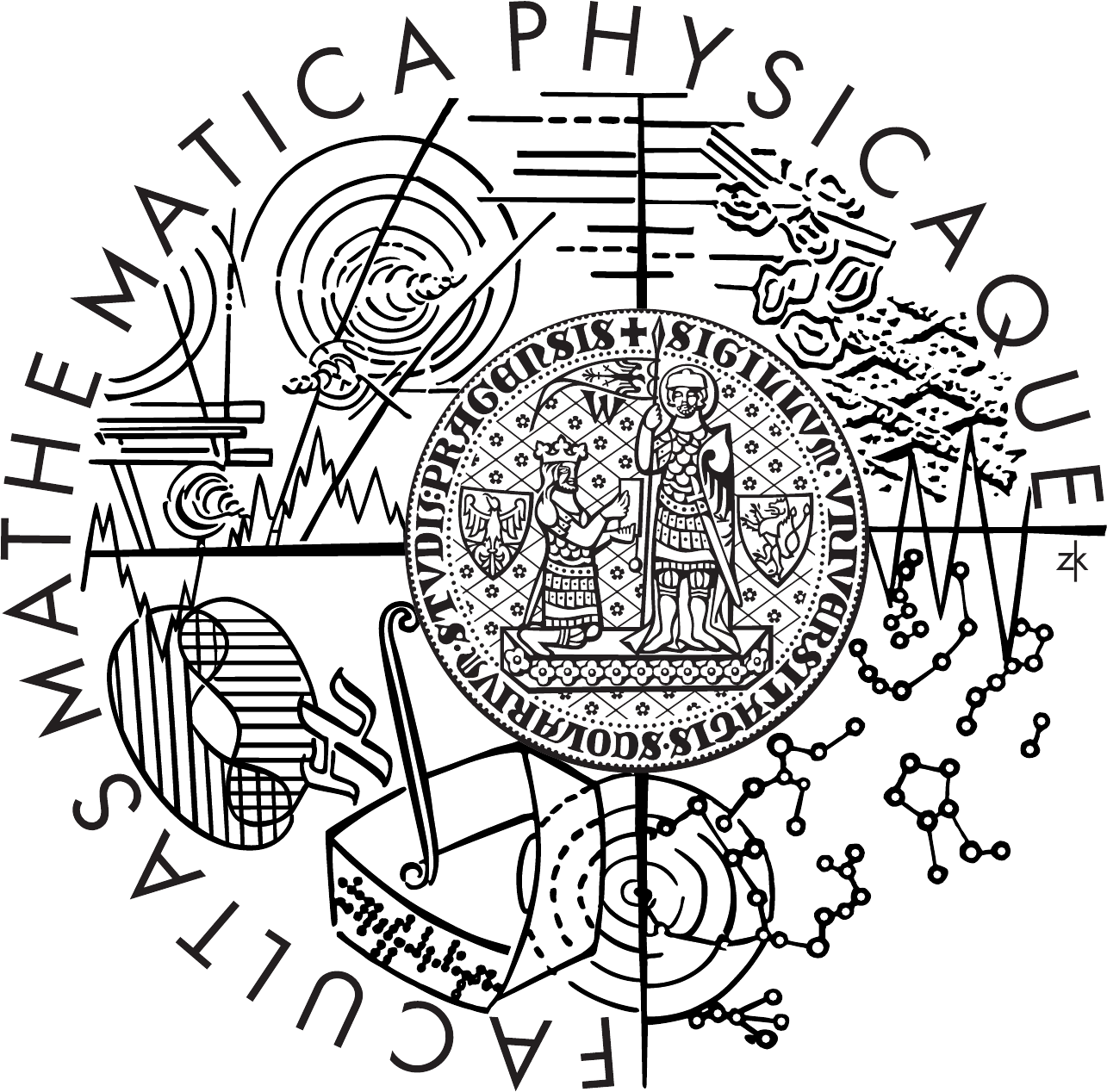}}}

\vfill
\vspace{5mm}

{\LARGE Vojtěch Kovařík}

\vspace{15mm}

{\LARGE\bfseries Absolute and non-absolute $\mc F$-Borel spaces}

\vfill

Department of Mathematical Analysis

\vfill

{\small
\begin{tabular}{rl}

Supervisor of the doctoral thesis: & Prof. RNDr. Ondřej Kalenda, Ph.D., DSc.\\
\noalign{\vspace{2mm}}
Study programme: & Mathematics \\
\noalign{\vspace{2mm}}
Specialization: & Mathematical Analysis \\
\end{tabular}}

\vfill

Prague 2018

\end{center}

\newpage
\pagestyle{plain}
\pagenumbering{roman}
\setcounter{page}{1}


\openright

\noindent
I would like to thank my supervisor, Ondřej Kalenda, for his patience and for the time he invested into showing me how to do research better.
\newpage


\vglue 0pt plus 1fill

\noindent
I declare that I carried out this doctoral thesis independently, and only with the cited
sources, literature and other professional sources.

\medskip\noindent
I understand that my work relates to the rights and obligations under the Act No.
121/2000 Coll., the Copyright Act, as amended, in particular the fact that the Charles
University in Prague has the right to conclude a~license agreement on the use of this
work as a~school work pursuant to Section 60 paragraph 1 of the Copyright Act.

\vspace{10mm}

\hbox{\hbox to 0.5\hsize{%
In ........ date ............
\hss}\hbox to 0.5\hsize{%
Vojtěch Kovařík
\hss}}

\vspace{20mm}
\newpage


\vbox to 0.5\vsize{
\setlength\parindent{0mm}
\setlength\parskip{5mm}

Název práce:
Absolutně a neabsolutně $\mc F$-borelovské prostory

Autor:
Vojtěch Kovařík

Katedra:  
Katedra matematické analýzy

Vedoucí disertační práce:
Prof. Ondřej Kalenda, Ph.D., DSc., Katedra matematické analýzy

Abstrakt:
Zabýváme se $\mc F$-borelovskou složitostí topologických prostorů a~tím, jak \err{se} tato složitost liší v závislosti na prostoru, kam je daný topologický prostor vnořen. Obzvláště nás pak zajímá, kdy je tato složitost absolutní, tj. stejná ve všech kompaktifikacích.
Ukazujeme, že složitost metrizovatelných prostorů je absolutní. Dále odvozujeme postačující podmínku pro to, aby byl prostor absolutně $\fsd$. Studujeme vztah lokální a~globální složitosti a~závádíme různé reprezentace \mbox{$\mc F$-borelovských} množin.
Tyto nástroje používáme k důkazu několika výsledků, zejména pak k získání hierarchie prostorů, jejichž složitost je neabsolutní.

Klíčová slova: deskriptivní složitost, kompaktifikace, F-borelovská množina, absolutní složitost

\vss}\nobreak\vbox to 0.49\vsize{
\setlength\parindent{0mm}
\setlength\parskip{5mm}

Title: Absolute and non-absolute $\mc F$-Borel spaces

Author:
Vojtěch Kovařík

Department:
Department of Mathematical Analysis

Supervisor:
Prof. Ondřej Kalenda, Ph.D., DSc., Department of Mathematical Analysis

Abstract:
We investigate $\mc F$-Borel topological spaces. We focus on finding out how a~complexity of a~space depends on where the space is embedded. Of a~particular interest is the problem of determining whether a~complexity of given space $X$ is absolute (that is, the same in every compactification of $X$).
We show that the complexity of metrizable spaces is absolute and provide a~sufficient condition for a~topological space to be absolutely $\fsd$. We then investigate the~relation between local and global complexity. To improve our understanding of $\mc F$-Borel spaces, we introduce different ways of representing an $\mc F$-Borel set. We use these tools to construct a~hierarchy of $\mc F$-Borel spaces with non-absolute complexity, and to prove several other results.


\vss}

\newpage


\openright
\tableofcontents


\pagestyle{plain}
\pagenumbering{arabic}
\setcounter{page}{1}

\begin{bibunit}[alphanum]
\renewcommand{\thetheorem}{\arabic{theorem}}

\chapter{Introduction}

The present thesis consists of three chapters, each of them corresponding to a research
paper. All these articles contain original results which are ultimately meant to serve the goal of better understanding the descriptive properties of Banach spaces, and descriptive properties of topological spaces in general.

Of particular interest to this topic is the class of weakly compactly generated Banach spaces (WCG); this is because every WCG space is an $\fsd$-subset of its second dual endowed with the $\text{weak}^\star$ topology (\cite[Theorem 3.2]{talagrand1979espaces}).
A natural question is then to ask whether a WCG space $X$, equipped with the weak topology, is ``\emph{absolutely}`` $\fsd$ -- that is, $\fsd$ in every space which contains it, or equivalently, whether $(X,w)$ is $\fsd$ in every compactification.
It is not obvious what the answer to this question should be, because there exists a topological space which is $\fsd$ in some compactification, but \emph{not} absolutely $\fsd$ (\cite{talagrand1985choquet} or \cite[Remark 6.1]{argyros2008talagrand}).

Strictly larger than WCG spaces is the class of those Banach spaces which are $\mc K$-analytic in their weak topology.
A long-standing open problem posed by Talagrand (\cite[Problem 4.6.a]{talagrand1979espaces}) was to decide whether every such space is $\fsd$ in its second dual endowed with the weak$^\star$ topology.
This was answered in the negative by \cite{argyros2008talagrand}, but the authors also show that the answer would be positive under some additional assumptions. This topic was further studied in \cite{kampoukos2013certain,argyros2009reznichenko}.
A related question which we would like to address is whether a general weakly-$\mc K$-analytic Banach space might at least be $\fsd$ in \emph{some} compactification of $(X,w)$.

A more general version of the aforementioned questions is the following:
\begin{problem}\label{problem:Fsd_Banach_spaces}
\begin{enumerate}[(i)]
\item Which Banach spaces, endowed with the weak topology, are $\fsd$ in some compactification?
\item Which of these spaces are absolutely $\fsd$?
\end{enumerate}
\end{problem}

We now highlight the~most important contributions made by the present thesis. In Chapters II, III, and IV respectively, we present the following papers:
\begin{itemize}
\item Kalenda, Ondřej, and Kovařík, Vojtěch. "Absolute $\fsd$ spaces." \emph{Topology and its Applications} 233 (2018): 44-51.
\item Kovařík, Vojtěch. "Absolute $\mathcal F $-Borel classes." \emph{Fundamenta Mathematicae}, \emph{published electronically}.
\item Kovařík, Vojtěch. "Complexities and Representations of $\mc F$-Borel Spaces" \emph{Submitted, available at arXiv:1804.08367}.
\end{itemize}

In Chapter \ref{chapter:abs_fsd_spaces} we provide a sufficient condition for a topological space to be absolutely $\fsd$, which partially answers a problem of Frolík (Problem 1 in \cite{frolik1963descriptive}), who asked for a characterization of such spaces.
A corollary of this result is that hereditarily Lindelöf spaces which are $\fsd$ in some compactification are absolutely $\fsd$.
In particular, this shows that separable Banach spaces are absolutely $\fsd$ when endowed with the weak topology.

However, this turns out to be the only contribution the author was able to make with regards to Problem 1.
Rather than due to the lack of trying, this was caused by the absence of tools with which one could approach the problem of finding the complexity of a space in any compactification which is not both explicitly given and well behaved.
For example, we might be able to easily compute the complexity of $(X,w)$ in $(X,w\bd)$, but find it very difficult to compute its complexity in the Čech-Stone compactification $\beta (X,w)$. Even worse, we might be completely at loss when it comes to computing the complexity in an abstract compactification of $(X,w)$.
Another complication is that there are currently only two ``genuine'' examples of spaces with non-absolute complexity (\cite{talagrand1985choquet}, \cite[Remark 6.1]{argyros2008talagrand}), and the construction of these particular examples is rather complicated.

We instead study the problem in the general setting where $X$ is a Tychonoff topological space, contained in another space $Y$.
The class of $\mc F$-Borel sets is the smallest system containing closed sets and closed under countable unions and intersections. We enumerate its ``levels'' as $\mc F_0 :=$ closed sets, $\mc F_1 :=F_\sigma$-sets, $\mc F_2 := \fsd$-sets, and so forth, \errata{obtaining classes $\Fa$, $\alpha < \omega_1$, which are multiplicative for even $\alpha$ and additive classes for odd $\alpha$.}
We then ask what is the exact complexity of $X$ in $Y$, that is, what is the smallest $\alpha$ for which $X$ is an $\Fa$-subset of $Y$? And how does this complexity of $X$ vary between different spaces $Y$ containing $X$?
This is a problem which has been earlier studied either for different descriptive classes (\cite{holicky2003perfect,holicky2004internal,raja2002some,holicky2004fsigma}), or with attention restricted to those super-spaces $Y$ which are metrizable (\cite{marciszewski1997absolute,junnila1998characterizations,holicky2010descriptive}), but not in the full generality.

It is not hard to show that complexity of separable metrizable spaces is absolute (see
Theorem~\ref{theorem: separable metrizable} in Chapter \ref{chapter:brooms}).
This can be extended to hereditarily Lindelöf spaces, as observed by Jiří Spurný. An
elegant formal proof is given in Chapter \ref{chapter:representations} (Proposition \ref{proposition: hereditarily lindelof spaces are absolute}) using the~theory of simple represantations of $\mc F$-Borel sets developed therein (see Section \ref{section: simple representations}).

In Chapters \ref{chapter:brooms} and \ref{chapter:representations}, we further develop a stronger concept of ``regular'' representations of $\mc F$-Borel sets (Section \ref{sec:broom spaces and parametrization}, resp. \ref{section: regular representations}). These representations can be used to bound the absolute complexity of a space from above, for which (as far as the author is aware of) there was previously no known method.
These representations can be used to a more refined study of the Talagrand's example from  \cite{talagrand1985choquet}.
We show the existence of spaces $X_\alpha^\beta$ which are absolutely $\mc F_\beta$, but there exists some compactification in which the complexity of $X_\alpha^\beta$ is precisely $\Fa$ (Theorems \ref{chapter:brooms}.\ref{theorem: main theorem} and \ref{chapter:representations}.\ref{theorem:X_2_beta}).
In Chapter \ref{chapter:representations}, we further strengthen this result by constructing spaces $X_{[\alpha,\beta]}$ which are absolutely $\mc F_\beta$, but for each $\gamma \in [\alpha,\beta]$ there is a compactification in which the complexity of $X_{[\alpha,\beta]}$ is precisely $\mc F_\gamma$.
Finally, we study the relation between complexity and local complexity and conclude that the ``global'' complexity is determined by its local behavior.

Apart from the one result mentioned earlier, the questions surrounding the complexity of Banach spaces remain unresolved. However, we believe that the tools introduced in this work should make further investigations of this topic significantly easier.
There also remain some open problems regarding the general topological setting. Some of the more technical issues are listed in Chapter \ref{chapter:representations}, but we should also mention one question with slightly ``philosophical'' flavor:
There are only two proper examples of spaces with non-absolute complexity. Is it because spaces with non-absolute complexity are rare and anomalous? Or could it be that having non-absolute complexity is in fact common (but has not yet been proven for the ``natural'' examples, due to the difficulty of the task)?

\bigskip
In Chapter \ref{chapter:representations}, we build \err{on the ideas} from second half of Chapter~\ref{chapter:brooms}, \err{repeating and improving} most of the~proofs from Sections \ref{chapter:brooms}.4.3 and \ref{chapter:brooms}.5 \errata{(and also obtaining many completely new results).}
The reader should however be wary of the fact that a part of the notation is different between the two chapters (because the original notation turned out to be unsuitable for formulating the generalized results which came later). This in particular concerns the symbols $\Fa$ denoting the $\mc F$-Borel classes and the notation related to the ``infinite broom sets''.

\renewcommand{\thetheorem}{\arabic{section}.\arabic{theorem}}
\putbib[refs]
\end{bibunit}
\begin{bibunit}[alphanum]

\chapter{Absolute \texorpdfstring{$F_{\sigma\delta}$}{F sigma delta} Spaces}
 \label{chapter:abs_fsd_spaces}

\begin{center}
(with Ondřej Kalenda) \\
\vskip 5 mm
\textit{(Published in Topology and its Applications 233 (2018), 44–51.)}
\end{center}

\textbf{Abstract:}
We prove that hereditarily Lindelöf space which is $\fsd$ in some compactification is absolutely $\fsd$. In particular, this implies that any separable Banach space is absolutely $\fsd$ when equipped with the weak topology.

\section{Introduction}
Throughout the paper, all spaces will be Tychonoff. Central to the topic of our paper is the following definition:
\begin{definition}\label{definition: F sigma delta spaces}
Let $X$ be a Tychonoff topological space. We say that $X$ is an \emph{$\fsd$ space} if there exists a compactification $cX$ of $X$, such that $X\in\fsd(cX)$.

We say that $X$ is an \emph{absolute $\fsd$ space} (or that $X$ is \emph{absolutely $\fsd$}) if $X\in\fsd(cX)$ holds for every compactification $cX$ of $X$.
\end{definition}
Note that $X$ is absolutely $\fsd$ if and only if $X\in\fsd(Y)$ holds for every Tychonoff topological space $Y$ in which $X$ is embedded.

If, in the above definition, we replace the class $\fsd$ by $G_\delta$, we get the definition of the well known concept of Čech-completeness -- however, in such a case the situation is less complicated, because every Čech-complete space is automatically absolutely $G_\delta$. Internal characterization \err{of} Čech-complete spaces was given by Zdeněk Frolík, who also gave a characterization of $\fsd$ spaces in terms of complete sequences of covers (see Definition \ref{def: complete seq. of covers} below). He then asked for a description of those spaces which are \emph{absolutely} $\fsd$ (Problem 1 in \cite{frolik1963descriptive}), and this problem is still open.
 
However, Frolík did not know whether there actually exist non-absolute $\fsd$ spaces. This part of the problem was solved later by Talagrand, who found an example of such a space (\cite{talagrand1985choquet}). Thus, we formulate Frolík's problem as follows:
\begin{problem} \label{problem: 1}
Among all $\fsd$ spaces, describe those which are absolutely $\fsd$.
\end{problem}
If we are unable to completely determine the answer to Problem \ref{problem: 1}, the next best thing to do is to find a partial answer to Problem \ref{problem: 2} for as many spaces as possible.
\begin{problem} \label{problem: 2}
Let $X$ be a (possibly non-absolute) $\fsd$ space. Describe those compactifications of $X$ in which it is $\fsd$.
\end{problem}

In Section \ref{section:fsd topological spaces}, we give a partial answer to Problem \ref{problem: 2} by showing that if a $X$ is $\fsd$ in some compactification $cX$, it is automatically $\fsd$ in all larger compactifications (which is easy) and also in all compactification which are not much smaller than $cX$ (see Corollary \ref{corollary: countable quotients} for the details).

In Proposition \ref{proposition: fs disjoint covers and abs fsd} we give a partial answer to Problem \ref{problem: 1} by finding a sufficient condition for a space to be absolutely $\fsd$. This condition is similar in flavor to Frolík's characterization of $\fsd$ spaces. Applying this result, we get that hereditarily Lindelöf $\fsd$ spaces are absolutely $\fsd$ (Theorem \ref{theorem: countable network}) and that separable Banach spaces are absolutely $\fsd$ in the weak topology (Corollary \ref{corollary: banach spaces}).

In the rest of the introductory section we collect some known results and background information.

We could adapt Definition \ref{definition: F sigma delta spaces} for the lower classes of Borel hierarchy, where we have the following standard results. Their proof consists mostly of using the fact that continuous image of a compact space is compact.
\begin{remark}
Let $X$ be a topological space.
\begin{enumerate}
	\item $X$ is absolutely closed $\iff$ $X$ is compact.
	\item $X$ is absolutely $F_\sigma$ $\iff$ $X$ is $\sigma$-compact.
	\item $X$ is absolutely open $\iff$ $X$ is locally compact.
	\item $X$ is absolutely $G_\delta$ $\iff$ $X$ is Čech-complete.
\end{enumerate}
In the first two cases, $X$ being closed ($F_\sigma$) in some compactification automatically implies that $X$ is closed ($F_\sigma$) in every Tychonoff space where it is embedded. For open and $G_\delta$ spaces, we only get that $X$ is open ($G_\delta$) in those Tychonoff spaces where it is densely embedded.
\end{remark}

As shown in \cite{talagrand1985choquet}, not every $\fsd$ space is absolutely $\fsd$. This means that the class of $\fsd$ sets is the first one for which it makes sense to study Problem \ref{problem: 2}, which is one of the reasons for our interest in this particular class. However, Talagrand's is the only result of this kind (as far as the authors are aware of), and not much else is known about `topologically' absolute $\fsd$ spaces. In \cite{kovarik2018brooms}, topological absoluteness is studied for general $\mc F$-Borel classes, providing more examples based on Talagrand's construction and also some theoretical results.

Several authors have investigated slightly different notions of absoluteness for $\fsd$ spaces. Recall that in separable metrizable setting, $\fsd$ sets coincide with $\textrm{alg}\left(F\right)_{\sigma\delta}$ sets (where $\textrm{alg}\left(F\right)$ is the algebra generated by closed sets). As shown in \cite{holicky2003perfect}, the class of $\textrm{alg}\left(F\right)_{\sigma\delta}$ sets is absolute (in the sense that if a set is in $\textrm{alg}\left(F\right)_{\sigma\delta}(cX)$ for some compactification $cX$, it is automatically of this same class in every Tychonoff space where it is embedded).

In \cite{marciszewski1997absolute} and \cite{junnila1998characterizations}, the authors study metric spaces which are absolutely $\fsd$ `in a metric sense' - that is, $X\in\fsd (Y)$ for every \emph{metrizable} space $Y$ in which $X$ is embedded. In \cite{junnila1998characterizations}, the authors give a characterization of `metric absoluteness' for $\fsd$ spaces in terms of complete sequences of covers - namely that $X$ is absolutely $\fsd$ in the metric sense if and only if it has a complete sequence of $\sigma$-discrete closed covers.

Unfortunately, this result is not useful in our setting, because Talagrand's space is an example of non-metrizable space, which does have such a complete sequence, but it is not absolutely $\fsd$ (in our - topological - sense).

In \cite{kovarik2018brooms}, it is shown that if a metric space is separable, its complexity is automatically absolute even in the topological sense. For $\fsd$ spaces, this is a special case of Theorem \ref{theorem: countable network}.

\section{Compactifications}\label{section:compactifications}
We recall the standard definitions of compactifications and their partial ordering. By \emph{compactification} of a topological space $X$ we understand a pair $(cX,\varphi)$, where $cX$ is a compact space and $\varphi$ is a homeomorphic embedding of $X$ onto a dense subspace  of $cX$. Symbols $cX$, $dX$ etc. will always denote compactifications of $X$. 

Compactification $(cX,\varphi)$ is said to be \emph{larger} than $(dX,\psi)$, if there exists a continuous mapping $f : cX\rightarrow dX$, such that $\psi = f \circ \varphi$. We denote this as $ cX \succeq dX $. Recall that for a given $T_{3\slantfrac{1}{2}}$ topological space $X$, its compactifications are partially ordered by $\succeq$ and Stone-Čech compactification $\beta X$ is the largest one.

Often, we encounter a situation where $X\subset cX$ and the corresponding embedding is identity. In this case, we will simply write $cX$ instead of $(cX,\textrm{id}|_X)$.
Much more about this topic can be found in many books - for a more recent one, see for example \cite{freiwald2014introduction}.

In the introduction, we defined the notion of being an $\fsd$ space and an absolute $\fsd$ space. Having defined the partial order $\succeq$ on the class of compactifications of $X$, we note the basic facts related to Problem \ref{problem: 2}. The proof of this remark consists of using the fact that continuous preimage of an $\fsd$ set is an $\fsd$ set.
\begin{remark}\label{rem:absolute complexity}
For a topological space $X$, we have the following:
\begin{itemize}
	\item $X$ is an $\fsd$ space $\iff$ $X\in\fsd(\beta X)$;
	\item $X$ is an absolute $\fsd$ space $\iff$ $X\in\fsd(cX)$ for every $cX$;
	\item $X\in\fsd (dX)$, $cX\succeq dX$ $\implies$ $X\in\fsd(cX)$.
\end{itemize}
\end{remark}

\begin{notation}
Let $X$ be a topological space, suppose that two of its compactifications satisfy $dX\preceq cX$ and that $\varphi : cX\rightarrow dX$ is the mapping which witnesses this fact. We denote 
\[ \mathcal{F}\left( cX, dX\right):=\left\lbrace \varphi^{-1} \left(x\right)\big | \ x\in dX, \ \varphi^{-1} \left(x\right) \textrm{ is not a singleton} \right\rbrace . \]
\end{notation}
In this sense, every compactification $dX$ smaller than $cX$ corresponds to some disjoint system $\mathcal{F}$ of compact subsets of $cX\setminus X$. Conversely, some disjoint systems of compact subsets of $cX\setminus X$ correspond to quotient mappings, which correspond to compactifications smaller than $cX$. Not every such system $\mc F$ corresponds to a compactification, but surely every finite (disjoint, consisting of compact subsets of $cX\setminus X$) $\mc F$ does.

\section{\texorpdfstring{$\fsd$}{F sigma delta} Spaces}
In this section, we will list some results related to $\fsd$ spaces.
\subsection{Banach Spaces}
Unless otherwise specified, a Banach space $X$ (resp. its second dual), will always be equipped with weak (resp. $w^\star$) topology. In \cite{argyros2008talagrand}, a Banach space $X$ is said to be $F_{\sigma\delta}$ if it is is an $F_{\sigma\delta}$ subset of $X^{\star\star}$. Note that the space $X^{\star\star}$ is always $\sigma$-compact, so it is $F_\sigma$ in $\beta X^{\star\star}$. Consequently, any $\fsd$ Banach space is automatically an $\fsd$ space (in the sense of Definition \ref{definition: F sigma delta spaces}).

An important class of Banach spaces are the spaces which are weakly compactly generated (WCG). Recall that a Banach space $X$ is said to be WCG, if there exists a set $K\subset X$ which is weakly compact, such that $\textrm{span}(K)$ is dense in $(X,||\cdot||)$. Clearly all separable spaces and all reflexive spaces are WCG. The canonical example of non-separable non-reflexive WCG space is the space $c_0(\Gamma)$ for uncountable index set $\Gamma$. For more information about WCG spaces, see for example \cite{fabian2013functional}. The reason for our interest in WCG spaces is the following result (\cite[Theorem 3.2]{talagrand1979espaces}):
\begin{proposition}\label{proposition: WCG spaces are fsd}
Any WCG space is an $\fsd$ Banach space.
\end{proposition}
In fact, even every \errata{Banach} subspace of a WCG space is $\fsd$ in its second dual. Talagrand has found an example of an $\fsd$ Banach space which is not a subspace of a WCG space \cite{talagrand1979espaces}. This space belongs to a more general class of weakly $\mc K$-analytic spaces. A problem which had been open for a long time is whether every weakly $\mc K$-analytic space is an $\fsd$ Banach space. A counterexample has been found in \cite{argyros2008talagrand} (as well as some sufficient conditions for a weakly $\mc K$-analytic space to be an $\fsd$ Banach space). The problem which still remains unsolved is whether weakly $\mc K$-analytic spaces are topologically $\fsd$.

\subsection{Topological Spaces}\label{section:fsd topological spaces}

\begin{proposition}\label{proposition: G delta characterization}
Suppose that $X\in F_{\sigma\delta}\left(cX\right)$ and $ dX \preceq cX $. Then $X\in F_{\sigma\delta}\left(dX\right)$ holds if and only if there exists a sequence $(H_n)_n$ of $F_\sigma$ subsets of $cX$, such that
\[ \left(\forall F\in\mathcal{F}\left(cX,dX\right)\right)\left(\exists n\in\N\right): X\subset H_n \subset cX \setminus F. \]
\end{proposition}
\begin{proof}
Denote by $\varphi$ the map witnessing that $dX\preceq cX$. \\*
$\implies $: Assume that $X=\bigcap F_n$, where the sets $F_n$ are $F_\sigma$ in $dX $. Denote $H_n:=\varphi^{-1}(F_n)$. Clearly, $H_n\subset cX$ are $F_\sigma$ sets containing $X$. Let $F$ be $\mathcal{F}\left(cX,dX\right)$, that is, $F=\varphi^{-1}(y)$ for some $y\in dX\setminus X$. By the assumption, we have $X\subset F_n \subset dX\setminus \left\{y\right\}$ for some $n\in \mathbb N$. By definition of $\varphi$, we get the desired result:
\[ X=\varphi^{-1}\left(X\right)\subset \varphi^{-1}\left(F_n\right)\subset \varphi^{-1}\left(dX\setminus \left\{y\right\}\right)=\varphi^{-1}\left(X\right)\setminus \varphi^{-1}\left(y\right)=X\setminus F. \]
$\Longleftarrow$: Let the sequence of sets $H_n\subset cX$ be as in the proposition. We know that $X=\bigcap F_n$ for some $F_\sigma$ sets $F_n$. We now receive sets $F'_n:=\varphi(F_n)$ and $H'_n:=\varphi(H_n)$, $n\in\N$, all of which are $F_\sigma$ in $dX$. Clearly, we have
\[ X\subset\bigcap F'_n \cap \bigcap {H'_n}. \]
For the converse inclusion, let $y\in dX \setminus X$. If $\varphi^{-1}(y)$ is a singleton, we have $\varphi^{-1}(y) \subset cX\setminus F_n$ for some $n\in\N$, and therefore $y\notin \varphi(F_n)=F'_n$. If $\varphi^{-1}(y)$ is not a singleton, then $\varphi^{-1}(y)\in \mathcal{F}\left(cX,dX\right)$, so there exists some $n\in\N$, such that $H_n\subset X\setminus \varphi^{-1}(y)$. In this case, we have $y\notin \varphi(H_n)=H'_n$.
\end{proof}

Since any $\fsd$ space is Lindelöf, we can make use of the following lemma, which follows immediately from \cite[Lemma 14]{spurny2006solution}.
\begin{lemma}\label{lemma:fsd spaces are Fs separated}
Let $X$ be a Lindelöf subspace of a compact space $L$. Then for every compact set $K\subset L \setminus X$, there exists $H\in F_\sigma\left(L\right)$, such that $X\subset H\subset L\setminus K$.
\end{lemma}

Once we have Lemma \ref{lemma:fsd spaces are Fs separated}, Proposition \ref{proposition: G delta characterization} yields the following corollary, which gives a partial answer to Problem \ref{problem: 2}:
\begin{corollary}\label{corollary: countable quotients}
Suppose that $X$ is an $\fsd$ subspace of $cX$ and $dX\preceq cX$. Then $X$ is $\fsd$ in $dX$ as well, provided that the family $\mathcal{F}\left(cX,dX\right)$ is at most countable.
\end{corollary}
In particular, this implies that there exists no such thing as a "minimal compactification in which $X$ is $\fsd$" (unless, of course, $X$ is locally compact).

\section{Hereditarily Lindelöf Spaces}\label{section:countable network}
In this section, we present the following main result:
\begin{theorem}\label{theorem: countable network}
Every hereditarily Lindelöf $\fsd$ space is absolutely $\fsd$.
\end{theorem}
Note that every $\fsd$ space is Lindelöf (because it is $\mc K$-analytic), but the converse implication to Theorem \ref{theorem: countable network} does not hold - that is, not every absolutely $\fsd$ space is hereditarily Lindelöf. Indeed, any compact space which is not hereditarily normal is a counterexample. The proof of Theorem \ref{theorem: countable network} will require some technical preparation, but we can state an immediate corollary for Banach spaces:
\begin{corollary}\label{corollary: banach spaces}
Every separable Banach space is absolutely $\fsd$ (in the weak topology).
\end{corollary}
\begin{proof} 
By Proposition \ref{proposition: WCG spaces are fsd}, every separable Banach space $X$ is $\fsd$. The countable basis of the norm topology of $X$ is a network for the weak topology. The proposition follows from the fact that spaces with countable network are hereditarily Lindelöf.
\end{proof}

We will need the notion of complete sequence of covers:
\begin{definition}[Complete sequence of covers]\label{def: complete seq. of covers}
Let $X$ be a topological space. \emph{Filter} on $X$ is a family of subsets of $X$, which is closed with respect to supersets and finite intersections and does not contain the empty set. A point $x\in X$ is said to be an \emph{accumulation point} of a filter $\mc F$ on $X$, if each neighborhood of $x$ intersects each element of $\mc F$.

A sequence $\left( \mc F_n \right)_{n\in\mathbb N}$ of covers of $X$ is said to be \emph{complete}, if every filter which intersects each $\mc F_n$ has an accumulation point in $X$.
\end{definition}

The connection between this notion and our topic is explained by Proposition \ref{proposition: fsd iff complete sequence of covers}. Note that a cover of $X$ is said to be \emph{closed} (\emph{open}, $F_\sigma$, \emph{disjoint}) if it consists of sets which are closed  in $X$ (open, $F_\sigma$, disjoint). As a slight deviation from this terminology, a cover of $X$ is said to be \emph{countable} if it contains countably many elements.
\begin{proposition}\label{proposition: fsd iff complete sequence of covers}Any topological space $X$ satisfies
\begin{enumerate}
\item $X$ is Čech-complete $\iff$ $X$ has a complete sequence of open covers,
\item $\begin{aligned}[t]
X \textrm{ is } \fsd & \iff X \textrm{ has a complete sequence of countable closed covers} \\
& \iff 	X \textrm{ has a complete sequence of countable } F_\sigma \textrm{ covers},
\end{aligned}$
\item $X$ is $\mc K$-analytic $\iff$ $X$ has a complete sequence of countable covers.
\end{enumerate}
\end{proposition}
The equivalence between first and second part of $2.$ is easy, and follows from Lemma \ref{lemma:refinement}. The remaining assertion are due to Frolík (\cite{frolik1960generalizations}, \cite[Theorem 7]{frolik1963descriptive} and \cite[Theorem 9.3]{frolik1970survey}).  To get our main result, we will prove a statement which has a similar flavor:
\begin{proposition}\label{proposition: fs disjoint covers and abs fsd}
Any topological space with a complete sequence of countable disjoint $F_\sigma$ covers is absolutely $\fsd$.
\end{proposition}

We will need the following observation:
\begin{lemma}\label{lemma:refinement} Let $X$ be a topological space.
\begin{enumerate}
\item If $\left( \mc F_n \right)_{n\in\mathbb N}$ is a complete sequence of covers on $X$ and for each $n\in\mathbb N$, the cover $\mc F_n'$ is a refinement of $\mc F_n$, then the sequence $\left( \mc F_n' \right)_{n\in\mathbb N}$ is complete.
\item If $X$ has a complete sequence of countable closed (open, $F_\sigma$) covers, then it also has a complete sequence of countable closed (open, $F_\sigma$) covers $\left( \mc F_n \right)_{n\in\mathbb N}$, in which each $\mc F_{n+1}$ refines $\mc F_n$.
\end{enumerate}
\end{lemma}
\begin{proof}
The first part follows from the definition of complete sequence of covers of $X$. For the second part, let $\left( \mc F_n' \right)_{n\in\mathbb N}$ be a complete sequence of covers of $X$. We define the new sequence of covers as the refinement of $\left( \mc F_n' \right)_{n\in\mathbb N}$, setting $\mc F_1:=\mc F_1'$ and
\[ \mc F_{n+1}:=\mc F_{n+1}'\wedge\mc F_n := \left\{ F'\cap F|\ F'\in\mc F_n',\ F\in\mc F_{n+1} \right\}. \]
Clearly, the properties of being countable and closed (or $F_\sigma$) are preserved by this operation.
\end{proof}

The main reason for the use of complete sequences of covers is the following lemma:
\begin{lemma}\label{lemma:x is in X and c.s.of c.}
Let $\left( \mc F_n \right)_{n\in\mathbb N}$ be a complete sequence of covers of $X$ and $cX$ a compactification of $X$. If a sequence of sets $F_n\in\mc F_n$ satisfies $F_1\supset F_2\supset \dots$, then $\bigcap_{n \in\mathbb N}\overline{F_n}^{cX} \subset X$.
\end{lemma}
\begin{proof}
Fix $x\in \bigcap_{n \in\mathbb N}\overline{F_n}^{cX}$. We observe that the family
\[ \mc B:=\left\{ U \cap F_n|\ U\in\mc U\left(x\right),\ n\in\mathbb N \right\} \]
is, by hypothesis, formed by nonempty sets and closed under taking finite intersections, therefore it is a basis of some filter $\mc F$ (note that this is the only step where we use the monotonicity of $(F_n)_n$). Since every $F_n$ belongs to both $\mc F$ and $\mc F_n$, $\mc F$ must have some accumulation point $y$ in $X$. If $x$ and $y$ were distinct, they would have some neighborhoods $U$ and $V$ with disjoint closures. This would imply that $V\in\mc U\left (y\right )$, $U\supset U\cap F_1\in \mc F$ and $V\cap\overline{U}$, which contradicts the definition of $y$ being an accumulation point of $\mc F$. This means that $x$ is equal to $y$ and, in particular, $x$ belongs to $X$.
\end{proof}

The property of being hereditarily Lindelöf will be used in the following way:
\begin{lemma}\label{lemma:countable network and disjoint cover}
Every hereditarily Lindelöf $\fsd$ space $X$ has a complete sequence of countable disjoint $F_\sigma$ covers.
\end{lemma}
\begin{proof}
Let $(\mc F_n )_n$ be a complete sequence of countable closed covers of $X$ (which exists by Proposition \ref{proposition: fsd iff complete sequence of covers}). To get the desired result, it suffices to show that each $\mc F_n$ admits a disjoint $F_\sigma$ refinement (by $1.$ in Lemma \ref{lemma:refinement}).

Recall that in a hereditarily Lindelöf Tychonoff space, open sets are $F_\sigma$ and consequently, closed sets are $G_\delta$. Moreover, in any topological space, a countable cover by sets which are both $F_\sigma$ and $G_\delta$ has a disjoint countable refinement by sets of the same type (this is standard - simply enumerate the cover as $\{C_n|\ n\in\N\}$ and define the refinement as $\{\widetilde{C}_n|\ n\in\N\}$, where $\widetilde{C}_n:=C_n\setminus (C_1\cup\dots\cup C_{n-1})$). From these two observations, the existence of the desired refinements is immediate.
\end{proof}

In order to get Theorem \ref{theorem: countable network}, it remains to prove Proposition \ref{proposition: fs disjoint covers and abs fsd}:
\begin{proof}[Proof of Proposition \ref{proposition: fs disjoint covers and abs fsd}]
Let $(\mc D_n)_n$ be a complete sequence of countable disjoint $F_\sigma$ covers of $X$. Without loss of generality, we can assume (by Lemma \ref{lemma:refinement}) that each $\mc D_{n+1}$ refines $\mc D_n$. Also, let $cX$ be a compactification of $X$. Since $cX$ is fixed, all closures will automatically be taken in this compactification.

We enumerate each cover as $\mc D_n=\left\{D^n_m|\ m\in\N \right\}$ and write each of its elements as countable union of closed sets: $D^n_m=\bigcup_i D^n_{m,i}$. We set $\widetilde { \mc D_n } := \left\{ D^n_{m,i} |\ m,i\in\N \right\}$ and $\widetilde { \mc D } :=\bigcup_n \widetilde { \mc D_n }$.

It is clear that $X\subset \bigcap_n \bigcup \left\{ \overline{D}|\ D\in \widetilde {\mc D_n} \right\}$. Note that the set on the right hand side is $\fsd$. The equality does not, in general, hold, but we can modify the right hand side using Lemma \ref{lemma:fsd spaces are Fs separated}.

Indeed, suppose that $x$ belongs to $\bigcap_n \bigcup \left\{ \overline{D}|\ D\in \widetilde {\mc D_n} \right\}$, but not to $X$. By definition, such $x$ satisfies $x \in \overline{D^n_{m_n,i_n}} \subset \overline{D^n_{m_n}}$ for some sequences $(m_n)_n$ and $(i_n)_n$. However, $\left( \mc D_n \right)_n$ is a complete sequence of covers, so by Lemma \ref{lemma:x is in X and c.s.of c.} $D^1_{m_1}\supset D^2_{m_2}\supset \dots$ does \emph{not} hold. Since each $\mc D_{n+1}$ refines $\mc D_n$, the sets $D^{n+1}_{m_{n+1}}$ and $D^n_{m_n}$ must be disjoint for some $n$.

In particular, any such $x$ satisfies $x\in\overline D\cap \overline E$ for some disjoint sets $D,E\in\widetilde {\mc D}$. Since both $D$ and $E$ are closed in $X$, we have $\overline D \cap \overline E \subset cX \setminus X$. This means we can use Lemma \ref{lemma:fsd spaces are Fs separated} to obtain an $F_\sigma$ subset $H_{D,E}$ of $cX$ satisfying $X\subset H_{D,E}\subset cX\setminus \left( \overline D \cap \overline E \right)$. We claim that
\[ X = \bigcap_n \bigcup \left\{ \overline{D}|\ D\in \widetilde {\mc D_n} \right\} \cap \bigcap \left\{H_{D,E}|\  D,E\in\widetilde {\mc D} \textrm{ disjoint} \right\}. \]
By definition of $H_{D,E}$, the set on the right side contains $X$, and the  opposite inclusion follows from the observation above. Since $\widetilde {\mc D}$ is countable, the right hand side is $\fsd$. This proves that that $X\in\fsd \left(cX\right)$, which completes the proof (and also the whole section).
\end{proof}

\section{Conclusion}
We have shown that being hereditarily Lindelöf is a sufficient condition for an $\fsd$ space to be absolutely $\fsd$ - this is a fairly useful condition for applications.
The problem of finding the description of absolute $\fsd$ spaces remains yet unsolved, but we have gotten one step closer to the characterization: By Frolík's result, absolutely $\fsd$ space must have a complete sequence of countable $F_\sigma$ covers. If a space has such a sequence of covers which are also disjoint, then it must be absolutely $\fsd$. Therefore, if the desired characterization can be formulated in terms of complete sequences of covers, it must be something between these two conditions.

\putbib[refs]
\end{bibunit}
\begin{bibunit}[alphanum]

\chapter{Absolute \texorpdfstring{$\mc F$}{F}-Borel Classes}
 \label{chapter:brooms}
\begin{center}
\textit{(Published electronically in Fundamenta Mathematicae (2018), \\
DOI: 10.4064/fm412-10-2017.)}
\end{center}

\textbf{Abstract:}
We investigate and compare $\mc F$-Borel classes and absolute $\mc F$-Borel classes. We provide precise examples distinguishing these two hierarchies. We also show that for separable metrizable spaces, $\mc F$-Borel classes are automatically absolute.

\section{Introduction}\label{sec:introduction}
\emph{Borel} sets (in some topological space $X$) are the smallest family containing the open sets, which is closed under the operations of taking countable intersections, countable unions and complements (in $X$). The family of \emph{$\mc F$-Borel} (resp. \emph{$\mc G$-Borel}) sets is the smallest family containing all closed (resp. open) sets, which is closed under the operations of taking countable intersections and countable unions of its elements. In metrizable spaces the families of Borel, $\mc F$-Borel and $\mc G$-Borel sets coincide. In non-metrizable spaces, open set might not necessarily be a countable union of closed sets, so we need to make a distinction between Borel, $\mc F$-Borel and $\mc G$-Borel sets.

In the present paper, we investigate absolute $\mc F$-Borel classes. While Borel classes are absolute by \cite{holicky2003perfect} (see Proposition \ref{proposition: HS} below), by \cite{talagrand1985choquet} it is not the case for $\mc F$-Borel classes (see Theorem \ref{theorem: talagrand} below). We develop a method of estimating absolute Borel classes from above, which enables us to compute the exact complexity of Talagrand's examples and to provide further examples by modifying them. Our main results are Theorem \ref{theorem: main theorem} and Corollary \ref{corollary: non K alpha space which is abs K alpha+1}.

The paper is organized as follows: In the rest of the introductory section, we define the Borel, $\mc F$-Borel and $\mc G$-Borel hierarchies and recall some basic results. In Section \ref{sec:compactifications}, we recall the definitions and basic results concerning
compactifications and their ordering. In Theorem \ref{theorem: separable metrizable}, we
show that the complexity of separable metrizable spaces coincides with their absolute complexity (in any of the hierarchies). Section \ref{sec:subspaces} is devoted to showing that absolute complexity is inherited by closed subsets. In Section \ref{sec:sequences and brooms}, we study special sets of sequences of integers -- trees, sets which extend to closed discrete subsets of $\baire$ and the `broom sets' introduced by Talagrand. We study in detail the hierarchy of these sets using the notion of rank. In Section \ref{sec:broom spaces and parametrization}, we introduce the class of examples of spaces used by Talagrand. We investigate in detail their absolute complexity and in Section \ref{sec:absolute complexity of brooms}, we prove our main results.

Let us start by defining the basic notions. Throughout the paper, all the spaces will be Tychonoff. For a family of sets $\mc C$, we will denote by $\mc C_\sigma$ the collection of all countable unions of elements of $\mc C$ and by $\mc C_\delta$ the collection of all countable intersections of elements of $\mc C$.
\begin{definition}[Borel classes]\label{def:borel classes}
Let $X$ be a topological space. We define the Borel multiplicative classes $\mc M_\alpha(X)$ and Borel additive classes $\mc A_\alpha(X)$, $\alpha<\omega_1$, as follows:
\begin{itemize}
\item $\mc M_0(X)=\mc A_0(X):=$ the algebra generated by open subsets of $X$,
\item $\mc M_\alpha(X):=\left( \underset {\beta<\alpha} \bigcup \left( \mc M_\beta(X) \cup \mc A_\beta(X) \right) \right)_\delta$ for $1\leq \alpha <\omega_1$,
\item $\mc A_\alpha(X):=\left( \underset {\beta<\alpha} \bigcup \left( \mc M_\beta(X) \cup \mc A_\beta(X) \right) \right)_\sigma$ for $1\leq \alpha <\omega_1$.
\end{itemize}
\end{definition}
In any topological space $X$, we have the families of closed sets and open sets, denoted by $F(X)$ and $G(X)$, and we can continue with classes $F_\sigma(X)$, $G_\delta(X)$, $F_{\sigma\delta}(X)$ and so on. However, this notation quickly gets impractical, so we use the following notation.
\begin{definition}[$\mc F$-Borel and $\mc G$-Borel classes]
We define the hierarchy of $\mc F$-Borel sets on a topological space $X$ as follows:
\begin{itemize}
\item $\mc F_1(X):=$ closed subsets of $X$,
\item $\mc F_\alpha(X):=\left( \underset {\beta<\alpha} \bigcup \mc F_\beta(X) \right)_\sigma$ for $2\leq \alpha <\omega_1$ even,
\item $\mc F_\alpha:=\left( \underset {\beta<\alpha} \bigcup \mc F_\beta (X)\right)_\delta$ for $3\leq \alpha <\omega_1$ odd.
\end{itemize}
The sets of $\alpha$-th $\mc G$-Borel class, $\mc G_\alpha(X)$, are the complements of $\mc F_\alpha(X)$ sets.
\end{definition}
By a \emph{descriptive class of sets}, we will always understand one of the Borel, $\mc F$-Borel or $\mc G$-Borel classes.

\begin{remark}\label{rem:F,G and M,A}
$(i)$ In any topological space $X$, we have
\[ \mc F_1(X) = F(X) \subset \mc M_0(X). \]
It follows that each class $\mc F_n$, $n\in\N$, is contained in $\mc A_{n-1}$ or $\mc M_{n-1}$ (depending on the parity). The same holds for classes $\mc F_\alpha$, $\alpha \geq \omega$, except that the difference between ranks disappears.

$(ii)$ If $X$ is metrizable, the Borel, $\mc F$-Borel and $\mc G$-Borel classes are related to the standard Borel hierarchy \cite[ch.11]{kechris2012classical}. In particular, each $\mc F_\alpha(X)$ is equal to either $\Sigma^0_\alpha(X)$ or $\Pi^0_\alpha(X)$ (depending on the parity of $\alpha$) and we have
\[ F(X)\cup G(X) \subset \mc M_0(X)=\mc A_0(X) \subset F_\sigma(X)\cap G_\delta(X). \]
The relations between our descriptive classes can then be summarized as:
\begin{itemize}
\item $\mc F_n(X)=\Pi^0_n(X)=\mc M_{n-1}(X)$ holds if $1\leq n<\omega$ is odd,
\item $\mc F_n(X)=\Sigma^0_n(X)=\mc A_{n-1}(X)$ holds if $2\leq n<\omega$ is even,
\item $\mc F_\alpha(X)=\Sigma^0_\alpha(X)=\mc A_\alpha(X)$ holds if $\omega\leq \alpha <\omega_1$ is even,
\item $\mc F_\alpha(X)=\Pi^0_\alpha(X)=\mc M_\alpha(X)$ holds if $\omega\leq \alpha <\omega_1$ is odd.
\end{itemize}
Clearly, the $\mc G$-Borel classes satisfy the dual version of $(i)$ and $(ii)$.
\end{remark}
Note also that in compact spaces, closed sets are compact, so in this context the $\mc F$-Borel sets are sometimes called $\mc K$-Borel, $F_\sigma$ sets are called $K_\sigma$ and so on.

Let us define two notions central to the topic of our paper.
\begin{definition}\label{def:F sigma delta spaces}
Let $X$ be a (Tychonoff) topological space and $\mc C$ be a descriptive class of sets. We say that $X$ is a \emph{$\mc C$ space} if there exists a compactification $cX$ of $X$, such that $X\in\mc C(cX)$.

If $X\in\mc C(Y)$ holds for any Tychonoff topological space $Y$ in which $X$ is densely embedded, we say that $X$ is an \emph{absolute $\mc C$ space}\err{.} We call the class $\mc C$ absolute if every $\mc C$ space is an absolute $\mc C$ space.
\end{definition}

The basic relation between complexity and absolute complexity are noted in the following remark.
\begin{remark}\label{rem:1-5}
Consider the following statements:
\begin{enumerate}[(i)]
\item $X\in\mc C(cX)$ holds for some compactification $cX$;
\item $X\in\mc C(\beta X)$ holds for the Čech-Stone compactification;
\item $X\in\mc C(cX)$ holds for every compactification $cX$;
\item $X\in\mc C(Y)$ holds for every Tychonoff space where $Y$ is densely embedded;
\item $X\in\mc C(Y)$ holds for every Tychonoff space where $Y$ is embedded.
\end{enumerate}
Clearly, the implications $(v)\implies(iv)\implies(iii)\implies(ii)\implies(i)$ always hold. In the opposite direction, we always have $(i)\implies(ii)$ (this is standard, see Remark \ref{rem:abs_complexity}) and $(iii)\implies (iv)$. For Borel and $\mc F$-Borel classes, $(iv)$ is equivalent to $(v)$ (since these classes are closed under taking intersections with closed sets). For $\mc G$-Borel classes, $(iv)$ is never equivalent to $(v)$ (just take $Y$ \err{which is ``so large''} that $\overline{X}^Y$ is not $G_\delta$ in $Y$).
\end{remark}
The interesting part is therefore the relation between $(i)$ and $(iii)$ from Remark \ref{rem:1-5}. For the first two levels of $\mc F$-Borel and $\mc G$-Borel hierarchies, $(i)$ is equivalent to $(iii)$, that is, these classes are absolute. A more precise formulation is given in the following remark.

\begin{remark}\label{rem:F,G,Fs,Gd are absolute}
For a topological space $X$ we have
\begin{enumerate}[(i)]
	\item $X$ is an $\mc F_1$ space $\iff$ $X$ is absolutely $\mc F_1$ $\iff$ $X$ is compact;
	\item $X$ is an $\mc F_2$ space $\iff$ $X$ is absolutely $\mc F_2$ $\iff$ $X$ is $\sigma$-compact;
	\item $X$ is a $\mc G_1$ space $\iff$ $X$ is absolutely $\mc G_1$ $\iff$ $X$ is locally compact;
	\item $X$ is a $\mc G_2$ space $\iff$ $X$ is absolutely $\mc G_2$ $\iff$ $X$ is Čech-complete
\end{enumerate}
(for the proofs of $(iii)$ and $(iv)$ see \cite[Theorem 3.5.8 and Theorem 3.9.1]{engelking1989general}).
\end{remark}

The first counterexample was found by Talagrand, who showed that already the class of $\fsd$ sets is not absolute:
\begin{theorem}[\cite{talagrand1985choquet}]\label{theorem: talagrand}
There exists an $\mc F_3$ space $T$ and its compactification $K$, such that $T$ is not $\mc F$-Borel in $K$.
\end{theorem}
This not only shows that none of the classes $\mc F_\alpha$ for $\alpha\geq 4$ are absolute, but also that the difference between absolute and non-absolute complexity can be large. Indeed, the `non-absolute complexity' of $T$ is `$\mc F_3$', but its `absolute complexity' is `$\mc K$-analytic'. In \cite{kalenda2018absolute}, the authors give a sufficient condition for an $\mc F_3$ space to be absolutely $\mc F_3$, but a characterization of absolutely $\mc F_\alpha$ spaces is still unknown for all $\alpha\geq 3$.

In \cite{holicky2003perfect} and \cite{raja2002some}, the authors studied (among other things) absoluteness of Borel classes. In particular, the following result is relevant to our topic:
\begin{proposition}[{\cite[Corollary 14]{holicky2003perfect}}]\label{proposition: HS}
For every $1\leq \alpha < \omega_1$, the classes $\mc A_\alpha$ and $\mc M_\alpha$ are absolute.
\end{proposition}
Another approach is followed in \cite{marciszewski1997absolute} and \cite{junnila1998characterizations}, where the absoluteness is investigated in the metric setting (that is, for spaces which are of some class $\mc C$ when embedded into any \emph{metric} space). In \cite{junnila1998characterizations} a characterization of `metric-absolute' $\mc F_3$ spaces is given in terms of existence complete sequences of covers (a classical notion used by Frolík for characterization of Čech-complete spaces, see \cite{frolik1960generalizations}) which are $\sigma$-discrete. Unfortunatelly, this is not applicable to the topological version of absoluteness, because every countable cover is (trivially) $\sigma$-discrete, and any $\mc F_3$ space does have such a cover by \cite{frolik1963descriptive} -- even the Talagrand's non-absolute space.

\section{Compactifications and Their Ordering}\label{sec:compactifications}
By a \emph{compactification} of a topological space $X$ we understand a pair $(cX,\varphi)$, where $cX$ is a compact space and $\varphi$ is a homeomorphic embedding of $X$ onto a dense subspace  of $cX$. Symbols $cX$, $dX$ and so on will always denote compactifications of $X$.

Compactification $(cX,\varphi)$ is said to be \emph{larger} than $(dX,\psi)$, if there exists a continuous mapping $f : cX\rightarrow dX$, such that $\psi = f \circ \varphi$. We denote this as $ cX \succeq dX $. Recall that for a given $T_{3\slantfrac{1}{2}}$ topological space $X$, its compactifications are partially ordered by $\succeq$ and Stone-Čech compactification $\beta X$ is the largest one.

Often, we encounter a situation where $X\subset cX$ and the corresponding embedding is identity. In this case, we will simply write $cX$ instead of $(cX,\textrm{id}|_X)$.

Much more about this topic can be found in many books, see for example \cite{freiwald2014introduction}. The basic relation between the complexity of a space $X$ and the ordering of compactifications is the following:
\begin{remark}\label{rem:abs_complexity}
If $\mc C$ is a descriptive \err{class} of sets, we have
 \[ X\in\mc C(dX), \ cX\succeq dX \implies X\in\mc C(cX). \]
In particular, $X$ is a $\mc C$ space if and only if $X\in\mc C(\beta X)$.
\end{remark}
We will also make use of the following result about existence of small metrizable compactifications.
\begin{proposition}\label{proposition: metrizable sub-compactification}
Let $X$ be a separable metrizable space and $cX$ its compactification. Then $X$ has some metrizable compactification $dX$, such that $dX\preceq cX$.
\end{proposition}
This proposition is an exercise, so we only include a sketch of its proof:
\begin{proof}
We can assume that $cX\subset [0,1]^\kappa$ for some $\kappa$. Since $X$ has a countable base, there is a countable set of coordinates $I\subset \kappa$ such that the family of all $V \cap X$,
where $V$ is a basic open set $V\subset [0,1]^\kappa$ depending on coordinates in $I$ only is a base
of $X$. Then we can take $dX := \pi_I (cX)$, where $\pi_I$ denotes the projection.
\end{proof}

We put together several known results to obtain the following theorem (which we have not found anywhere in the literature):
\begin{theorem}\label{theorem: separable metrizable}
Let $\mc C$ be a descriptive class of sets. Then any separable metrizable $\mc C$ space is an absolute $\mc C$ space (so the first four conditions of Remark \ref{rem:1-5} are equivalent).

If $\mc C$ is one of the Borel or $\mc F$-Borel classes, this is further equivalent to $X$ being in $\mc C(Y)$ for every Tychonoff space in which $Y$ is embedded.
\end{theorem}
\begin{proof}
The statement for the classes $\mc A_\alpha$ and $\mc M_\alpha$ follows from the general result of \cite{holicky2003perfect} (see Proposition \ref{proposition: HS} above). Let us continue by proving the result for the $\mc F$-Borel classes. For the first two levels, it follows from Remark \ref{rem:F,G,Fs,Gd are absolute}.
Suppose that $X\in \mc F_\alpha(\beta(X))$ for some $\alpha\ge 3$ and let $cX$ be any compactification of $X$.

By Proposition \ref{proposition: metrizable sub-compactification} we can find a smaller metrizable compactification $dX$. Let $\mc D$ be the Borel class corresponding to $\mc F_\alpha$ (see Remark \ref{rem:F,G and M,A}). Then $X\in \mc D(\beta X)$. Since $\mc D$ is absolute by Proposition \ref{proposition: HS}, we get $X\in\mc D(dX)$. Since $dX$ is metrizable, we deduce that $X\in\mc F_\alpha(dX)$, so $X\in \mc F_\alpha(cX)$ holds by Remark \ref{rem:abs_complexity}.

The proof for $\mc G$-Borel classes is analogous.

\end{proof}

\section{Hereditarity of Absolute Complexity} \label{sec:subspaces}
In this section, we show that absolute complexity is hereditary with respect to closed subspaces. To do this, we first need to be able to extend compactifications of subspaces to compactifications of the original space. We start with a topological lemma:

\begin{lemma}\label{lemma: quotient of K}
Suppose that $K$ is compact, $F\subset K$ is closed and $f:F\rightarrow L$ is a continuous surjective mapping. Define a mapping $q: K \to L\cup (K\setminus F)$ as $q(x)=x$ for $x\in K\setminus F$ and $q(x)=f(x)$ for $x\in F$. Then $\left(K\setminus F\right)\cup L$ equipped with the quotient topology (induced by $q$) is a compact space.
\end{lemma}
\begin{proof}
Since $f$ is surjective, the space $(K\setminus F)\cup L$ coincides with the adjunction space $K\cup_f L$ determined by $K$, $L$ and $f$ (for definition, see \cite[below Ex. 2.4.12]{engelking1989general}). Since $K\cup_f L$ is a quotient of a compact space $K\oplus L$ with respect to a closed equivalence relation (\cite[below Ex. 2.4.12]{engelking1989general}), we get by Alexandroff theorem (\cite[Thm 3.2.11]{engelking1989general}) that the space $K\cup_f L$ is compact (in particular, it is Hausdorff).
\end{proof}
This gives us a way to extend a compactification, provided that we already know how to extend some bigger compactification:
\begin{proposition}\label{proposition: extending a compactification}
Let $Y$ be a topological space and suppose that for a closed $X\subset Y$ we have $cX\preceq \overline{X}^{dY}$ for some compactifications $cX$ and $dY$. Then $cX$ is equivalent to \errata{$\overline{X}^{cY}$} for some compactification $cY$ of $Y$.
\end{proposition}
\begin{proof}
Let $X$, $Y$, $cX$ and $dY$ be as above and denote by $f:\overline{X}^{dY}\rightarrow cX$ the mapping which witnesses that $cX$ is smaller than $\overline{X}^{dY}$. We will assume that $X\subset cX$ and $Y\subset dY$, which means that $f|_X=\textrm{id}_X$.

By Lemma \ref{lemma: quotient of K} (with $K:=dY$, $F:=\overline{X}^{dY}$ and $L:=cX$), the space $cY$, defined by the formula
\[ cY:= \left(dY\setminus \overline{X}^{dY} \right)\cup cX =  \left(K\setminus F\right)\cup L, \]
is compact. To show that $cY$ is a compactification of $Y$, we need to prove that the mapping $q: dY\rightarrow cY$, defined as $f$ extended to $dY\setminus \overline{X}^{dY}$ by identity, is a homeomorphism when restricted to $Y$. The continuity of $q|_Y$ follows from the continuity of $q$. The restriction is injective, because \errata{$q$} is injective on $\left( dY \setminus \overline{X}^{dY}\right) \cup X$ (and $X$ is closed in $Y$). The definition of $q$ and the fact that $X$ is closed in $Y$ imply that the inverse mapping is continuous on $Y\setminus X$. To get the continuity of $(q|_Y)^{-1}$ at the points of $X$, let $x\in X$ and $U$ be a neighborhood of $x$ in $K$. We want to prove that $q(U\cap Y)$ is a neighborhood of $q(x)$ in $q(Y)$. To see this, observe that $C:=dY\setminus U$ is compact, therefore $q(C)$ is closed, its complement $V:=cY \setminus q(C)$ is open and it satisfies $V\cap q(Y) = q(U\cap Y)$.

Lastly, we note that for $y\in cY\setminus \overline{q(X)}^{cY}$, $q^{-1}(y)$ is in $dY\setminus \overline{X}^{dY}$. This gives the second identity on the following line, and completes the proof:
\[ \overline{X}^{cY} = \overline{q(X)}^{cY} = q\left( \overline{X}^{dY} \right) = f\left( \overline{X}^{dY} \right)=cX. \]
\end{proof}

Recall that a subspace $X$ of $Y$ is said to be \emph{$C^\star$-embedded} in $Y$ if any bounded continuous real function on $X$ can be continuously extended to $Y$.
\begin{corollary}\label{corollary: extending all compactifications}
If a closed set $X$ is $C^\star$-embedded in $Y$, then each compactification $cX$ is of the form $cX=\overline{X}^{cY}$ for some compactification $cY$.
\end{corollary}
\begin{proof}
Recall that the Čech-Stone compactification $\beta X$ is characterized by the property that each bounded continuous function from $X$ has a unique continuous extension to $\beta X$. Using this characterization, it is a standard exercise to show that $X$ is $C^\star$-embedded in $Y$ if and only if $\overline{X}^{\beta Y}=\beta X$. It follows that for $X$ and $Y$ as above, any compactification $cX$ is smaller than $\overline{X}^{\beta Y}$, and the result follows from Proposition \ref{proposition: extending a compactification}.
\end{proof}
We use these results to get the following corollary:
\begin{proposition}\label{proposition: absoluteness is hereditary}
Let $\mc C$ be one of the classes $\mc F_\alpha$, $\mc M_\alpha$ or $\mc A_\alpha$ for some $\alpha<\omega_1$.
\begin{enumerate}[(i)]
\item Any closed subspace of a $\mc C$ space is a $\mc C$ space;
\item Any closed subspace of an absolute $\mc C$ space is an absolute $\mc C$ space.
\end{enumerate}
\end{proposition}
\begin{proof}
$(i)$: This is trivial, since if $X$ is a closed subspace of $Y$ and $Y$ is of the class $\mc C$ in some compactification $cY$, then $X$ is of the same class in the compactification $\overline{X}^{cY}$.

$(ii)$: For Borel classes, which are absolute, the result follows from the first part -- therefore, we only need to prove the statement for the $\mc F$-Borel classes. Let  $\alpha<\omega_1$ and assume that $X$ is a closed subspace of some $\mc F_\alpha$ space $Y$.

Firstly, if $Y$ is an $\mc F_\alpha$ subset of $\beta Y$, it is $\mc K$-analytic and, in particular, Lindelöf. Since we also assume that $Y$ is \err{Tychonoff}, we get that $Y$ is normal \cite[Propositions 3.4 and 3.3]{kkakol2011descriptive}. By Tietze's theorem, $X$ is $C^\star$-embedded in $Y$, which means that every compactification of $X$ can be extended to a compactification of $Y$ (Corollary \ref{corollary: extending all compactifications}). We conclude the proof as in $(i)$.
\end{proof}

Note that the second part of Proposition \ref{proposition: absoluteness is hereditary} does not hold, in a very strong sense, if we replace the closed subspaces by $\fsd$ subspaces. Indeed, the space $T$ from Theorem \ref{theorem: talagrand} is an $\fsd$ subspace of a compact space, but it is not even absolutely $\mc F$-Borel (that is, $\mc F$-Borel in every compactification). Whether anything positive can be said about $F_\sigma$ subspaces of absolute $\mc F_\alpha$ spaces is unknown.

\section{Special Sets of Sequences}\label{sec:sequences and brooms}
In this chapter, we define `broom sets' -- a special type \err{of} sets of finite sequences, which have a special `tree' structure and all of their infinite extensions are closed and discrete. As in \cite{talagrand1985choquet}, these sets will then be used to construct spaces with special complexity properties.

\subsection{Trees and a Rank on Them}
First, we introduce the (mostly standard) notation which will be used in the sequel.
\begin{notation}[Integer sequences of finite and infinite length]\label{not:sequences}
We denote
\begin{itemize}
\item $\omega^\omega:=$ infinite sequences of non-negative integers $:=\left\{\sigma : \omega \rightarrow \omega\right\}$,
\item $\omega^{<\omega}:=$ finite sequences of non-neg. integers $:=\left\{s : n \rightarrow \omega |\ n\in\omega \right\}$.
\end{itemize}
Suppose that $s\in\omega^{<\omega}$ and $\sigma\in\baire$. We can represent $\sigma$ as $(\sigma(0),\sigma(1),\dots)$ and $s$ as $(s(0),s(1),\dots,s(n-1))$ for some $n\in\omega$. We denote the \emph{length} of $s$ as $|s|=\textrm{dom}(s)=n$, and set $|\sigma|=\omega$. If for some $t\in \omega^{<\omega} \cup \baire$ we have $|t|\geq |s|$ and $t|_{|s|}=s$, we say that $t$ \emph{extends} $s$, denoted as $\err{s}\sqsubset t$. We say that  $u,v\in\omega^{<\omega}$ are non-comparable, denoting as $u\perp v$, when neither $u\sqsubset v$ nor $u\sqsupset v$ holds.

Unless we say otherwise, $\omega^\omega$ will be endowed with the standard product topology, whose basis consists of sets $\mc N (s):=\left\{ \sigma \sqsupset s | \ \sigma\in\baire \right\}$, $s \in \omega^{<\omega}$. 

For $n\in\omega$, we denote by $(n)$ the corresponding sequence of length $1$. By $s\ext t$ we denote the concatenation of a sequences $s\in\omega^{<\omega}$ and $t\in\omega^{<\omega}\cup \baire$. We will also use this notation to extend finite sequences by integers and sets \err{of sequences}, using the convention $s \hat{\ } k:=s \hat{\ } (k)$ for $k\in\omega$ and $s \hat \ T:=\left\{s\hat\ t\ |\ t\in T\right\}$ for $T\subset \omega^{<\omega} \cup \baire$.
\end{notation}

\begin{definition}[Trees on $\omega$]
A \emph{tree} (on $\omega$) is a set $T\subset \omega^{<\omega}$ which satisfies
\[ (\forall s,t\in\omega^{<\omega}): s\sqsubset t \ \& \ t\in T \implies s\in T. \]
By $\textrm{Tr}$ we denote the space of all trees on $\omega$. For $S\subset \omega^{<\omega}$ we denote by $\mathrm{cl}_\mathrm{Tr}(S):=\left\{u\in\omega^{<\omega} | \ \exists s\in S: s\sqsupset u \right\}$ the smallest tree containing $S$. Recall that the empty sequence $\emptyset$ can be thought of as the `root' of each tree, since it is contained in any nonempty tree.

If each of the initial segments $\sigma|n$, $n\in\omega$, of some $\sigma\in\baire$ belongs to $T$, we say that $\sigma$ is an infinite branch of $T$. By $\textrm{WF}$ we denote the space of all trees which have no infinite branches (the `well founded' trees).
\end{definition}

We define a variant of the standard derivative on trees and the corresponding $\Pi^1_1$-rank. We do not actually use the fact that the rank introduced in Definition \ref{def:derivative} is a $\Pi^1_1$-rank, but more details about such ranks and derivatives can be found in \cite[ch. 34 D,E]{kechris2012classical}.
\begin{definition}\label{def:derivative}
Let $T\in \textrm{Tr}$. We define its `infinite branching' \emph{derivative} $D_i(T)$ as
\[ D_i(T):=\left\{ t\in T| \ t\sqsubset s \textrm{ holds for infinitely many }s\in T \right\}. \]
This operation can be iterated in the standard way:
\begin{align*}
D^0_i(T) & :=  T, \\
D^{\alpha +1}_i \left( T \right) & :=  D_i\left(D_i^\alpha\left(T\right)\right) \textrm{ for successor ordinals},\\
D_i^\lambda \left( T \right) & := \underset {\alpha < \lambda } \bigcap D_i^\alpha \left( T \right) \textrm{ for limit ordinals}.
\end{align*}
This allows us to define a rank $r_i$ on $\textrm{Tr}$ as $r_i(\emptyset):=-1$, $r_i(T):=0$ for finite trees and
\[ r_i(T):=\min \{ \alpha<\omega_1| \ D_i^{\alpha+1}\left( T \right) = \emptyset \}. \]
If no such $\alpha<\omega_1$ exists, we set $r_i(T):=\omega_1$.
\end{definition}

Note that on well founded trees, the derivative introduced in Definition \ref{def:derivative} behaves the same way as the derivative from \cite[Exercise 21.24]{kechris2012classical}, but it leaves any infinite branches untouched. It follows from the definition that the rank $r_i$ of a tree $T$ is countable if and only if $T$ is well founded. We consider this approach better suited for our setting.

While this definition of derivative and rank is the most natural for trees, we can extend it to all subsets $S$ of $\omega^{<\omega}$ by setting $D_i(S):=D_i(\mathrm{cl}_\mathrm{Tr}(S))$ and defining the rest of the notions using this extended `derivative'. Clearly, if $S\in Tr$, this coincides with the original definition.

\begin{lemma}\label{lemma: rank and finite covers}
No set $S\subset \omega^{<\omega}$ which satisfies $r_i(S)<\omega_1$ can be covered by finitely many sets of rank strictly smaller than $r_i(S)$.
\end{lemma}
\begin{proof}
Note that if $T_1$ and $T_2$ are trees, $T_1\cup T_2$ is a tree as well, and we have $D_i^\alpha(T_1\cup T_2)=D_i^\alpha(T_1)\cup D_i^\alpha(T_2)$ for any $\alpha<\omega_1$. Moreover, for any $S_1,S_2\subset \omega^{<\omega}$ we have $\mathrm{cl}_\mathrm{Tr}(S_1\cup S_2)=\mathrm{cl}_\mathrm{Tr}(S_1)\cup \mathrm{cl}_\mathrm{Tr}(S_2)$. This yields the formula
\[ r_i(S_1\cup\dots\cup S_k)=\max \{ r_i(S_1),\dots,r_i(S_k) \} \]
for any $S_i\subset \omega^{<\omega}$, which gives the result.
\end{proof}

\subsection{Sets which Extend into Closed Discrete Subsets of \texorpdfstring{$\baire$}{the Baire space}}
\begin{definition}\label{def:D}
A set $A\subset \baire$ is said to be an infinite extension of $B\subset \omega^{<\omega}$, if there exists a bijection $\varphi : B\rightarrow A$, such that $(\forall s\in B) : \varphi (s)\sqsupset s$. If $\tilde A$ has the same property as $A$, except that we have $\tilde A \subset \omega^{<\omega}$, it is said to be a (finite) extension of $B$.

By $\mc D\subset \mc P\left(\omega^{<\omega}\right)$ we denote the system of those sets $D\subset \omega^{<\omega}$ which satisfy
\begin{itemize}
\item[(i)] $(\forall s,t\in D)$ : $s\neq t \implies s \perp t$;
\item[(ii)] every infinite extension of $D$ is closed and discrete in $\baire$.
\end{itemize}  
\end{definition}

\begin{example}[Broom sets of the first class]\label{example: brooms}
Let $h\in \omega^{<\omega}$ and $s_n\in\omega^{<\omega}$ for $n\in\omega$ be finite sequences and suppose that $(f_n)_{n\in\omega}$ is an injective sequence of elements of $\omega$. Using these sequences, we define a `broom' set $B$ as
\[ B:=\left\{ h\hat{\ }f_n\hat{\ }s_n|\ n\in\omega\right\}. \]
The intuition is that $h$ is the `handle' of the broom, $s_n$ are the `bristles', and the sequence $(f_n)_n$ causes these bristles to `fork out' from the end of the handle. Note that the injectivity of $(f_n)_n$ guarantees that $B\in \mc D$ (see Proposition \ref{proposition: properties of D}).
\end{example}

\begin{proposition}[Properties of $\mc D$]\label{proposition: properties of D}
The system $\mc D$ has the following properties (where $D$ is an arbitrary element of $\mc D$):
\begin{enumerate}[(i)]
\item $\{ (n) | \ n\in\omega \} \in \mc D$;
\item $(\forall E\subset \omega^{<\omega}): E\subset D \implies E\in\mc D$;
\item $(\forall E\subset \omega^{<\omega}): E$ is a finite extension of $D \implies E\in\mc D$;
\item $\left(\forall h\in\omega^{<\omega}\right): h\ext D\in\mc D$;
\item For each sequence $\left(D_n\right)_n$ of elements of $\mc D$ and each one-to-one enumeration  $\left\{d_n| \ n\in\omega\right\}$ of $D$, we have $\bigcup _{n\in\omega} d_n\hat{\ }D_n \in \mc D$.
\end{enumerate}
\end{proposition}
\begin{proof}
The first four properties are obvious. To prove \err{(v)}, let $D$, $D_n$, $E:=\bigcup _{n\in\omega} d_n\hat{\ }D_n$ be as above. \err{To} show that $E$ satisfies (i) from Definition \ref{def:D}, suppose that $e=d_n\ext e',\ f=d_m\ext f' \in D$ are comparable. Then $d_n$ and $d_m$ are also comparable, which means that $m=n$ and $e',f'\in D_n=D_m$ are comparable. It follows that $e'=f'$ and hence $e=f$.

For $A\subset \omega^{<\omega}$, the condition (ii) from Definition \ref{def:D} is equivalent to
\begin{equation}\label{eq: equivalence}
(\forall \sigma \in \baire)(\exists m\in\omega): \  \left\{ s\in A | \ s \textrm{ is comparable with } \sigma|m \right\} \textrm{ is finite}.
\end{equation}
Indeed, suppose \eqref{eq: equivalence} holds and let $H$ be an infinite extension of $A$. For $\sigma\in\baire$ let $m\in\omega$ be as in \eqref{eq: equivalence}. Any $\mu\in H\cap \mc N(\sigma|m)$ is an extension of some $s\in A$, which must be comparable with $\sigma|m$. It follows that the neighborhood $\mc N(\sigma|m)$ of $\sigma$ only contains finitely many elements of $H$, which means that $\sigma$ is not a cluster point of $H$ and $H$ is discrete.
In the opposite direction, suppose that the condition does not hold. Then we can find a one-to-one sequence $(s_m)_{m\in\omega}$ of elements of $A$, such that $s_m \sqsupset \sigma|m$ for each $m\in\omega$. Then $\sigma$ is a cluster point of (any) infinite extension of $\{ s_m|\ m\in\omega \}$.

If $A$ also satisfies (i) from Definition \ref{def:D}, \eqref{eq: equivalence} is clearly equivalent to
\begin{equation}\label{eq: equivalence with (i)}
(\forall \sigma \in \baire)(\exists m\in\omega): \  \sigma|m \textrm{ is comparable with at most one } s\in A.
\end{equation}
To see that $E$ satisfies \eqref{eq: equivalence with (i)}, let $\sigma\in \baire$. Since $D\in \mc D$, let $m$ be the integer obtained by applying \eqref{eq: equivalence with (i)} to this set. If no $d\in D$ is comparable with $\sigma|m$, then neither is any $e\in E$ and we are done.

Otherwise, let $d_{n_0}$ be the only element of $D$ comparable with $\sigma|m$ and let  $m_0$ be the integer obtained by applying \eqref{eq: equivalence with (i)} to $d_{n_0}\ext D_{n_0}$. Then $m':=\max\{m_0,m\}$ witnesses that $E$ satisfies \eqref{eq: equivalence with (i)}: Indeed, any $e=d_n\ext e'$ comparable with $\sigma|{m'}$ must have $d_n$ comparable with $\sigma|{m'}$. In particular, $d_n$ is comparable with $\sigma|m$. It follows that $d_n=d_{n_0}$ and $e\in d_{n_0}\ext D_{n_0}$. Therefore, $e$ must be comparable with $\sigma|{m_0}$ and since there is only one such element, the proof is finished.
\end{proof}

\subsection{Broom Sets}
We are now ready to define the broom sets and give some of their properties.
\begin{definition}[Broom sets of higher classes]\label{def:brooms}
We define the hierarchy of broom sets $\mc B_\alpha$, $\alpha < \omega_1$. We set $\mc B_0:=\left\{ \left\{ s \right\} | \ s\in\omega^{<\omega} \right\}$ and, for $1<\alpha<\omega_1$ we define $\mc B_\alpha$ as the union of $\bigcup_{\beta < \alpha}\mc B_\beta$ and the collection of all sets $B\subset \omega^{<\omega}$ of the form
\[ B= \underset {n\in \omega}\bigcup h\hat{\ }f_n\hat{\ }B_n \]
for some finite sequence $h\in\omega^{<\omega}$ (`handle'), broom sets $B_n\in \bigcup_{\beta < \alpha}\mc B_\beta$ (`bristles') and one-to-one sequence $(f_n)_{n\in\omega} \in \baire$ (`forking sequence'). We also denote as $\mc B_{\omega_1}:=\underset {\alpha<\omega_1} \bigcup\mc B_\alpha$ the collection of all these broom classes.
\end{definition}

In Section \ref{sec:broom spaces and parametrization}, we will need a result of the type `if $B$ is an element of $\mc B_{\alpha+1}\setminus \mc B_\alpha$, then $B$ cannot be covered by finitely many brooms of class $\alpha$'. \err{Since this is not true in general}, we restrict ourselves to the following subcollection $\widetilde{\mc B}_{\alpha+1} \subset \mc B_{\alpha+1}\setminus \mc B_\alpha$ of `$(\alpha+1)$-brooms whose bristles are of the highest allowed class' for which the result holds (as we will see later):
\begin{definition} \label{def:tilde brooms}
We set $\widetilde {\mc B}_0:=\mc B_0$. For successor ordinals we define $\widetilde {\mc B}_{\alpha +1}$ as the collection of all sets $B\subset \omega^{<\omega}$ of the form
\[ B= \underset {n\in \omega}\bigcup h\hat{\ }f_n\hat{\ }B_n, \]
where $h\in\omega^{<\omega}$, $B_n\in \widetilde{\mc B}_\alpha$ and $(f_n)_n$ is the (1-to-1) forking sequence. For a limit ordinal $\alpha<\omega_1$, the family $\widetilde {\mc B}_{\alpha}$ is defined in the same way, except that the bristles $B_n$ belong to $\widetilde {\mc B}_{\beta_n}$ for some $\beta_n<\alpha$, where $\sup_n \beta_n=\alpha$.
\end{definition}

\begin{remark}[Basic properties of broom sets]\label{rem:basic broom properties}
It is clear that for every broom $B\in\mc B_\alpha$, any finite extension of $B$ is also a broom of class $\alpha$, and so is the set $h\ext B$ for any $h\in\omega^{<\omega}$. Moreover, it follows from Proposition \ref{proposition: properties of D} that $\mc B_{\omega_1} \subset \mc D$. These properties also hold for collections $\widetilde {\mc B}_\alpha$, $\alpha<\omega_1$.

Actually, broom sets from Definition \ref{def:brooms} are exactly those elements of $\mc D$ for which every `branching point' is a `point of infinite branching'. However, we do not need this description of broom sets, so we will not give its proof here.
\end{remark}

Broom sets are connected with rank $r_i$ in the following way:
\begin{proposition}\label{proposition: rank of B}
For $B\in\mc B_{\omega_1}$ and $\alpha<\omega_1$, the following holds:
\begin{enumerate}[(i)]
\item $B\in \mc B_\alpha \implies r_i(B)\leq \alpha$,
\item $B\in \widetilde{\mc B}_\alpha \implies r_i(B) = \alpha$.
\end{enumerate}
\end{proposition}
\begin{proof}
Firstly, we make the following two observations:
\begin{equation*}
\begin{split}
(\forall B\subset \omega^{<\omega}): r_i(B)=-1 & \iff B=\emptyset, \\
(\forall B \in\mc B_{\omega_1}): r_i(B)=0 & \iff B \in\mc B_0 = \widetilde {\mc B}_0. 
\end{split}
\end{equation*}
We use this as a starting point for transfinite induction. We now prove $(i)$ and $(ii)$ separately, distinguishing also between the cases of successor $\alpha$ and limit $\alpha$.

$(i)$: Suppose that the statement holds for $\alpha<\omega_1$ and let $B\in {\mc B}_{\alpha+1}$. By definition, $B=\bigcup_n h\hat{\ }f_n\hat{\ }B_n$ holds for some handle $h\in\omega^{<\omega}$, forking sequence $(f_n)_n$ and $B_n\in\mc B_\alpha$. For each $n\in\omega$ we have $r_i(B_n)\leq \alpha$, which gives $D_i^{\alpha+1}(B_n)=\emptyset$ and hence $D_i^{\alpha+1}(f_n\hat{\ }B_n)=\emptyset$. It follows that
\[ D_i^{\alpha+1}(\bigcup_n h\hat{\ }f_n\hat{\ }B_n) \subset \mathrm{cl}_\mathrm{Tr}(h), \]
and
\[ D_i^{\alpha+2}(\bigcup_n h\hat{\ }f_n\hat{\ }B_n) = \emptyset ,\]
which, by definition \err{of} $r_i$, means that $r_i(B)\leq \alpha+1$.

For the case of limit ordinal $\alpha<\omega_1$, assume that the statement holds for every $\beta<\alpha$, and let $B\in {\mc B}_{\alpha}$. As above, $B=\bigcup_n h\hat{\ }f_n\hat{\ }B_n$ holds for some $h$, $(f_n)_n$ and $B_n$, where $B_n\in\mc B_{\beta_n}$, $\beta_n + 1 <\alpha$. We successively deduce that $r_i(B_n)$ is at most $\beta_n$, hence $D_i^{\alpha}(B_n) \subset D_i^{\beta_n+1}(B_n)$ holds and we have $D_i^{\alpha}(f_n\hat{\ }B_n)=\emptyset$.
Since $(f_n)_n$ is a forking sequence, then for any $s\in f_{n_0}\hat{\ }B_{n_0}$ which is different from the empty sequence, we get the following equivalence:
\begin{equation}\label{eq: r(B)<=alpha}
h\hat{\ }s\in D_i^{\alpha}(\bigcup_n h\hat{\ }f_n\hat{\ }B_n) \iff
 \left( \exists n_0\in \omega \right): s\in D_i^{\alpha}(f_{n_0}\hat{\ }B_{n_0}).
\end{equation}
From \eqref{eq: r(B)<=alpha} we conclude that $D_i^{\alpha}(B)\subset \mathrm{cl}_\mathrm{Tr}(\{h\})$ is finite, which means that $r_i(B)\leq \alpha$.

$(ii)$: Suppose that the second part of the proposition holds for some $\alpha<\omega$ and let $B=\bigcup_n h\hat{\ }f_n\hat{\ }B_n\in\widetilde{\mc B}_{\alpha+1}$. For each $n\in\omega$, $B_n$ is an element of $\widetilde{\mc B}_\alpha$, which by the induction hypothesis means that $r_i(B_n) = \alpha$, and thus there exists some sequence $s_n\in D_i^\alpha(B)$. In particular, we have $h\hat{\ }f_n\hat{\ }s_n\in D_i^\alpha(h\hat{\ }f_n\hat{\ }B_n)$. We conclude that $h\in D_i^{\alpha+1}(B)$ and $r_i(B)\geq \alpha+1$.

Assume now that $\alpha<\omega_1$ is a limit ordinal and that $(ii)$ holds for all $\beta<\alpha$. Let $B=\bigcup_n h\hat{\ }f_n\hat{\ }B_n\in\widetilde{\mc B}_{\alpha}$, where $B_n\in\widetilde{\mc B}_{\beta_n}$, $\sup_n \beta_n=\alpha$. For each $\beta <\alpha$, there is some $n\in\omega$ with $\beta_n\geq \beta$. Since $D_i^{\beta_n}(B_n)\neq\emptyset$, we have
\[ D_i^\beta (B)\supset D_i^\beta(h\hat{\ }f_n\hat{\ }B_n)\supset D_i^{\beta_n}(h\hat{\ }f_n\hat{\ }B_n) \supset h\ext f_n\ext D_i^{\beta_n}(B_n) \ni h\hat{\ }f_n\ext s \sqsupset h \]
for some $s\in\omega^{<\omega}$. This implies that $h\in D_i^\beta(B)$ holds for each $\beta<\alpha$, which gives $h\in \underset{\beta<\alpha} \bigcap D_i^\beta(B) = D_i^\alpha(B)$ and $r_i(B)\geq \alpha$.
\end{proof}

From Proposition \ref{proposition: rank of B}, it follows that no \err{$\widetilde {\mc B_{<\alpha}}$} set can be covered by finitely many sets from $\bigcup_{\beta<\alpha} \mc B_\beta$. We will prove a slightly stronger result later (Lemma \ref{lemma: E, H and finite unions}).

\section{Broom Spaces and Their \texorpdfstring{$\mc F_\alpha$}{F alpha} Parametrization}\label{sec:broom spaces and parametrization}
We will study a certain collection of one-point Lindelöfications of the discrete uncountable space which give the space a special structure. These spaces were used by (among others) Talagrand, who constructed a space which is non-absolutely $\fsd$ (\cite{talagrand1985choquet}). Our goal is to compute the absolute complexity of spaces of this type.

For an arbitrary system $\mc E\subset \mc P \left( \baire \right)$ which contains only countable sets, we define the space $X_{\mc E}$ as the set $\baire\cup\{\infty\}$ with the following topology: each $\sigma\in\baire$ is an isolated point and a neighborhood subbasis of $\infty$ \err{in $X_{\mc E}$} consists of all sets of the form $\{\infty \}\cup (\baire\setminus E)$, $E\in\mc E$, \errata{and $\{\infty \}\cup (\baire\setminus \{\sigma\})$, $\sigma \in \baire$.}

The outline of Section \ref{sec:broom spaces and parametrization} is the following: We note in Example \ref{example: Y(T)} that if $X_{\mc E} \subset cX_{\mc E}$, then there is a natural formula for an $\fsd$ set $Y_1\subset cX_{\mc E}$ containing $X_{\mc E}$. Then in Definition \ref{def:Y_alpha}, we generalize the formula and define sets \err{$Y_\alpha\supset X_{\mc E}$}. These are useful for two reasons: the first is that their definition is `absolute' (it does not depend on the choice of $cX_{\mc E}$). The other is that it is easy to compute (an upper bound on) their complexity, which we do in Lemma \ref{lemma: complexity of Y}.

Our goal is to show that $X_{\mc E}=Y_\alpha$ holds for a suitable $\alpha$. In Section \ref{sec:auxiliary}, we prepare the tools for establishing a connection between the complexity of $\mc E$ and the value of `the $\alpha$ that works'. In Proposition \ref{proposition: XA are absolutely K alpha}, we use these results to prove that $X_{\mc E}=Y_\alpha$ holds in every compactification. Applying Lemma \ref{lemma: complexity of Y}, we obtain the absolute complexity of $X_{\mc E}$. We finish by giving two corollaries of this proposition.

Rather than working with a general $\mc E\subset \mc P(\baire)$, we will be interested in the following systems:

\begin{definition}[Infinite broom sets\footnote{In \cite{talagrand1985choquet} these systems are denoted as $\mc A_\alpha$, but we have to use a different letter to distinguish between infinite brooms and additive Borel classes. The letter $\mc E$ is supposed to stand for ``extension" (of $\mc B$). Note also that Talagrand's collections $\mc A_n$, $n\in\omega$, correspond to our $\mc E_{n+1}$ (for $\alpha\geq \omega$ the enumeration is the same).}]\label{def:infinite brooms}
For $\alpha\leq\omega_1$ we define
\[ \mc E_\alpha:=\left\{ E\subset \baire | \ E \textrm { is an infinite extension of some } B\in\mc B_\alpha \right\}.\footnote{Where $\mc B_\alpha$ is the $\alpha$-th collection of broom sets, introduced in Definition \ref{def:brooms}.} \]
\end{definition}
Since $\mc B_{\omega_1}\subset \mc D$, every $E\in \mc E_{\omega_1}$ is closed and discrete in $\baire$, and the following result holds:

\begin{proposition}[\cite{talagrand1985choquet}]
If each $E\in\mc E$ is closed discrete in $\baire$, then $X_{\mc E}\in\fsd\left(\beta X_{\mc E}\right)$.
\end{proposition}

In the remainder of Section \ref{sec:broom spaces and parametrization}, $X$ will stand for the space $X_{\mc E}$ for some family $\mc E \subset \mc E_{\omega_1}$ and $cX$ will be a fixed compactification of $X$. All closures will automatically be taken in $cX$. The choice of $\mc E$ is important, because it gives the space $X$ the following property:

\begin{lemma}\label{lemma: x in cX -- X and H}
For every $x\in cX \setminus X$ there exists $H\subset X$ satisfying
\begin{enumerate}[(i)]
\item $H$ can be covered by finite union of elements of $\mc E$;
\item $\left( \forall F\subset X \right) : x\in\overline{F} \implies x\in \overline{F\cap H}$.
\end{enumerate}
\end{lemma}
\begin{proof}
Since $x\neq\infty$, there is some $U\in\mc U(x)$ with $\infty \notin \overline{U}$. Because $U$ is a neighborhood of $x$, we get
\[ \left( \forall F\subset X \right) : x\in\overline{F} \implies x\in\overline{F\cap U}  .\]
Therefore, we define $H$ as $H:=U\cap (X \cap \bigcup \mc E)$.
By definition of $X_\mc E$, sets of the form $\{\infty\}\cup (\baire \setminus E)$, $E\in\mc E$, and $\{\infty\}\cup (\baire \setminus \{\sigma\})$, $\sigma \in\baire$ form a~neighborhood subbasis of $\infty$ in $X$.
In particular, $U \setminus H = U \cap (X \setminus \bigcup \mc E)$ must be finite, and thus $H$ satisfies (ii).
Moreover, $U\cap \bigcup \mc E$ can be covered by finitely elements of $\mc E$ and points from $X \setminus \bigcup \mc E$. This implies that $H$ can be covered by finitely many elements of $\mc E$, which gives $(i)$.
\end{proof}

\subsection{\texorpdfstring{$\mc F$}{F}-Borel sets \texorpdfstring{$Y_\alpha$}{Y alpha} Containing \texorpdfstring{$X_{\mc E}$}{X E}}
First, we start with motivation for our definitions:
\begin{example}[$\fsd$ and $F_{\sigma\delta\sigma\delta}$ sets containing $X_{\mc E}$]\label{example: Y(T)}
When trying to compute the complexity of $X$ in $cX$, the following canonical $K_{\sigma\delta}$ candidate comes to mind\footnote{Recall that the notions related to sequences are defined in Notation \ref{not:sequences}. In particular, $\mc N(s)$ denotes the standard Baire interval.}
\[ X \subset \{\infty\} \cup \underset {n\in\omega} \bigcap \underset {s\in \omega^n} \bigcup \overline{\mc N(s)}=:Y_1. \]
The two sets might not necessarily be identical, but the inclusion above always holds. At the cost of increasing the complexity of the right hand side to $F_{\sigma\delta\sigma\delta}$, we can define a smaller set, which will still contain $X$:
\[ X \subset \{\infty\} \cup \underset {n\in\omega} \bigcap \  \underset {s\in \omega^{n}} \bigcup \ \underset {k\in\omega} \bigcap \ \underset {t\in \omega^k} \bigcup \ \overline{\mc N(s\hat{\ }t)}=:Y_2. \]
\end{example}
In the same manner as in Example \ref{example: Y(T)}, we could define sets $Y_1\supset Y_2 \supset Y_3 \supset \dots \supset X$. Unfortunately, the notation from Example \ref{example: Y(T)} is impractical if we continue any further, and moreover, it is unclear how to extend the definition to infinite ordinals. We solve this problem by introducing `admissible mappings' and a more general definition below.

\begin{definition}[Admissible mappings]\label{def:admissible mapping}
Let $T\in\textrm{Tr}$ be a tree and $\varphi : T\rightarrow \omega^{<\omega}$ a mapping from $T$ to the space of finite sequences on $\omega$. We will say that $\varphi$ is \emph{admissible}, if it satisfies
\begin{enumerate}[(i)]
\item $(\forall s,t\in T): s \sqsubset t \implies \varphi (s) \sqsubset \varphi(t)$
\item $(\forall t\in T): \left| \varphi (t) \right| = t(0)+t(1)+\dots+t(|t|-1)$.
\end{enumerate}
For $S\subset T$ we denote by $\widetilde{\varphi}(S):=\mathrm{cl}_\mathrm{Tr}(\varphi (S))$ the tree generated by $\varphi (S)$.
\end{definition}

\begin{notation} \label{not:trees of height alpha}
For each limit ordinal $\alpha<\omega_1$ we fix a bijection $\pi_\alpha : \omega \rightarrow \alpha$. If $\alpha=\beta+1$ is a successor ordinal, we set $\pi_\alpha(n):=\beta$ for each $n\in\omega$. For $\alpha=0$, we define $T_0:=\{\emptyset\}$ to be the tree which only contains the root. For $\alpha\geq 1$, we define $T_\alpha$ (`maximal trees of height $\alpha$') as
\[ T_\alpha := \{\emptyset\} \cup \underset {n\in\omega} \bigcup n\hat{\ }T_{\pi_\alpha (n)}. \]
\end{notation}

\begin{definition} \label{def:Y_alpha}
For $\alpha<\omega_1$ we define
\[ Y_\alpha:=\left\{ x\in cX | \ \left(\exists \varphi : T_\alpha\rightarrow \omega^{<\omega} \textrm{ adm.} \right)\left( \forall t\in T_\alpha \right) : x\in \overline {\mc N \left( \varphi (t)\right)} \right\}\cup \left\{ \infty \right\} .\]
\end{definition}
For $\alpha=1$, it is not hard to see that this definition coincides with the one given in Example \ref{example: Y(T)}. This follows from the equivalence
\begin{equation*}
\begin{split}
x \in \underset {n\in\omega} \bigcap \underset {s\in \omega^n} \bigcup \overline{\mc N(s)} \iff & (\forall n\in\omega)(\exists s_n\in \omega^{<\omega},\ |s_n|=n): x\in\overline { \mc N(s_n) } \\
\overset{(\star)} \iff & (\exists \varphi : T_1 \rightarrow \omega^{<\omega} \textrm{ adm.})(\forall t\in T_1) :x\in\overline{\mc N\left(\varphi(t)\right)} \\
\overset{\textrm{def}} \iff & x\in Y_1,
\end{split}
\end{equation*}
where $(\star)$ holds because the formula $\varphi(\emptyset):=\emptyset$, $\varphi((n)):=s_n$ defines an admissible mapping on $T_1=\{\emptyset\}\cup\{ (n)| \ n\in\omega \}$ and also any admissible mapping on $T_1$ is defined in this way for some sequences  $s_n$, $n\in\omega$. The definition also coincides with the one given above for $\alpha=2$ -- this follows from the proof of Lemma \ref{lemma: complexity of Y}. We now state the main properties of sets $Y_\alpha$.

\begin{lemma}
$Y_\alpha\supset X$ holds for every $\alpha<\omega_1$.
\end{lemma}
\begin{proof}
Let $\sigma\in \baire$. For $t\in T_\alpha$ we set $\varphi(t):=\sigma|_{t(0)+\dots+t(|t|-1)}$, and observe that this is an admissible mapping which witnesses that $\sigma\in Y_\alpha$.
\end{proof}

\begin{lemma}[Complexity of $Y_\alpha$]\label{lemma: complexity of Y}
For any $\alpha=\lambda+m<\omega_1$ (where $m\in\omega$ and $\lambda$ is either zero or a limit ordinal), $Y_\alpha$ belongs to $\mc F_{\lambda + 2m+1}(cX)$.
\end{lemma}
\begin{proof}
For each $h\in\omega^{<\omega}$ and $\alpha<\omega_1$, we will show that the set $Y^h_\alpha$ defined as
\[ Y^h_\alpha:=\left\{ x\in cX| \ (\exists \varphi : T_\alpha \rightarrow \omega^{<\omega} \textrm{ adm.})(\forall t\in T_\alpha):x\in\overline{\mc N(h\hat{\ }\varphi(t)} \right\} \]
belongs to $\mc F_{\lambda + 2m+1}(cX)$.
The only admissible mapping from $T_0$ is the mapping $\varphi:\emptyset\mapsto \emptyset$, which means that \errata{$Y^h_0=\overline{\mc N(h)}\in \mc F_1$.} In particular, the claim holds for $\alpha=0$.

Let $1\leq\alpha<\omega_1$. First, we prove the following series of equivalences
\begin{equation*}
\begin{split}
x \in Y^h_\alpha \overset{(a)} \iff & (\exists \varphi : T_\alpha \rightarrow \omega^{<\omega} \textrm{ adm.})(\forall t\in T_\alpha): \\
	& x\in \overline{\mc N \left( h\hat{\ }\varphi(t)\right) } \\
\overset{(b)} \iff & (\exists \varphi : T_\alpha \rightarrow \omega^{<\omega} \textrm{ adm.})(\forall n\in\omega)(\forall t\in T_{\pi_\alpha(n)}): \\
	& x\in\overline{\mc N\left(h\hat{\ }\varphi(n\hat{\ }t)\right)}  \\
\overset{(c)} \iff & (\forall n\in\omega)(\exists s_n\in\omega^ n)(\exists \varphi_n : T_{\pi_\alpha(n)} \rightarrow \omega^{<\omega} \textrm{ adm.})(\forall t\in T_{\pi_\alpha(n)}): \\
	& x\in\overline{\mc N\left(h\hat{\ }s_n\hat{\ }\varphi_n(t)\right)} \\
\overset{(d)} \iff & x\in\bigcap_{n\in\omega}\bigcup_{s_n\in\omega^n}Y^{h\hat{\ }s_n}_{\pi_\alpha(n)},
\end{split}
\end{equation*}
from which the equation \eqref{eq: Y_alpha} follows:
\begin{equation}\label{eq: Y_alpha}
Y^h_\alpha=\bigcap_{n\in\omega}\bigcup_{s_n\in\omega^n}Y^{h\hat{\ }s_n}_{\pi_\alpha(n)}.
\end{equation}
The equivalences $(a)$ and $(d)$ are simply the definition of $Y^h_\alpha$, and $(b)$ follows from the definition of trees $T_\alpha$. The nontrivial part is $(c)$. To prove the implication `$\implies$', observe that if $\varphi$ is admissible, then
\[ (\forall n\in\omega)(\forall t\in T_{\pi_\alpha(n)})(\exists \varphi_n(t)\in \omega^{<\omega}): \varphi(n\ext t)=\varphi\left((n)\right)\ext \varphi_n(t),  \]
the mappings $\varphi_n:T_{\pi_\alpha(n)}\rightarrow \omega^{<\omega}$ defined by this formula are admissible and $|\varphi((n))|=n$. The implication from right to left follows from the fact that whenever $\varphi_n:T_{\pi_\alpha(n)}\rightarrow \omega^{<\omega}$ are admissible mappings and $s_n\in\omega^n$, the mapping $\varphi$ defined by formula $\varphi(\emptyset):=\emptyset$, $\varphi(n\ext t):=s_n\ext \varphi_n(t)$ is admissible as well.

We now finish the proof. If $\alpha$ is a successor ordinal, then $\pi_\alpha(n)=\alpha-1$ holds for all $n\in\omega$. Therefore, we can rewrite \eqref{eq: Y_alpha} as
\begin{equation*}
Y^h_\alpha = \bigcap_{n\in\omega}\bigcup_{s_n\in\omega^n}Y^{h\hat{\ }s_n}_{\pi_\alpha(n)} 
= \bigcap_{n\in\omega}\bigcup_{s_n\in\omega^n}Y^{h\hat{\ }s_n}_{\alpha-1}
\end{equation*}
and observe that each `successor' step increases the complexity by $(\cdot )_{\sigma\delta}$. Lastly, assume that $\alpha$ is a limit ordinal and $Y^s_{\alpha'}\in \underset {\beta<\alpha} \bigcup \mc F_\beta$ holds for all $s\in\omega^{<\omega}$ and $\alpha'<\alpha$. Since for limit $\alpha$, we have $\pi_\alpha(n)<\alpha$ for each $n\in\omega$, \eqref{eq: Y_alpha} gives
\begin{equation*}
Y^h_\alpha = \bigcap_{n\in\omega}\bigcup_{s_n\in\omega^n}Y^{h\hat{\ }s_n}_{\pi_\alpha(n)} 
\in \left( \underset {\beta<\alpha} \bigcup \mc F_\beta  \right)_{\sigma\delta}=\mc F_{\alpha+1},
\end{equation*}
which is what we wanted to prove.
\end{proof}

\subsection{Some Auxiliary Results} \label{sec:auxiliary}
In this section, we give a few tools which will be required to obtain our main result. Later, we will need to show that for a suitable $\alpha$, $Y_\alpha \setminus X$ is empty. In order to do this, we first explore some properties of those $x\in Y_\alpha$ which are in $X$, and of those $x\in Y_\alpha$ which do not belong to $X$. We address these two possibilities separately in Lemma \ref{lemma: Y(T) -- IF case} and Lemma \ref{lemma: Y(T) -- WF case} and show that they are related with the properties of the corresponding admissible mappings\footnote{By such a mapping we mean an admissible mapping which witnesses that $x\in Y_\alpha$ (for details, see Definitions \ref{def:admissible mapping} and \ref{def:Y_alpha}).}.
\begin{lemma}[The non-WF case]\label{lemma: Y(T) -- IF case}
For any $\sigma\in\baire$ and any increasing sequence of integers $(n_k)_{k\in\omega}$, we have $\bigcap_k\overline{\mc N(\sigma|n_k)}^{cX} = \{\sigma,\errata{\infty\}.}$ In particular, we have
\[ \bigcap_{k\in\omega} \overline{\mc N(\sigma|n_k)}^{cX}\subset X . \]
\end{lemma}
\begin{proof}
Suppose that $\sigma$ and $n_k$ are as above and the intersection $\bigcap_k \overline{\mc N(\sigma|n_k)}$ contains some \errata{$x\notin \{ \sigma, \infty \}$.} Clearly, we have $x\in cX\setminus X$. By Lemma \ref{lemma: x in cX -- X and H}, there is some discrete (in $\baire$) set $H\subset \baire$, such that for each $k\in\omega$, $x$ belongs to $\overline{\mc N(\sigma|n_k)\cap H}$. Since the sets $\mc N(\sigma|n_k)$ form a neighborhood basis of $\sigma$, there exists $n_0\in\omega$, such that $\mc N(s_{n_0})\cap H$ is a singleton. In particular, this means that $\mc N(s_{n_0})\cap H$ is closed in $cX$, and we have $x\in \overline{\mc N(s_{n_0})\cap H} = \mc N(s_{n_0})\cap H \subset H \subset X$ -- a contradiction.
\end{proof}

On the other hand, Lemma \errata{\ref{lemma: Y(T) -- WF case}} shows that if the assumptions of Lemma \ref{lemma: Y(T) -- IF case} are not satisfied, we can find a broom set\footnote{Recall that the collections $\widetilde{B}_\alpha$ were introduced in Definition \ref{def:tilde brooms}.} of class $\alpha$ which corresponds to $x$. We remark that Lemma \ref{lemma: Y(T) -- WF case} is a key step on the way to Proposition \ref{proposition: XA are absolutely K alpha} -- in particular, it is the step in which the admissible mappings are an extremely useful tool.
\begin{lemma}[Well founded case]\label{lemma: Y(T) -- WF case}
Let $x\in Y_\alpha$ and  suppose that $\varphi$ is any admissible mapping witnessing this fact. If $\widetilde{\varphi} (T_\alpha)\in\textrm{WF}$, then there exists $B \in \widetilde{ \mc B} _\alpha$ with $B\subset \widetilde{\varphi} (T_\alpha)$.
\end{lemma}
\begin{proof}
We prove the statement by induction over $\alpha$. For $\alpha=0$ we have $T_0=\{\emptyset\}$ and $\varphi(\emptyset)\in\omega^{<\omega} = \widetilde{\mc B}_0$.

Consider $\alpha > 0$ and assume that the statement holds for every $\beta < \alpha$. We have $\{(n)|\ n\in \omega\} \subset T_{\alpha}$. By the second defining property of admissible mappings, we get $\varphi((m))\neq \varphi((n))$ for distinct $m,n\in\omega$, which means that the tree $T:=\widetilde{\varphi} (\{(n)| \ n\in\omega\})$ is infinite. Since $T \subset \widetilde{\varphi} (T_{\alpha}) \in \textrm{WF}$, we use König's lemma to deduce that $T$ contains some $h\in\omega^{<\omega}$ with infinite set of successors $S:=\textrm{succ}_T(h)$.

For each $s\in S$ we choose one $n_s\in\omega$ with $\varphi(n_s)\sqsupset s$. In the previous paragraph, we have observed that $\left(s(|h|)\right)_{s\in S}$ is a forking sequence. From the first property of admissible mappings, we get
\begin{equation} \label{eq: WF case}
(\forall s\in S)(\forall t \in T_{\alpha} \textrm{ s.t. } t(0)=n_s): \varphi(t)\sqsupset s.
\end{equation}
In particular, for each $s\in S$ we have
\[ \varphi ( \{ t\in T_{\alpha}| \ t(0)=n_s \} ) = s\hat{\ }S_s \textrm{ for some }S_s\subset \omega^{<\omega}. \]
By definition\footnote{Recall that $T_\alpha$ is an $\omega$-ary tree of height $\alpha$, defined in Notation \ref{not:trees of height alpha}. } of $T_{\alpha}$, we have $\{ t\in \omega^{<\omega}| \ n_s\hat{\ }t\in T_{\alpha} \}=T_{\pi_\alpha (n_s)}$. By \eqref{eq: WF case}, each $\varphi (n_s \hat{\ } t)$ for $t\in T_{\pi_\alpha (n_s)}$ is of the form $\varphi(n_s\hat{\ }t)=\varphi(n_s)\hat{\ }\varphi_s(t)$ for some $\varphi_s(t)\in\omega^{<\omega}$, and, since $\varphi$ is admissible, the mapping $\varphi_s(t):T_{\pi_\alpha (n_s)}\rightarrow\omega^{<\omega}$ defined by this formula is admissible as well.

It follows from the induction hypothesis that $\mathrm{cl}_\mathrm{Tr}\left( \varphi_s (T_{\pi_\alpha (n)})  \right)$ contains some $\widetilde{\mc B}_{\pi_\alpha (n)}$-set $D_s$. Let $C_s \subset \varphi_s (T_{\pi_\alpha (n)})$ be some finite extension of $D_s$. By Remark \ref{rem:basic broom properties}, $C_s$ belongs to $\widetilde{\mc B}_{\pi_\alpha (n)}$. Finally, since
\[ s\ext S_s=\varphi(n_s)\hat{\ }\varphi_s(T_{\pi_\alpha(n_s)})\supset \varphi(n_s)\ext C_s \in \widetilde{\mc B}_{\pi_\alpha (n)} , \]
there also exists some $B_s\subset S_s$ with $B_s\in \widetilde{\mc B}_{\pi_\alpha (n)}$ (the set $B_s$ can have the same bristles and forking sequence as $C_s$, but $s$ might be shorter than $\varphi(n_s)$, so $B_s$ might have a longer handle than $C_s$). We finish the proof by observing that
\begin{equation}\label{eq: Phi contains Balpha}
\widetilde{\varphi} (T_\alpha) \supset \varphi \left( \underset{s\in S} \bigcup \left\{ t\in T_{\alpha}| \ t(0)=n_s \right\} \right) \err{=} \underset{s\in S} \bigcup s\hat{\ }S_s \supset \underset{s\in S} \bigcup h\hat{\ }s(|h|)\hat{\ }B_s.
\end{equation}
Because the set $S$ is infinite, the rightmost set in \eqref{eq: Phi contains Balpha} is, by definition, an element of $\widetilde{ \mc B} _\alpha$.
\end{proof}

We will eventually want to show that the assumptions of Lemma \ref{lemma: Y(T) -- WF case} can only be satisfied if the family $\mc E$ is sufficiently rich. To this end, we need to extend Lemma \ref{lemma: rank and finite covers} to a situation where a set $H\subset\baire$ cannot be covered by finitely many `infinite broom sets'\footnote{That is, the infinite extensions of broom sets, introduced in Definition \ref{def:infinite brooms}.}:
\begin{lemma}\label{lemma: E, H and finite unions}
Suppose that a set $H\subset \baire$ contains an infinite extension of some $B\subset \omega^{<\omega}$ which satisfies $r_i(B)\geq \alpha$.\footnote{Recall that $r_i$ denotes the `infinite branching rank' introduced in Definition \ref{def:derivative}.} Then $H$ cannot be covered by finitely many elements of $\bigcup_{\beta<\alpha} \mc E_\beta$.
\end{lemma}
\begin{proof}
Let $\alpha<\omega_1$, $B$ and $H$ be as above. For contradiction, assume that $H\subset \bigcup_{j=0}^k E_j $ holds for some $E_j\in \mc E_{\beta_j}$, $\beta_j<\alpha$.

By definition of $\mc E_{\beta_j}$, each $E_j$ is an infinite extension of some $B_j\in \mc B_{\beta_j}$ -- in other words, there are bijections $\psi_j:B_j \rightarrow E_j$ satisfying $\psi_j(s)\sqsupset s$ for each $s\in B_j$. Similarly, $H$ contains some infinite extension of $B$, which means there exists an injective mapping $\psi : B\rightarrow H$ satisfying $\psi(s)\sqsupset s$ for each $s\in B$. Clearly, we have
\begin{equation}\label{eq: psi B is covered by psi B_j}
\psi(B)\subset H \subset \bigcup_{j=0}^k E_j = \bigcup_{j=0}^k \psi_j(B_j) .
\end{equation}

Firstly, observe that if $t\in B$ and $s\in B_j$ satisfy $\psi_j(s)=\psi(t)$, the sequences $s$ and $t$ must be comparable. This means that for each $j\leq k$, the following formula correctly defines a mapping $\varphi_j : B_j\rightarrow \omega^{<\omega}$:
$$
\varphi_j(s):=
\begin{cases}
s,  & \textrm{ for } \psi_j(s) \notin \psi(B), \\
\psi^{-1}(\psi_j(s)), & \textrm{ for } s\sqsubset \psi^{-1}(\psi_j(s)), \\
s, & \textrm{ for } s\sqsupset \psi^{-1}(\psi_j(s)),
\end{cases}
$$
which satisfies $\varphi_j(s)\sqsupset s$ for each $s\in B_j$. Since no two elements of $B_j$ are comparable, it also follows that $\varphi_j$ is injective and thus $\varphi(B_j)$ is an extension of $B_j$. By Remark \ref{rem:basic broom properties}, this implies that $\varphi_j(B_j)\in\mc B_{\beta_j}$.

By \eqref{eq: psi B is covered by psi B_j}, for $t\in B$ there exists some $j\leq k$ and $s\in B_j$, such that $\psi_j(s)=\psi(t)$. It follows from the definition of $\psi_j$ that
\[ t=\psi^{-1}(\psi(t))=\psi^{-1}(\psi_j(s))\sqsubset \varphi_j(s) .\] In other words, we get
\[ B\subset \mathrm{cl}_\mathrm{Tr}\left( \bigcup_{j=0}^k \varphi_j(B_j) \right) .\]
We observe that this leads to a contradiction:\errata{}
\begin{align*}
\alpha & \err{\leq} r_i(B) \leq r_i \left( \mathrm{cl}_\mathrm{Tr}\left( \bigcup_{j=0}^k \varphi_j(B_j) \right) \right) = r_i\left( \bigcup_{j=0}^k \varphi_j(B_j) \right) = \\
	& = \max_j \ r_i \left( \varphi_j(B_j) \right) \leq \max_j \beta_j < \alpha.
\end{align*}
\end{proof}

\subsection{Absolute Complexity of \texorpdfstring{$X_{\mc E}$}{X E}} \label{sec:absolute complexity of brooms}
In this section, we prove an upper bound on the absolute complexity of spaces $X_{\mc E}$ for any collection of infinite broom sets $\mc E$. Talagrand's earlier result will then imply that, in some cases, this bound is sharp.
\begin{proposition}[Complexity of $X_{\mc E}$]\label{proposition: XA are absolutely K alpha}
For any integer $m\in\omega$ and limit ordinal $\lambda < \omega_1$, we have
\begin{enumerate}[1)]
\item $\mc E\subset \mc E_m$ $\implies$ $X_{\mc E}$ is absolutely $\mc F_{2m+1}$;
\item $\mc E\subset \underset {\beta<\lambda} \bigcup \mc E_\beta$ $\implies$ $X_{\mc E}$ is absolutely $\mc F_{\lambda+1}$;
\item $\mc E\subset \mc E_{\lambda+m}$ $\implies$ $X_{\mc E}$ is absolutely $\mc F_{\lambda+2m+3}$.
\end{enumerate}
\end{proposition}
\begin{proof}
Suppose that $\lambda$ and $m$ are as above, $\mc E$ is a family of discrete subsets of $\baire$ and $cX$ is a compactification of the space $X:=X_{\mc E}$. Denote $\alpha_1 := m$, $\alpha_2:=\lambda$ and $\alpha_3:=\lambda+m+1$. Note that we are only going to use the following property of $\mc E$ and $\alpha_i$:
\begin{equation} \label{eq: main proposition cases}
\begin{split}
 & \textrm{$\mc E$ satisfies the hypothesis of } i) \textrm{ for } i=1 \ \ \ \ \ \implies \mc E\subset \bigcup_{\beta<{\alpha_i+1}} \mc E_\beta, \\
 & \textrm{$\mc E$ satisfies the hypothesis of } i) \textrm{ for } i\in\{2,3\} \implies \mc E\subset \bigcup_{\beta<{\alpha_i}} \mc E_\beta.
\end{split}
\end{equation}
We will show that $X=Y_{\alpha_i}$ holds in each of these cases\footnote{Where $Y_\alpha$ is given by Definition \ref{def:Y_alpha}.}. Once we have this identity, the conclusion immediately follows from Lemma \ref{lemma: complexity of Y}.

Let $i\in\{1,2,3\}$. Suppose that there exists $x\in Y_{\alpha_i}\setminus X$ and let $\varphi : T_{\alpha_i} \rightarrow \omega^{<\omega}$ be an admissible mapping\footnote{\label{footnote: tilde admissible mapping}Recall that admissible mappings were introduced in Definition \ref{def:admissible mapping}. The `tilde version' is defined as $\widetilde{\varphi}(S):=\mathrm{cl}_\mathrm{Tr}(\varphi (S))$. } witnessing that $x\in Y_{\alpha_i}$. Moreover, let $H\subset \baire$ be a set satisfying the conclusion of Lemma \ref{lemma: x in cX -- X and H}.

Recall that by definition of $Y_{\alpha_i}$, we have
\begin{equation} \label{eq: def of Y alpha}
\left(\forall t\in T_{\alpha_i}\right): x\in\ \overline{ \mc N(\varphi(t))}.
\end{equation}
By Lemma \ref{lemma: Y(T) -- IF case}, $\widetilde{\varphi} (T_{\alpha_i})$ contains no infinite branches. Therefore, we can apply Lemma \ref{lemma: Y(T) -- WF case} to obtain a $\widetilde{\mc B}_{\alpha_i}$-set\footnote{$\widetilde {\mc B}_\alpha$ is the `disjoint version' of broom collection $\mc B_\alpha$ (see Definitions \ref{def:brooms} and \ref{def:tilde brooms}).} $B\subset \widetilde{\varphi}( T_{\alpha_i})$. Since, by \eqref{eq: def of Y alpha}, $x$ belongs to $\overline{ \mc N(s)}$ for each $s\in B$, we conclude that $x\in\overline{H\cap \mc N(s)}$ holds for every $s\in B$ (Lemma \ref{lemma: x in cX -- X and H}, $(ii)$). In particular, all the intersections $H\cap \mc N(s)$ must be infinite. Note that for the next part of the proof, we will only need the intersections to be non-empty. The fact that they are infinite will be used in the last part of the proof.

We can now conclude the proof of $2)$ and $3)$. Let $i\in\{2,3\}$ and assume that the hypothesis of $i)$ holds. Since $H\cap \mc N(s)$ is non-empty for each $s\in B$, it follows that $H$ contains an infinite extension of $B$. Since $r_i(B)=\alpha_i$ holds by Lemma \ref{proposition: rank of B}, Lemma \ref{lemma: E, H and finite unions} yields that $H$ cannot be covered by finitely many elements of $\bigcup_{\beta<{\alpha_i}} \mc E_\beta$. By \eqref{eq: main proposition cases}, we have $\mc E\subset \bigcup_{\beta<{\alpha_i}} \mc E_\beta$. This contradicts the first part of Lemma \ref{lemma: x in cX -- X and H} (which claims that $H$ \emph{can} be covered by finitely many elements of $\mc E$).

For the conclusion of the proof of $1)$, let $i=1$. \err{By \eqref{eq: main proposition cases}, we have $\mc E\subset  \bigcup_{\beta<\alpha_i+1} \mc E_\beta$.} Compare this situation with the setting from the previous paragraph -- it is clear that if we show that $H$ in fact contains an infinite extension of some $\widetilde B$ with $r_i(\widetilde B)\geq \alpha_i+1$, we can replace $B$ by $\widetilde B$. We can then apply the same proof which worked for $2)$ and $3)$.

To find $\widetilde B$, enumerate $B$ as $B=\{ s_n | \ n\in\omega \} $. Since each $H\cap \mc N(s_n)$ is infinite, there exist distinct sequences $\sigma^k_n\in H$, $k\in\omega$, satisfying $s_n\sqsubset\sigma^k_n$. We use these sequences to obtain $s_n^k\in\omega^{<\omega}$ for $n,k\in\omega$, which satisfy $s_n\sqsubset s^k_n\sqsubset \sigma^k_n$ and $k\neq l$ $\implies$ $s^k_n\neq s^l_n$. We denote $\widetilde{B}:=\{ s^k_n | \ n,k\in\omega \}$. Clearly, this set satisfies $D_i(\widetilde B)\supset B$ and $H$ contains some infinite extension of $\widetilde B$. Because $r_i(B)=\alpha_i$ is finite, we have $r_i(\widetilde{B}) \geq 1+\alpha_i \overset{\alpha_i<\omega}{=} \alpha_i +1 > \alpha_i = r_i(B)$, which completes the proof.
\end{proof}

Applying Proposition \ref{proposition: XA are absolutely K alpha} to the construction from \cite{talagrand1985choquet} immediately yields the following corollary:
\begin{theorem}\label{theorem: main theorem}
For every even $4\leq \alpha < \omega_1$, there exists a space \errata{$\mathbf T_\alpha$} with properties
\begin{enumerate}[(i)]
\item $\mathbf T_\alpha$ is an $F_{\sigma\delta}$ space;
\item $\mathbf T_\alpha$ is not an absolute $\mc F_\alpha$ space;
\item $\mathbf T_\alpha$ is an absolute $\mc F_{\alpha+1}$ space.
\end{enumerate}
\end{theorem}
\begin{proof}
The existence of spaces $\mathbf T_\alpha$ which satisfy $(i)$ and $(ii)$ follows from \cite{talagrand1985choquet}. In this paper, Talagrand constructed a family $\mc A_T\subset \mc E_{\omega_1}$, which is both `rich enough' and almost-disjoint (that is, the intersection of any two families is finite). This is accomplished by a smart choice of bristles and forking sequences on the `highest level of each broom'.

He then shows that for a suitable $\widetilde \alpha$, any space $X_{\mc E}$ with $\mc E\supset \mc A_T \cap \mc E_{\widetilde \alpha}$, where $\mc E$ is almost-disjoint, satisfies $(i)$ and $(ii)$. In particular, this holds for $\mathbf T_\alpha := X_{\mc A_T \cap \mc E_{\widetilde \alpha}}$. Note however that Talagrand was interested in the `maximal' version of the construction, which is the space $T:= X_{\mc A_T}$. We, on the other hand, will use the `intermediate steps' of his construction.

For us, the details of the construction are not relevant -- the only properties we need are $(i)$, $(ii)$ and the correspondence between $\alpha$ and $\widetilde \alpha$.\footnote{Note that in \cite{talagrand1985choquet}, the author uses slightly different (but equivalent) definition of broom families $\mc E_\alpha$, which shifts their numbering for finite $\alpha$ by $1$.} This correspondence is such that the broom class only increases with countable intersections (that is, odd steps of $\mc F_\alpha$), which translates to our notation as follows:
\begin{itemize}
\item $\alpha=2n\in\omega$, $\alpha\geq 4$ $\implies \mathbf T_{\alpha}=X_{\mc E}$ for some ${\mc E}\subset \mc E_n$;
\item $\alpha<\omega_1$ is a limit ordinal $\implies \mathbf T_\alpha=X_{\mc E}$ for some ${\mc E}\subset \underset{\beta<\alpha} \bigcup \mc E_\beta$;
\item $\alpha=\lambda+2n+2$ for $n\in\omega$ and limit $\lambda$ $\implies \mathbf T_\alpha=X_{\mc E}$ for some ${\mc E}\subset \mc E_{\lambda+n}$.
\end{itemize}
In all the cases, it follows from Proposition \ref{proposition: XA are absolutely K alpha}  that $X_{\mc E}$ is absolutely $\mc F_{\alpha+1}$.
\end{proof}

If $\beta<\omega_1$ is the least ordinal for which $X$ is an (absolute) $\mc F_\beta$ space, we say that the \emph{(absolute) complexity} of $X$ is $\mc F_\beta$. Using this notation, Theorem \ref{theorem: main theorem} can be rephrased as ``for every odd $5\leq \alpha<\omega_1$, there exists an $\fsd$ space whose absolute complexity is $\mc F_\alpha$''. By modifying the spaces from Theorem \ref{theorem: main theorem}, we obtain the following result:

\begin{corollary}\label{corollary: non K alpha space which is abs K alpha+1}
For every two countable ordinals $3 \leq \alpha \leq \beta$, $\beta$ odd, there exists a space $X_\alpha^\beta$, such that
\begin{enumerate}[(i)]
\item the complexity of $X_\alpha^\beta$ is $\mc F_\alpha$;
\item the absolute complexity of $X_\alpha^\beta$ is $\mc F_\beta$.
\end{enumerate}
\end{corollary}
\begin{proof}
If $\alpha=\beta$, this is an immediate consequence of absoluteness of Borel classes. To make this claim more precise, let $P$ be some uncountable Polish space. By Remark \ref{rem:F,G and M,A} $(ii)$, $\mc F_\alpha(P)$ corresponds to some Borel class $\mc C(P)$. We define $Z_\alpha$ as one of the subspaces of $P$, which are of the class $\mc C(P)$, but not of the `dual' class -- that is, if $\mc C(P)=\Sigma^0_\gamma(P)$ for some $\gamma<\omega_1$, then $Z_\alpha\in \Sigma^0_\gamma(P)\setminus \Pi^0_\gamma (P)$ (and vice versa for $\mc C(P)=\Pi^0_\gamma(P)$). For the existence of such a set, see for example Corollary 3.6.8 in \cite{srivastava2008course}.

Since $Z_\alpha$ is separable and metrizable, Theorem \ref{theorem: separable metrizable} guarantees that it is an absolute $\mc F_\alpha$ space. However, it is not of the class $\mc F_{\alpha'}$ for any $\alpha'<\alpha$, because that would by Remark \ref{rem:F,G and M,A} imply that $Z_\alpha$ belongs to $\mc F_{\alpha'} (P)\subset \Pi^0_\gamma(P)$. This shows that if $\alpha=\beta$, the space $X_\alpha^\beta:=Z_\alpha$ has the desired properties.

If $\alpha<\beta$, we define \errata{$X_\alpha^\beta$} as the topological sum of $Z_\alpha$ and the space $\mathbf T_{\beta-1}$ from Theorem \ref{theorem: main theorem}. Since $\mathbf T_{\beta-1}\in \fsd(\beta \mathbf T_{\beta-1})=\mc F_3(\beta \mathbf T_{\beta-1})$ and $\alpha \geq 3$, we get that $X_\alpha^\beta$ is of the class $\mc F_\alpha$ in the topological sum  $\beta X_\alpha^\beta = \beta Z_\alpha \oplus \beta \mathbf T_{\beta-1}$. By previous paragraph, it is of no lower class in $\beta X_\alpha^\beta$. By Theorem \ref{theorem: main theorem} (and the previous paragraph), $X_\alpha^\beta$ also has the correct absolute complexity.
\end{proof}

To get a complete picture of possible combinations of complexity and absolute complexity for $\mc F$-Borel spaces, it remains to answer the following questions:
\begin{problem}\label{problem: X alpha beta}
\begin{enumerate}[(i)]
\item Let $4\leq\beta<\omega_1$ be an \emph{even} ordinal number. Does there exist a space $X^\beta$, whose absolute complexity is $\mc F_\beta$ and complexity is $\mc F_\alpha$ for some $\alpha<\beta$? What is the lowest possible value of $\alpha$?
\item If a space $X$ is $\mc F$-Borel in every compactification, is it necessarily absolutely $\mc F_\beta$ for some $\beta<\omega_1$?
\end{enumerate}
\end{problem}
One could expect that Corollary \ref{corollary: non K alpha space which is abs K alpha+1} might also hold for even $\beta$, so that we should be able to find an $\fsd$ space answering the first part of Problem \ref{problem: X alpha beta} in positive.

\section*{Acknowledgment}
I would like to thank my supervisor, Ondřej Kalenda, for numerous very helpful
suggestions and fruitful consultations regarding this paper. This work was supported by the research grant GAUK No. 915.

\putbib[refs]
\end{bibunit}
\begin{bibunit}[alphanum]

\chapter{Complexities and Representations of \texorpdfstring{$\mc F$}{F}-Borel Spaces}
 \label{chapter:representations}

\begin{center}
\textit{(Submitted, available at arXiv:1804.08367.)}
\end{center}

\textbf{Abstract:}
We investigate the~$\mc F$-Borel complexity of topological spaces in their different compactifications.
We provide a~simple proof of the~fact that a~space can have \errata{anywhere between 1 and $\omega_1$-many} different complexities in different compactifications.
We also develop a~theory of representations of $\mc F$-Borel sets,
and show how to apply this theory to prove that the~complexity of hereditarily Lindelöf spaces is absolute (that is, it is the~same in every compactification).
We use these representations to characterize the~complexities attainable by a~specific class of topological spaces.
This provides an alternative proof of the~first result,
and implies the~existence of a~space with non-absolute additive complexity.
We discuss the~method used by Talagrand to construct the~first example of a~space with non-absolute complexity, hopefully providing an explanation which is more accessible than the~original one.
We also discuss the~relation of complexity and local complexity, and show how to construct amalgamation-like compactifications.

\section{Introduction}\label{section: intro}


We investigate complexity of $\mc F$-Borel sets, that is, of the sets from the smallest system containing closed sets and stable under taking countable unions and intersections.
This is of particular interest because the~$\mc F\textnormal{-Borel}$ classes are not absolute (unlike the Borel classes - see \cite{talagrand1985choquet} and, for example, \cite{holicky2003perfect}).
In particular, we investigate which values of $\mc F$-Borel complexity may a given space assume.
We further study various representations of $\mc F$-Borel sets.

Let us start by basic definitions and exact formulation of problems.
All topological spaces in this work will be Tychonoff (unless the contrary is explicitly specified -- see Section~\ref{section: amalgamation spaces}). 
For a~family of sets $\mc C$, we will denote by $\mc C_\sigma$ the~collection of all countable unions of elements of $\mc C$ and by $\mc C_\delta$ the~collection of all countable intersections of elements of $\mc C$.
If $\alpha$ is an ordinal and $\mc C_\beta$, $\beta \in [0,\alpha)$, are collections of sets, we denote $\mc C_{<\alpha} := \bigcup_{\beta < \alpha} \mc C_\beta$.

In any topological space $Y$, we have the~family $\mc F$ of closed sets, $\mc F_\sigma$ sets, $\mc F_{\sigma\delta}$ sets and so on. Since this notation quickly gets impractical \errata{(and isn't uniquely defined on limit steps)}, we use the~following definition:

\begin{definition}[$\mc F$-Borel\footnote{The $\mc F$-Borel sets are also sometimes called $\mc F$-Borelian, to denote the~fact that we only take the~system generated by countable unions and intersections of closed sets, rather than the~whole \emph{$\sigma$-algebra} generated by closed sets. Since there is little space for confusion in this paper, we stick to the~terminology of Definition~\ref{definition: F Borel sets}.} hierarchy]\label{definition: F Borel sets}
We define the~hierarchy of $\mc F$-Borel sets in a~topological space $Y$ as
\begin{itemize}
\item $\mc F_0(Y) := \mc F(Y) := $ closed subsets of $Y$,
\item $\mc F_\alpha(Y) := \left( \mc F_{<\alpha} (Y) \right)_\sigma$ for $0 < \alpha <\omega_1$ odd,
\item $\mc F_\alpha(Y) := \left( \mc F_{<\alpha} (Y) \right)_\delta$ for $0 < \alpha <\omega_1$ even.\footnote{Recall that every ordinal $\alpha$ can be uniquely written as $\alpha = \lambda + m$, where $\lambda$ is a~limit ordinal or 0 and $m\in\omega$. An ordinal is said to be even, resp. odd, when the~corresponding $m$ is even, resp. odd.}
\end{itemize}
\end{definition}

\noindent The classes stable under countable unions (resp. intersections) -- that is, $\Fa$ for $\alpha$ odd (resp. even) -- are called additive (resp. multiplicative).
Next, we introduce the~class of $\mc K$-analytic spaces, and a~related class of Suslin-$\mc F$ sets, which contains all $\mc F$-Borel sets (\cite[Part\,I, Cor.\,2.3.3]{rogers1980analytic}):

\begin{definition}[Suslin-$\mc F$ sets and $\mc K$-analytic spaces]
Let $Y$ be a~topological space and $S\subset Y$. $S$ is said to be \emph{Suslin-$\mc F$} in $Y$ if it is of the~form
\[ S = \bigcup_{n_0,n_1,\dots} \bigcap_{k\in\omega} F_{n_0,\dots,n_k} \]
for closed sets $F_{n_0,\dots,n_k} \subset Y$, $n_i\in\omega$, $i,k\in\omega$. By $\mc F_{\omega_1} (Y)$ we denote the~collection of Suslin-$\mc F$ subsets of $Y$.

A topological space is \emph{$\mc K$-analytic} if it is a~Suslin-$\mc F$ subset of some compact space.
\end{definition}

\noindent (Typically, $\mc K$-analytic spaces are defined differently, but our definition equivalent by \cite[Part\,I,\,Thm.\,2.5.2]{rogers1980analytic}.)

Note that if $Y$ is compact, then $\mc F_0$-subsets of $Y$ are compact, $\mc F_1$-subsets are $\sigma$-compact, $\mc F_2$-subsets are $\mc K_{\sigma\delta}$, $\Fa$-subsets for $\alpha<\omega_1$ are ``$\mc K_\alpha$'', and Suslin-$\mc F$ subsets of $Y$ are $\mc K$-analytic.

By complexity of a~topological space $X$ in some $Y$ we mean the~following:

\begin{definition}[Complexity of $X$ in $Y$]
Suppose that $X$ is Suslin-$\mc F$ subset of a~topological space $Y$. By \emph{complexity} of $X$ in $Y$ we mean
\begin{align*}
\textnormal{Compl}\,(X,Y):=
	& \text{ the~smallest } \alpha\leq\omega_1 \text{ for which } X\in \Fa(Y), \text{ that is} \\
	& \text{ the~unique } \alpha\leq\omega_1 \text{ for which } X\in \Fa(Y)\setminus \mc F_{<\alpha}(Y) .
\end{align*}
\end{definition}

\noindent Recall that a~$\mc K$-analytic space is Suslin-$\mc F$ in every space which contains it (\cite[Theorem 3.1]{hansell1992descriptive}).
Consequently, $\mc K$-analytic spaces are precisely those $X$ for which $\textnormal{Compl}\,(X,Y)$ is defined for every $Y\supset X$, and this is further equivalent to $\textnormal{Compl}\,(X,Y)$ being defined for \emph{some} space $Y$ which is compact.

Central to this work is the~following notion of complexities attainable in different spaces $Y$ containing a~given space. To avoid pathological cases such as $\Compl{X}{X}=0$ being true for arbitrarily ``ugly'' $X$, we have to restrict our attention to ``sufficiently nice'' spaces $Y$.

\begin{definition}[Attainable complexities]\label{definition:att_compl}
The \emph{set of complexities attained by a~space $X$} is defined as
\begin{align*}
\textnormal{Compl}(X) := \left\{ \alpha \leq \omega_1 \, | \ \alpha = \textnormal{Compl}\,(X,cX) \text{ for some compactification $cX$ } \right\} \\
= \left\{ \alpha \leq \omega_1 \, | \ \alpha = \textnormal{Compl}\,(X,Y) \text{ for some $Y\supset X$ s.t. } \overline{X}^Y \text{ is compact} \right\} .
\end{align*}
\end{definition}

(The~second identity in Definition \ref{definition:att_compl} holds because $X\subset K \subset Y$ implies $\Compl{X}{K} \leq \Compl{X}{Y}$.)
If $X$ satisfies $\textnormal{Compl}(X) \subset [0,\alpha]$ for some $\alpha \leq \omega_1$, we say that $X$ is an \emph{absolute $\Fa$ space}. The smallest such ordinal $\alpha$ is called the \emph{absolute complexity of $X$} -- we clearly have $\alpha = \sup \textnormal{Compl}(X)$.
Sometimes, $X$ for which $\textnormal{Compl}(X)$ contains $\alpha$ is said to be an \emph{$\Fa$ space}.
Finally, if $\textnormal{Compl}(X)$ is empty or a~singleton, the~complexity of $X$ said to be \emph{absolute}. Otherwise, the~complexity of $X$ is \emph{non-absolute}.

The goal of this paper is to investigate the~following two problems:

\begin{problem}[$X$ of a~given complexity]\label{problem: X with Compl = C}
Let $C\subset [0,\omega_1]$. Is there a~topological space $X$ with $\textnormal{Compl}(X)=C$?
\end{problem}

\begin{problem}[Complexity of an arbitrary $X$]\label{problem: Compl of X}
For any topological space $X$, describe $\textnormal{Compl}(X)$.
\end{problem}

Inspired by behavior of complexity in Polish spaces, one might suspect that any $\Fa$ space is automatically absolutely $\Fa$. However, this quite reasonably sounding conjecture does \emph{not} hold -- there exists an $\fsd$ space $\mathbf{T}$ which is not absolutely $\Fa$ for any $\alpha<\omega_1$ (\cite{talagrand1985choquet}; $\mathbf{T}$ is also defined here, in Definition~\ref{definition: AD topology}).
The following proposition summarizes the~results related to Problem~\ref{problem: X with Compl = C} known so far.

\begin{proposition}[Attainable complexities: state of the~art]\label{proposition: basic attainable complexities}\mbox{}
\begin{enumerate}[(1)]
\item A topological space $X$ is $\mc K$-analytic if and only if $\textnormal{Compl}(X)$ is non-empty.\label{case:K-analytic}
\item For a~topological space $X$, the~following propositions are equivalent:
	\begin{enumerate}[(i)]
		\item $X$ is compact (resp. $\sigma$-compact);
		\item $\textnormal{Compl}(X)=\{0\}$ (resp. $\{1\}$);
		\item $\textnormal{Compl}(X)$ contains $0$ (resp. $1$).
	\end{enumerate} \label{case: compact and sigma compact}
\item The complexity of any separable metrizable space is absolute. In particular, for every $\alpha\leq \omega_1$, there exists a~space satisfying $\textnormal{Compl}(X)=\{\alpha\}$. \label{case: polish spaces and absolute complexity}
\item For every $2 \leq \alpha \le \beta \le \omega_1$, $\beta$ even, there exists a~topological space $X=X_\alpha^\beta$ satisfying \label{case: X alpha}
\[ \left\{ \alpha, \beta \right\} \subset \textnormal{Compl}\left(X \right) \subset [\alpha, \beta] .\]
\end{enumerate}
\end{proposition}

\begin{proof}[Proof of Proposition \ref{proposition: basic attainable complexities}]
\eqref{case: compact and sigma compact} is trivial, since continuous images of compact sets are compact.
Regarding \eqref{case:K-analytic}, we have already mentioned that any Suslin-$\mc F$ subset of a compact space is $\mc K$-analytic, and any $\mc K$-analytic space is Suslin-$\mc F$ in every space which contains it. Since $(\textnormal{Suslin-}\mc F)(Y) = \mc F_{\omega_1}(Y) \supset \bigcup_{\alpha < \omega_1} \Fa(Y)$ holds for any $Y$, \eqref{case:K-analytic} follows.

The first part of \eqref{case: polish spaces and absolute complexity} is by no means obvious -- we are dealing with the~complexity of $X$ in \emph{all} compactifications, not only those which are metrizable. Nonetheless, it's proof is fairly elementary; see for example \cite[Theorem 2.3]{kovarik2018brooms} (and Remark~\ref{remark: F and F tilde}) or Proposition \ref{proposition: hereditarily lindelof spaces are absolute}.
The ``in particular'' part follows from the~fact that the~Borel hierarchy in Polish spaces is strictly increasing (see, for example, the~existence of universal sets in \cite{srivastava2008course}).

\eqref{case: X alpha}:
By \eqref{case: polish spaces and absolute complexity} of this proposition, there is some space $X_\alpha^\alpha$ satisfying
\[ \textnormal{Compl}(X^\alpha_\alpha) = \{ \alpha \} .\]
By \cite{talagrand1985choquet}, there exists a~space $\mathbf{T}$ which is $\fsd$ in $\beta \mathbf{T}$, but we have $\mathbf{T}\notin \mc F_{<\omega_1} (c\mathbf{T})$ for some compactification $c\mathbf{T}$ of $\mathbf{T}$. Since such a~space is $\mc K$-analytic, it satisfies 
\[ \left\{ 2, \omega_1 \right\} \subset \textnormal{Compl}\left(\mathbf{T} \right) \subset [2, \omega_1] .\]
For $\beta = \omega_1$, the~topological sum $X := X^\alpha_\alpha \oplus \mathbf{T}$ clearly has the~desired properties.

For even $\beta<\omega_1$, take the~topological sum $X := X^\alpha_\alpha \oplus X^\beta_2 $, where $X^\beta_2$ is some space satisfying
\begin{equation}\label{equation: X 2 beta}
\left\{ 2, \beta \right\} \subset \textnormal{Compl}\left(\err{X_2^\beta} \right) \subset [2, \beta ] .
\end{equation}
\errata{As explained in Remark \ref{remark: F and F tilde}, the existence of such $X^\beta_2$ follows from \cite[Theorem~5.14]{kovarik2018brooms}. For a more general independent proof, see Corollary \ref{corollary: T alpha for odd}.}
\end{proof}

\errata{The following remark explains the differences in notation between this paper and \cite{kovarik2018brooms}, as well as how the results of \cite{kovarik2018brooms} should be interpreted in the language of the present paper.
\begin{remark}[Relation to \cite{kovarik2018brooms}]\label{remark: F and F tilde}
Formally, the results from \cite{kovarik2018brooms} ``do not say much'' about the classes $\Fa$. However, the proofs presented therein can be used to obtain \eqref{case: polish spaces and absolute complexity} and \eqref{case: X alpha} from Proposition \ref{proposition: basic attainable complexities}.
\end{remark}

\noindent The class of $\mc F$-Borel sets is the smallest class which contains closed sets and is closed under taking countable unions and intersections. When constructing the hierarchy corresponding to this class, two choices need to be made. First, should the classes corresponding to limit ordinals be additive or multiplicative? Second, should closed sets be of the 1st or 0th class? Moreover, the notation would be more elegant if either all even classes are multiplicative, or all even classes are additive. Since limit ordinals are even, the class of closed sets is multiplicative and the second choice is merely a cosmetic one, this leaves us with two natural definitions: A) limit classes are additive and closed sets are of the 1st class and B) limit classes are multiplicative and closed sets are of the 0th class.

The classes $\Fa$ from Definition \ref{definition: F Borel sets} correspond to B), while the classes used in \cite{kovarik2018brooms} correspond to A) --- for the purposes of this discussion, we denote these as $\widetilde {\mc F}_\alpha$. First reason for this change was that the newer notation allowed for more elegant formulations (compare for example Proposition \ref{chapter:brooms}.\ref{proposition: XA are absolutely K alpha} with Proposition \ref{chapter:representations}.\ref{proposition: basic broom complexities}\,(ii)). More importantly, the more general construction presented in the current paper naturally worked well with the classes $\Fa$, but not with $\widetilde{\mc F}_\alpha$ --- in particular, it was unclear to the author how to construct a space whose absolute complexity is exactly $\widetilde{\mc F}_\lambda$ for limit $\lambda$. Unfortunately, this only became apparent to the author after the first paper has already been published.

Since the finite classes satisfy $\mc F_n = \widetilde{\mc F}_{n+1}$ we easily see that any $\alpha<\omega_1$ satisfies $\Fa \subset \widetilde{\mc F}_{\alpha+1}$, $\widetilde{\mc F}_\alpha \subset \mc F_{\alpha+1}$ and for limit $\lambda$ we have $\mc F_{<\lambda} = \widetilde {\mc F}_{<\lambda}$. In this sense, the results from \cite{kovarik2018brooms} can be translated into the language of the present paper, but the ``+1 shift'' makes the translated results weaker than Proposition \ref{proposition: her Lind spaces absolute}.
Nonetheless, we claim that \eqref{case: polish spaces and absolute complexity} and \eqref{case: X alpha} from Proposition \ref{proposition: basic attainable complexities} ``can be derived from the proofs in \cite{kovarik2018brooms} in a straightforward manner''. Rather than trying to pass this claim for a proof of Proposition \ref{proposition: basic attainable complexities}, we mention this for the purposes of judging the novelty of the present paper.
Indeed, while the required modifications are rather elementary in the case of \eqref{case: polish spaces and absolute complexity}, a non-trivial effort would be needed to derive \eqref{case: X alpha} from \cite{kovarik2018brooms}. However, these propositions are special cases of Proposition \ref{proposition: her Lind spaces absolute} and Theorem \ref{theorem:X_2_beta}. A reader interested in the proofs is therefore invited to read these latter results instead.}

\section{Main Results and Open Problems}\label{section:overview}

The first contribution of this paper is showing that \eqref{case: X alpha} from Proposition~\ref{proposition: basic attainable complexities} holds not only for even ordinals, but for general $\beta \in [2,\omega_1]$:

\begin{theorem}[Non-absolute space of additive complexity]
 \label{theorem:X_2_beta}
For every $2 \leq \alpha \le \beta \le \omega_1$, there exists a~topological space $X=X_\alpha^\beta$ satisfying
\begin{equation}\label{eq:X_a^b}
\left\{ \alpha, \beta \right\} \subset \textnormal{Compl}\left(X \right) \subset [\alpha, \beta] .
\end{equation}
\end{theorem}

\begin{proof}
As in the~proof of Proposition~\ref{proposition: basic attainable complexities}\,\eqref{case: X alpha}, we can assume that $\alpha=2$.
\errata{The existence of $X_2^\beta$ follows from Corollary~\ref{corollary: T alpha for odd}.}
\end{proof}

The second contribution to Problem~\ref{problem: X with Compl = C} made by this paper is the~following generalization of Theorem~\ref{theorem:X_2_beta}, which resolves the~uncertainty about the~set $(\alpha,\beta) \cap \textnormal{Compl}(X_\alpha^\beta)$ by showing that there is such a~space $X_{[\alpha,\beta]}$ for which the set of attainable complexities is the~whole interval $[\alpha,\beta]$.

\begin{theorem}[Attainable complexities]\label{theorem: summary}
For every closed interval $I\subset [2,\omega_1]$, there is a~space $X$ with
\[ \textnormal{Compl}\left(X\right) = I .\]
\end{theorem}

\errata{\noindent We present one proof of Theorem~\ref{theorem: summary} in Section~\ref{section: complexity of brooms}, where it is shown that appropriately chosen Talagrand's broom spaces satisfy the conclusion (Theorem~\ref{theorem: complexity of talagrands brooms}).
Alternatively, Section~\ref{section: topological sums} shows how to obtain Theorem~\ref{theorem: summary} in a straightforward manner by combining the spaces $X_\alpha^\beta$ from Theorem~\ref{theorem:X_2_beta}.
While both of these proofs currently rely on Section~\ref{section: complexity of brooms}, another source of spaces satisfying \eqref{eq:X_a^b} might be found in the future, and then the two approaches would be truly independent.}

Regarding Problem~\ref{problem: Compl of X}, we (weakly) conjecture that the~set of attainable complexities always has the~following properties:

\begin{conjecture}[$\textnormal{Compl}(\cdot)$ is always a~closed interval]\label{conjecture: complexities form closed interval}
For $\mc K$-analytic space $X$,
\begin{enumerate}[(i)]
	\item the~set $\textnormal{Compl}(X)$ is always an interval,
	\item the~set $\textnormal{Compl}(X)$ is always closed in $[0,\omega_1]$.
\end{enumerate}
\end{conjecture}

\noindent If Conjecture~\ref{conjecture: complexities form closed interval} holds, Theorem~\ref{theorem: summary} would actually be a~\emph{complete} solution of Problem~\ref{problem: X with Compl = C}:

\begin{conjecture}[Solution of Problem~\ref{problem: X with Compl = C}]
For any space $X$, exactly one of the~following options holds:
\begin{enumerate}[(1)]
\item $X$ is not $\mc K$-analytic and $\textnormal{Compl}(X)=\emptyset$.
\item $X$ is compact  and $\textnormal{Compl}(X) = \{0\}$.
\item $X$ is $\sigma$-compact and $\textnormal{Compl}(X) = \{1\}$.
\item $\textnormal{Compl}(X) = [\alpha,\beta]$ holds for some $2\leq \alpha \leq \beta \leq \omega_1$.
\end{enumerate}
Moreover, any of the~possibilities above is true for some $X$.
\end{conjecture}

Regarding Problem~\ref{problem: Compl of X}, \errata{in Proposition~\ref{proposition: hereditarily lindelof spaces are absolute}} we prove the~following generalization of \eqref{case: polish spaces and absolute complexity} from Proposition~\ref{proposition: basic attainable complexities}:

\begin{proposition*}\label{proposition: her Lind spaces absolute}
The complexity of any hereditarily Lindelöf space is absolute.
\end{proposition*}

\noindent (This result was previously unpublished, but the~core observation behind it is due to J. Spurný and P. Holický -- our contribution is enabling a~simple formal proof.)

We also mention an open problem related to Problem~\ref{problem: Compl of X} which is not further discussed in this paper:
So far, the~only known examples of spaces with non-absolute complexity are based on Talagrand's broom space $\mathbf{T}$ (resp. on some other Talagrand compacta -- see \cite{argyros2008talagrand}). An interesting question is therefore whether having absolute complexity is the~``typical case'' for a~topological space (and $\mathbf{T}$ is an anomaly), or whether there in fact exist many spaces with non-absolute complexity (but somehow this is difficult to prove).
As a~specific example, recall that for some Banach spaces, the~unit ball $B_X$ with weak topology is $\fsd$ the~bi-dual unit ball $B_{X^{\star\star}}$ with $w^\star$ topology. Is $(B_X,w)$ absolutely $\fsd$?

We now give a~brief overview of the~contents of the~paper.
Section~\ref{section: background} contains some preliminary concepts -- compactifications, the~(mostly standard) terminology describing the~trees on $\omega$, and several elementary results regarding derivatives on such trees.
Section~\ref{section: topological sums} is devoted to providing a~simple proof of Theorem~\ref{theorem: summary} by treating the~spaces $X_\alpha^\beta$ from Theorem~\ref{theorem:X_2_beta} as atomic and taking their generalized topological sums (``zoom spaces'', introduced in Section~\ref{section: zoom spaces}).

In Section~\ref{section: simple representations}, we introduce the~concept of a~\emph{simple $\Fa$-representation}, and show how this concept can be used to give an elegant proof of Proposition~\ref{proposition: hereditarily lindelof spaces are absolute}.
We also investigate the~concept of local complexity and its connection with the~standard complexity. As a~side-product, we prove that a~``typical'' topological space $X$ cannot have a~``universal $\Fa$-representation'', even when its complexity is absolute.

In Section~\ref{section: regular representations}, we introduce the~concept of a~\emph{regular $\Fa$-representation}. As the~name suggests, this is a~(formally) stronger notion than that of a~simple $\Fa$-representation. In particular, it allows us to take any two spaces $X\subset Y$, and construct ``$\Fa$-envelopes'' of $X$ in $Y$, looking for an ordinal $\alpha$ for which the~$\Fa$-envelope of $X$ will be equal to $X$.
We also justify the~concept of a~regular $\Fa$-representation by showing that if $X$ is an $\Fa$ subset of $Y$, then it does have a~regular $\Fa$-representation in $Y$.

In Section~\ref{section: complexity of brooms}, we study the~class of the~so-called ``Talagrand's broom spaces'' -- spaces based on Talagrand's example $\mathbf{T}$.
We investigate the~class of broom-space compactifications which are ``amalgamation-like''.
As an application of regular $\Fa$-representations, we compute which complexities are attainable by these spaces. In particular, this \errata{provides a specific example whose existence proves} Theorem~\ref{theorem: summary}.
What makes this method valuable is that, unlike the~simple approach from Section~\ref{section: topological sums}, it holds a~promise of being applicable not only to broom spaces, but to some other topological spaces as well.

It should be noted that many parts of this article can be read independently of each other, as explained by the~following remark.

\begin{remark}[How to use this paper]\label{remark:how_to_use}
It is \emph{not} necessary to read Sections~\ref{section: topological sums},~\ref{section: simple representations} and~\ref{section: regular representations} in order to read the~subsequent material.
\end{remark}

\noindent Sections~\ref{section: simple representations},~\ref{section: regular representations} and~\ref{section: complexity of brooms} all rely on the notation introduced in Section~\ref{section: sequences}. This notation is, however, not required in Section~\ref{section: topological sums}.
Section~\ref{section: topological sums} deals with topological sums, and can be read independently of the~content of any of the~following sections.
Moreover, Sections~\ref{section: simple representations},~\ref{section: regular representations} and~\ref{section: complexity of brooms} are only loosely related, and can be read independently of each other with one exception: Section~\ref{section: broom space properties} relies on one result from Section~\ref{section: Suslin scheme rank}. But if the~reader is willing to use this result as a~``black box'', this dependency can be ignored.

\section{Preliminaries}\label{section: background}

This section reviews some preliminary results.
The article assumes familiarity with the~concept of compactifications, whose basic overview can be found in Section~\ref{section: compactifications}.
Section~\ref{section: sequences} introduces the~notation used to deal with sequences, trees, and derivatives on trees.


\subsection{Compactifications and Their Ordering}\label{section: compactifications}

By a~\emph{compactification} of a~topological space $X$ we understand a~pair $(cX,\varphi)$, where $cX$ is a~compact space and $\varphi$ is a~homeomorphic embedding of $X$ onto a~dense subspace  of $cX$. Symbols $cX$, $dX$ and so on will always denote compactifications of $X$.

Compactification $(cX,\varphi)$ is said to be \emph{larger} than $(dX,\psi)$, if there exists a~continuous mapping $f : cX\rightarrow dX$, such that $\psi = f \circ \varphi$. We denote this as $ cX \succeq dX $.
Recall that for a~given $T_{3\slantfrac{1}{2}}$ topological space $X$, its compactifications are partially ordered by $\succeq$ and Stone-Čech compactification $\beta X$ is the~largest one.
Sometimes, there also exists the~smallest compactification $\alpha X$, called \emph{one-point compactification} or \emph{Alexandroff compactification}, which only consists of a~single additional point. 

In this paper, we will always assume that $cX \supset X$ and that the~corresponding embedding is identity. In particular, we will simply write $cX$ instead of $(cX,\textrm{id}|_X)$.

Much more about this topic can be found in many books, see for example \cite{freiwald2014introduction}. The basic relation between the~complexity of a~space $X$ and the~ordering of compactifications is the~following observation:

\begin{remark}[Larger compactification means smaller complexity]\label{remark: absolute complexity}
For any $\alpha\leq \omega_1$, we have
 \[ X\in\Fa (dX), \ cX\succeq dX \implies X\in\Fa (cX). \]
\end{remark}

\subsection{Trees and Derivatives on Trees}\label{section: sequences}

We now introduce the notation needed by Sections~\ref{section: simple representations},~\ref{section: regular representations} and~\ref{section: complexity of brooms}.
We start with sequences in $\omega$:

\begin{notation}[Finite and infinite sequences in $\omega$]\label{notation: sequences}
We denote
\begin{itemize}
\item $\omega^\omega:=$ infinite sequences of non-negative integers $:=\left\{\sigma : \omega \rightarrow \omega\right\}$,
\item $\seq:=$ finite sequences of non-neg. integers $:=\left\{s : n \rightarrow \omega |\ n\in\omega \right\}$.
\end{itemize}
\end{notation}

Suppose that $s\in\seq$ and $\sigma\in\baire$. We can represent $\sigma$ as $(\sigma(0),\sigma(1),\dots)$ and $s$ as $(s(0),s(1),\dots,s(n-1))$ for some $n\in\omega$.
We denote the~\emph{length} of $s$ as $|s|=\textrm{dom}(s)=n$, and set $|\sigma|=\omega$.
If for some $t\in \seq \cup \baire$ we have $|t|\geq |s|$ and $t|_{|s|}=s$, we say that $t$ \emph{extends} $s$, denoted as $t\sqsubset s$.
We say that  $u,v\in\seq$ are non-comparable, denoting as $u\perp v$, when neither $u\sqsubset v$ nor $u\sqsupset v$ holds.


Unless we say otherwise, $\omega^\omega$ will be endowed with the~standard product topology $\tau_p$, whose basis consists of sets $\mc N (s):=\left\{ \sigma \sqsupset s | \ \sigma\in\baire \right\}$, $s \in \seq$.

For $n\in\omega$, we denote by $(n)$ the~corresponding sequence of length $1$.
By $s\ext t$ we denote the~\emph{concatenation} of a~sequences $s\in\seq$ and $t\in\seq\cup \baire$.
We will also use this notation to extend finite sequences by integers and sets, using the~convention $t \hat{\ } k:=t \ext (k)$ for $k\in\omega$ and $t \hat \ S:=\left\{t\ext s\, | \ s\in S\right\}$ for $S\subset \seq \cup \baire$.

\bigskip
Next, we introduce the~notation related to trees.

\begin{definition}[Trees on $\omega$]
A \emph{tree} (on $\omega$) is a~set $T\subset \seq$ which satisfies
\[ (\forall s,t\in\seq): s\sqsubset t \ \& \ t\in T \implies s\in T. \]
\end{definition}

\noindent By $\textrm{Tr}$ we denote the~space of all trees on $\omega$.
For $S\subset \seq \cup \baire$ we denote by
\[ \cltr{S} := \left\{u\in\seq | \ \exists s\in S: s\sqsupset u \right\} \]
the smallest tree ``corresponding'' to $S$ (for $S\subset\seq$, $\mathrm{cl}_\mathrm{Tr}(S)$ is the~smallest tree containing $S$).
Recall that the~empty sequence $\emptyset$ can be thought of as the~`root' of each tree, since it is contained in any nonempty tree.
For any $t \in T$, we denote the~set of \emph{immediate successors of $t$ in $T$} as
\begin{equation}\label{equation: immediate successor in T}
\ims{T}{t} := \{ t\ext n | \ n\in\omega , \, t\ext n \in T \} .
\end{equation}
A \emph{leaf} of a~tree $T$ is a~sequence $t\in T$ with no immediate successors in $T$, and the~set of all \emph{leaves of $T$} is denoted as $l(T)$. Finally, we denote as $T^t$ the~tree corresponding to $t$ in $T$:
\begin{equation}\label{equation: T of t}
T^t := \{ t' \in \seq | \ t\ext t' \in T \} .
\end{equation}
For each non-empty tree $T$, we clearly have
\[ T = \{\emptyset\} \cup \bigcup \{ m \ext T^{(m)} \, | \ (m) \in \ims{T}{\emptyset} \} .\]

If each of the~initial segments $\sigma|n$, $n\in\omega$, of some $\sigma\in\baire$ belongs to $T$, we say that $\sigma$ is an infinite branch of $T$.
By $\textrm{WF}$ we denote the~space of all trees which have no infinite branches (the \emph{well-founded} trees). By $\textnormal{IF}$, we denote the~complement of $\textnormal{WF}$ in $\Tr$ (the \emph{ill-founded} trees).

\bigskip
Useful notion for studying trees is the~concept of a~derivative:

\begin{definition}[Derivative on trees]
A mapping $D:\Tr\rightarrow\Tr$ is a~\emph{derivative} on trees if it satisfies $D(T)\subset T$ for every $T\in \Tr$.
\end{definition}

Any derivative on trees admits a~natural extension to subsets of $\seq \cup \baire$ by the~formula $D(S) := D( \cltr{S} )$.
We define \emph{iterated derivatives} $D^\alpha$ in the~standard way:
\begin{align*}
D^0(S) & := \cltr{S}, \\
D^{\alpha +1} \left( S \right) & :=  D\left(D^\alpha\left(S\right)\right) \textrm{ for successor ordinals},\\
D^\lambda \left( S \right) & := \underset {\alpha < \lambda } \bigcap D^\alpha \left( S \right) \textrm{ for limit ordinals}.
\end{align*}

Clearly, the~iterated derivatives of any $S$ either keep getting smaller and smaller, eventually reaching $\emptyset$, or there is some $\alpha$ for which the~iterated derivatives no longer change anything, giving $D^\alpha(T)=D^{\alpha+1}(T)\neq \emptyset$.
Since each $T\in\Tr$ is countable, it suffices to consider $\alpha<\omega_1$.
This allows us to define a~\emph{rank} corresponding $D$ (on subsets of $\seq\cup \baire$):
\[ r(S) := \begin{cases}
-1 & \text{for } S=\emptyset ;\\
\min \{ \alpha<\omega_1| \ D^{\alpha+1}\left( S \right) = \emptyset \}
	& \text{for $S\neq \emptyset$, if there is some $\alpha<\omega_1$} \\
	& \text{ such that } D^{\alpha+1}\left( S \right) = \emptyset ;\\
\omega_1 & \text{if no $\alpha$ as above exists.}
\end{cases} \]

Note that when $S$ is non-empty and $r(S)<\omega_1$, then $r(S)$ is the~highest ordinal for which $D^{r(S)}(S)$ is non-empty (or equivalently, for which $D^{r(S)}(S)$ contains the~empty sequence).
For more details and examples regarding ranks, see for example \cite[ch. 34 D,E]{kechris2012classical}.
We will be particularly interested in the~following three derivatives on trees:

\begin{definition}[Examples of derivatives]\label{definition: derivatives}
For $T\in \textrm{Tr}$, we denote
\begin{align*}
D_l (T) \ := \big\{ t\in T| \ & T \text{ contains some extension $s\neq t$ of $t$} \big\} ,\\
D_i (T) \ := \big\{ t\in T| \ & T \text{ contains infinitely many extensions of $t$} \big\} ,\\
\D (T)  := \big\{ t\in T| \ & T \text{ contains infinitely many incomparable} \\
								& \text{ extensions of $s$ of different length} \big\} .
\end{align*}
We use the~appropriate subscripts to denote the~corresponding iterated derivatives and ranks.
\end{definition}
Using transfinite induction, we obtain the~following recursive formula for the~leaf-rank 
\begin{equation}\label{equation: leaf rank formula}
\left( \forall t\in T \right) :
r_l(T^t) = \sup \{ r_l(T^s)+1 \,| \ s\in \ims{T}{t} \} 
\end{equation}
(with the~convention that a~supremum over an empty set is 0).
It is straightforward to check that the~leaf-rank $r_l(T^t)$ of $T^t$ is the~highest ordinal for which $t$ belongs to the~iterated leaf-derivative $D^{r_l(T)}_l(T)$ of $T$.
In particular, leaves of $T$ are precisely those $t\in T$ which satisfy $r_l(T^t)=0$.

For a~general $T\in \Tr$, $r_l(T)$ is countable if and only if $r_i(T)$ is countable, and this happens if and only if $T$ is well founded.
Note that on well founded trees, $D_i$ behaves the~same way as the~derivative from \errata{example 3) from \cite[34.D]{kechris2012classical}}, but it leaves any infinite branches untouched.

The following trees serve as examples of trees of rank $\alpha$ for both $r_l$ and $r_i$:

\begin{notation}[``Maximal'' trees of height $\alpha$] \label{notation: trees of height alpha}
For each limit $\alpha<\omega_1$, fix a~bijection $\pi_\alpha : \omega \rightarrow \alpha$.
We set\errata{}
\begin{align*}
T_0 			& := \{\emptyset\} ,\\
T_{\alpha}	& := \{\emptyset\} \cup \underset {\err{j}\in\omega} \bigcup \err{j}\hat{\ }T_{\alpha-1} & \text{for countable successor ordinals,} \\
T_\alpha		& := \{\emptyset\} \cup \underset {\err{j}\in\omega} \bigcup \err{j}\hat{\ }T_{\pi_\alpha (\err{j})}	& \text{for countable limit ordinals,} \\
T_{\omega_1}	& := \seq .
\end{align*}
\end{notation}

In particular, we have $T_k = \omega^{\leq k}$ for $k\in\omega$ and $T_\omega = \{\emptyset\} \, \cup \, (0) \, \cup \, 1\ext \omega^{\leq 1} \, \cup \, 2\ext \omega^{\leq 2} \, \cup \, \dots $ .
Denoting $\pi_\alpha(\err{j}):=\alpha-1$ for successor $\alpha$, we can write $T_\alpha$ as $T_\alpha = \{\emptyset\} \cup \underset {\err{j}\in\omega} \bigcup \err{j}\hat{\ }T_{\pi_\alpha (\err{j})}$ for both limit and successor ordinals.
A straightforward induction over $\alpha$ yields $r_l(T_\alpha) = r_i(T_\alpha) = \alpha$ for each $\alpha \leq\omega_1$.

Every ordinal $\alpha$ can be uniquely written as $\alpha = \lambda + 2n + i$, where $\lambda$ is $0$ or a~limit ordinal, $n\in\omega$ and $i\in\{0,1\}$. We denote $\alpha' := \lambda + n$.
In Section~\ref{section: R T sets}, we will need trees which serve as intermediate steps between $T_\alpha$ and $T_{\alpha+1}$. These will consist of trees $\{\emptyset\} \cup \err{j}\ext T_\alpha$, $\err{j}\in \omega$. While we could call these trees something like `$T_{\alpha+\frac 1 2}$', we instead re-enumerate them as
\begin{align} \label{equation: canonical trees}
T^c_\alpha 		& := T_{\alpha'} & \text{for even $\alpha\leq\omega_1$} \nonumber ,\\
T^c_{\alpha,\err{j}} 	& := \{\emptyset\} \cup \err{j}\ext T_{\alpha'}
				& \text{for $\err{j}\in\omega$ and odd $\alpha<\omega_1$} ,\\
T^c_{\alpha} 	& := \{\emptyset\} \cup 1\ext T_{\alpha'} = T^c_{\alpha,1}
					& \text{for odd $\alpha<\omega_1$} . \nonumber
\end{align}
\section{Attaining Complexities via Topological Sums} \label{section: topological sums}

Suppose we have topological spaces $X$ and $Y$ satisfying
\[ \{\alpha,\beta\}\subset\textnormal{Compl}(X) \subset [\alpha,\beta] \ \ \ \text{ and } \ \ \ \{\alpha,\gamma\} \subset \textnormal{Compl}(Y) \subset [\alpha,\gamma ] \]
for some $\alpha \leq \beta \leq \gamma$.
It is straightforward to verify that the~topological sum $X\oplus Y$ satisfies
\[ \{ \alpha,\beta,\gamma\} \subset \textnormal{Compl}(X\oplus Y) = [\alpha,\gamma] .\]
In this section, we extend this observation to topological sums of infinitely many spaces.
We shall do so with the~help of the~concept of ``zoom'' spaces, of which topological sums are a~special case.
Moreover, we will be able to put together uncountably many spaces -- not exactly as a~topological sum, but in a~very similar manner -- and show that this gives the~existence of a~space with $\textnormal{Compl}(X)=[2,\omega_1]$.

In Section~\ref{section: zoom spaces}, we introduce and investigate zoom spaces. In Section~\ref{section:top_sums_as_zoom}, we apply this theory to topological sums, thus proving Theorem~\ref{theorem: summary}. Apart from that, the~results of Section~\ref{section: topological sums} are independent on the~other parts of the~paper.

\subsection{Zoom spaces} \label{section: zoom spaces}
Let $Y$ be a~topological space. Throughout this section, we shall denote by $I_Y$ the~set of all isolated points of $Y$ and by $Y' = Y \setminus I_Y$ the~set of all non-isolated points of $Y$.
\begin{definition}[Zoom space]\label{definition: zoom space}
Let $Y$ be a~topological space and $\mc X=(X_i)_{i\in I_Y}$ a~collection of non-empty topological spaces.
We define the~\emph{zoom space $Z(Y,\mc X)$ of $Y$ with respect to $\mc X$} as the~disjoint union\errata{\footnote{\errata{By saying that $Y'\cup \bigcup \mc X$ is a disjoint union we mean that the corresponding spaces are assumed to be disjoint (or changed accordingly). We avoid the term `topological sum' since it is typically associated with a particular topology.}}} $Y' \cup \bigcup \mc X$. The basis of topology of $Z(Y,\mc X)$ consists of open subsets of $X_i$, $i\in I_Y$, and of all sets of the~form
\[ V_U := (U \setminus I_Y ) \cup \bigcup \{ X_i  | \ i \in U \cap I_Y \}
		, \ \ \ \ U\subset Y \text{ open in }Y .\]
\end{definition}

\noindent \errata{Note that $V_{U_1}\cap V_{U_2} = V_{U_1\cap U_2}$, so these sets indeed form a basis.}
Moreover, the~definition above works even if the~indexing set $I_Y$ of the~collection $\mc X$ is replaced by some $I\subset I_Y$.
Lastly, any collection $(X_i)_{i\in I}$ indexed by such a~subset can be extended to $\mc X \cup (\{x\})_{x\in I_Y\setminus I}$ and the~corresponding zoom space is identical to $Z(Y,\mc X)$.

The following notation for families of sets shall be used throughout the~paper:

\begin{notation}[Collection of compactifications]\label{notation: c A}
Let $\mc X=(X_i)_{i\in I}$ be a~collection of topological spaces and suppose that for every $i\in I$, $c X_i$ is a~topological space containing $X_i$. We denote
\[ c\mc X:= \{cX_i | \ i \in I \} .\]
\end{notation}

The basic properties of zoom spaces are as follows:

\begin{proposition}[Basic properties of zoom spaces]\label{proposition: properties of Z}
Let $Z(Y,\mc X)$ be a~zoom space.
\begin{enumerate}[(i)]
\item $Z(Y,\mc X)$ is a~Tychonoff space and $Y$ is its quotient. Each $X_i$ is homeomorphic to the~clopen subset $X_i \subset Z(Y,\mc X)$.  For any selector $s$ of $\mc X$, $Y$ is homeomorphic to the~closed set $Y_s := Y' \cup s(I_Y) \subset Z(Y,\mc X)$. 	\label{case: Y as subspace of Z(Y,Z)}
\item If  $Y$ is a~dense subspace of $\widetilde Y$ and each $X_i$ is a~(dense) subspace of $\widetilde X_i$, then $Z(Y,\mc X)$ is a~(dense) subspace of $Z(\widetilde Y,\widetilde{\mc X})$. \label{case: subspaces}
\item If $Y$ and all $X_i$ are compact (Lindelöf), then $Z(Y,\mc X)$ is compact (Lindelöf). \label{case: Z is compact}
\item If $cY$ and $cX_i$, $i\in I_Y$, are compactifications, then $Z(cY, c\mc X)$ is a~compactification of $Z(Y,\mc X)$. \label{case: compactifications of Z(Y,Z)}
\end{enumerate}
\end{proposition}

\begin{proof}
$(i)$: Clearly, Definition~\ref{definition: zoom space} correctly defines a~topology on $Z(Y,\mc X)$ by giving its basis.
By definition of this topology, a~set $U\subset X_i$ is open in $Z(Y,\mc X)$ if and only if it is open in $X_i$, so the~topology on each $X_i$ is preserved. For every $i\in I_Y$, $Y\setminus \{i\}$ is open in $Y$. It follows that each $Z(Y,\mc X)\setminus X_i = V_{Y\setminus \{ i \}}$ is open and hence each $X_i$ is clopen in $Z(Y,\mc X)$.

To see that $Y$ is a~quotient of $Z(Y,\mc X)$, it is enough to observe that the~following mapping $q : Z(Y,\mc X) \rightarrow Y$  is continuous, surjective and open:
\[ q(x) :=
\begin{cases}
x, 	& \textrm{for } x\in Y'	\\
i, 	& \textrm{for } x\in X_i,\ i\in I_Y .
\end{cases} \]
Indeed, $q$ is clearly continuous (since $q^{-1}(U)=V_U$) and surjective.  Moreover, it maps basic open sets in $Z(Y,\mc X)$ onto open sets: we have $q(V)=\{i\} \in I_Y$ for non-empty $V\subset X_i$ and $q(V_U)=U$ for $U\subset Y$.

Let $s:I_Y\rightarrow \bigcup \mc X$ be a~selector of $\mc X$ and denote by $f_s$ the~following restriction of $q$:
\[ f_s := q|_{Y'\cup s(I_Y)} : Y' \cup s(I_Y) \rightarrow Y .\]
Clearly, this restriction of $q$ is injective and continuous. The restriction is an open mapping, because the~image of $V\subset X_i$ is either empty (when $s(i)\notin V$) or equal to $\{i\}$ and for $U\subset Y$ we have $f_s(V_U) = q(V_U) = U$. In particular, the~range of $f_s$ is $Y$, so $Y_s := Y'\cup s(I_Y)$ and $Y$ are homeomorphic.

The topology on $Z(Y,\mc X)$ is easily seen to be Hausdorff -- indeed, if $x,y\in Z(Y,\mc X)$ are distinct, then either both of them belong to some $X_i$ (which is Hausdorff and open in $Z(Y,\mc X)$), or one of them belongs to $X_i$ and the~other to $Z(Y,\mc X)\setminus X_i$ (so we separate $x$ from $y$ by $X_i$ and its complement), or both belong to $Y'$. In the~last case, we use the~fact that $Y$ is Hausdorff to separate $x$ and $y$ in $Y$ by open subsets $U, U'$ of $Y$ and note that $V_U$ and $V_{U'}$ are open sets separating $x$ and $y$ in $Z(Y,\mc X)$.

To see that $Z(Y,\mc X)$ is Tychonoff, let $F\subset Z(Y,\mc X)$ be closed and $x\in Z(Y,\mc X) \setminus F$. If $x\in Y'$, we find and open subset $U$ of $Y$ such that $x\in V_U \subset Z(Y,\mc X) \setminus F$. We separate $x$ from $Y\setminus U$ in $Y$ by a~continuous function $f:Y\rightarrow [0,1]$ and note that $f\circ q$ is a~continuous function separating $x$ from $F$ in $Z(Y,\mc X)$.

If $x\in X_i$, we find continuous $f:X_i\rightarrow [0,1]$ which separates $x$ from $F\cap X_i$. Since $X_i$ is clopen in $Z(Y,\mc X)$, $f$ can be extended into a~function $\tilde f$ which separates $x$ from $F$ in $Z(Y,\mc X)$.

$(ii)$ immediately follows from the~definition of $Z(\cdot,\cdot)$ and its topology.

$(iii)$: Assume that $Y$ and all $X_i$ are Lindelöf and let $\mc V$ be an open cover of $Z(Y,\mc X)$. Without loss of generality, we can assume that $\mc V$ consists only of basic open sets, that is, we have $\mc V = \mc V_0 \cup \bigcup_{i\in I} \mc V_i$ where each $\mc V_i$, $i\in I$, only contains open subsets of $X_i$ and $\mc V_0$ only contains sets of the~form $V_U$ for $U\subset Y$ open in $Y$.
Denote
\[ \mc U_0:=\{ U \subset Y | \ V_U \in \mc V_0 \} . \]
By definition of topology on $Z(Y,\mc X)$, $\mc U_0$ is a~cover \err{of} $Y'$.
Since $Y$ is Lindelöf and $Y'$ is closed, there is a~countable  $\mc U^{\star}_0 \subset \mc U_0$ such that $\bigcup \mc U^{\star}_0 \supset Y'$. The set
\[ I^{\star}:= Y\setminus \bigcup \mc U^{\star}_0 = I_Y \setminus \bigcup \mc U^{\star}_0 \]
is closed and discrete in $Y$, and hence countable. Since each $X_i$ is Lindelöf, each $\mc V_i$ has a~countable subcover $\mc V^{\star}_i$. The system $\mc V^{\star} := \mc V_0 \cup \bigcup_{i\in I^{\star}} \mc V^{\star}_i$ is then a~countable subcover of $\mc V$.

The proof of the~compact case is the~same.

$(iv)$: This follows from $(ii)$ and $(iii)$.
\end{proof}

Recall that Definition~\ref{definition: zoom space} introduced the~sets of the~form $V_F$ for $F\subset W$:
\[ V_F := F \cup \bigcup \{ X_i | \ i\in I_W, \, F \ni i \} \subset Z(W,\mc Z) .\]
The following lemma investigates the~complexity of these sets.

\begin{lemma}[Complexity of basic sets]\label{lemma: complexity of V F sets}
For a~zoom space $Z(W,\mc Z)$, $\mc Z=(Z_i)_{i}$, the~sets $V_F$ satisfy
\begin{equation}\label{equation: complexity of V F sets}
\left( \forall F \subset W \right) : F\in \Fa (W) \implies V_F \in \Fa(Z(W,\mc Z)) .
\end{equation}
\end{lemma}

\begin{proof}
Indeed, if $F$ is closed, then $Z(W,\mc Z)\setminus V_F = V_{W\setminus F}$ is open by definition of topology on $Z(W,\mc Z)$, so $V_F$ is closed. Moreover, $\bigcup_{E\in \mc E} V_E = V_{\bigcup {\mc E}}$ and $\bigcap_{E\in \mc E} V_E = V_{\bigcap {\mc E}}$ holds for any $\mc E\subset \mc P(W)$. This implies the~result for higher $\mc F$-Borel classes (by transfinite induction) and for $F\in\mc F_{\omega_1}$ (by showing that if $F$ is Suslin in $W$, then $V_F$ is Suslin in $Z(W,\mc Z)$).
\end{proof}

Before computing the~complexity of a~general zoom space, we need to following technical lemma:

\begin{lemma}[Upper bound on zoom space complexity]\label{lemma: complexity in zoom spaces}
Let $Z(W,\mc Z)$, $\mc Z=(Z_i)_{i\in I}$, be a~zoom space and $\mc C = (C_i)_{i\in I}$ a~system satisfying $C_i \subset Z_i$ for each $i\in I$. Then for any $\alpha \leq \omega_1$, we have
\begin{equation}\label{equation: complexity in zoom spaces}
\left( \left( \forall i \in I \right) : C_i \in \Fa (Z_i) \right) \implies Z(W,\mc C) \in \Fa(Z(W,\mc Z)).
\end{equation}
\end{lemma}

\begin{proof}
For $\alpha=0$, $Z(W,\mc C)$ is closed in $Z(W,\mc Z)$ since $Z_i \setminus C_i$ are basic open sets and
\[ Z(W,\mc C) = Z(W,\mc Z) \setminus \bigcup_{i\in I} (Z_i \setminus C_i) .\]
For $0<\alpha<\omega_1$, \eqref{equation: complexity in zoom spaces} follows by transfinite induction, because we have
\[ \left( \left( \forall i\in I \right) : C_i=\bigcup_n C_i^n \right) \implies
Z\left(W,\mc C\right) = Z(W, ( \bigcup_n C_i^n )_i ) =
 \bigcup_n Z\left(W, (C_i^n )_i \right) \]
and the~analogous formula holds when each $C_i$ satisfies $C_i=\bigcap_n C_i^n$.

For $\alpha=\omega_1$, the~following formula shows that $Z(W,\mc C)$ is Suslin in $Z(W, \mc Z)$:
\begin{align*}
& \left( \left( \forall i\in I \right) : C_i = \bigcup_{\sigma\in\baire} \bigcap_n C_i^{\sigma|n}
 \text{, where each } C_i^{\sigma|n} \text{ is closed in } Z_i  \right) \implies \\
& Z\left(W,\mc C\right) = Z\left(W, \left( \bigcup_{\sigma\in\baire} \bigcap_n C_i^{\sigma|n} \right)_i \right) = \bigcup_{\sigma\in\baire} \bigcap_n Z\left(W, (C_i^{\sigma |n} )_i \right) \\
& \in \mc F_{\omega_1} \left( Z(W,\mc Z)\right) ,
\end{align*}
where at the~last line we have used \eqref{equation: complexity in zoom spaces} with $\alpha=0$.
\end{proof}

The following result then states that the~complexity of zoom spaces can be retrieved from the~complexity of its parts.

\begin{proposition}[Complexities attained by a~zoom space]\label{proposition: complexity of Z(Y,Z)}
Let $Z(Y,\mc X)$ be a~zoom space. If $Y$ is a~dense subspace of $\widetilde Y$ and each $X_i$ is a~subspace of $\widetilde X_i$, then the~spaces $Z(Y, {\mc X}) \subset Z(\widetilde Y, \widetilde {\mc X})$ satisfy
\[ \textnormal{Compl} \left( Z(Y,\mc X), Z(\widetilde Y, \widetilde {\mc X}) \right) = 
	\max\left\{ \textnormal{Compl} \left( Y, \widetilde Y \right), \ 
	\sup_{i\in I_Y} \textnormal{Compl} \left( X_i, \widetilde X_i \right) \right\} ,\]
whenever at least one of the~sides is defined.
\end{proposition}

\begin{proof}
"$\geq$": Let $Y$, $\widetilde Y$, $\mc X$ and $\widetilde {\mc X}$ be as in the~statement. Assume that the~LHS is defined and $Z(Y,\mc X) \in \Fa ( Z(\widetilde Y, \widetilde {\mc X}) )$ holds for some $\alpha \leq \omega_1$.
Let $s$ be a~selector of $\mc X$. Since $\widetilde Y$ and $Y$ have the~same isolated points, $s$ is also a~selector of $\widetilde{\mc X}$. Using the~notation from Proposition~\ref{proposition: properties of Z}\,\eqref{case: Y as subspace of Z(Y,Z)}, we have
\begin{align*}
Y_s \subset \widetilde Y_s \subset Z(\widetilde Y, \widetilde {\mc X}) && \ \& \ && Y_s   & = \widetilde Y_s \cap Z(Y,\mc X) \subset Z(\widetilde Y, \widetilde {\mc X}) \\
X_i \subset \widetilde X_i \subset Z(\widetilde Y, \widetilde {\mc X}) && \ \& \ && X_i & = \widetilde X_i \cap Z(Y,\mc X) \subset Z(\widetilde Y, \widetilde {\mc X}) ,
\end{align*}
where $\widetilde Y_s$ and $\widetilde X_i$ are closed subsets of $Z(\widetilde Y,\widetilde{\mc X})$.
In particular we get (by definition of the~subspace topology) that each $X_i$ satisfies $X_i \in \Fa ( \widetilde X_i)$ and $Y_s$ satisfies $Y_s \in \Fa ( \widetilde Y_s)$. Since $\widetilde Y_s$ is homeomorphic to $\widetilde Y$, we get $Y \in \Fa (\widetilde Y)$. It follows that the~RHS is defined and at most equal to the~LHS.

"$\leq$":
Let $Y$, $\widetilde Y$, $\mc X$ and $\widetilde {\mc X}$ be as in the~proposition. Assume that the~RHS is defined and we have $Y \in \Fa (\widetilde Y)$ and $X_i \in \Fa ( \widetilde X_i)$, $i\in I_Y$, for some $\alpha \leq \omega_1$.
Rewriting the~space $Z(Y,\mc X)$ as
\begin{equation} \label{equation: Z as intersection in tilde Z}
Z := Z\left(Y,\mc X\right) = Z\left(Y,\widetilde {\mc X}\right) \cap Z\left(\widetilde Y,\mc X\right)
 \subset Z\left(\widetilde Y,\widetilde {\mc X} \right) =: \widetilde Z ,
\end{equation}
we obtain an upper bound on its complexity:
\[ \textnormal{Compl}\left( Z, \widetilde Z \right) \le \max \left\{
	\textnormal{Compl}\left( Z\left(Y, \widetilde{\mc X}\right), \widetilde Z \right) ,
	\textnormal{Compl}\left( Z\left(\widetilde Y,\mc X\right), \widetilde Z \right)
 \right\} .\]
Working in the~zoom space $Z(\widetilde Y,\widetilde{ \mc X})$, we have $Z(Y,\widetilde{\mc X})=V_Y$. Applying \eqref{equation: complexity of V F sets} (with $W:=\widetilde Y$, $\mc Z := \widetilde {\mc X}$ and $\alpha$), we get
\[ \textnormal{Compl}\left( Z\left(Y, \widetilde{\mc X}\right), \widetilde Z \right) =
 \textnormal{Compl}\left( V_Y, Z(\widetilde Y, \widetilde{\mc X}) \right)
 \overset{\eqref{equation: complexity of V F sets}}{\le}
 \textnormal{Compl}\left( Y, \widetilde Y \right) \leq \alpha .\]
Application of \eqref{equation: complexity in zoom spaces} (with $W:= \widetilde Y$, $\mc Z := \widetilde{\mc X}$, $\mc C := \mc X$ and $\alpha$) implies that $\textnormal{Compl}\left( Z\left(\widetilde Y, {\mc X}\right), \widetilde Z \right) \leq \alpha $. It follows that the~LHS is defined and no greater than the~RHS.
\end{proof}

As a~corollary of the~proof of Proposition~\ref{proposition: properties of Z}, we obtain the~following result:

\begin{corollary}[Stability under the~zoom space operation]\label{corollary: zoom spaces are K analytic}
A zoom space $Z(Y,\mc X)$ is compact ($\sigma$-compact, $\mc K$-analytic) if and only if $Y$ \err{and} each $X_i$ is compact ($\sigma$-compact, $\mc K$-analytic).
\end{corollary}
\begin{proof}
$\Longrightarrow$: By Proposition~\ref{proposition: properties of Z} \eqref{case: Y as subspace of Z(Y,Z)}, $Z(Y,\mc X)$ contains $Y$ and each $X_i$ as a~closed subspace -- this implies the~result.

$\Longleftarrow$: The compact case is the~same as Proposition~\ref{proposition: properties of Z} \eqref{case: Z is compact}. If $Y$ and each $X_i$ is $\sigma$-compact, then each of these spaces is $F_\sigma$ in its Čech-stone compactification. By Proposition~\ref{proposition: complexity of Z(Y,Z)}, $Z(Y,\mc X)$ is an $F_\sigma$ subset of $Z(\beta Y,\beta \mc X)$. The lat\err{t}er space is compact by the~compact case of this proposition, which implies that $Z(Y,\mc X)$ is $\sigma$-compact.

Consider the~case where $Y$ and each $X_i$ is $\mc K$-analytic. In the~proof of Proposition~\ref{proposition: complexity of Z(Y,Z)}, we have shown that both $Z\left(Y,\beta {\mc X}\right)$ and $Z\left(\beta Y,\mc X\right)$ are Suslin in $Z(\beta Y, \beta \mc X)$. It follows that both these sets are $\mc K$-analytic. By \eqref{equation: Z as intersection in tilde Z}, $Z(Y,\mc X)$ is an intersection of two $\mc K$-analytic sets, which implies that it it itself $\mc K$-analytic.
\end{proof}

\subsection{Topological sums as zoom spaces}\label{section:top_sums_as_zoom}
In this section, we give some natural examples of zoom spaces, and prove the~first part of our main results.

\begin{example}[Topological sums as zoom spaces] \label{example: application to topological sums}
Whenever $\mc X=(X_i)_{i\in I}$ is a~collection of topological spaces, the~topological sum $\bigoplus \mc X$ is homeomorphic to the~zoom space $Z(I,\mc X)$ of the~discrete space $I$.
\end{example}
\begin{proof}
This follows from the~fact that each $X_i$ is clopen in $Z(I,\mc X)$ and that in this particular case, we have $I_I=I$ and hence $Z(I,\mc X)$ is a~disjoint union of $X_i$, $i\in I$.
\end{proof}

In the~particular case of countable sums, we can use Proposition~\ref{proposition: complexity of Z(Y,Z)} to fully describe the~complexity of the~resulting space:
\begin{proposition}[Complexities attainable by a~topological sum]\label{proposition: complexities of topological sums}\mbox{}
\begin{enumerate}[(i)]
\item If $X_k$, $k\leq n$, are compact, then $\overset{n}{\underset{k=0}{\bigoplus}} X_k$ is compact.
\item If $X_k$, $k\in \omega$, are $\sigma$-compact, then $\underset{k\in \omega}{\bigoplus} X_k$ is $\sigma$-compact.
\item If $X_k$, $k\in\omega$, are $\mc K$-analytic and at least one $X_k$ is not compact, then
\[ \textnormal{Compl}\left( \bigoplus_{k\in\omega} X_k \right) = \left\{ \sup f | \ f\in \prod_{k\in\omega} \textnormal{Compl}(X_k) \right \} . \]
\end{enumerate}
\end{proposition}
\begin{proof}
$(i)$ and $(ii)$ are obvious. To get $(iii)$, note that by Example~\ref{example: application to topological sums}, the~topological sum $\bigoplus_{k} X_k$ can be rewritten as
\[ X:=\bigoplus_{k\in\omega} X_k = Z(\omega,(X_k)_{k}) .\]
We will prove each of the~two inclusions between $\textnormal{Compl}(X)$ and $\{ \sup f | \ f\in \Pi_k \textnormal{Compl}(X_k) \}$.
\errata{To simplify the notation, we assume that $X_k$ are disjoint.}

"$\supset$": Let $f$ be a~selector for $(\textnormal{Compl}(X_k))_{k}$. By definition of $\textnormal{Compl}(\cdot)$, there exist compactifications  $cX_k$, $k\in\omega$, such that $\Compl{X_k}{cX_k}=f(k)$.
Recall that $K:=\omega+1$ is a (one-point) compactification of $\omega$.
Applying Proposition~\ref{proposition: properties of Z}\,\eqref{case: compactifications of Z(Y,Z)}, we obtain a~compactification of $X$:
\[ cX:=Z(K,(cX_k)_k) .\]
Since $\omega$ is $\sigma$-compact \errata{but not compact}, its complexity in $K$ is 1. We have
\[ \sup_{k\in\omega} \Compl{X_k}{cX_k} \geq 1 \]
(since not every $X_k$ is compact). By Proposition $\ref{proposition: complexity of Z(Y,Z)}$, the~complexity of $X$ in $cX$ is the~maximum of $\textnormal{Compl}(\omega,K)$ and this supremum, so we get
\[ \sup f = \sup_{k\in\omega} \ \Compl{X_k}{cX_k} = \Compl{X}{cX} \in \textnormal{Compl}(X) . \]

"$\subset$": Let $cX$ be a~compactification of $X$ and denote $Y:= \bigcup_{k\in\omega} \overline{X_k}^{cX}$. Define $f:\omega \rightarrow [1,\omega_1]$ by $f(k) := \Compl{X_k}{\overline{X_k}^{cX}}$ and set $\alpha:=\sup f$. Note that each $X_k$ is an $\Fa$-subset of $Y$. Since $Y$ is $\sigma$-compact and $\alpha\geq 1$, it is enough to show that $X$ is an $\Fa$-subset of $Y$.

If $\Fa$ is an additive class, $X=\bigcup X_k$ is a~countable union of $\Fa$ sets, so we have $X\in \Fa (Y)$. If $\Fa$ is a~multiplicative class, we have to proceed more carefully. Let $k,l\in \omega$ be distinct. We set $K^l_k := \overline{X_k}^Y\cap \overline{X_l}^Y$ and note that by definition of topological sum, $K^l_k$ is a~compact disjoint with $X_l$. Since $X_l$ is $\mc K$-analytic, and hence Lindelöf, we can apply \cite[Lemma 14]{spurny2006solution} to obtain a~set $E^l_k\in F_\sigma(Y)$ satisfying
\begin{equation} \label{equation: separation of E from K}
X_l \subset E^l_k \subset \overline{X_l}^Y \setminus K^l_k  = \overline{X_l}^Y \setminus \overline{X_k}^Y .
\end{equation}
Denote $E_k := \bigcup_{l \neq k} E^l_k \in F_\sigma(Y)$. We finish the~proof by starting with \eqref{equation: separation of E from K} and taking the~union over $l\in\omega \setminus \{k\}$ \dots
\begin{align*}
\bigcup_{l \neq k} X_l \ & \subset & \bigcup_{l \neq k} E^l_k = E_k & & \subset &
\ \ \bigcup_{l \neq k} \overline{X_l}^Y \setminus \overline{X_k}^Y = Y \setminus \overline{X_k}^Y \\
&&&&& \dots \text{then adding } X_k \dots \\
X = X_k \cup \bigcup_{l \neq k} X_l \ & \subset & X_k \cup E_k \ \ & & \subset &
\ \ X_k \cup \left( Y \setminus \overline{X_k}^Y \right) = Y \setminus \left( \overline{X_k}^Y \setminus X \right) \\
&&&&& \dots \text{and intersecting over } k\in\omega \\
X\, \ & \subset & \bigcap_{k\in\omega} \Big( X_k \cup E_k & \Big) & \subset & \ \ Y \setminus \bigcup_{k\in\omega} \left( \overline{X_k}^Y \setminus X \right) = X .
 \end{align*}
This proves that $X$ is an intersection of $\Fa$-subsets of $Y$. Since $\Fa$ is multiplicative, this completes the~proof.
\end{proof}

\errata{Assuming that we already have Theorem~\ref{theorem:X_2_beta}, a slight modification of the proof of Proposition~\ref{proposition: complexities of topological sums} shows} that for any closed interval $I\subset [2,\omega_1]$, there exists a~space satisfying $\textnormal{Compl}(X)=I$.
\err{On the one hand, some extra work will be required to deal with the case of uncountable $I$. On the other hand, this also makes our work easier because it makes no sense to worry about the upper bound on complexity -- it suffices to prove that $X$ is $\mc K$-analytic, and that's that.}

\begin{proof}[Proof of Theorem~\ref{theorem: summary}]
Let $I=[\alpha, \beta] \subset [2,\omega_1]$ be a~closed interval.
For $\gamma \in [\alpha,\beta]$, let $X_\alpha^\gamma$ be a~topological space from Theorem~\ref{theorem:X_2_beta} satisfying
\[ \left\{ \alpha,\gamma \right\} \subset \textnormal{Compl}\left( X_\alpha^\gamma \right) \subset [\alpha,\gamma] \]
and denote $\mc X := (X_\alpha^\gamma)_{\gamma\in I}$.
Clearly, the family $\mc X$ satisfies
\[ \left\{ \sup f | \ f \text{ is a~selector of } \left( \textnormal{Compl}(X) \right)_{X\in \mc X} \right\} = I .\]

If $\alpha=\beta = \omega_1$, we simply set $X := X^{\omega_1}_{\omega_1}$.

For $\beta < \omega_1$, the family $\mc X$ is countable. This means that Proposition~\ref{proposition: complexities of topological sums}\, $(iii)$ applies, and implies that the~space $X := \bigoplus \mc X$ has the~desired properties.


Finally, suppose that $\err{\alpha<}\beta=\omega_1$.
\errata{Let $T$ be some space with a~single non-isolated point, such that $T$ is $F_{\sigma\delta}$ in $K := \beta T$ but isn't $\sigma$-compact (see e.g. \cite{talagrand1985choquet} for the existence of such a~$T$).
Clearly, $T$ must have uncountably many isolated points (otherwise it would be $\sigma$-compact).
In particular, we can assume that the~set $I$ is contained in $I_T$.

To summarize, we have a space $T$ and its compactification $K$ such that $X := Z(T,\mc X)$ is well defined and $\textnormal{Compl}(T,K) = 2$.}
By Corollary~\ref{corollary: zoom spaces are K analytic}, $X$ is $\mc K$-analytic but not $\sigma$-compact, which implies that $\textnormal{Compl}(X)\subset [2,\omega_1]$.

To finish the~proof, it is enough to prove the~inclusion
\[ \left\{ \sup f | \ f \text{ is a~selector of } \left( \textnormal{Compl}(X) \right)_{X\in\mc X} \right\} \subset \textnormal{Compl}( Z(T,\mc X) ) .\]
This can be done analogously to the~``$\supset$''-part of the~proof of Proposition~\errata{\ref{proposition: complexities of topological sums}, the only difference being that $\omega$ is replaced by $T$}.
\end{proof}

\section{Simple Representations}\label{section: simple representations}

We begin with an example of different ways of describing and representing of $\fsd$-sets.
This example will be rather trivial; however, it will serve as a motivation for studying higher classes of the~$\mc F$-Borel hierarchy.

Suppose that $Y$ is a topological space and $X\in \mc F_2(Y)$.

By definition of $\mc F_2(Y)$, $X$ can be described in terms of $\mc F_1$-subsets of $Y$ as $X=\bigcap_m X_m$, where $X_m \in \mc F_1(Y)$.
And by definition of $\mc F_1(Y)$, each $X_m$ can be rewritten as $X_m := \bigcup_n X_{m,n}$, where $X_{m,n}\in \mc F_0(Y)$.
We could say that the~fact that $X$ is an $\mc F_2$-subset of $Y$ is represented by the~collection $\mc H := (X_{m,n})_{(m,n)\in\omega^2}$.

For $(n_0,n_1,\dots,n_k) \in \seq$, we can set $C(n_0,n_1,\dots,n_k) := X_{0,n_0} \cap X_{1,n_1} \cap \dots \cap X_{k,n_k}$.
We then have $s\sqsubset t \implies C(t) \subset C(s)$ and for each $m\in \omega$, $\{ C(s) | \ s\in \omega^m \}$ covers $X$.
The fact that $X$ is an $\mc F_2$-subset of $Y$ is then represented by the~collection $\mc C := (C(s))_{s}$, in the~sense that $X= \bigcap_{m\in\omega} \bigcup_{s\in \omega^m} C(s)$.
We denote this as $X = \mathbf{R}_2(\mc C)$.

In both cases, $X$ is represented in terms of closed subsets of $Y$. We can also represent $X$ by collections of subsets of $X$, by working with $\mc H_X := (X_{m,n} \cap X)_{(m,n)\in\omega^2}$ and $\mc C_X := (C(s)\cap X )_{s}$ instead.
$X$ is then represented by the~collection $\mc H' := \overline{\mc H_X}^Y$, resp. as $X = \mathbf{R}_2 (\overline{\mc C_X}^Y) =: \mathbf{R}_2 (\mc C_X, Y)$ (where for a collection $\mc A$ of subsets of $Y$, $\overline{\mc A}^Y$ denotes the~collection of the~corresponding closures in $Y$).

We call the~former representation a \emph{simple $\mc F_2$-representation of $X$ in $Y$}, and the~later representation a \emph{regular $\mc F_2$-representation of $X$ in $Y$}.

The advantage of representing $X$ by a collection of subsets of $X$ is that this collection might be able to represent $X$ in more than just one space $Y$.
This allows us to define a \emph{universal} representation of $X$ (either simple or regular) -- a collection which represents $X$ in every space which contains it.
It turns out that only very special topological spaces admit universal representations, but in Section~\ref{section: broom space properties}, we will see that it is the~case for the~class of the~so-called broom spaces.

This section is organized as follows: In Section~\ref{section:simple_rep_def}, we introduce the~concept of a simple representation and discuss its main properties. In Section~\ref{section:her_lind}, we give a sample application of this concept by providing a simple proof of the~fact that the~complexity of hereditarily Lindelöf spaces is absolute. In Section~\ref{section: universal representation}, we investigate the~concept of local complexity and its relation to standard complexity and universal representations (resp. their non-existence).
Reading this section might make it easier to understand the~intuition behind Section~\ref{section: regular representations}, where regular representations are introduced. However, the~results of this section are in no way formally required by the~subsequent parts of this paper.

\subsection{Definition of Simple Representations}\label{section:simple_rep_def}

This section will make an extensive use of the~notions introduced in Section~\ref{section: sequences}, in particular those related to the~leaf-derivative on trees.
We will repeatedly use collections of sets which are indexed by leaves of some tree:

\begin{definition}[Leaf-scheme]\label{definition: leaf scheme}
A collection $\mc H = (\err{\mc H(t)})_{t\in l(T)}$ \err{(of subsets of a topological space $Y$)}, where $T \in \textnormal{WF}$, is said to be a~\emph{leaf-scheme} \errata{(in $Y$)}.
\end{definition}

\noindent When we need to specify the~indexing set, we will say that $\mc H$ is an \emph{$l(T)$-scheme}.
Saying that a~$l(T)$-scheme is \emph{closed in $Y$}, for a~topological space $Y$, shall mean that $\mc H(t)$ is closed in $Y$ for each $t\in l(T)$.

Any $l(T)$-scheme $\mc H$ can be viewed as a~mapping with domain $l(T)$. $\mc H$ has a~natural extension to $T$, defined by the~following recursive formula:
\begin{equation} \label{equation: H t extensions}
\mc H(t) := \begin{cases}
\mc H(t) & \text{when } t \text{ is a~leaf of }T \\
\bigcup \{ \mc H(s) | \ s\in \ims{T}{t} \} & \text{when } r_l(T^t) > 0 \text{ is odd} \\
\bigcap \{ \mc H(s) | \ s\in \ims{T}{t} \} & \text{when } r_l(T^t) > 0 \text{ is even}.
\end{cases}
\end{equation}

It should be fairly obvious that closed leaf-schemes in any topological space $Y$ \emph{naturally} correspond to $\mc F$-Borel subsets of $Y$, in the~sense that a~set $X\subset Y$ is $\mc F$-Borel in $Y$ if and only if it is of the~form $X=\mc H(\emptyset)$ for some closed leaf-scheme $\mc H$ in $Y$.
For now, we only prove one direction of this claim (the other will be shown in Proposition~\ref{proposition: existence of simple representations}).

\begin{lemma}[Complexity of sets corresponding to leaf-schemes]\label{lemma: complexity of S T}
\errata{Let $\alpha<\omega_1$. For any closed $l(T)$-scheme $\mc H$ in a~topological space $Y$, we have
\begin{equation}\label{equation: complexity of H t}
\left( \forall t\in T \right) : r_l(T^t)=\alpha \implies \mc H(t) \in \Fa (Y) .
\end{equation}}
In particular, if $\mc H$ is a~closed $l(T_\alpha)$-scheme in $Y$, then $\mc H(\emptyset) \in \Fa(Y)$.
\end{lemma}

\begin{proof}
The ``in particular'' part is a~special case of the~first part because $\err{r_l(T_\alpha^\emptyset)=}$ $r_l(T_\alpha)=\alpha$.
Therefore, it suffices to use induction over $\alpha$ to prove \eqref{equation: complexity of H t}.
For $\alpha = 0$, the~$t$-s for which $r_l(T^t)=0$ are precisely the~leaves of $T$. For these we have $\mc H(t) \in \mc F_0(Y)$ because $\mc H$ is closed in $Y$.

$<\alpha \mapsto \alpha$:
Suppose that \eqref{equation: complexity of H t} holds for every $\beta < \alpha$, and $t\in T$ satisfies $r_l(T^t)= \alpha$.
By \eqref{equation: leaf rank formula}, we have $r_l(T^s) < r_l(T^t) = \alpha$ for every $s\in \ims{T}{t}$.
It then follows from the~induction hypothesis that each $s\in \ims{T}{s}$, $\mc H(s)$ belongs to $\mc F_{<\alpha}(Y)$.
When $\alpha$ is odd, we have
\[ \mc H(t) \overset{\alpha \text{ is}}{\underset{\text{odd}}=}
\bigcup \left\{ \mc H(s) \, | \ s\in\ims{T}{s} \right\} \in (\mc F_{<\alpha}(Y))_\sigma
\overset{\alpha \text{ is}}{\underset{\text{odd}}=} \Fa(Y) .\]
When $\alpha$ is even, we have 
\[ \mc H(t) \overset{\alpha \text{ is}}{\underset{\text{even}}=}
\bigcap \left\{ \mc H(s) \, | \ s\in\ims{T}{s} \right\} \in (\mc F_{<\alpha}(Y))_\delta
\overset{\alpha \text{ is}}{\underset{\text{even}}=} \Fa(Y) .\]
\end{proof}

\bigskip
Lemma~\ref{lemma: complexity of S T} motivates the~following definition of (simple) $\Fa$-representation:

\begin{definition}[Simple representations]\label{definition: simple representation}
Let $X\subset Y$ be topological spaces, and $\mc H$ an $l(T)$-scheme in $Y$. $\mc H$ is said to be a
\begin{itemize}
\item \emph{simple representation} of $X$ in $Y$ if it satisfies $\mc H(\emptyset)=X$;
\item \emph{simple $\mc F$-Borel-representation} of $X$ in $Y$ if it is both closed in $Y$ and a~simple representation of $X$ in $Y$;
\item \emph{simple $\Fa$-representation} of $X$ in $Y$, for some $\alpha<\omega_1$, if it is a~simple $\mc F$-Borel-representation of $X$ in $Y$ and we have $r_l(T)\leq \alpha$.
\end{itemize}
\end{definition}

\noindent Unless we need to emphasize that we are not talking about the~regular representations from the~upcoming Section~\ref{section: regular representations}, we omit the~word `simple'.
We note that the~following stronger version of Lemma~\ref{lemma: complexity of S T} also holds, which allows for a~different, but equivalent, definition of a~simple $\Fa$-representation:

\begin{remark}[Alternative definition of $\Fa$-representations]
In Lemma~\ref{lemma: complexity of S T}, we even get $\mc H(\emptyset) \in \mc F_{r_i(T)}(Y)$.
\end{remark}

\noindent Consider the~derivative $D_i$ from Definition~\ref{definition: derivatives}  which cuts away sequences which only have finitely many extensions, and the~corresponding rank $r_i$.
Using the~fact that the~class of closed sets (and of $\Fa$-sets as well) is stable under finite unions and intersections, we can prove that when a~tree $T$ is finite and we have $\mc H(t)\in \Fa(Y)$ for each $t\in l(T)$, then $\mc H(\emptyset) \in \Fa(Y)$.
It follows that the~conclusion of Lemma~\ref{lemma: complexity of S T} also holds when we  replace `$r_l$' by `$r_i$'.
In particular, we could equally well define $\Fa$-representations as those $\mc F$-Borel representations which are indexed by trees which satisfy $r_i(T)\leq \alpha$.
We will mostly work with trees on which the~ranks $r_l$ and $r_i$ coincide, so the~distinction will be unimportant.

\bigskip

Before proceeding further, we remark that as long as the~intersection with $X$ is preserved, cutting away some pieces of a~representing leaf-scheme still yields a~representation of $X$:

\begin{lemma}[Equivalent representations]\label{lemma: leaf scheme making smaller}
Let $X\subset Y$ be topological spaces, $\mc H$ and $\mc H'$ two $l(T)$-schemes in $Y$, and suppose that $\mc H$ is a~representation of $X$ in $Y$.

If $\mc H'$ satisfies $\mc H(t) \cap X \subset \mc H'(t) \subset \mc H(t)$ for every $t\in l(T)$, then $\mc H'$ is also a~representation of $X$ in $Y$.
\end{lemma}

\begin{proof}
Let $X$, $Y$, $T$, $\mc H$ and $\mc H'$ be as above.
By the~assumptions of the~lemma, we have the~following formula for every $t\in l(T)$:
\begin{equation}\label{equation: H and H'}
\mc H(t) \cap X \subset \mc H'(t) \ \ \ \  \& \ \ \ \ \mc H'(t) \subset \mc H(t).
\end{equation}
Using \eqref{equation: H t extensions} and induction over $r_l(T^t)$, we get \eqref{equation: H and H'} for every $t\in T$.
Applying \eqref{equation: H and H'} to $t=\emptyset$, we get the~desired conclusion of the~lemma:
\[ X = \mc H(\emptyset) = \mc H(\emptyset) \cap X \subset \mc H' (\emptyset) \subset
\mc H(\emptyset) = X .\]
\end{proof}

The existence of simple representations is guaranteed by the~following result:

\begin{proposition}[Existence of $\mc F$-Borel representations]\label{proposition: existence of simple representations}
For topological spaces $X\subset Y$ and $\alpha<\omega_1$, the~following conditions are equivalent.
\begin{enumerate}[(i)]
\item $X \in \Fa(Y)$
\item $X$ has a~simple $\Fa$-representation in $Y$; that is, $\mc H(\emptyset) = X$ holds for some closed $l(T)$-scheme $\mc H$ in $Y$, where $r_l(T)\leq \alpha$.
\item $\mc H(\emptyset) = X$ holds for some closed $l(T_\alpha)$-scheme $\mc H$ in $Y$.
	\newcounter{nameOfYourChoice}
	\setcounter{nameOfYourChoice}{\value{enumi}}
\end{enumerate}
Denoting $\overline{\mc H}^Y\!\! := \left( \overline{\mc H(t) \cap Y}^Y\!\! \right)_{t\in l(T)}$ for any $l(T)$-scheme $\mc H$, this is further equivalent to:
\begin{enumerate}[(i)]
	\setcounter{enumi}{\value{nameOfYourChoice}}
\item $\overline{\mc H}^Y\!\!(\emptyset) = X$ holds for some $l(T)$-scheme $\mc H$ in $X$, where $r_l(T)\leq \alpha$.
\item $\overline{\mc H}^Y\!\!(\emptyset) = X$ holds for some $l(T_\alpha)$-scheme $\mc H$ in $Y$.
\end{enumerate}
\end{proposition}

\begin{proof}
Since the~leaf-schemes $\overline{\mc H}^Y$ in (iv) and (v) are closed in $Y$, it follows from Lemma~\ref{lemma: complexity of S T} that any of the~conditions (ii)-(v) implies (i).
The implications (iii)$\implies$(ii) and (v)$\implies$(iv) are trivial.
By Lemma~\ref{lemma: leaf scheme making smaller}, we have (ii)$\implies$(iv) and (iii)$\implies$(v) (set $\mc H_X(t) := \mc H(t) \cap X$ and $\mc H'(t) := \overline{\mc H_X}^Y\!\!$ in Lemma~\ref{lemma: leaf scheme making smaller}).
Therefore, it remains to prove (i)$\implies$(iii).
We will first prove (i)$\implies$(ii), and then note how to modify the~construction such that it gives (i)$\implies$(iii).

(i)$\implies$(ii):
We proceed by induction over $\alpha$.
For $\alpha=0$, let $X$ be a~closed subset of $Y$. We set $\mc H(\emptyset) := X$, and observe that such $\mc H$ is a~closed $l(S)$-scheme in $Y$, where $S = \{\emptyset\}$.

$\alpha \mapsto \alpha+1$ for odd $\alpha+1$:
Suppose that $X = \bigcup_{n\in \omega} X_n$, where $X_n \in \Fa(Y)$.
For each $n\in \omega$, we apply the~induction hypothesis to obtain a~closed $l(S_n)$-scheme $\mc H_n$ such that $\mc H_n(\emptyset) = X_n$ and $r_l(S_n) = \alpha$.
For $n\in\omega$ and $s\in l(S_n)$, we set $\mc H (n\ext s) := \mc H_n(s)$. This defines a~leaf-scheme indexed by the~leaves of the~tree $S := \{\emptyset\} \cup \bigcup_n n\ext S_n$, and the~formula $\mc H(n\ext s) = \mc H_n(s)$ obviously holds for every $n\ext s \in S$.
By \eqref{equation: leaf rank formula}, we have $r_l(S) = \sup_n r_l(S_n) +1 = \alpha+1$.
Since $\alpha+1$ is odd, it follows that
\[ \mc H(\emptyset) = \bigcup_n \mc H(\emptyset \ext n) = \bigcup_n \mc H(n \ext \emptyset ) = \bigcup_n \mc H_n(\emptyset) = \bigcup_n X_n = X .\]

$<\alpha \mapsto \alpha$ for even $\alpha$:
Suppose that $X = \bigcap_{n\in \omega} X_n$, where $X_n \in \mc F_{\alpha_n}(Y)$ holds for some $\alpha_n < \alpha$.
Since $\mc F_\beta (Y) \subset \mc F_\gamma(Y)$ holds whenever $\beta$ is smaller than $\gamma$, we can assume that either $\alpha$ is a~successor and we have $\alpha_n = \alpha-1$ for each $n$, or $\alpha$ is limit and we have $\sup_n \alpha_n = \alpha$. In either case, we have $\sup_n ( \alpha_n +1 ) = \alpha$.

The remainder of the~proof proceeds as in the~previous case -- applying the~induction hypothesis, defining an $l(S)$-scheme as $\mc H(n\ext s) := \mc H_n(s)$, and observing that $\mc H(\emptyset) = \bigcap_n X_n$ (because $r_l(S) = \alpha$ is even).

(i)$\implies$(iii):
We prove by induction that in each induction-step from the~proof of '(i)$\implies$(ii)', the~tree $S$ can be of the~form $S=T_\alpha$.
Recall here the~Notation~\ref{notation: trees of height alpha} (which introduces these trees).

In the~initial step of the~induction, we have $S=\{\emptyset\}=T_0$ for free.
Similarly for the~odd successor step, we have $S=T_{\alpha+1}$ provided that each $S_n$ is equal to $T_\alpha$.

Consider the~``$<\alpha \mapsto \alpha$ for even $\alpha$'' step.
By the~remark directly below Notation~\ref{notation: trees of height alpha}, we have $T = \{\emptyset\} \cup \bigcup_n n\ext T_{\pi_\alpha(n)}$. Therefore, it suffices to show that we can assume $\alpha_n \leq \pi_\alpha(n)$ for each $n\in\omega$.
By replacing the~formula $X = \bigcap_n X_n$ by
\[ X = X_0 \cap \overline{X}^Y \cap X_1 \cap \overline{X}^Y \cap X_2 \cap \dots ,\]
we can guarantee that $\alpha_n = 0$ holds for infinitely many $n\in\omega$.
Re-enumerating the~sets $X_n$, we can ensure that $\alpha_n \leq \pi_\alpha(n)$ holds for each $n\in\omega$.
 \end{proof}

The conditions (iv) and (v) motivate the~definition of a~universal simple representation:

\begin{definition}[Universal representation]\label{definition: universal simple representation}
A leaf-scheme $\mc H$ in $X$ is said to be a~\emph{universal simple $\mc F$-Borel-representation of $X$} if for every $Y\supset X$, $\overline{\mc H}^Y\!\!$ is a~simple $\mc F$-Borel-representation of $X$ in $Y$.
\end{definition}

Clearly, compact and $\sigma$-compact spaces have universal representations.
In Section~\ref{section: broom space properties} \errata{(Proposition~\ref{proposition: basic broom complexities})}, we show an example of more complex spaces with universal representation.
In Section~\ref{section: universal representation} \errata{(Proposition~\ref{proposition: non-existence of universal repre})}, we show that these spaces are rather exceptional, because sufficiently topologically complicated spaces do not admit universal  representations.

\subsection{The Complexity of Hereditarily Lindelöf Spaces}\label{section:her_lind}

We note that simple representations can be used not only with the $\mc F$-Borel class, but also be with other descriptive classes:

\begin{remark}\label{remark: simple C representations}
Simple representations based on $\mc C$-sets.
\end{remark}

\noindent Let $Y$ be a topological space and $\mc C(Y)$ a family of subsets of $Y$.
Analogously to Definition \ref{definition: F Borel sets}, we define
\begin{itemize}
\item $\mc C_0(Y) := \mc C(Y)$,
\item $\mc C_\alpha(Y) := \left( \mc C_{<\alpha} (Y) \right)_\sigma$ for $0 < \alpha <\omega_1$ odd,
\item $\mc C_\alpha(Y) := \left( \mc C_{<\alpha} (Y) \right)_\delta$ for $0 < \alpha <\omega_1$ even.
\end{itemize}

Analogously to \ref{definition: simple representation}, an $l(T)$-scheme $\mc H$ which represents $X$ in $Y$ is said to be a \emph{simple $\mc C_\alpha$-representation of $X$ in $Y$} when the tree $T$ satisfies $r_l(T)\leq \alpha$ and we have $\mc H(t) \in \mc C(Y)$ for each $t\in l(T)$.

Clearly, the equivalence (i)$\iff$(ii)$\iff$(iii) in Proposition~\ref{proposition: existence of simple representations} holds for $\mc C_\alpha$-sets and $\mc C_\alpha$-representations as well. (The proof literally consists of replacing `$\mc F$' by `$\mc C$' in Lemma \ref{lemma: complexity of S T} and in the proof of Proposition~\ref{proposition: existence of simple representations}.)

\bigskip

To prove the main result of this subsection, we need the following separation lemma, which is an immediate result of \cite[Lemma\,14]{spurny2006solution}.

\begin{lemma}[$F_\sigma$-separation for Lindelöf spaces]\label{lemma: fsd spaces are Fs separated}
Let $L$ be a Lindelöf subspace of a compact space $C$. Then for every compact set $K\subset C \setminus L$, there exists $H\in \mc F_\sigma\left(C\right)$, such that $L\subset H\subset C\setminus K$.
\end{lemma}

Recall that a topological space is said to be \emph{hereditarily Lindelöf} if its every subspace is Lindelöf. For example, separable metrizable spaces (or more generally, spaces with countable weight) are hereditarily Lindelöf. The following proposition shows that the $\mc F$-Borel complexity of such spaces is absolute:

\begin{proposition}[J. Spurný and P. Holický]\label{proposition: hereditarily lindelof spaces are absolute}
For a hereditarily Lindelöf space and $\alpha \leq \omega_1$,
the following statements are equivalent:
\begin{enumerate}[(i)]
\item $X \in \Fa(Y)$ holds for every $Y\supset X$;
\item $X \in \Fa(cX)$ holds for some compactification $cX$.
\end{enumerate}
\end{proposition}

In \cite{kalenda2018absolute} the author of the present paper, together with his supervisor, proved Proposition~\ref{proposition: hereditarily lindelof spaces are absolute} for the class of $\fsd$ sets ($\alpha=2$). J. Spurný later remarked that it should be possible to give a much simpler proof using the fact that the classes originating from the \emph{algebra} generated by closed sets are absolute (\cite{holicky2003perfect}). Indeed, this turned out to be true, and straightforward to generalize for $\mc F_n$, $n\in\omega$. However, the proof in the general case of $\Fa$, $\alpha \in \omega_1$, would not be very elegant, and can be made much easier by the use of simple representations.

\begin{proof}
For $\alpha=0$, $\alpha=1$, and $\alpha=\omega_1$, this follows from the fact that compact, $\sigma$-compact and $\mc K$-analytic spaces are absolute.

To prove (ii)$\implies$(i), it suffices to show that (ii) implies that $X$ is $\Fa$ in every compactification (by, for example, \cite[Remark\,1.5]{kovarik2018brooms}).

Let $\alpha \in [2,\omega_1)$, suppose that $X\in \Fa(cX)$ holds for some compactification $cX$, and let $dX$ be another compactification of $X$.
Since we have
\[ \mc F(cX) \subset \{ F\cap G | \ F\subset cX \text{ is closed, $G\subset cX$ is open} \}  =: (\mc F\land \mc G)(cX) ,\]
it follows that $X \in (\mc F\land \mc G)_\alpha (cX)$.
The classes $(\mc F\land \mc G)_\alpha$ are absolute by \cite[Corollary 14]{holicky2003perfect}, so we have $X \in (\mc F\land \mc G)_\alpha(dX)$.
By Remark \ref{remark: simple C representations}, $X$ has a simple $(\mc F\land \mc G)_\alpha$-representation in $dX$ -- that is, there exists an $l(T_\alpha)$-scheme $\mc H$ in $dX$, such that $X = \mc H(\emptyset)$, and we have $\mc H(t)\in (\mc F \land \mc G)(dX)$ for each $t\in l(T_\alpha)$.

Let $t\in l(T_\alpha)$. \err{We can express} $\mc H(t)$ as $\mc H(t) = F_t \cap G_t$, where $F_t$ is closed in $dX$ and $G_t$ is open in $dX$.
Since $X$ is hereditarily Lindelöf, the intersection $X\cap G_t$ is Lindelöf.
Applying Lemma \ref{lemma: fsd spaces are Fs separated} to ``$L$''$=X\cap G_t$ and ``$K$''$=cX\setminus G_t$, we get a set $H_t \in \mc F_\sigma(dX)$ s.t. $X \cap G_t \subset H_t \subset G_t$.
Denote $\mc H'(t) := F_t \cap H_t$.
$\mc H'$ is an $l(T_\alpha)$-scheme in $dX$, such for each $t\in l(T_\alpha)$, we have $\mc H'(t) \in \mc F_\sigma(dX)$ and $X\cap \mc H(t) \subset \mc H'(t) \subset \mc H(t)$.

It follows that $\mc H'(\emptyset) \in (\mc F_\sigma)_\alpha(dX)$ (Remark \ref{remark: simple C representations}) and $\mc H'(\emptyset) = \mc H(\emptyset) = X$ (Lemma \ref{lemma: leaf scheme making smaller}).
Since $(\mc F_\sigma)_1(\cdot)$ coincides with $\mc F_1(\cdot)$, it follows that $(\mc F_\sigma)_\alpha(dX) = \Fa(dX)$.
This shows that $X \in \Fa(dX)$.
\end{proof}

\subsection{Local Complexity and Spaces with No Universal Representation}
	\label{section: universal representation}

It is evident that only a~space which is absolutely $\Fa$ stands a~chance of having a~universal (simple) $\Fa$-representation. We will show that even among such spaces, having a~universal representation is very rare.
(Obviously, this does not include spaces which are $F_\sigma$ in some compactification -- such spaces are $\sigma$-compact, and the~representation $X = \bigcup_n K_n$ works no matter where $X$ is embedded.)
To obtain this result, we first discuss the~notion of local complexity.

Let $X$ be a~$\mc K$-analytic space, $Y$ a~topological space containing $X$, and $y\in Y$.
Consider \err{the collection $\mc U(y)$ of all open neighborhoods of $y$ in $Y$, partially ordered by inclusion}, and the~corresponding sets of the~form $\overline{U}^Y\!\! \cap X$, \err{$U\in \mc U(y)$}.
\err{Since $\overline{U}^Y\!\! \cap X$ is a closed subset of $X$}, we have $\err{\alpha_U :=} \, \Compl{\overline{U}^Y \!\! \cap X}{Y} \leq \Compl{X}{Y}$.
\err{Similarly, if $U\subset V$ are two neighborhoods of $y$, then we have $\alpha_U \leq \alpha_V$.}
It follows that \errata{the net $(\alpha_U)_{U\in\mc U(y)}$ converges to some $\alpha \in [0,\omega_1]$.}
And since the~ordinal numbers are well-ordered, the~limit is attained for some neighborhood $V$ of $y$ (and thus for all $U\subset V$ as well).
We shall call this limit the~local complexity of $X$ in $Y$ at $y$:

\begin{definition}[Local complexity] \label{definition: local complexity}
Let $X\in \mc F_{\omega_1}(Y)$ and $y\in Y$. The \emph{local complexity of $X$ in $Y$ at $y$} is the~ordinal
\[ \textnormal{Compl}_y(X,Y) := \min \left\{ \Compl{\overline{U}^Y\!\!\cap X}{Y} | \ 
 U \textnormal{ is a~neighborhood of $y$ in $Y$} \right\}. \]
We define the \emph{local complexity of $X$ in $Y$} as
\begin{equation}\label{equation:local_complexity}
\textnormal{Compl}_\textnormal{loc}(X,Y) := \sup_{x\in X} \textnormal{Compl}_x(X,Y) .
\end{equation}
We also set $\textnormal{Compl}_{\overline{\textnormal{loc}}}\,(X,Y) := \sup_{y\in Y} \textnormal{Compl}_y(X,Y)$.
\end{definition}

\noindent Note that for $x\in X$, we can express the~local complexity of $X$ in $Y$ at $x$ as
\[ \textnormal{Compl}_x(X,Y) = \min \left\{ \Compl{\overline{U}^X\!\!}{Y} | \ 
 U \textnormal{ is a~neighborhood of $x$ in $X$} \right\}. \]
 Also, if the local complexity of $X$ in $Y$ is a non-limit ordinal, then the supremum in \eqref{equation:local_complexity} is in fact a maximum.

We have the~following relation between the~local and ``global'' complexity:

\begin{lemma}[Complexity and local complexity]\label{lemma: local and global complexity}
For any $\mc K$-analytic space $X$ and any $Y\supset X$, we have 
\begin{equation}\label{equation:loc_and_glob_compl}
\textnormal{Compl}_{\textnormal{loc}}\,(X,Y) \leq \textnormal{Compl}_{\overline{\textnormal{loc}}}\,(X,Y) \leq \textnormal{Compl}(X,Y) \leq \textnormal{Compl}_\textnormal{loc}(X,Y) + 1 .
\end{equation}
Moreover, if either $\textnormal{Compl}_\textnormal{loc}(X,Y)$ is odd, or $\textnormal{Compl}(X,Y)$ is even, we get
\[ \textnormal{Compl}_{\textnormal{loc}}\,(X,Y) = \textnormal{Compl}_{\overline{\textnormal{loc}}}\,(X,Y) = \textnormal{Compl}(X,Y) .\]
When $Y$ is compact, we have
\[ \textnormal{Compl}(X,Y) = \textnormal{Compl}_{\overline{\textnormal{loc}}}\,(X,Y) = \max_{y\in Y} \textnormal{Compl}_y(X,Y) .\]
\end{lemma}

\begin{proof}
Let $X$ be a $\mc K$-analytic space and $Y\supset X$, and denote $\gamma := \textnormal{Compl}_{\textnormal{loc}}\,(X,Y)$ and $\alpha := \textnormal{Compl}(X,Y)$.

The first inequality in \eqref{equation:loc_and_glob_compl} is trivial.
The second follows from the~fact that the~complexity of $\overline{U}^Y\!\! \cap X$ in $Y$ is at most the~maximum of the~complexity of $\overline{U}^Y\!\!$ in $Y$ (that is, $0$) and $\Compl{X}{Y}$.
The last one immediately follows from the~following claim:

\begin{claim}
$X$ can be written as a countable union of $\mc F_\gamma$-subsets of $Y$.
\end{claim}

To prove the~claim, assume that for each $x\in X$ there exists an open neighborhood $U_x$ of $x$ in $X$ which satisfies $\overline{U_x}^X \in \mc F_{\gamma}(Y)$.
Since $X$ is $\mc K$-analytic, it is in particular Lindelöf (\cite[Proposition\,3.4]{kkakol2011descriptive}).
Let $\{U_{x_n} | \ n\in\omega \}$ be a~countable subcover of $\{ U_x | \ x\in X \}$.
It follows that $X = \bigcup_{n\in\omega} U_{x_n} = \bigcup_{n\in\omega} \overline{U_{x_n}}^X$ is a~countable union of $\mc F_{\gamma}(Y)$ sets, and hence $X \in (\mc F_\gamma(Y) )_\sigma $.

Next, we prove the ``moreover'' part of the lemma.
When $\gamma$ is odd, we have $(\mc F_\gamma(Y))_\sigma = \mc F_\gamma$, so the claim yields $\textnormal{Compl}(X,Y) \leq \gamma$ (and the conclusion follows from \errata{\eqref{equation:loc_and_glob_compl}}).

Suppose that $\alpha$ is even and assume for contradiction that $\gamma < \alpha$.
The claim implies that $X$ belongs to the additive $\mc F$-Borel class $(\mc F_\gamma(Y) )_\sigma$. Since we have $\gamma < \alpha$ and the class $\Fa(Y)$ is multiplicative, it follows that $(\mc F_\gamma(Y) )_\sigma \subsetneq \Fa(Y)$. This contradicts the definition of $\textnormal{Compl}(X,Y)$ (which states that $\alpha$ is the smallest ordinal satisfying $X\in \Fa(Y)$).

Lastly, suppose that $Y$ is compact and denote $\eta := \textnormal{Compl}_{\overline{\textnormal{loc}}}\,(X,Y)$.
For each $y\in Y$, let $U_y\subset Y$ be an open neighborhood of $y$ in $Y$ s.t. $\overline{U}^Y\cap X$ is an $\mc F_\eta$-subset of $Y$.
Since $Y$ is compact, we have $Y = U_{y_0} \cup \dots \cup U_{y_k}$ for some $y_i \in Y$, $k\in\omega$.
It follows that $X$ can be written as $X = \bigcup_{i=0}^k (\overline{U_{y_i} }^Y \cap X)$ of finitely many $\mc F_\eta$-subsets of $Y$.
This shows that $X\in \mc F_\eta(Y)$, and thus $\textnormal{Compl}(X,Y) \leq \gamma$. Since the supremum in $\textnormal{Compl}_{\overline{\textnormal{loc}}}\,(X,Y)$ is equal to the supremum of finitely many numbers $\textnormal{Compl}_{y_i}\,(X,Y)$, it is in fact a maximum.
\end{proof}

Note that \eqref{equation:loc_and_glob_compl} is optimal in the sense that when $\textnormal{Compl}(X,Y)$ is odd, it could be equal to either of the numbers $\textnormal{Compl}_{\textnormal{loc}}\,(X,Y)$ and $\textnormal{Compl}_\textnormal{loc}(X,Y) + 1$.
Indeed, consider the examples $X_1 := (0,1) \subset [0,1] =: Y$ and $X_2 = \Q \cap [0,1] \subset [0,1] = Y$.

In Section~\ref{section: broom space properties}, we construct absolute $\Fa$ spaces $\mathbf{T}_\alpha$, each of which has a~special point $\infty \in \mathbf{T}_\alpha$. For any $Y\supset \mathbf{T}_\alpha$, the~local complexity of $\mathbf{T}_\alpha$ in $Y$ is 0 at every $x\in \mathbf{T}_\alpha \setminus \infty$ and 1 at every $y \in Y\setminus \mathbf{T}_\alpha$. (This shall immediately follow from the~fact that the only non-isolated point of $X$ is $\infty$ and any closed subset of $\mathbf{T}_\alpha$ which does not contain $\infty$ is at most countable.)
This in particular shows that in general, it is not useful to consider the~supremum of $\sup_y \textnormal{Compl}_y(X,Y)$ over $Y\setminus X$.

We will use the~following local variant of the~notion of descriptive ``hardness'':

\begin{definition}[local $\Fa$-hardness]\label{definition: Fa hardness}
Let $X\subset Y$ be topological spaces s.t. $X\in \mc F_{\omega_1}(Y)$.
For $\alpha\leq \omega_1$, $X$ is said to be
\begin{itemize}
\item \emph{$\Fa$-hard in $Y$ at $y$}, for some $y\in Y$, if we have $\textnormal{Compl}_y\left(X,Y\right) \geq \alpha$;
\item \emph{locally $\Fa$-hard in $Y$} if it is $\Fa$-hard in $Y$ at every $x\in X$.
\end{itemize}
\end{definition}

\noindent Note that $X$ is locally $\Fa$-hard in $Y$ if and only if we have $\overline{U}^X \notin \mc F_{<\alpha}(Y)$ for every open subset $U$ of $X$ (or equivalently, when $F\notin \mc F_{<\alpha}(Y)$ holds for every regularly closed\footnote{Recall that a~set $F$ is regularly closed if it is equal to the~closure of its interior.} subset of $X$).

Note that if $X$ is locally $\Fa$-hard in $Y$, then it is $\Fa$-hard in $Y$ at every point of $\overline{X}^Y\setminus X$.

\begin{lemma}[Density of ``too nice'' sets]\label{lemma: Fa hardness and density}
Let $X\subset Y$ be topological spaces. Suppose $X\subset H \in \mc F_{<\alpha}(Y)$ holds for some $H$ and $\alpha \leq \omega_1$.
\begin{enumerate}[(i)]
\item If $X$ is dense in $Y$ and locally $\Fa$-hard in $Y$, then $H\setminus X$ is dense in $Y$.
\item More generally, if $X$ is $\Fa$-hard in $Y$ at some $y\in Y$, then $y\in \overline{H\setminus X}^Y$.
\end{enumerate}
\end{lemma}

In particular, (i) says that if a~space $X$ is locally hard in $Y$, then any set $H\supset X$ of lower complexity than $X$ must be much bigger than $X$, in the~sense that the~closure of $H\setminus X$ is the~same as that of $X$.

\begin{proof}
(i): By density of $X$ in $Y$ (and the~remark just above this lemma), (i) follows from (ii).

(ii): Suppose we have $X\subset H \subset Y$, $H \in \mc F_{<\alpha}(Y)$, and that $X$ is $\Fa$-hard in $Y$ at some $y\in Y$.
Let $U$ be an open neighborhood of $y$ in $Y$.
Let $V$ be some open neighborhood of $y$ in $Y$ which satisfies $\overline{V}^Y \subset U$.
Since $\textnormal{Compl}_y\left(X,Y\right) \geq \alpha$, we have $\overline{V}^Y \cap X \notin \mc F_{<\alpha}(Y)$.
On the~other hand, the~assumptions imply that
\begin{equation}\label{equation: Fa hardness and density}
\overline{V}^Y \!\! \cap X \subset \overline{V}^Y \!\! \cap H \in \mc F_{<\alpha}(Y) .
\end{equation}
Since the~two sets in \eqref{equation: Fa hardness and density} have different complexities, they must be distinct. In particular, the~intersection $\overline{V}^Y \!\! \cap (H\setminus X) \subset U \cap (H\setminus X)$ is non-empty.
Since $U$ was arbitrary, this shows that $y$ is in the~closure of $H\setminus X$.
\end{proof}

Before proceeding further, make the~following simple observation:

\begin{lemma}[$G_\delta$-separation of $\mc F$-Borel sets]\label{lemma: G delta separation}
Let $X$ be \err{an}~$\mc F$-Borel subset of $Y$. Then for every $y\in Y\setminus X$, there is a~$G_\delta$ subset $G$ of $Y$ which satisfies $\err{y}\in G \subset Y \setminus X$.
\end{lemma}

\begin{proof}
We proceed by transfinite induction over the~complexity of $X$ in $Y$.
When $X$ is closed in $Y$, the~conclusion holds even with the open set $G:= Y \setminus X$.

Suppose that $X=\bigcup_n X_n$ and $y\in Y \setminus X$. If each $X_n$ satisfies the~conclusion -- that is, if we have $y\in G_n \subset Y \setminus X_n$ for some $G_\delta$ subsets $G_n$ of $Y$ -- then the~set $\bigcap_n G_n$ is $G_\delta$ in $Y$ and we have
\[ y \in \bigcap_n G_n \subset \bigcap_n (Y \setminus X_n) = Y \setminus \bigcup_n X_n = Y \setminus X .\]

Suppose that $X = \bigcap_n X_n$ and let $y\in Y \setminus X$. Let $n_0\in\omega$ be s.t. $y \err{\in} Y \setminus X_n$. By the induction hypothesis, there is some $G_\delta$ set $G$ which satisfies the~conclusion for $X_{n_0}$. Clearly, $G$ satisfies the~conclusion for $X$ as well.
\end{proof}

We are now ready to prove that the~existence of universal representation is rare.
Recall that a point $x$ is a~cluster point of a sequence $(x_n)_n$ when each neighborhood of $x$ contains $x_n$ for infinitely many $n$-s.

\begin{proposition}[Non-existence of universal representations]
	\label{proposition: non-existence of universal repre}
Let $X$ be a~$\mc K$-analytic space which is not $\sigma$-compact.
Suppose that $X$ has a~compactification $cX$, such that one of the~following conditions holds for some \emph{even} $\alpha < \omega_1$:
\begin{enumerate}[(a)]
\item $X$ is locally $\Fa$-hard in $cX$;
\item there is a~sequence $(x_n)_{n\in\omega}$ in $X$ with no cluster point in $X$, such that $X$ is $\Fa$-hard in $cX$ at each $x_n$;
\item $X$ is $\Fa$-hard in $cX$ at some point of $cX\setminus X$.
\end{enumerate}
Then $X$ does not have a~universal simple $\Fa$-representation.
\end{proposition}

\begin{proof}
Note that both (a) and (b) imply (c).
Indeed, by the~remark between Definition~\ref{definition: Fa hardness} and Lemma~\ref{lemma: Fa hardness and density}, (a) implies that $X$ is $\Fa$-hard at each point of the~(non-empty) set $cX\setminus X$.
If (b) holds, then $(x_n)_n$ must have some cluster point $y$ in $cX\setminus X$. And since each neighborhood of $y$ is a~neighborhood of some $x_n$, $X$ must be $\Fa$-hard in $cX$ at $y$:
\[ \textnormal{Compl}_y(X,cX) \geq \inf_{n\in\omega} \textnormal{Compl}_{x_n}(X,cX) \geq \alpha .\]

It remains to show that (c) implies that the~universal $\Fa$-representation of $X$ does not exist.
Suppose that $X$ is $\Fa$-hard in $cX$ at some $y\in cX\setminus X$.
By Lemma~\ref{lemma: G delta separation}, there is some $G_\delta$ subset $G$ of \errata{$cX$}, such that $y\in G \subset cX \setminus X$.
Since $cX$ is \err{compact}, we can assume that $G = \bigcap_n G_n$, where $G_n$ are open subsets of $cX$ which satisfy $\overline{G_{n+1}} \subset G_n$.

Let $\mc H$ be a~``candidate for a~universal simple representation of $X$'' -- that is, an $l(T_\alpha)$-scheme in $X$ satisfying $X = \overline{\mc H}^{cX}\!\!(\emptyset)$.
We shall prove that $\mc H$ is not a~universal representation by constructing a~compactification $dX$ in which $\overline{\mc H}^{dX}\!\!$ is not a representation of $X$ (that is, we have $\overline{\mc H}^{dX}\!\!(\emptyset) \neq X$).

Since $\alpha$ is even, we get $\overline{\mc H}^{cX}\!\!(\emptyset) = \bigcap_n \overline{\mc H}^{cX}\!\!(n)$.
By Lemma~\ref{lemma: complexity of S T} (or more precisely, by \eqref{equation: complexity of H t} from the~proof of Lemma~\ref{lemma: complexity of S T}), the~sets $\overline{\mc H}^{cX}\!\!(n)$ satisfy $X \subset \overline{\mc H}^{cX}\!\!(n) \in \mc F_{<\alpha}(\errata{cX).}$
In this setting, we can use Lemma~\ref{lemma: Fa hardness and density}\,(ii) to get $y\in  \overline{\overline{\mc H}^{cX}\!\!(n) \setminus X}^{cX}$.
In particular, there exists some $x_n \in \left( \overline{\mc H}^{cX}\!\!(n) \setminus X \right) \cap G_n$.

Set $K := \overline{\{ x_n | \ n\in\omega\} }^{cX}$.
Since we have $\bigcap_n \overline{G_n}^{cX} = \bigcap_n G_n \subset cX \setminus X$, it follows that $K$ is a~compact subset of $cX$ which is disjoint with $X$.
Let $dX := cX/_K$ be the~compact space obtained by ``gluing together'' the~points of $K$ (formally, we define an equivalence $\sim$ on $cX$ as
\[ x \sim y \iff ( x=y ) \lor ( x,y \in K ) ,\]
and define $dX$ as the~corresponding quotient of $cX$).
We can identify $dX$ with the~set $\{[K]\} \cup cX \setminus X$, where a~set containing $[K]$ is open in $dX$ if and only if it is of the~form $\{[K]\} \cup U \setminus K$ for some open subset $U\supset K$ of $cX$.
Since $K$ is disjoint from $X$, $dX$ is a~compactification of $X$.

We finish the~proof by showing $[K] \in \overline{\mc H}^{dX}\!(\emptyset) \setminus X$.
Let $n\in\omega$.
The topology of $dX$ is such that for every $A\subset X$, we have $\overline{A}^{cX} \ni x_n \implies \overline{A}^{dX} \ni [K]$. In particular, we obtain the~following implication for any leaf $t\in l(T_\alpha)$:
\begin{equation*}
\overline{\mc H}^{cX}\!(t) = \overline{\mc H(t)}^{cX} \ni x_n \implies
\overline{\mc H}^{dX}\!(t) = \overline{\mc H(t)}^{dX} \ni [K] .
\end{equation*}
It follows (by \eqref{equation: H t extensions} and induction over $r_l(T^t_\alpha)$) that the~following implication holds for every $t\in T_\alpha$
\begin{equation*}
\overline{\mc H}^{cX}\!(t) \ni x_n \implies \overline{\mc H}^{dX}\!(t) \ni [K] .
\end{equation*}
And since $\overline{\mc H}^{cX}\!\!(n)$ contains $x_n$, we get $[K] \in \overline{\mc H}^{dX}\!(n)$.
Because $n$ was arbitrary, this concludes the~proof:
\[ \overline{\mc H}^{dX}\!(\emptyset)
\overset{\alpha \text{ is}}{\underset{\text{even}}=}
\bigcap_n \overline{\mc H}^{dX}\!(n) \ni [K] .\]
\end{proof}

In particular, we obtain the~following corollary for Banach spaces:

\begin{corollary}[No universal representations for Banach spaces]\label{corollary: Banach spaces have no universal representation}
A Banach space $X$ is either reflexive (and therefore weakly $\sigma$-compact), or $(X,w)$ has no universal $\fsd$ representation.
\end{corollary}

\begin{proof}
Let $X$ be a~non-reflexive Banach space.
By Proposition~\ref{proposition: existence of simple representations}\,(a), it is enough to note that $X$ is locally $\fsd$-hard in some compactification.
Since $(X^{\star\star},w^\star)$ is $\sigma$-compact, being locally $\fsd$-hard in $(X^{\star\star},w^\star)$ implies being locally $\fsd$-hard in, for example, $\beta (X^{\star\star},w^\star)$.
By the~remark between Definition~\ref{definition: Fa hardness} and Lemma~\ref{lemma: Fa hardness and density}, it suffices to show that no regularly closed subset of $(X,w)$ is $F_\sigma$ in $(X^{\star\star},w^\star)$.
In this particular setting, this means showing that weakly regularly closed sets in $X$ are not weakly $\sigma$-compact.

For contradiction, assume that there is a~sequence of weakly compact sets $K_n \subset X$, such that $\bigcup_n K_n$ has got a~non-empty interior in weak topology.
In particular, the~interior contains a~norm-open set.
Consequently, Baire category theorem yields $n_0\in \omega$ such that $K_{n_0}$ contains an open ball $U_X(x,\epsilon)$.
However, this would imply that the~closed ball $B_X(x,\frac{\epsilon}{2})$ -- a~weakly closed subset of $K_{n_0}$ -- is weakly compact.
This is only true for reflexive spaces.
\end{proof}
\section{Regular Representations}\label{section: regular representations}

In this section, we introduce representations which are more structured than the~simple representations from Section~\ref{section: simple representations}.
In Section~\ref{section:Suslin_schemes}, we recall the~notion of a Suslin scheme and some related results. Section~\ref{section: R T sets} introduces the~concept of a regular representation and investigates its basic properties. In Section~\ref{section: Suslin scheme rank}, we give an alternative description of regular representations. This yields a criterion for estimating $\Compl{X}{Y}$, which will in particular be useful in Section~\ref{section: complexity of brooms}. Section~\ref{section: existence of regular representation} is optional, and justifies the~concept of regular representation\err{s} by proving their existence.

\subsection{Suslin Schemes}\label{section:Suslin_schemes}

A tool relevant to $\mc K$-analytic sets and $\mc F$-Borel complexities is the~notion of complete sequence of covers.

\begin{definition}[Complete sequence of covers]\label{definition: complete seq. of covers}
Let $X$ be a topological space.
\emph{Filter} on $X$ is a family of subsets of $X$, which is closed with respect to supersets and finite intersections and does not contain the~empty set.
A point $x\in X$ is said to be an \emph{accumulation point} of a filter $\mc F$ on $X$, if each neighborhood of $x$ intersects each element of $\mc F$.

A sequence $\left( \mc C_n \right)_{n\in\mathbb N}$ of covers of $X$ is said to be \emph{complete}, if every filter which intersects each $\mc C_n$ has an accumulation point in $X$.
\end{definition}

A notion related to Suslin sets is that of a Suslin scheme:

\begin{definition}[Suslin schemes and sets]\label{definition: Suslin schemes}
By a \emph{Suslin scheme} we understand a family $\mc C =\{ C(s) | \ s\in\seq \}$ which satisfies the~following monotonicity condition:
$$\left(\forall s,t\in\seq\right) : t\sqsupset s \implies C(t)\subset C(s) .$$
A \emph{Suslin operation} is the~mapping
\[ \mc A : \mc C \mapsto \mc A(\mc C) := \bigcup_{\sigma\in\baire} \bigcap_{n\in\omega} C(\sigma|n) .\]

\noindent Let $X\subset Y$ be topological spaces and $\mc C$ a Suslin scheme. We say that $\mc C$
\begin{itemize}
\errata{\item is a \emph{Suslin scheme in $Y$} if $\mc C \subset \mc P(Y)$;}
\item is \emph{closed} in $Y$ if $\mc C \subset \mc F_0 (Y)$;
\item \emph{\errata{$\mc A$-covers} $X$} if it satisfies $\mc A(\mc C) \supset X$;
\item is a Suslin scheme \emph{on $X$} if it is a~Suslin scheme in $X$ which \err{$\mc A$-covers} $X$;
\item is \emph{complete on $X$} if it is Suslin scheme on $X$ and $\left( \mc C_n \right)$ is a complete sequence of covers of $X$, where
\[ \mc C_{n} := \left\{ C(s) | \ s\in\omega^n \right\} .\]
\end{itemize}
\end{definition}

Note that a subset of $Y$ is Suslin-$\mc F$ in $Y$ if it is the~image under Suslin operation of some closed Suslin scheme in $Y$.
The existence of complete Suslin schemes is guaranteed by the~following result.

\begin{proposition}[Existence of complete Suslin schemes]
\label{proposition: K analytic spaces have complete Suslin schemes}
On any $\mc K$-analytic space, there is a complete Suslin scheme.
\end{proposition}

\begin{proof}
By \cite[Theorem 9.3]{frolik1970survey}, any $\mc K$-analytic space has a complete sequence of countable covers. Enumerate the~$n$-th cover as $\mc C_n = \{ C_n^k | \ k\in\omega \}$. Denoting $C(\emptyset) := X$ and $C(s) := C_0^{s(0)} \cap C_1^{s(1)} \cap \dots \cap C_k^{s(k)}$, we obtain a Suslin scheme satisfying $\mc A(\mc C)=X$. Moreover, $\mc C$ is easily seen to be complete on $X$.
\end{proof}

The main reason for our interest in complete Suslin schemes is the~following property:

\begin{lemma}\label{lemma: complete suslin schemes in super spaces}
Let $\mc C$ be a Suslin scheme on $X$. If $\mc C$ is complete on $X$, then $\mc A(\overline{\mc C}^Y\!\!) = X$ holds for every $Y\supset X$.
\end{lemma}

\begin{proof}
This holds, for example, by \cite[Lemma 4.7]{kalenda2018absolute}.
\end{proof}

\subsection{Definition of a Regular Representation} \label{section: R T sets}

The basic tool for the~construction of regular representations will be the~following mappings from trees to sequences:

\begin{definition}[Admissible mappings]\label{definition: admissible mapping}
Let $T\in\Tr$ be a tree and $\varphi : T\rightarrow \seq$ a mapping from $T$ to the~space of finite sequences on $\omega$. We will say that $\varphi$ is \emph{admissible}, if it satisfies
\begin{enumerate}[(i)]
\item $(\forall s,t\in T): s \sqsubset t \implies \varphi (s) \sqsubset \varphi(t)$
\item $(\forall t = (t_0,\dots,t_k)\in T): \left| \varphi (t) \right| = t(0)+\dots+t(k)$.
\end{enumerate}
\end{definition}

The next lemma completely describes how these mappings look like.
(Recall from Section~\ref{section: sequences} that each non-empty tree $T$ can be rewritten as $T = \{ \emptyset \} \cup \bigcup \{ m\ext T^{(m)} \, | \ m\in \ims{T}{\emptyset} \}$, where $T^{(m)}$ is defined as $\{ t' \in \seq \, | \ m\ext t' \in T \}$.)

\begin{lemma}[Construction of admissible mappings] \label{lemma: construction of admissible mappings}$\ $ 
\begin{enumerate}[(i)]
\item The only admissible mapping with range $\{\emptyset\}$ is the~mapping $\varphi : \emptyset \mapsto \emptyset$.
\item A mapping $\varphi : T \rightarrow \seq$ is admissible if and only if it is a restriction of some admissible mapping $\psi : \seq \rightarrow \seq $.
\item For $\{\emptyset\}\neq T\in \Tr$, a mapping $\varphi : T\rightarrow \seq$ is admissible if and only if it is defined by the~formula
\[ \varphi(t) =
\begin{cases}
\emptyset, 					& \textrm{for } t = \emptyset \\
s_m\ext \varphi_m(t'), 	& \textrm{for } t=m\ext t' \ \ \ (\text{where } (m)\in \ims{T}{\emptyset}, t'\in T^{(m)})
\end{cases} \]
for some sequences $s_m\in\omega^m$ and admissible mappings $\varphi_m : T^{(m)}\rightarrow \seq$.
\end{enumerate}
\end{lemma}
\begin{proof}
$(i)$: By definition of an admissible mapping, we have $|\varphi(\emptyset)|=0$, so $\varphi(\emptyset)$ must be equal to the~empty sequence $\emptyset$.

$(ii)$: Let $T\in \Tr$. From the~definition of an admissible mapping, we see that if $\psi: \seq \rightarrow \seq$ is admissible, then $\psi|_T : T \rightarrow \seq$ is also admissible. Moreover, any admissible mapping $\varphi: T \rightarrow \seq$ is of this form. To see this, fix any $\sigma\in \baire$.
For every $s\in \seq $, there exists some $t\in T$, such that $s=t\ext h$ holds for some $h\in\seq$.
Denoting by $h=(h(0),\dots,h(k))\in\seq$ the~sequence corresponding to the~longest such $t\in T$, we set $\psi(s):=\varphi(t)\ext \sigma|{h(0)+\dots+h(k)}$. Clearly, $\psi$ is admissible and coincides with $\varphi$ on $T$.

$(iii)$:
Denote $M := \ims{T}{\emptyset}$.
Firstly, note that if $\varphi:T\rightarrow \seq$ is admissible, then $\varphi(\emptyset)=\emptyset$, for every $m\in M$ we have $s_m := \varphi(m) \in \omega^m$ and for other elements of the~tree, we have $\varphi(t)=\varphi(m\ext t')\sqsupset s_m$ for some $m\in M$ and $t'\in T^{(m)}$. Denoting $\varphi(m\ext t') =: s_m \ext \varphi_m(t')$, we get mappings $\varphi_m : T^{(m)}\rightarrow \seq$ which are easily seen to be admissible.

On the~other hand, it is not hard to see that if mappings $\varphi_m : T^{(m)} \rightarrow \seq$ are admissible and the~sequences $s_m$ for $m\in M$ are of the~correct length, then the~formula in $(iii)$ defines an admissible mapping.
\end{proof}

Before defining the~desired regular representations, we introduce the~following $\mathbf{R}_T$-sets (of which regular representations will be a special case, for particular choices of $\mc C$ and $T$):

\begin{definition}[Regular representation] \label{definition: R T sets}
Let $T$ be a tree and $\mc C$ a Suslin scheme. We define
\[ \mathbf{R}_T(\mc C):=\left\{ x \in C(\emptyset) \ | \ \left(\exists \varphi : T\rightarrow \seq \textrm{ admissible} \right)\left( \forall t\in T \right) : x\in C \left(\varphi (t)\right) \right\} . \]
Let $X\subset Y$ be topological spaces. \err{If $\mc C$ and $T$ additionally satisfy $\mc C \subset \mc P(Y)$ and $X = \mathbf{R}_T(\mc C)$, then the pair $(\mc C,T)$ is said to be a \emph{regular representation of $X$ in $Y$}.}
\end{definition}

\noindent When $x$ satisfies $x \in \bigcap_T C(\varphi(t))$ for some admissible $\varphi:T\rightarrow \seq$, we shall say that \emph{$\varphi$ witnesses that $x$ belongs to $\mathbf{R}_T(\mc C)$}.
We also need a technical version of the~notation above. For $h\in\seq$, we set
\[ \mathbf{R}^h_T(\mc C):=\left\{ x \in C(\emptyset) \ | \ \left(\exists \varphi : T\rightarrow \seq \textrm{ admissible} \right)\left( \forall t\in T \right) : x\in C \left(h\ext\varphi (t)\right) \right\} .\]

The technical properties of $\mathbf{R}_T$-sets are summarized by the~following lemma:

\begin{lemma}[Basic properties of $\mathbf{R}_T$-sets] \label{lemma: basic properties of R T sets}
Let $\mc C$ be a Suslin scheme, $S,T\in\Tr$ and $h\in\seq$.
\begin{enumerate}[(i)]
\item Suppose that there is a mapping $f:T \rightarrow S$ satisfying \label{case: R T and embeddings}
	\begin{itemize}
	\item $\left( \forall t_0, t_1 \in T \right): t_0 \sqsubset t_1 \implies f(t_0) \sqsubset f(t_1)$,
	\item $\left(\forall t = (t(0), \dots, t(k)) \in T \right) : f(t)(0)+\dots + f(t)(k) \geq t(0) + \dots + t(k)$.
	\end{itemize}
Then $\mathbf{R}_S(\mc C) \subset \mathbf{R}_T(\mc C)$.
\item In particular, we have $\mathbf{R}_S(\mc C) \subset \mathbf{R}_T(\mc C)$ whenever $S\supset T$. \label{case: R T and smaller trees}
\item We have $\mathbf{R}^h_{\{\emptyset\}}(\mc C) = C(h)$ and the~following recursive formula: \label{case: R T formula}
\begin{equation*}
 \mathbf{R}_T^h(\mc C) \ = \ \bigcap_{(m) \in \ims{T}{\emptyset}} \bigcup_{s_m\in\omega^m} \mathbf{R}_{T^{(m)}}^{h\ext s_m}(\mc C) .
\end{equation*}
\item If $T\in\textrm{IF}$ has a branch with infinitely many non-zero elements, then $\mathbf{R}_T(\cdot)$ coincides with the~Suslin operation $\mc A(\cdot)$. \label{case: R T and IF}
\end{enumerate}
\end{lemma}

We can view the~mapping $f:T\rightarrow S$ from Lemma~\ref{lemma: basic properties of R T sets}\,\eqref{case: R T and embeddings} as an \errata{``embedding''} of the~tree $T$ into $S$.
In this light, \eqref{case: R T and embeddings} implies that ``larger tree means smaller $\mathbf{R}_{(\cdot)}(\mc C)$'' and for ``equivalent'' trees, the~corresponding sets $\mathbf{R}_{(\cdot)}(\mc C)$ coincide.

\begin{proof}
Let $\mc C$, $T$ and $h$ be as above. For simplicity of notation, we will assume that $h=\emptyset$ and therefore we will only work with the~sets $\mathbf{R}_T(\mc C)$. However, the~proof in the~general case is exactly the~same as for $h=\emptyset$.

\eqref{case: R T and embeddings}: Let $f:T\rightarrow S$ be as in the~statement and suppose that $\varphi : S \rightarrow \seq$ witnesses that $x\in \mathbf{R}_S(\mc C)$.
Using the~properties of $f$, we get that the~following formula defines an admissible mapping:
\[ \psi : t\in T \mapsto \varphi(f(t)) |_{t(0)+\dots + t(|t|-1)} . \]
Using monotonicity of $\mc C$, we prove that $\psi$ witnesses that $x$ belongs to $\mathbf{R}_T(\mc C)$:
\begin{equation*} \begin{split}
x \in 		& \bigcap_{s\in S} C(\varphi(s)) \overset{f(T)\subset S}\subset
				\bigcap_{t\in T} C \left(\varphi(f(t))\right)
				\overset{\textnormal{mon.}}{\underset{\text{of }\mc C}\subset}
				\bigcap_{t\in T} C \left(\varphi(f(t)) |_{\sum t(k)}\right)
				\overset{\textnormal{def.}}{\underset{\text{of }\psi}\subset}				
				\bigcap_{t\in T} C(\psi(t)) .
\end{split} \end{equation*}

\errata{\eqref{case: R T and smaller trees} is a special case of \eqref{case: R T and embeddings} with $f(t)=t$.}

(iii): The identity $\mathbf{R}_{\{\emptyset\}}(\mc C) = C(\emptyset)$ follows from Lemma~\ref{lemma: construction of admissible mappings}~$(i)$.
For $\{\emptyset\} \subsetneq T$, we can rewrite $T$ as
\[  T = \{\emptyset\}\cup \bigcup \{ m\ext T^{(m)} | \ m\in M\} ,\]
where $M := \ims{T}{\emptyset}$.

Let $x\in C(\emptyset)$.
By definition of $\mathbf{R}_T$, we have $x \in \mathbf{R}_T(\mc C)$ if and only if there exists a witnessing mapping -- that is, an admissible mapping $\varphi : T \rightarrow \seq$ which satisfies $x \in \bigcap_T C(\varphi(t))$.
Since $\mc C$ is monotone, this is further equivalent to the~existence of admissible $\varphi : T \rightarrow \seq$ s.t.\errata{}
\[ x \in \bigcap_T C(\varphi(t)) = \bigcap_{\err{m\in}M} \bigcap_{\err{t} \in T^{(m)}} C(\varphi(m\ext t)) .\]
By Lemma~\ref{lemma: construction of admissible mappings}\,(iii), this is true if and only if there are sequences $s_m \in \omega^m$, $m\in M$, and admissible mappings $\varphi_m : T^{(m)} \rightarrow \seq$, $m\in M$, for which $x$ belongs to $\bigcap_{T^m} C(s_m \ext \varphi_m (t'))$.
This is precisely when, for each $m\in M$, there exist some mapping witnessing that $x \in \mathbf{R}^{s_m}_{T^{(m)}}(\mc C)$ holds for some $s_m \in \omega^m$.

\errata{\eqref{case: R T and IF}:} Let $T$ be as in the~statement and $\mc C$ be a Suslin scheme. To show that $\mc A(\mc C) \subset \mathbf{R}_T(\mc C)$, it suffices to show that $\mc A(\mc C) \subset \mathbf{R}_{\seq}(\mc C)$ (by \eqref{case: R T and smaller trees}).
Let $x \in \bigcap_{m\in\omega} C(\sigma|m) \subset \mc A(\mc C)$ and define $\varphi : \seq \rightarrow \seq $ as
\[ \varphi(\, (t(0),\dots,t(k)) \,):= \sigma |_{t(0)+\dots+t(k)} .\]
Since $\varphi$ is admissible, it witnesses that $x$ belongs to $\mathbf{R}_{\seq}(\mc C)$.

For the~reverse inclusion, let $\nu$ be a branch of $T$ with $\sum_k \nu(k) = \infty$ and $\varphi:T\rightarrow \seq$ be an admissible mapping. The sequence $\varphi(\nu|_0) \sqsubset \varphi(\nu|_1) \sqsubset \varphi(\nu|_2) \sqsubset \dots$ eventually grows to an arbitrary length, which means that there exists $\sigma\in\baire$ such that for each $m\in\omega$ there is $k\in\omega$ such that $\sigma|m \sqsubset \varphi(\nu|_k)$. In particular, we have
\[ \bigcap_{t\in T} C(\varphi(t)) \subset \bigcap_{k\in\omega} C(\varphi(\nu|_k)) \subset \bigcap_{m\in\omega} C(\sigma|m) \subset \mc A (\mc C) . \]
\end{proof}

Since our goal is to study $\mc F$-Borel sets, we will focus on $\mathbf{R}_T$-sets corresponding to closed Suslin schemes. The following notation is handy for that purpose.

\begin{definition}[$\mathbf{R}_T(\mc C,Y)$-sets]
For any Suslin scheme $\mc C$, tree $T$, and a topological space $Y$, we define
\[ \mathbf{R}_T(\mc C,Y):=\mathbf{R}_T(\bar{\mc C}^Y) \textrm{, where } \bar{\mc C}^Y
 := \left\{ \overline{C(s)\cap Y}^Y | \ s\in\seq \right\} . \]
\end{definition}

We have the~following complexity estimate for $\mathbf{R}_T$-sets:

\begin{lemma}[Complexity of $\mathbf{R}_T$-sets] \label{lemma: complexity of R T}
Let $\mc C$ be a closed Suslin scheme in a topological space $Y$, and $T\in \textnormal{WF}$.
\err{Expressing} the~leaf-rank of $T$ as $r_l(T) = \lambda+n$ (where $\lambda$ is a limit ordinal or $0$ and $n\in\omega$), we have:
\begin{enumerate}[(i)]
\item $\mathbf{R}_T(\mc C) \in \mc F_{\lambda+2n}(Y)$;
\item If $\ims{T}{\emptyset}$ is finite, then $\mathbf{R}_T(\mc C) \in \mc F_{\lambda+2n-1}(Y)$.
\end{enumerate}
\end{lemma}

\noindent (The expression in (ii) is well-defined because if $\ims{T}{\emptyset}$ is finite, then $r_l(T)$ cannot be a limit ordinal.)

\begin{proof}
We proceed by induction over $r_l(T)$.
For $r_l(T) = 0$, we necessarily have $T= \{\emptyset\}$, which yields
\[ \mathbf{R}_T(\mc C) = \mathbf{R}_{\{\emptyset\}} (\mc C)
\overset{L\ref{lemma: construction of admissible mappings}(i)}= C(\emptyset) \in \mc F_0(Y) .\]

Let $\mc C$ be a closed Suslin scheme in $Y$, $T\in \textnormal{WF}$ a tree s.t. \errata{$r_l(T) = \lambda + n$.}
Suppose, as an induction hypothesis, that for every closed Suslin scheme $\mc D$ in $Y$ and every $S\in \Tr$ with $r_l(S) < \lambda + n$,
we have $\mathbf{R}_S(\mc D) \in \mc F_\beta (Y)$ for some \emph{even} ordinal $\beta_S<\lambda + 2n$.
Denoting $M:=\ims{T}{\emptyset}$, Lemma~\ref{lemma: basic properties of R T sets}\,\eqref{case: R T formula} yields
\begin{equation}\label{equation: R T formula}
\mathbf{R}_T(\mc C) = \bigcap_{m\in M} \bigcup_{s\in \omega^m} \mathbf{R}^s_{T^{(m)}}(\mc C) .
\end{equation}
Let $m\in M$.
We have $r_l(T^{(m)}) < \lambda + n$, and $\mathbf{R}^s_{T^{(m)}}(\mc C)$ can be rewritten as $\mathbf{R}_{T^{(m)}}(\mc D)$, where $\mc D := \left( C(s\ext u) \right)_{u\in \seq}$ is a closed Suslin scheme in $Y$.
Using the~induction hypothesis, we get $\mathbf{R}^s_{T^{(m)}}(\mc C) \in \mc F_{\beta_m}(Y) = \mc F_{<\beta_m+1}(Y)$ for some even $\beta_m < \lambda +2n$.
Since $\beta_m +1$ is odd, we get
\[ \bigcup_{s\in\omega^m} \mathbf{R}^s_{T^{(m)}}(\mc C) \in (\mc F_{<\beta_m +1}(Y))_\sigma = \mc F_{\beta_m+1}(Y) \subset \mc F_{<\lambda +2n}(Y) .\]

(i): It follows that $\mathbf{R}_T(\mc C) \in (\mc F_{<\lambda +2n}(Y))_\delta = \mc F_{\lambda+2n}(Y)$.

(ii): When the~set $M$ is finite, we necessarily have $n\geq 1$.
It follows that each $\bigcup_{s\in\omega^m} \mathbf{R}^s_{T^{(m)}}(\mc C)$ belongs to $\mc F_{<\lambda +2n}(Y) = \mc F_{\lambda +2n-1}(Y)$.
Since the~intersection in \eqref{equation: R T formula} is finite, it follows that $\mathbf{R}_T(\mc C)$ belongs to $\mc F_{\lambda +2n-1}(Y)$ as well.
\end{proof}

It follows from Lemma~\ref{lemma: complexity of R T} (and Lemma~\ref{lemma: basic properties of R T sets}\,\eqref{case: R T and IF}) that for each $T\in \Tr$, there exists $\alpha_T \leq \omega_1$ such that $\mathbf{R}_T(\mc C,Y) \in \mc F_{\alpha_T}(Y)$ holds for every $\mc C$ and $Y$.
In this sense, each $\mathbf{R}_T(\cdot , Y)$ can be understood as an operator which maps Suslin schemes to $\mc F_{\alpha_T}$-subsets of $Y$.

It follows from \errata{\eqref{case: R T and embeddings}} in Lemma~\ref{lemma: basic properties of R T sets} that for the~purposes of studying the~$\mathbf{R}_T$-sets, many trees are in fact equivalent.
We therefore focus our attention on the~``canonical'' trees $T^c_\alpha$, $\alpha \leq \omega_1$.
In Proposition~\ref{proposition: existence of regular F alpha representations}, we show that this is sufficient -- every $\Fa$-subset $X$ of $Y$ can be written as
$X = \mathbf{R}_{T^c_\alpha}(\mc C,Y)$ for some $\mc C$.

\begin{definition}[$\mathbf{R}_\alpha$-sets]\label{definition: R alpha sets}
For a Suslin scheme $\mc C$ and $\alpha\leq \omega_1$, we denote
\[ \mathbf{R}_\alpha(\mc C):=\mathbf{R}_{T^\textrm{c}_\alpha}(\mc C) .\]
We define $\mathbf{R}_\alpha(\mc C,Y)$, $\mathbf{R}^h_\alpha(\mc C)$ and $\mathbf{R}^h_\alpha(\mc C,Y)$ analogously.
\end{definition}

\errata{
For finite $\alpha$, the $\mathbf{R_\alpha}$-sets can be expressed in a~simpler way:

\begin{remark}[$\mathbf{R}_\alpha$-sets for $\alpha\in\omega$]\label{rem:R_n}
For any Suslin scheme $\mc C$ and $k\in\omega$, we have
\begin{align*}
& \mathbf{R}_0(\mc C) = C(\emptyset), \
\mathbf{R}_2(\mc C) = \bigcap_{m\in\omega} \bigcup_{s\in\omega^m} C(s), \
\mathbf{R}_4(\mc C) =
\bigcap_{m\in\omega} \bigcup_{s\in\omega^m} \bigcap_{n\in\omega} \bigcup_{t\in\omega^n}C(s\ext t), & \dots \\
& \mathbf{R}_{2k}(\mc C) =
\bigcap_{m_1\in\omega} \bigcup_{|s_1| = m_1} \cdots \bigcap_{m_k\in\omega} \bigcup_{|s_k| =m_k}C(s_1 \ext \dots \ext s_k)
, \ \dots \ .
\end{align*}
To get the formula for $\mathbf{R}_{2k+1}(\mc C)$, set $m_1 = 1$ in the expression for $\mathbf{R}_{2(k+1)}(\mc C)$.
\end{remark}

\begin{proof}
For $k=0$, the result is immediate (and also follows from the $\{\emptyset\}$-case of Lemma~\ref{lemma: basic properties of R T sets}\,\eqref{case: R T formula} with $h=\emptyset$).
The ``2k'' part of the remark for $0<k<\omega$ follows by induction, using Lemma~\ref{lemma: basic properties of R T sets}\,\eqref{case: R T formula} with $T=T_k$ (since $T^c_{2k} = T_k$ and $\ims{T_k}{\emptyset} = \{(j) | \, j\in\omega \}$).
For the ``2k+1'' part for $0\leq k <\omega$, we have $T^c_{2k+1} = \{\emptyset\} \cup 1 \ext T_k$, and the result once again follows from Lemma~\ref{lemma: basic properties of R T sets}\,\eqref{case: R T formula} (using the already proven ``2k'' part).
\end{proof}

\noindent As a corollary of Remark~\ref{rem:R_n}, we get that $\mathbf{R}_k(\cdot,\cdot)$-sets belong to the class $\mc F_k$. The general version of this result is 
}
a special case of Lemma~\ref{lemma: complexity of R T}:

\begin{proposition}[$\mathbf{R}_\alpha$-sets are $\Fa$]\label{proposition: complexity of R alpha sets}
For any Suslin scheme $\mc C$ and $\alpha \leq \omega_1$, we have $\mathbf{R}_\alpha (\mc C, Y) \in \Fa(Y)$ for any topological space $Y$.
\end{proposition}

\begin{proof}
For $\alpha<\omega_1$, this follows from Lemma~\ref{lemma: complexity of R T}. For $\alpha=\omega_1$, this holds by \eqref{case: R T and IF} from Lemma~\ref{lemma: basic properties of R T sets}.
\end{proof}

Let $X\subset Y$ be a~topological space and $\mc C$ a Suslin scheme in $X$, and suppose that $\mc A(\mc C) \supset X$.
These assumptions in particular imply $\mc A(\mc C) = X$ and $C(\emptyset) = X$.
By definition of $\mathbf{R}_0(\cdot)$ (and Lemma~\ref{lemma: construction of admissible mappings}\,(i)), we have $\mathbf{R}_0(\mc C) = C(\emptyset)$.
Consequently, the~corresponding closed Suslin scheme $\overline{\mc C}^Y$ in $Y$ satisfies $\mc A(\overline{\mc C}^Y\!) \supset X$ and $\mathbf{R}_0(\mc C,Y) = \overline{C(\emptyset)}^Y = \overline{X}^Y$.
As a particular case of \eqref{case: R T and embeddings} from Lemma~\ref{lemma: basic properties of R T sets}, we get
\begin{align} \label{equation: R alpha as approximations of X}
\overline{X}^Y = \overline{C(\emptyset)}^Y = \mathbf{R}_0(\mc C,Y)
\supset \mathbf{R}_1(\mc C,Y) \supset \dots \\
\dots \supset \mathbf{R}_\alpha(\mc C,Y) \supset \dots
\supset \mathbf{R}_{\omega_1}(\mc C,Y) = \mc A(\overline{\mc C}^Y\!) \supset X .
\nonumber
\end{align}

In the~sense of \eqref{equation: R alpha as approximations of X}, the~sets $\mathbf{R}_\alpha(\mc C,Y)$, $\alpha\leq \omega_1$, can be viewed as approximations of $X$ in $Y$, or as its $\Fa$-envelopes in $Y$.
The set $\mathbf{R}_0(\mc C,Y)$ is the~nicest (in other words, closed) approximation, but also the~biggest. We can get a more accurate approximation by increasing $\alpha$, but this comes at a cost of increased complexity.
We might get lucky with the~choice of $\mc C$, and get $\mathbf{R}_\alpha(\mc C,Y) = X$ at some point (and thus also $\mathbf{R}_\beta (\mc C,Y)$ for every $\beta \in [\alpha,\omega_1]$). In that case, we would say that $\mc C$ is a regular $\Fa$-representation of $X$ in $Y$.
This is the~case when $\mc C$ is complete, because then $X = \mc A(\overline{\mc C}^Y\!)$ (Lemma~\ref{lemma: complete suslin schemes in super spaces}). But even then, there is no reason to expect that $\alpha = \Compl{X}{Y}$.

We have seen that a regular representation of $X$ in $Y$ can also be specified by picking a Suslin scheme $\mc C$ in $X$ and providing an ordinal $\alpha$ for which this $\mc C$ ``works''. This motivates the~following definition:

\begin{definition}[Regular $\Fa$-representations]\label{definition: regular representations}
A Suslin scheme $\mc C$ in $X$ is said to be a
\begin{itemize}
\item \emph{regular $\Fa$-representation of $X$ in $Y$}, for some $Y\supset X$ and $\alpha\leq\omega_1$, if it satisfies $X = \mathbf{R}_\alpha (\mc C,Y)$;
\item \emph{universal regular $\Fa$-representation of $X$} if it is a regular $\Fa$-representation of $X$ in every $Y\supset X$.
\end{itemize}
\end{definition}

The existence of regular $\Fa$-representations is studied in Section~\ref{section: existence of regular representation}.
However, as this work is not a detective novel, we can go ahead and spoil the~surprise right away: any $\Fa$-subset $X$ of $Y$ has some regular $\Fa$-representation (Proposition~\ref{proposition: existence of regular F alpha representations}). And when either $X$ or $Y$ is $\mc K$-analytic, this representation can be made complete ``without loss of generality'' (Theorem~\ref{theorem: existence of regular representation}).

Clearly, the~condition of having a universal regular $\Fa$-representation is even stronger than being an absolute $\Fa$ space (at least formally).
In Section~\ref{section: broom space properties} we study a class of spaces (the so-called ``broom spaces'') which do admit universal regular representations. However, these spaces are rather simple from the~topological point of view -- they only have a single point which is not isolated.
As we have seen is Section~\ref{section: universal representation}, this is rather an exception.

\bigskip

Next, we describe $\mathbf{R}_\alpha$-sets in terms of a ``Suslin scheme rank'' on trees.

\subsection{Regular Representations and Suslin Scheme Ranks}\label{section: Suslin scheme rank}

For a Suslin scheme $\mc C$ in $X$ and a point $y\in Y \supset X$ we denote
\begin{equation} \label{equation: S C of y}
S_{\mc C}(y)  := \{ s\in \seq | \ \overline{C(s)}^Y \ni y  \} .
\end{equation}
Note that $S_{\mc C}(y)$ is always a tree (by the~monotonicity property of Suslin schemes).
The \emph{rank of $y$ corresponding to $\mc C$} is defined as $r_{\mc C}(y) := \rank (S_{\mc C}(y) )$.
The following lemma ensures that as long as we investigate only points from $Y\setminus \mc A (\overline{\mc C}^Y\!\!)$, we can assume that $r_{\mc C}$ is countable.

\begin{lemma}[$S_{\mc C}(\cdot)$ and $\textnormal{WF}$]\label{lemma: IF trees}
Let $X\subset Y$ be topological spaces, $\alpha\leq \omega_1$, and $\mc C$ a Suslin scheme in $Y$.
\begin{enumerate}[(i)]
\item For every $y\in Y$, we have $y \in \mc A( \overline{\mc C}^Y\!) \iff
	S_{\mc C}(y) \in \textnormal{IF} \iff r_i (S_{\mc C}(y) ) = \omega_1$. \label{case: Suslin operation and S of y}
\item If $\mc A (\overline{\mc C}^Y\!\!) = X$, we have $S_{\mc C}(y) \in \textnormal{WF}$ for every $y\in \mathbf{R}_\alpha(\mc C,Y) \setminus X$. \label{case: complete C and WF}
\end{enumerate}
\end{lemma}

By Lemma~\ref{lemma: complete suslin schemes in super spaces}, the~assumptions of (ii) are in particular satisfied when $\mc C$ is complete on $X$.

\begin{proof}
The first equivalence in \eqref{case: Suslin operation and S of y} is immediate once we rewrite $\mc A( \overline{\mc C}^Y\!)$ as
\[ \mc A( \overline{\mc C}^Y\!) = \bigcup_\sigma \bigcap_n \overline{C(\sigma|n)}^Y .\]
The second equivalence follows from the~fact the~derivative $D_i$ leaves a branch of a tree ``untouched'' if and only if the~branch is infinite.

\eqref{case: complete C and WF}: We have $\mc A( \overline{\mc C}^Y\!) \subset X$.
By \eqref{case: Suslin operation and S of y}, the~tree $S_{\mc C}(y)$ corresponding to $y \in Y\setminus X$ must not be ill-founded.
\end{proof}

We have seen that $\mathbf{R}_\alpha(\mc C,Y) \supset \mc A(\overline{\mc C}^Y\!) \supset X$ holds for every $\alpha\leq \omega_1$. In other words, $\mathbf{R}_\alpha(\mc C,Y)$ always contains $X$, and then maybe some extra points from $Y\setminus X$.
It is our goal in this subsection to characterize this remainder of $\mathbf{R}_\alpha(\mc C,Y)$ in $Y$.
We claim that for even $\alpha$, this remainder can be written as
\[ \mathbf{R}_\alpha(\mc C,Y) \setminus X = \{ y\in Y \setminus X | \ r_{\mc C}(y) \geq \alpha' \} . \footnote{Recall that for $\alpha=\lambda+2n+i$ we have $\alpha ' := \lambda +n$ (where $\lambda$ is limit or $0$, $n\in\omega$ , and $i\in\{0,1\}$).} \]
Unfortunately, we haven't found any such succinct formulation for odd $\alpha$. We will therefore describe the~remainder in terms of $\D$-derivatives of $S_{\mc C}(y)$.
\errata{With the exception of Corollary~\ref{corollary: sufficient condition for Fa}, the remainder of Section~\ref{section: Suslin scheme rank} is dedicated to the proof of this result.}

\begin{proposition}[Description of regular representations via $\D$]\label{proposition: regular representations and D}
Let $\alpha < \omega_1$ and let $\mc C$ be a Suslin scheme on $X \subset Y$ satisfying $\mc A(\overline{\mc C}^Y\!\!) = X$.
\begin{enumerate}[(i)]
\item For even $\alpha$, we have $\mathbf{R}_\alpha (\mc C ,Y) \setminus X =
	\{ y\in Y\setminus X | \ \D^{\alpha'}(S_{\mc C} (y)) \neq \emptyset \}$
\item For odd $\alpha$, we have $\mathbf{R}_{\alpha} (\mc C ,Y) \setminus X =
	\{ y\in Y \setminus X| \ \D^{\alpha'}(S_{\mc C} (y)) \supsetneq \{\emptyset\} \} $.
 \item Moreover, we have $\mathbf{R}_{T_{\alpha',i}^{\textnormal{c}}} (\mc C ,Y) \setminus X =
 	\{ y\in Y \setminus X| \ \D^{\alpha'}(S_{\mc C} (y)) \cap \omega^i \neq \emptyset \} $ for any $i\in\omega$.
\end{enumerate}
\end{proposition}

\noindent (Recall that for $\alpha=\omega_1$ and complete $\mc C$ on $X$, we have $\mathbf{R}_{\omega_1}(\mc C ,Y) = X$.)

For even $\alpha$ we have $\alpha'=(\alpha+1)'$ and the~set $\D^{\alpha'}(S_{\mc C}(y))$ is always a tree. Therefore, the~case (i), resp. (ii) and (iii), of the~proposition says that the~remainder is equal to those $y\in Y\setminus X$ for which $\D^{\alpha'}(S_{\mc C}(y))$ contains a sequence of length 0 (that is, the~empty sequence $\emptyset$), resp. some sequence of length 1 (equivalently, any $s\neq \emptyset$), resp. some sequence of length $i$.

This gives the~following criterion for bounding the~complexity of $X$ in $Y$ from above (actually, even for showing that $X$ admits a regular $\Fa$-representation in $Y$):

\begin{corollary}[Sufficient condition for being $\Fa$]
	\label{corollary: sufficient condition for Fa}
A space $X\subset Y$ satisfies $X\in \Fa(Y)$ for some $\alpha<\omega_1$, provided that there is a Suslin scheme $\mc C$ on $X$, s.t. $\mc A(\overline{\mc C}^Y\!\!) = X$ and one of the~following holds:
\begin{enumerate}[(i)]
\item $\alpha$ is even and $\D^{\alpha'} (S_{\mc C}(y))$ is empty for every $y \in Y\setminus X$;
\item $\alpha$ is odd and $\D^{\alpha'} (S_{\mc C}(y))$ is either empty or it only contains the~empty sequence for every $y \in Y\setminus X$;
\item $\alpha$ is odd and there is $i\in\omega$ s.t. $\D^{\alpha'} (S_{\mc C}(y))$ only contains sequences of length $\leq i$ (that is, it is a subset of $\omega^{\leq i}$) for every $y \in Y\setminus X$.
\end{enumerate}
\end{corollary}

We now proceed to give the~proof of Proposition~\ref{proposition: regular representations and D}.
For the~``$\supset$'' inclusion, it suffices to construct a suitable admissible mapping.

\begin{proof}[Proof of Proposition~\ref{proposition: regular representations and D}, the~``$\supset$'' part]
Recall that for any derivative $D$ on $\seq$, the~corresponding rank $r$ is defined such that for every $S\subset \seq$, $r(S)$ is the~highest ordinal for which $D^{r(S)}(S)$ is non-empty.

(i): Let $\mc C$ be a Suslin scheme in $X\subset Y$, $\alpha<\omega_1$ even, and $y\in Y$ s.t. $r_{\mc C}(y) = \rank (S_{\mc C}(y)) \geq \alpha'$.
(This direction of the~proof does not require $\mc C$ to satisfy $\mc A(\overline{\mc C}^Y\!\!) \errata{= X.}$)
Recall that $\mathbf{R}_\alpha(\mc C,Y) \overset{\text{def.}}= \mathbf{R}_{T_{\alpha'}}(\mc C,Y)$ and, as noted below Notation~\ref{notation: trees of height alpha}, $r_l(T_{\alpha'})=\alpha' \leq r_{\mc C}(y)$.

We will prove that for every $T\in \Tr$ and $y\in Y$ satisfying $r_l(T) \leq r_{\mc C}(y)$, there exists an admissible mapping witnessing that $y$ belongs to $\mathbf{R}_T(\mc C,Y)$.
In particular, we will use induction over $|t|\in \omega$, to construct a mapping $\varphi : T \rightarrow \seq$ which satisfies (a) and (b) for every $t\in T$:
\begin{enumerate}[(a)]
\item If $t=t'\ext m$ for some $t'\in T$ and $m\in\omega$, then we have $\varphi(t) \sqsupset \varphi(t')$ and $|\varphi(t)| = |\varphi(t')| + m$.
\item For every $\gamma<\omega_1$, we have $t\in D^\gamma_l(T) \implies \varphi(t) \in \D^\gamma (S_{\mc C}(y))$.
\end{enumerate}
By (a), the~resulting mapping will be admissible. By (b) with $\gamma=0$, we have $\varphi: T \rightarrow S_{\mc C}(y)$, which proves that $\varphi$ witnesses that $y\in \mathbf{R}_T(\mc C,Y)$.

$|t|=0$: The only sequence of length $0$ is the~empty sequence, so (a) is trivially satisfied by setting $\varphi(\emptyset) := \emptyset$. Any $\gamma$ satisfying $\emptyset\in D^\gamma_l(T)$ is smaller than $r_l(T)$ by definition, and we have $r_l(T) \leq \rank(S_{\mc C}(y)) = r_{\mc C}(y)$. It follows that $\varphi(\emptyset)$ belongs to $\D^\gamma(S_{\mc C}(y))$:
\[ \D^\gamma(S_{\mc C}(y)) \supset \D^{r_l(T)}(S_{\mc C}(y))
\supset \D^{\rank(S_{\mc C}(y))}(S_{\mc C}(y)) \ni \emptyset \]
(where the~last set is non-empty by definition of $\rank$).

Let $t=t'\ext m\in T$ (where $t'\in T$ and $m\in\omega$) and suppose we already have $\varphi(t')$ satisfying (a) and (b).
Let $\gamma_t$ be the~highest ordinal for which $t\in D^{\gamma_t}_l(T)$. By definition of $D_l$, we have $t' \in D^{\gamma_t+1}_l(T)$.
Consequently, the~induction hypothesis gives $\varphi(t') \in \D( \D^{\gamma_t} (S_{\mc C}(y)) )$.
By definition of $\D$, there is some $s_m \in \D^{\gamma_t} (S_{\mc C}(y))$ for which $|s_m| \geq |\varphi(t')| + m$.
We set $\varphi(t) = \varphi(t'\ext m) := s_m |_{|\varphi(t')|+m}$. Clearly, $\varphi(t)$ satisfies (a) and (b).

(ii), (iii): (ii) is a special case of (iii), so it suffices to prove (iii).
Let $\mc C$ be a Suslin scheme on $X\subset Y$, $i\in\omega$, and $\alpha<\omega_1$ an even ordinal (so that $\alpha+1$ is odd -- recall that $(\alpha+1)'=\alpha'$).
Let $y\in Y$ be s.t. $\D^{\alpha'}(S_{\mc C}(y))$ contains some sequence $s_y$ of length $i$.

Since $T^c_{\alpha+1,i} \overset{\text{def.}}= \{\emptyset\} \cup i\ext T^c_\alpha$, an application of Lemma~\ref{lemma: basic properties of R T sets}\,\eqref{case: R T formula} yields the~first identity in the~following formula:
\[ \mathbf{R}_{T^c_{\alpha,i}} (\mc C,Y)
= \bigcup_{s\in\omega^i} \mathbf{R}^s_{T^c_\alpha} (\mc C,Y)
\supset \mathbf{R}^{s_y}_{T^c_\alpha} (\mc C,Y)
= \mathbf{R}^{s_y}_{\alpha} (\mc C,Y)
= \mathbf{R}_{\alpha} \left( \mc D,Y\right) ,\]
where $\mc D := \left( C(s_y \ext t) \right)_{t\in\seq}$.
Clearly, $S_{\mc D}(y)$ contains all sequences $t$ for which $s_y \ext t$ belongs to $S_{\mc C}(y)$.
In particular, $\D^{\alpha'}(S_{\mc D}(y))$ is non-empty, so $y$, $\alpha$ and $\mc D$ satisfy all assumptions used in the~proof of case (i). Therefore, $y$ belongs to $\mathbf{R}_{\alpha} \left( \mc D,Y\right) \subset \mathbf{R}_{T^c_{\alpha,i}} (\mc C,Y)$, which concludes the~proof. 
\end{proof}

For the~``$\subset$'' inclusion, we start with the~following lemma:

\begin{lemma}[Points from the~remainder force big $S_{\mc C}(\cdot)$]\label{lemma: admissible mappings and rank}
Let $\mc C$ be a Suslin scheme on $X$, $T\in\Tr$, and $\varphi$ witnesses that $y\in \mathbf{R}_T(\mc C,Y)$. If $\mc A(\overline{\mc C}^Y\!\!) = X$ and $y\notin X$, then $S_{\mc C}(y) \supset \varphi( D^\gamma_i (T) )$ holds for every $\gamma<\omega_1$.
\end{lemma}

\begin{proof}
Let $\mc C$, $T$, $y$ and $\varphi$ be as in the~lemma.
Since $y\in \bigcap_T \overline{C(\varphi(t))}^Y$, we have $\varphi (T) \subset S_{\mc C}(y)$ by definition of $S_{\mc C}(y)$. Since $S_{\mc C}(y)$ is a tree (by monotonicity of $\mc C$), we even get
\[ \D^\gamma(\varphi(T)) \subset \cltr{\varphi(T)} \subset S_{\mc C}(y) \]
for every $\gamma<\omega_1$.
Consequently, it suffices to prove the~following implication for every $\gamma<\omega_1$
\begin{equation}\label{equation: varphi induction hypothesis}
t\in D^\gamma_i (T) \implies \varphi(t) \in \D^\gamma(\varphi(T)) .
\end{equation}

We proceed by induction over $\gamma$.
For $\gamma=0$, \eqref{equation: varphi induction hypothesis} holds trivially:
\[ t\in D^0_i (T) = \cltr{T} = T \implies \varphi(t) \in \varphi(T) \subset \cltr{\varphi(T)} = \D^0 (\varphi(T)) .\]

$\gamma\mapsto \gamma+1$:
Let $t\in D^{\gamma+1}_i(T) = D_i ( D^\gamma_i(T) )$.
By definition of $D_i$, there exist infinitely many $n\in\omega$ for which $t\ext n \in D^\gamma_i(T)$.
By the~induction hypothesis, we have $\varphi(t\ext n) \in \D^\gamma (\varphi(T))$ for any such $n$.
Consider the~tree
\[ S := \cltr{\{\varphi(t\ext n) | \ n\in \omega, \ \varphi(t\ext n) \in \D^\gamma (\varphi(T)) \} } \subset \D^\gamma (\varphi(T)) \subset S_{\mc C}(y) .\]
Since $\varphi$ is admissible, we have $\varphi(t) \sqsubset \varphi(t\ext n)$ and $|\varphi(t\ext n)| = | \varphi(t)| + n$ for every $n\in \omega$.
It follows that $S$ is infinite -- since $S_{\mc C}(y)\in \textnormal{WF}$ holds by Lemma~\ref{lemma: IF trees} -- well-founded.
By König's lemma, there is some $s\in\seq$ and an infinite set $M\subset \omega$ such that $\{ s\ext m | \ m\in M \} \subset S$.
For every $m\in M$, we denote by $n_m$ the~integer for which $\varphi(t\ext n) \sqsupset s\ext m$.
The set $\{ \varphi(t\ext n_m) | \ m\in M \}$ then witnesses that $S$ contains infinitely many incomparable extensions of $\varphi(t)$. It follows that
\[ \varphi(t) \in\D(S) \subset \D (\D^\gamma(\varphi(T))) = \D^{\gamma+1}(\varphi(T)) .\]

Let $\lambda<\omega_1$ be limit and suppose that \eqref{equation: varphi induction hypothesis} holds for every $\gamma<\lambda$. By definition of $\D^\lambda$ and $D_i^\lambda$, \eqref{equation: varphi induction hypothesis} holds for $\lambda$ as well:\errata{}
\begin{align*}
t\in D^\lambda_i (T) = \bigcap_{\gamma<\lambda} D^\gamma_i (T)
\overset{\eqref{equation: varphi induction hypothesis}}{\implies}
\forall \gamma<\lambda : \varphi(t) \in \D^\gamma(\varphi(T)) \implies \\
\implies \varphi \err{(t)} \in \bigcap_{\gamma<\lambda} \D^\gamma (\varphi(T)) = \D^\lambda (\varphi(T)) .
\end{align*}
\end{proof}

We can now finish the~proof of Proposition~\ref{proposition: regular representations and D}:

\begin{proof}[Proof of Proposition~\ref{proposition: regular representations and D}, the~``$\subset$'' part]
Let $\alpha<\omega_1$ and let $\mc C$ be a Suslin scheme satisfying $\mc A(\overline{\mc C}^Y\!\!) = X$.
Recall that by definition of $\mathbf{R}_\alpha(\cdot)$, we have $\mathbf{R}_\alpha(\cdot) = \mathbf{R}_{T^c_\alpha}(\cdot)$.

(i): Suppose that $\alpha$ is even.
The tree $T^c_\alpha \overset{\text{def.}}= T_{\alpha'}$ is constructed to satisfy $r_i(T_{\alpha'}) = \alpha'$, so we have $D^{\alpha'}_i (T_{\alpha'}) \neq \emptyset$.
In particular, $D^{\alpha'}_i (T_{\alpha'})$ contains the~empty sequence. By Lemma~\ref{lemma: admissible mappings and rank}, $\D^{\alpha'} (S_{\mc C}(y))$ contains $\varphi(\emptyset)$. In particular, $\D^{\alpha'} (S_{\mc C}(y))$ is non-empty, which proves the~inclusion ``$\subset$'' in (i).

(ii): Suppose that $\alpha$ is odd.
The tree $T^c_\alpha$ is defined as $T^c_\alpha \overset{\text{def.}}= \{\emptyset\} \cup 1\ext T_{\alpha'}$.
Since $\emptyset\in D^{\alpha'}_i (T_{\alpha'})$, it follows that $D^{\alpha'}_i (1\ext T_{\alpha'})$ contains the~sequence $1\ext \emptyset = (1)$.
By Lemma~\ref{lemma: admissible mappings and rank}, $\D^{\alpha'} (S_{\mc C}(y))$ contains $\varphi( (1) )$.
Since $\varphi$ is admissible, we have $|\varphi((1))| = 1$, which implies $\emptyset\neq \varphi((1))$.
This shows that $\D^{\alpha'} (S_{\mc C}(y))$ contains some other sequence than $\emptyset$, and therefore proves the~inclusion ``$\subset$'' in (ii).

The proof of (iii) is the~same as the~proof of (ii), except that we get $(i)\in \D^{\alpha'} (S_{\mc C}(y))$ instead of $(1)\in \D^{\alpha'} (S_{\mc C}(y))$, and the~admissibility of $\varphi$ gives $|\varphi((i))| = i$.
\end{proof}
\subsection{Existence of a~Regular Representation for \texorpdfstring{$\mc F$}{F}-Borel Sets} \label{section: existence of regular representation}

The goal of this section is to prove the~following theorem:

\begin{theorem}[Existence of regular $\Fa$-representations]\label{theorem: existence of regular representation}
Let $X$ be a~$\mc K$-analytic space and $\alpha\leq \omega_1$. For $Y \supset X$, the~following conditions are equivalent:
\begin{enumerate}[(1)]
\item $X\in\Fa(Y)$; \label{case: X Fa in Y}
\item $X=\mathbf{R}_\alpha(\mc C)$ for some closed Suslin scheme in $Y$;
	\label{case: R(C) with closed C}
\item $X$ has a~regular $\Fa$-representation in $Y$ -- that is, $X=\mathbf{R}_\alpha(\mc C,Y)$ for some Suslin scheme in $X$;
	\label{case: regular representation}
\item $X=\mathbf{R}_\alpha(\mc C,Y)$ for some complete Suslin scheme on $X$.
	\label{case: regular representation by complete C}
\end{enumerate}
\end{theorem}

First, we take a closer look at the~implications in Theorem~\ref{theorem: existence of regular representation}. The implications $\eqref{case: regular representation by complete C} \implies \eqref{case: regular representation} \implies \eqref{case: R(C) with closed C}$ are trivial,
and the~implication $\eqref{case: R(C) with closed C} \implies \eqref{case: X Fa in Y}$ follow\err{s} from Proposition~\ref{proposition: complexity of R alpha sets}.
The implication $\eqref{case: X Fa in Y} \implies \eqref{case: R(C) with closed C}$ is the~hardest one, and follows from Proposition~\ref{proposition: existence of regular F alpha representations}.
The implication $\eqref{case: R(C) with closed C} \implies \eqref{case: regular representation}$ is straightforward -- when $X=\mathbf{R}_\alpha(\mc C)$ holds for some closed Suslin scheme $\mc C$ in $Y$, we define a~Suslin scheme on $X$ as $\mc C' := (C(s) \cap X)_s$, and note that $\mathbf{R}_\alpha(\mc C',Y) = \bigcup_\varphi \bigcap_t \overline{C(\varphi(t))\cap X}^Y$ satisfies
\begin{align*}
\mathbf{R}_\alpha(\mc C) = \mathbf{R}_\alpha(\mc C) \cap X \overset{\text{def.}}=
 \bigcup_\varphi \bigcap_t C(\varphi(t)) \cap X
 \subset  \\
 \subset \bigcup_\varphi \bigcap_t \overline{C(\varphi(t))}^Y
 \overset{\mc C \text{ is}}{\underset{\text{closed}}=} \bigcup_\varphi \bigcap_t C(\varphi(t))
 \overset{\text{def.}}= \mathbf{R}_\alpha(\mc C)
\end{align*}
 (where the~unions are taken over all admissible mappings $\varphi : T^c_{\alpha} \rightarrow \seq$, and the~intersections over $t\in T^c_\alpha$).
Lastly, any regular representation can be made complete by Lemma~\ref{lemma: completing representation} (provided that $X$ is $\mc K$-analytic), which proves the~implication $\eqref{case: regular representation} \implies \eqref{case: regular representation by complete C}$.

Next, we aim to prove the~implication $\eqref{case: X Fa in Y} \implies \eqref{case: R(C) with closed C}$.
The following lemma shows that \eqref{case: X Fa in Y} and \eqref{case: R(C) with closed C} are equivalent for $\alpha=2$. More importantly, it will be used in the~induction step when proving the~equivalence for general $\alpha$.

\begin{lemma}[Existence of regular representations for $\fsd$ sets]\label{lemma: indexing representation by Suslin scheme}
For any $\alpha<\omega_1$ and $X\in \left(\Fa(Y)\right)_{\sigma\delta}$, there exist a~Suslin scheme \errata{$\mc X = \{X_s| \ s\in\seq\}$} in $Y$ satisfying
\begin{enumerate}[(i)]
\item $X=\bigcap_{m\in\omega} \bigcup_{s\in\omega^m} X_s$,
\item $\left( \forall s\in\seq \right) : X_s \in \Fa(Y)$,
\item $\mc X$ \err{$\mc A$-}covers $X$.
\end{enumerate}
\end{lemma}

\err{\noindent By Remark~\ref{rem:R_n}, $\alpha=0$ indeed gives a regular $\mc F_2$-representation of $X$ in $Y$.}

\begin{proof}
Let $\alpha$ and $X$ be as above.
Since $X$ belongs to $\left(\Fa(Y)\right)_{\sigma\delta}$, there exists some countable families $\mc P_m\subset \mc \Fa(Y)$, such that $X=\bigcap_{m\in\omega} \bigcup \mc P_m$. 
Using the~notation $\mc A \land \mc B = \{ A \cap B | \ A\in \mc A, B\in \mc B \}$, we set
\[ \mc R_m:= \mc P_0 \land \mc P_1 \land \dots \land \mc P_m .\]
Clearly, we have $X=\bigcap_{m\in\omega} \bigcup \mc R_m$ and $\mc R_m \subset \Fa(Y)$. Since $\mc P_m$ are all countable, we can enumerate them as $\mc P_m =: \{ P^m_n | \ n\in\omega \}$. We then get
\[ \mc R_m = \{ P^0_{n_0} \cap \dots \cap P^m_{n_m} | \ n_i \in \omega, \ i\leq m \} = \{ P^0_{s(0)} \cap \dots P^m_{s(m)} | \ s\in \omega^{m+1} \} . \]
Denoting $X_s := P^0_{s(0)} \cap \dots \cap P^m_{s(m)}$ for $s=(s(0),\dots,s(m))\in\seq$, we get a~Suslin scheme $\mc X$ in $Y$ which satisfies the~conditions $(i)$ and $(ii)$.
Moreover, any $\mc P_m$ is a~cover of $X$, so $\mc X$ \err{$\mc A$-}covers $X$.
\end{proof}

The following proposition shows that the~equivalence between \eqref{case: X Fa in Y} and \eqref{case: R(C) with closed C} holds for $\alpha=4$. While the~proposition is not required for the~proof of Theorem~\ref{theorem: existence of regular representation}, we include it to demonstrate the~method of the~proof in a~less abstract setting. The proof also introduces the~notation necessary to get the~general result.

\begin{proposition}[Existence of regular representations for $\mc F_{\sigma\delta\sigma\delta}$ sets] \label{proposition: representation for F 4 sets}
A set $X\subset Y$ is $\mc F_{\sigma\delta\sigma\delta}$ in $Y$ if and only if there is a~closed Suslin Scheme $\mc C$ in $Y$ s.t.
\[ X = \bigcap_{m\in\omega} \bigcup_{s\in \omega^m} \bigcap_{n\in\omega} \bigcup_{t\in \omega^n} C(s\ext t). \]
\end{proposition}

\err{\noindent By Remark~\ref{rem:R_n}, this indeed gives a regular $\mc F_4$-representation of $X$ in $Y$.}

\begin{proof}
Let $X$ be a~subset of $Y$. Since the~implication from right to left is immediate, we need to prove that if $X$ is an $\mc F_{\sigma\delta\sigma\delta}$ subset of $Y$, it has a~representation with the~desired properties.

By Lemma~\ref{lemma: indexing representation by Suslin scheme}, there exists a~Suslin scheme $\mc X = \{X_s | \ s\in\seq \} \subset F_{\sigma\delta}(Y)$ \err{$\mc A$-}covering $X$ such that
\[ X = \bigcap_{m\in\omega} \bigcup_{s\in \omega^m} X_s . \]
Using Lemma~\ref{lemma: indexing representation by Suslin scheme} once more on each $X_s$, we obtain Suslin schemes $\mc X_s = \{ X_s^t | \ t\in\seq \}$ such that
	\begin{itemize}
	\item for each $s$, we have $ X_s = \bigcap_{n\in\omega} \bigcup_{t\in \omega^n} X_s^t $,
	\item for each $s$ and $t$, $X_s^t$ is closed in $Y$,
	\item for each $s$,  $\{ X_s^t | \ t\in\seq \}$ \err{$\mc A$-}covers $X_s$.
	\end{itemize}
It follows that $X$ can be written as
\[ X = \bigcap_{m\in\omega} \bigcup_{s\in \omega^m} \bigcap_{n\in\omega} \bigcup_{t\in \omega^n} X_s^t. \]

The idea behind the~proof is the~following: If it was the~case that $X_s^t=X_u^v$ held whenever $s\ext t = u \ext v$, we could define $C(s\ext t)$ as $X_s^t$, and the~proof would be finished.
Unfortunately, there is no reason these sets should have such a~property. However, by careful refining and re-indexing of the~collection $\{ X_s^t | \ s,t\in\seq \}$, we will be able to construct a~new family $\{ C_s^t | \ s,t\in\seq \}$ for which the~above condition will hold. And since the~condition above holds, we will simply denote the~sets as $C(s\ext t)$ instead of $C_s^t$.

We now introduce the~technical notation required for the~definition of $\mc C$. For a~less formal overview of what the~notation is about, see Figure~\ref{figure: notation for representation theorem}.
\begin{notation}\label{notation:sequences_of_sequences}
\begin{itemize}
	\item $S_0:=\{\emptyset\}$, $S_m:=\omega^2 \times \omega^3 \times \dots \times \omega^{m+1}$ for $m\geq 		1$ and $S:=\bigcup_{m\in\omega} S_m$,
	\item Elements of \err{$S_m$} will be denoted as
		\begin{equation*} \begin{split}
		\vec s_m 	& = \left(s_1,s_2,\dots,s_m\right)  \\
				& = \left( \left(s_1^0, s_1^1\right), \left(s_2^0, s_2^1, s_2^2\right), \dots, \left(s_m^0, s_m^1, \dots, s_m^m\right) \right) ,
		\end{split} \end{equation*}
		where $s_k \in \omega^{k+1}$ and
		$s^l_k\in\omega$. By length $|\vec s_m|$ of $\vec s_m$ we will understand the~number of sequences 			it contains -- in this case `$m$'. If there is no need to specify the~length, we will denote an 			element of $S$ simply as $\vec s$.
	\item $\pi_k:\omega^k\rightarrow \omega$ will be an arbitrary fixed bijection (for $k\geq 2$) and we 			set $\pi:=\bigcup_{k=2}^\infty \pi_k$.
	\item We define $\varrho_0: \emptyset \in S_0 \mapsto \emptyset$, for $m\geq 1$ we set
		$$\varrho_m : \vec s_m \in S_m \mapsto \left( \pi_2 (s_1), \pi_3 (s_2), \dots, \pi_{m+1} (s_m) 			\right) \in \omega^m$$
		and we put these mappings together as $\varrho:=\bigcup_{m\in\omega} \varrho_m : S \rightarrow \seq 		$.
	\item For $k\in\omega$ we define the~mapping $\Delta_k$ as
		$$\Delta_k : \vec s \in \bigcup_{m\geq k} S_m \mapsto \left( s_1^1, s^2_2,  \dots, s^k_k 			\right) \in \omega^k$$
		(where $\Delta_0$ maps any $\vec s \in S$ to the~empty sequence $\emptyset$).
	\item We also define a~mapping $\xi_k$ for $k\in\omega$:
		$$\xi_k : \vec s \in \bigcup_{m\geq k} S_m \mapsto \left( s^k_{k+1}, s^k_{k+2},
		\dots,s_{|\vec s|}^k \right) \in \omega^{|\vec s|-k}.$$
\end{itemize}

\begin{figure}
\definecolor{Gray}{gray}{0.9}
\definecolor{DarkGray}{gray}{0.6}
\definecolor{Cyan}{rgb}{0.88,1,1}
\begin{tabular}{r l l l l l l l l l}
$\omega^0 \ni$ & \!\!\!\!``$s_0$''=	& \ \cellcolor{DarkGray}$\emptyset$ & & &		& & & $\longrightarrow \emptyset$ & \\
$\omega^2 \ni$ & 	$s_1=$	& $(s^0_1$ 	& \cellcolor{DarkGray} $s^1_1)$	& &	& & & $\longrightarrow \pi(s_1)$ & $\in \omega$ \\
$\omega^3 \ni$ & 	$s_2=$	& $(s^0_2$	& $s^1_2$	& \cellcolor{DarkGray} $s^2_2)$ & & & & $\longrightarrow \pi(s_2)$ & $\in \omega$ \\
$\omega^4 \ni$ & 	$s_3=$	& $(s^0_3$	& $s^1_3$	& \cellcolor{Gray}$s^2_3$	& $s^3_3)$	& & & $\longrightarrow \pi(s_3)$ & $\in \omega$ \\
$\omega^5 \ni$ & 	$s_4=$	& $(s^0_4$	& $s^1_4$	& \cellcolor{Gray}$s^2_4$	& $s^3_4$	& $s^4_4)$ & & $\longrightarrow \pi(s_4)$ & $\in \omega$ \\
& \dots \ \ \ \ 	&		&			& \cellcolor{Gray} \dots	&			& & & \ \ \dots & \\
$\omega^{m+1} \ni$ & $s_m=$ 	& $(s^0_m$	& $s^1_m$	& \cellcolor{Gray} $s^2_m$	& $s^3_m$	& ...	& $s^m_m)$ & $\longrightarrow \pi(s_m)$ & $\in \omega$ \\
			\\
$S_m \ni$ & $\vec s_{m}$ \ \ \ &		&			&			&			& & & $\longrightarrow \varrho(\vec s_{m})$ & $\in \omega^m$ 
\end{tabular}
\caption{An illustration of the~notation from the~proof of Proposition~\ref{proposition: representation for F 4 sets}.
Each sequence $s_k\in\omega^{k+1}$ is mapped to an integer $\pi(s_k)$  via a~bijection.
This induces a~mapping of the~``sequence of sequences'' $\vec s_m$ to a~sequence $\varrho(\vec s_m)$ of length $|\vec s_m|=m$.
The diagonal sequence highlighted in dark gray corresponds to $\Delta_2(\vec s_{m})$. Note that the~letter $\Delta$ stands for `diagonal' and the~lower index ($2$ in this case) corresponds to the~length of this sequence. The part of the~column highlighted in light gray corresponds to the~sequence $\xi_2(\vec s_{m})$.
\label{figure: notation for representation theorem}}
\end{figure}

\end{notation}

Suppose that $w\in \omega^m$ satisfies $w=\varrho(\vec s)$ for some $\vec s \in S_m$. Without yet claiming that this correctly defines a~Suslin scheme, we define $C(w)$ as
\begin{equation} \begin{split} \label{equation: definition of C}
C\left(\varrho\left(\vec s\right)\right) \ 
:= \ & \bigcap_{k=0}^{|\vec s|} X^{\xi_k(\vec s)}_{\Delta_k(\vec s)} \ 
 = \ X^{\xi_0(\vec s)}_{\Delta_0(\vec s)} \cap X^{\xi_1(\vec s)}_{\Delta_1(\vec s)} \cap \dots \cap X^{\xi_{|\vec s|}(\vec s)}_{\Delta_{|\vec s|}(\vec s)} \\
 = \ & X_\emptyset^{s_1^0 s_2^0 \dots s_m^0} \cap X_{s_1^1}^{s_2^1 s_3^1 \dots s_m^1}
 	\cap X_{s_1^1 s_2^2}^{s_3^2 s_4^2 \dots s_m^2} \cap \dots \cap X^\emptyset_{s_1^1 s_2^2 \dots s_m^m} .
\end{split} \end{equation}

In order to show that $\mc C$ indeed does have the~desired properties, we first note the~following properties of $\varrho$, $\Delta_k$ and $\xi_k$ (all of which immediately follow from the~corresponding definitions).
\begin{lemma} \label{lemma: properties of auxiliary functions}
The functions $\varrho$, $\Delta_k$ and $\xi_k$ have the~following properties:
\begin{enumerate}[(a)]
\item $\left(\forall \vec s \in S \right) \left(\forall k\leq |\vec s|\right) :
	\left|\Delta_k(\vec s)\ext \xi_k(\vec s)\right| = \left| \vec s \right| = \left| \varrho(\vec s) \right|$.
\item $\varrho : S \rightarrow \seq$ is a~bijection.
\item $\left(\forall u,v \in \seq \right): v\sqsupset u \iff \varrho^{-1}(v) \sqsupset \varrho^{-1}(u)$.
\item $\left( \forall \vec s, \vec t \in S \right) \left( \forall k \leq |\vec s| \right) : \vec t \sqsupset \vec s \implies \Delta_k(\vec t) = \Delta_k(\vec s) \ \& \ \xi_k(\vec t) \sqsupset \xi_k(\vec s) $.
\end{enumerate}
\end{lemma}
\errata{\textbf{$\mc C$ is a Suslin scheme:}} From $(b)$ it follows that each $C(\varrho(\vec s))$ is well defined and that $C(w)$ is defined for every $w\in\seq$.
To see that $\mc C$ is a~Suslin scheme, we need to verify its monotonicity. For $u\sqsubset v\in \seq$ we have
\begin{equation*} \begin{split}
C(v) =
	\ & \bigcap_{k=0}^{|\varrho^{-1}(v)|} X^{\xi_k(\varrho^{-1}(v))}_{\Delta_k(\varrho^{-1}(v))} 
	\ \overset{(c)}\subset \ \bigcap_{k=0}^{|\varrho^{-1}(u)|} X^{\xi_k(\varrho^{-1}(v)) }_{\Delta_k(\varrho^{-1}(v))} \\
	\overset{(c),(d)}{=} & \ \bigcap_{k=0}^{|\varrho^{-1}(u)|} X^{\xi_k(\varrho^{-1}(v))}_{\Delta_k(\varrho^{-1}(u))}
	\ \subset \ \bigcap_{k=0}^{|\varrho^{-1}(u)|} X^{\xi_k(\varrho^{-1}(u))}_{\Delta_k(\varrho^{-1}(u))} = C(u) ,
\end{split} \end{equation*}
where the~first and last identities are just the~definition of $C(\cdot)$ and the~last inclusion holds because we have $\xi_k(\varrho^{-1}(v)) \sqsupset \xi_k(\varrho^{-1}(u))$ (by (c) and (d)) and $\mc X_{\Delta_k(\varrho^{-1}(u))} = (X^t_{\Delta_k(\varrho^{-1}(u))})_t$ is a~Suslin scheme (and hence monotone).

Next, we will show that
$$ X \subset \bigcup_{\sigma\in\baire} \bigcap_{m\in\omega} C(\sigma|m) = \mc A(\mc C)
 = \mathbf{R}_{\omega_1}(\mc C) \subset \mathbf{R}_4(\mc C)
 = \bigcap_{m\in\omega} \bigcup_{s\in \omega^m} \bigcap_{n\in\omega} \bigcup_{t\in \omega^n} C(s\ext t)
 .$$
Going from right to left, the~identities (resp. inclusions) above follow from: definition of $\mathbf{R}_4$, Lemma~\ref{lemma: basic properties of R T sets}\,\eqref{case: R T and embeddings}, Lemma~\ref{lemma: basic properties of R T sets}\,\eqref{case: R T and IF}, and definition of $\mc A(\cdot)$; it remains to prove the~first inclusion.

\errata{\textbf{$X$ is $\mc A$-covered by $\mc C$:}} Let $x\in X$.
Our goal is to produce $\mu \in \baire$ such that $x\in C(\mu|m)$ holds for each $m\in\omega$. We shall do this by finding a~sequence $\vec s_0 \sqsubset \vec s_1 \sqsubset \dots $, $\vec s_m \in S_m$, for which which $x\in \bigcap_m C(\varrho(\vec s_m))$. Once we have done this, $(c)$ ensures that $\mu|m := \varrho (\vec s_m)$ correctly defines the~desired sequence $\mu$.

To this end, observe that since $\{X_s | \ s\in\seq \}$ $\mc A$-covers $X$, there exists a~sequence $\sigma\in \baire$ such that $x\in \bigcap_k X_{\sigma|k}$. Similarly, $\{X_{\sigma|k}^t | \ t\in\omega \}$ $\mc A$-covers $X_{\sigma|k}$ for any $k\in\omega$, so there are sequences $\nu_k \in \baire$, $k\in\omega$, for which $x\in \bigcap_m X_{\sigma|k}^{\nu_k|m}$. For $m\in\omega$ we denote
\[ s_m := \left( \nu_0(m-1), \nu_1(m-2), \dots, \nu_k (m-1-k), \dots, \nu_{m-1}(0), \sigma (m-1) \right) \]
and define $\vec s_m := (s_1, s_2, \dots , s_m )$. For $k\leq m$ we have
\begin{align*}
\Delta_k(\vec s_m) = & ( s_1^1, \dots , s_k^k ) = ( \sigma(0), \sigma(1), \dots, \sigma(k-1) ) = \sigma|k ,\\
\xi_k(\vec s_m) = & ( s_{k+1}^k, s_{k+2}^k, \dots, s_m^k ) \\
 = & \left( \nu_k(k+1-1-k), \nu_k(k+2-1-k), \dots , \nu_k(m-1-k) \right) \\
 = & \left( \nu_k(0), \nu_k(1), \dots , \nu_k(m-k-1) \right) = \nu_k|m-k.
\end{align*}
It follows that $X_{\Delta_k(\vec s_m)}^{\xi_k(\vec s_m)} = X_{\sigma|k}^{\nu_k|m-k} \ni x$ holds for every $k\leq m$. As a~consequence, we get $C(\varrho(\vec s_m)) = \bigcap_{k=0}^m X_{\sigma|k}^{\nu_k|m-k} \ni x$ for every $m\in\omega$, which shows that $\mc C$ $\mc A$-covers $X$.

In order to finish the~proof of the~theorem, it remains to show that
\[ \bigcap_{m\in\omega} \bigcup_{u\in \omega^m} \bigcap_{n\in\omega} \bigcup_{v\in \omega^n} C(u\ext v)
\subset \bigcap_{m\in\omega} \bigcup_{s\in \omega^m} \bigcap_{n\in\omega} \bigcup_{t\in \omega^n} X_s^t = X
. \]
To get this inclusion, it is enough to prove that
\[ \left( \forall u \in \seq \right) \left( \exists s_u \in \omega^{|u|} \right)
	\left( \forall v \in \seq \right) \left( \exists t_{uv} \in \omega^{|v|} \right) : C(u\ext v) \subset X_{s_u}^{t_{uv}} . \]
We claim that a~suitable choice is $s_u := \Delta_{|u|}(\varrho^{-1}(u))$ and $t_{uv}:= \xi_{|u|}(\varrho^{-1}(u\ext v))$. Indeed, for any $u\in\omega^m$ and $v\in\omega^n$ we have
\[ C(u\ext v) \overset{\textrm{def.}}=
	\bigcap_{k=0}^{m+n} X^{\xi_k(\varrho^{-1}(u\ext v))}_{\Delta_k(\varrho^{-1}(u\ext v))}
	\subset X^{\xi_m(\varrho^{-1}(u\ext v))}_{\Delta_m(\varrho^{-1}(u\ext v))}
	\overset{(d)}= X^{\xi_m(\varrho^{-1}(u\ext v))}_{\Delta_m(\varrho^{-1}(u))} = X_{s_u}^{t_{uv}} . \]
\end{proof}

The following proposition proves the~implication \eqref{case: X Fa in Y}$\implies$\eqref{case: R(C) with closed C} of Theorem~\ref{theorem: existence of regular representation}.

\begin{proposition}[Existence of regular representations]
	\label{proposition: existence of regular F alpha representations}
Let $X\subset Y$ and $\alpha\leq \omega_1$. Then $X\in\Fa(Y)$ if and only if $X=\mathbf{R}_\alpha(\mc C)$ holds for some Suslin scheme $\mc C$ which is closed in $Y$ and $\mc A$-covers $X$.
\end{proposition}
\begin{proof}
The implication ``$\Longleftarrow$'' follows from Proposition~\ref{proposition: complexity of R alpha sets}, so it remains to prove the~implication ``$\implies$''.

For $\alpha=0$, $X$ is closed and any Suslin scheme $\mc C$ satisfies $\mathbf{R}_0(\mc C) = C(\emptyset)$. Consequently, it suffices to set $C(u):=\overline{X}^Y$ for every $s\in \seq$.
We already have the~statement for $\alpha=2$ (Lemma~\ref{lemma: indexing representation by Suslin scheme}), $\alpha=4$ (Proposition~\ref{proposition: representation for F 4 sets}) and $\alpha=\omega_1$ (\errata{\eqref{case: R T and IF}} from Lemma~\ref{lemma: basic properties of R T sets}).

Suppose that the~statement holds for an even ordinal $\alpha<\omega_1$ and let $X\in \mc F_{\alpha+1}(Y)$. By induction hypothesis we have
\[ X= \bigcup_{m\in\omega} X_m = \bigcup_{m\in\omega} \mathbf{R}_\alpha (\mc C_m) \]
for some Suslin schemes $\mc C_m$, $m\in\omega$, which are closed in $Y$ and $\mc A$-cover $X_m$.
We define $\mc C$ as $C(\emptyset):= \overline{X}^Y$ and, for $m\in\omega$ and $t\in \seq$, $C(m\ext t):=\overline{X}^Y \cap C_m(t)$.
Clearly, $\mc C$ is a~Suslin scheme which is closed in $Y$ and $\mc A$-covers $X$.
Recall that for odd $\alpha+1$ we have $T^\textrm{c}_{\alpha+1}=\{\emptyset\} \cup 1 \ext T_\alpha^\textrm{c}$.
For $x\in Y$ we have\errata{}
\begin{equation*} \begin{split}
x \in X & \iff \left(\exists m\in \omega \right) : x \in \mathbf{R}_\alpha(\mc C_m) \\
		& \iff \left(\exists m \in \omega \right) \left(\exists \varphi : T^\textrm{c}_\alpha \rightarrow \seq \textrm{ adm.}\right) : \\
		& \ \ \ \ \ \ \ \ \ \ \ x \in \bigcap_{t'\in T^\textrm{c}_\alpha} C_m(\varphi(t')) = \bigcap_{t'\in T^\textrm{c}_\alpha} C(m\ext \varphi(t')) \\
		& \iff \left(\exists \psi : 1\ext T^\textrm{c}_\alpha \rightarrow \seq \textrm{ adm.}\right) : \\
		& \ \ \ \ \ \ \ \ \ \ \ x \in \bigcap_{1\ext t'\in 1\ext T^\textrm{c}_\alpha} C(\psi(\err{1} \ext t')) \ \& \ x\in \overline{X}^Y = C(\emptyset) \\
		& \iff \left(\exists \psi : T^\textrm{c}_{\alpha+1} \rightarrow \seq \textrm{ adm.}\right) : x \in \bigcap_{t\in T^\textrm{c}_{\alpha+1}} C(\psi(t)) \\
		& \iff x \in \mathbf{R}_{\alpha+1}(\mc C).
\end{split} \end{equation*}

Finally, suppose that $X\in \Fa(Y)$ holds for $0<\alpha<\omega_1$ even and that the~statement holds for all $\beta<\alpha$.

If $\alpha$ is a~successor ordinal, we can use Lemma~\ref{lemma: indexing representation by Suslin scheme} to obtain a~Suslin scheme $\{ X_s | \ s\in\seq \}\subset \mc F_{\alpha-2}(Y)$ satisfying the~conditions (i)-(iii) from the~lemma.
For each $m\in\omega$, we set $\alpha_m:=\alpha-2$ and define
\begin{equation}\label{equation:T}
T:=\{\emptyset \} \cup \bigcup_{m\in\omega} m\ext T_{\alpha_m}^\textrm{c} .
\end{equation}
Since $T = T^c_\alpha$, we trivially have $\mathbf{R}_T(\cdot) = \mathbf{R}_\alpha(\cdot)$ and the~definition of $\alpha_m$ is chosen to satisfy the~following:
\begin{equation}\label{equation: X s is in F alpha m}
\left( \forall s \in \seq \right) : X_s \in \mc F_{\alpha_{|s|}}(Y) .
\end{equation}

For limit $\alpha$, we can repeat the~proof of Lemma~\ref{lemma: indexing representation by Suslin scheme} to obtain a~sequence $(\alpha_m)_m$ and a~Suslin scheme $\{ X_s | \ s\in\seq \}$ which satisfies (i)-(iii) from Lemma~\ref{lemma: indexing representation by Suslin scheme}, \eqref{equation: X s is in F alpha m} and $\alpha_m < \alpha$ for each $m\in\omega$.
Moreover, we can assume without loss of generality that $\sup_m \alpha_m = \alpha$.
Let $T$ be the~tree defined by \eqref{equation:T}. 
We claim that the~trees $T^\textrm{c}_\alpha$ and $T$ are equivalent in the~sense of \eqref{case: R T and embeddings} from Lemma~\ref{lemma: basic properties of R T sets}, so that we have $\mathbf{R}_T(\cdot ) = \mathbf{R}_\alpha(\cdot)$ even when $\alpha$ is limit.

Indeed, the~sequence $(\pi_\alpha(n))_n$, defined in Notation~\ref{notation: trees of height alpha}, converges to $\alpha$ with increasing $n$, so for every $m\in\omega$ there exists $n_m\geq m$ such that $\pi_\alpha(n_m) \geq \alpha_{|m|}$. Since any canonical tree can be embedded (in the~sense of Lemma~\ref{lemma: basic properties of R T sets}\,\eqref{case: R T and embeddings}) into any canonical tree with a~higher index, $m\ext T_{\alpha_m}^\textrm{c}$ can be embedded by some $f_m$ into $n_m \ext T_{\pi_\alpha (n_m)}^\textrm{c}$. Defining $f_\emptyset : \emptyset \mapsto \emptyset$ and $f:=f_\emptyset \cup \bigcup_{m\in\omega} f_m$, we have obtained the~desired embedding of $T$ into $T^\textrm{c}_\alpha$. The embedding of $T^\textrm{c}_\alpha$ into $T$ can be obtained by the~same method.

Since $\alpha_{|s|} < \alpha$ for every $s\in \seq$, we can apply the~induction hypothesis  and obtain a~Suslin scheme $\mc C_s$ which is closed in $Y$, $\mc A$-covers $X_s$, and satisfies $X_s = \mathbf{R}_{\alpha_m}\left(\mc C_s\right)$.
This in particular gives
\[ X = \bigcap_{m\in\omega} \bigcup_{s\in \omega^m} X_s
	= \bigcap_{m\in\omega} \bigcup_{s\in \omega^m} \mathbf{R}_{\alpha_m}\left(\mc C_s\right) . \]
\errata{Using Notation~\ref{notation:sequences_of_sequences}, we define $\mc C = (C(w))_{w\in\baire}$ as}
\[ C\left(\varrho(\vec s)\right) := \bigcap_{k=0}^{|\vec s|} C_{\Delta_k (\vec s)}\left( \xi_k(\vec s) \right) . \]
\errata{Note that this is essentially the same definition as \eqref{equation: definition of C}.
In particular, the arguments presented below Lemma~\ref{lemma: properties of auxiliary functions} in the proof of Proposition~\ref{proposition: representation for F 4 sets}, can be repeated word by word (with $X_{\Delta_k (\vec s)}^{\xi_k(\vec s)}$ replaced by $C_{\Delta_k (\vec s)}\left( \xi_k(\vec s) \right)$ ) to show that $\mc C$ is a Suslin scheme which $\mc A$-covers $X$. Clearly, $\mc C$ is closed in $Y$.}

To finish the~proof, it remains to show that $\mathbf{R}_\alpha(\mc C)\subset X$, which reduces to proving the~inclusion in the~following formula:
\[ \mathbf{R}_\alpha(\mc C) = \mathbf{R}_T(\mc C) \subset \bigcap_{m\in\omega} \bigcup_{s\in\omega^m} \mathbf{R}_{\alpha_m}(\mc C_s) = X . \]
Using the~definition of $\mathbf{R}_T(\cdot)$ and the~fact that $T^m := \{ t' \in \seq | \ m\ext t' \in T \} = T^\textrm{c}_{\alpha_m}$ for any $m\in\omega$, we get
\begin{equation*} \begin{split}
x \in \mathbf{R}_T(\mc C) \iff & \left( \exists \varphi : T \rightarrow \seq \text{ adm.}\right) : x \in \bigcap_{t\in T} C(\varphi(t)) \\
\underset{(iii)}{\overset{L\ref{lemma: construction of admissible mappings}}\iff} & \left( \forall m\in \omega \right) \left( \exists s\in\omega^m \right) \left( \exists \varphi_m : T^m \rightarrow \seq \text{ adm.} \right) : \\ 
	& \ \ x \in \bigcap_{t'\in T^m} C(s\ext \varphi_m(t')) \\
\iff & \left( \forall m\in\omega \right) \left( \exists s\in \omega^m \right) \left( \exists \varphi_m : T^\textrm{c}_{\alpha_m} \rightarrow \seq \text{ adm.} \right) : \\
	& \ \ x \in \bigcap_{t'\in T^\textrm{c}_{\alpha_m}} C(s\ext \varphi_m(t'))		\\
\iff & \left( \forall m\in\omega \right) \left( \exists s\in\omega^m \right) : 
		x \in \mathbf{R}^{s}_{\alpha_m} (\mc C) \\
\iff & x \in \bigcap_{m\in\omega} \bigcup_{s\in\omega^m} \mathbf{R}^{s}_{\alpha_m}(\mc C) .
\end{split} \end{equation*}
We claim that for each $s\in\omega^m$ there exists $\tilde{s}\in \omega^m$ such that $\mathbf{R}^s_{\alpha_m}(\mc C) \subset \mathbf{R}_{\alpha_m}(\mc C_{\tilde{s}})$. Once we prove this claim, we finish the~proof by observing that
\[ \mathbf{R}_T(\mc C) = \bigcap_{m\in\omega} \bigcup_{s\in\omega^m} \mathbf{R}^{s}_{\alpha_m}(\mc C)
	\subset \bigcap_{m\in\omega} \bigcup_{s\in\omega^m} \mathbf{R}_{\alpha_m}(\mc C_{\tilde{s}})
	\subset \bigcap_{m\in\omega} \bigcup_{\tilde s\in\omega^m} \mathbf{R}_{\alpha_m}(\mc C_{\tilde s}) = X . \]
Let $s\in\omega^m$. To prove the~claim, observe first that for any $v\in\seq$ we have
\begin{equation} \label{equation: final inequality in F alpha representation proposition} \begin{split}
 C(s\ext v) & =
 	\bigcap_{k=0}^{m+|v|}C_{\Delta_k(\rho^{-1}(s\ext v))}\left(\xi_k(\rho^{-1}(s\ext v))\right)
	\subset C_{\Delta_m(\rho^{-1}(s\ext v))}\left(\xi_m(\rho^{-1}(s\ext v))\right) \\
 & \overset{L\ref{lemma: properties of auxiliary functions}}{\underset{(d)}=} C_{\Delta_m(\rho^{-1}(s))}\left(\xi_m(\rho^{-1}(s\ext v))\right) 
	= C_{\tilde{s}}\left(\xi_m(\rho^{-1}(s\ext v))\right) 
\end{split} \end{equation}
(where we denoted $\tilde{s}:=\Delta_m(\rho^{-1}(s))\in\omega^m$).

Let $x\in \mathbf{R}^s_{\alpha_m}(\mc C)$. By definition of $\mathbf{R}^s_{\alpha_m}(\mc C)$, $x$ belongs to $\bigcap_{t\in T^\textrm{c}_{\alpha_m}} C(s\ext \varphi(t))$ for some admissible $\varphi : T^\textrm{c}_{\alpha_m} \rightarrow \seq$. Taking $v=\varphi(t)$ in the~computation above, we get
\begin{equation}\label{eq:ex_of_reg_repr_last}
x\in \bigcap_{t\in T^\textrm{c}_{\alpha_m}} C(s\ext \varphi(t)) \implies
x \in \bigcap_{t\in T^\textrm{c}_{\alpha_m}} C_{\tilde{s}}\left(\xi_m(\rho^{-1}(s\ext \varphi(t)))\right) .
\end{equation}
Since $s\in\omega^m$, it follows from Lemma~\ref{lemma: properties of auxiliary functions} that the~mapping $t\in T^\textrm{c}_{\alpha_m} \mapsto \xi_m\left((\rho^{-1}(s\ext \varphi(t))\right)$ is admissible. It follows that the~intersection on the~right side of \errata{\eqref{eq:ex_of_reg_repr_last}} is contained in $\mathbf{R}_{\alpha_m}(\mc C_{\tilde{s}})$, so \errata{$\mathbf{R}^s_{\alpha_m}(\mc C) \subset \mathbf{R}_{\alpha_m}(\mc C_{\tilde{s}})$} and the~proof is complete.
\end{proof}

To prove Theorem~\ref{theorem: existence of regular representation}, it remains the~prove that every regular representation can be ``made complete without loss of generality'', in sense of the~following proposition.
When we say that a~regular $\Fa$-representation $\mc C$ of $X$ in $Y$ is \emph{complete}, we mean that $\mc C$ is complete on $X$.
A regular $\Fa$-representation $\mc C$ is a~\emph{refinement} of a~regular $\Fa$-representation $\mc D$ when for every $s\in\seq$, there is some $t\in \seq$ with $|t|=|s|$ s.t. $C(s) \subset D(s)$.

\begin{lemma}[Completing regular representations]
	\label{lemma: completing representation}
Let $X$ be $\mc K$-analytic space and $Y\supset X$. For any regular $\Fa$-representation $\mc D$ of $X$ in $Y$, there is a~complete regular $\Fa$-representation $\mc C$ of $X$ in $Y$ which refines $\mc D$.
\end{lemma}

\begin{proof}
Suppose that a~$\mc K$-analytic space $X$ satisfies $X\in\Fa(Y)$ for some $\alpha\leq \omega_1$.
By Proposition~\ref{proposition: K analytic spaces have complete Suslin schemes}, there exists some Suslin scheme $\mc R$ which is complete on $X$.
Let $\mc D$ be a~regular $\Fa$-representation of $X$ in $Y$. We fix an arbitrary bijection $\pi : \omega^2 \rightarrow \omega$ and denote by $\varrho : \bigcup_{m\in\omega} \omega^m \times \omega^m \rightarrow \seq$,\errata{}
\begin{align*}
\err{(s,t) = \Big( \big(s(0),\dots,s(k)\big), \big(t(0),\dots,t(k)\big) \Big) \mapsto \Big( \pi\big(s(0),t(0)\big),\dots,\pi\big(s(k),t(k)\big) \Big) ,}
\end{align*}
the bijection between pairs of sequences and sequences induced by $\pi$.
Finally, we define a~complete Suslin scheme $\mc C$ as
\[ C(u) := D(\varrho^{-1}_1(u)) \cap R(\varrho^{-1}_2(u)) \ \ \text{ (or equivalently, }
 C(\varrho(s,t)) := D(s) \cap R(t) ).\]
Denoting $\mc C_m:=\{ C_s | \ s\in\omega^m \}$, we obtain a~sequence of covers $(\mc C_m)_m$ which is a~refinement of the~complete sequence $(\mc R_m)_m$.
By a~straightforward application of the~definition of completeness, $(\mc C_m)_m$ is complete on $X$ as well.
In particular, Lemma~\ref{lemma: complete suslin schemes in super spaces} yields
\[ \mathbf{R}_\alpha(\mc C,Y) = \mathbf{R}_\alpha(\overline{\mc C}^Y\!\!)
\supset \mathbf{R}_{\omega_1}(\overline{\mc C}^Y\!\!)
= \mc A(\overline{\mc C}^Y\!\!) = X .\]
For the~converse inclusion, let $x\in \mathbf{R}_\alpha (\mc C,Y)$ and suppose that $\varphi$ is an admissible mapping witnessing\footnote{Recall that an admissible mapping $\varphi$ witnesses that $x$ belongs to $\mathbf{R}_\alpha(\mc C,Y)$ if and only if its domain is the~canonical tree $T^c_\alpha$ and we have $x \in \bigcap_{T^c_\alpha} \overline{C(\varphi(t))}^Y$.} that $x$ belongs to $\mathbf{R}_\alpha (\mc C,Y)$.
The mapping $\psi := \varrho^{-1}_1 \circ \varphi$ is clearly admissible as well, and for any $u\in\seq$, we have
\[ C( \varphi(u)) = D( \varrho^{-1}_1 (\varphi(u)))
 \cap R( \varrho^{-1}_2 ( \varphi(u)))
 \subset D( \varrho^{-1}_1 (\varphi(u)))
 = D( \psi (u) ) .\]
In particular, $\psi$ witnesses that $x$ belongs to $\mathbf{R}_\alpha (\mc D,Y) = X$.
\end{proof}
\section{Complexity of Talagrand's Broom Spaces} \label{section: complexity of brooms}

In this section, we study the so-called ``broom spaces'', based on the non-absolute $\fsd$ space $\mathbf{T}$ from \cite{talagrand1985choquet} (also defined in Definition \ref{definition: AD topology} here).
In \cite{kovarik2018brooms} the author shows that for each \emph{even} $\alpha$, there is a broom space $T$ with
\begin{equation}\label{eq:T_and_X_2_alpha}
\{2,\alpha\} \subset \textnormal{Compl}(T) \subset [2,\alpha]
\end{equation}
\err{(see \cite[Theorem\,5.14]{kovarik2018brooms} and Remark~\ref{remark: F and F tilde})}.
We improve this result in two ways, by providing conditions under which the set $\textnormal{Compl}(T)$ of complexities attainable by a broom space $T$ is \emph{equal to the whole interval} $[2,\alpha]$ for \emph{any} $\alpha \in [2,\omega_1]$.

\errata{Upon closer examination of Proposition~\ref{proposition: complexities of topological sums} and the~proof of Theorem \ref{theorem: summary} that follows, we can see that \eqref{eq:T_and_X_2_alpha} already yields the existence of spaces satisfying
\[ [2,\alpha] \cap \{ \textnormal{even ordinals} \} \subset \textnormal{Compl}(X)
\subset [2,\alpha] \]
for any $2\leq \alpha \leq \omega_1$.}
The more interesting part will therefore be the methods used to obtain the proofs.
First of all, to show that a space $T$ attains many different complexities, we need to find many different compactifications of $T$.
In each such $cT$, we need to bound $\Compl{T}{cT}$ from below and from above. To obtain the lower estimate, we describe a refinement of a method from \cite{talagrand1985choquet}. For the upper bound, we use a criterion from Section \ref{section: Suslin scheme rank} (Corollary~\ref{corollary: sufficient condition for Fa}).

This section is organized as follows: In Section \ref{section: broom sets}, we describe the hierarchy of broom \emph{sets}, which are then used in Section \ref{section: broom space properties} to define corresponding hierarchy of topological spaces, called broom spaces. Section \ref{section: broom space properties} also explores the basic complexity results related to broom spaces.
Section \ref{section: Talagrand lemma} reformulates and refines some tools from \cite{talagrand1985choquet}, useful for obtaining complexity lower bounds.
In Section \ref{section: amalgamation spaces}, we study compactifications of broom spaces in the abstract setting of amalgamation spaces. We conclude with Sections \ref{section: cT} and \ref{section: dT}, where all of these results are combined in order to compute the complexities attainable by broom spaces.

We note that while reading the subsections in the presented order should make some of the notions more intuitive, it is not always strictly required; Sections \ref{section: Talagrand lemma} and \ref{section: amalgamation spaces} are completely independent, and their only relation to Section \ref{section: broom space properties} is that broom spaces constitute an example which one might wish to use to get an intuitive understanding of the presented abstract results.

\subsection{Broom sets}\label{section: broom sets}

This section introduces a special sets of finite sequences on $\omega$, called \emph{finite broom sets}, and related sets of infinite sequences on $\omega$, called \emph{infinite broom sets}.
We will heavily rely on the notation introduced in Section~\ref{section: sequences}.

\begin{definition}[Finite broom sets]\label{definition: finite brooms}
A \emph{forking sequence} is a sequence $( f_n )_{n\in\omega}$ of elements of $\seq$ such that for distinct $m,n\in\omega$ we have $f_n(0)\neq f_m(0)$.

We denote $\mc B_{0} := \left\{ \, \{ \emptyset \} \, \right\}$. For $\alpha\in [1,\omega_1]$ we inductively define the hierarchy of (finite) broom sets as
\begin{equation*}
\mc B_{\alpha} := 
\begin{cases}
\ \mc B_{<\alpha} \cup \{ h\ext B | \ B \in \mc B_{\alpha-1} , \ h\in \seq \} & \text{ for } \alpha \text{ odd} \\
\ \mc B_{<\alpha} \cup \{ \bigcup_n f_n \ext B_n | \ B_n \in \mc B_{<\alpha}, \, (f_n)_n \text{ is a forking seq.} \} & \text{ for } \alpha \text{ even.}
\end{cases}
\end{equation*}
By $h_B$ we denote the \emph{handle} of a finite broom set $B$, that is, the longest sequence common to all $s\in B$.
\end{definition}

Note that $h_B$ is either the empty sequence or the sequence $h$ from Definition \ref{definition: finite brooms}, depending on whether $B$ belongs to $\mc B_\alpha \setminus \mc B_{<\alpha}$ for even or odd $\alpha$.

For odd $\alpha$, every $B\in \mc B_\alpha$ can be written as $B = h_B \ext B'$, where $B' \in \mc B_{\alpha-1}$.
Moreover, we have $h_0 \ext B' \in \mc B_\alpha$ for any $h_0 \in \seq$.

We claim that the hierarchy of finite broom sets is strictly increasing and stabilizes at the first uncountable step (the `strictly' part being the only non-trivial one):
\begin{equation} \label{equation: broom hierarchy}
\mc B_{0} \subsetneq \mc B_{1} \subsetneq \mc B_{2} \subsetneq \dots \subsetneq \mc B_{\omega_1} = \mc B_{<\omega_1} .
\end{equation}
Indeed, for odd $\alpha$, if $h \in \seq$ is not the empty sequence and $\alpha$ is the smallest ordinal for which $B\in \mc B_{\alpha-1}$, then it follows from the above definition that $h\ext B \in \mc B_\alpha \setminus \mc B_{<\alpha}$.
For $\alpha$, if $(f_n)_n$ is a forking sequence and $\alpha$ is the smallest ordinal satisfying $\{ B_n | \ n\in \omega \} \subset \mc B_{<\alpha}$, then $\bigcup_n f_n\ext B_n$ belongs to $\mc B_\alpha \setminus \mc B_{<\alpha}$.

By extending each sequence from a finite broom set into countably many infinite sequences (in a particular way), we obtain an \emph{infinite broom set} of the corresponding rank:

\begin{definition}[Infinite broom sets]\label{definition: infinite brooms}
A countable subset $A$ of $\baire$ is \emph{an infinite broom set} if there exists some $B\in \mc B_{\omega_1}$ such that the following formula holds:
\begin{align}\label{equation: A extension}
\left( \exists \{ \nu^s_n | \, s\in B,n\in \omega \} \subset \baire \right)
\left( \exists \text{ forking sequences } (f^s_n)_n, \ s\in B \right) : \\ 
A = \{ s\ext f^s_n \ext \nu^s_n | \ s\in B, n\in\omega \} . \nonumber
\end{align}
If $A$ and $B$ satisfy \eqref{equation: A extension}, we say that $A$ is a \emph{broom-extension} of $B$.
For $\alpha\in[0,\omega_1]$, the family of $\mc A_\alpha$ consists of all broom-extensions of elements of $\mc B_\alpha$.
\end{definition}

When a broom extension $A$ of $B$ belongs to some $\mc A \subset \mc A_{\omega_1}$, we will simply say that $A$ is an $\mc A$-extension of $B$.
Note that the set $B$ of which $A$ is a broom extension is uniquely determined, so $A\in \mc A_\alpha$ holds if and only if $B\in \mc B_\alpha$.
Using the definition of broom-extension and \eqref{equation: broom hierarchy}, we get
\[ \mc A_{0} \subsetneq \mc A_{1} \subsetneq \mc A_{2} \subsetneq \dots \subsetneq \mc A_{\omega_1} = \mc A_{<\omega_1} .\]

\errata{In relation to Definition \ref{definition: infinite brooms}, we claim that for $B\in\mc B_\alpha$} we have
\begin{equation} \label{equation: B tilde is B 2 + alpha}
\widetilde B := \{ s\ext f^s_n | \ s\in B, n\in\omega \} \in \mc B_{2+\alpha} .
\end{equation}
(Recall that for infinite $\alpha$, we have $2+\alpha = \alpha$.)
It follows from \eqref{equation: B tilde is B 2 + alpha} that each $\mc A_\alpha$ is actually the family of ``normal infinite extensions'' of certain $\mc B_{2+\alpha}$-brooms.
We define $\mc A_\alpha$ as in Definition \ref{definition: infinite brooms} to simplify notation later (the previous paper of the author, \cite{kovarik2018brooms}, uses the corresponding ``alternative'' numbering).

We omit the proof of \eqref{equation: B tilde is B 2 + alpha}, since it is a simple application of transfinite induction, and the above remark is not needed anywhere in the remainder of this paper.


The following result estimates the rank of broom sets.

\begin{lemma}[Rank of broom sets]\label{lemma: rank of broom sets}
Let $\alpha\in [0,\omega_1]$.
\begin{enumerate}[1)]
\item Let $B \in \mc B_\alpha$. Then
	\begin{enumerate}[(i)]
	\item $\D^{\alpha'}(B)$ is finite and
	\item if $\alpha$ is even, then $\D^{\alpha'}(B)$ is either empty or equal to $\{\emptyset\}$.
	\end{enumerate}
\item Let $A \in \mc A_\alpha$. Then
	\begin{enumerate}[(i)]
	\item $\D^{\alpha'}( \D(A) )$ is finite and
	\item if $\alpha$ is even, then $\D^{\alpha'}( \D(A) )$ is either empty or equal to $\{\emptyset\}$.
	\end{enumerate}
\end{enumerate}
\end{lemma}

\begin{proof}
$1)$: The proof of this part is essentially the same as the proof of Proposition 4.11(i) from \cite{kovarik2018brooms}, where an analogous result is proven for derivative
\[ D_{i}(B) := \{ s\in \seq | \ \textnormal{cl}_{\Tr}(B) \text{ contains infinitely many extensions of } s \} .\]
For any $B\in \mc B_{\omega_1}$, the initial segment $\textnormal{cl}_{\Tr}(h_B)$ is finite. Moreover, for $B \in \mc B_\alpha \setminus \mc B_{<\alpha}$, $\alpha$ even, we have $\textnormal{cl}_{\Tr}(h_B) = \{ \emptyset \}$.
Therefore -- since we clearly have $\D^{\alpha'}(B) \subset D_{i}^{\alpha'}(B)$ -- it suffices to prove 
\begin{equation} \label{equation: B derivative induction hypothesis}
\left( \forall B \in \mc B_{\omega_1} \right) : B \in \mc B_\alpha \implies D^{\alpha'}_i (B) \subset \textnormal{cl}_{\Tr}(h_B) .
\end{equation}

To show that \eqref{equation: B derivative induction hypothesis} holds for $\alpha=0$ and $\alpha=1$, note that every $B\in \mc B_1$ is of the form $B = \{ h \} = \{ h_B \}$ for some $h\in \seq$. Since $1'=0'=0$, such $B$ satisfies
\[ D_i^{1'} ( B ) = D_i^{0'} ( B ) = D_i^{0} ( B ) = \cltr{B} = \cltr{ \{h_B\} } .\]

Let $\alpha = \lambda + 2n$, where $n\geq 1$, and suppose that \eqref{equation: B derivative induction hypothesis} (and hence also $1)$ from the statement) holds for every $\beta < \alpha$. We will prove that \eqref{equation: B derivative induction hypothesis} holds for $\alpha$ and $\alpha+1$ as well.
Let $B\in \mc B_\alpha \setminus \mc B_{<\alpha}$. By definition of $\mc B_\alpha$, we have $B = \bigcup_n h_B \ext f_n \ext B_n$ for some forking sequence $(f_n)_n$, broom sets $B_n \in \mc B_{<\alpha}$ and $h_B = \emptyset$.
Clearly, \errata{$\alpha-1= \lambda + 2n \err{-1}$ is an odd ordinal, so we have $B \in \mc B_{\lambda+2n-1}$. Since $\lambda+2n-1$ is strictly smaller than $\alpha$ and we have
\[ (\lambda+2n-1)' = (\lambda+2(n-1)+1)' = \lambda + n -1
 = (\lambda + 2n)' - 1 = \alpha'-1 ,\]}
the induction hypothesis implies that $D_i^{\alpha'-1} (f_n\ext B_n )$ is finite.
Consequently, the $\alpha'$-th derivative $D_i^{\alpha'} ( f_n \ext B_n )$ is empty. It follows that the longest sequence that might be contained in
\[ D_i^{\alpha'} ( B ) = D_i^{\alpha'} ( \bigcup_n h_B \ext f_n \ext B_n ) \]
is $h_B$, which gives \eqref{equation: B derivative induction hypothesis}.

For $B\in \mc B_{\alpha+1} \setminus \mc B_\alpha$, the proof is the same, except that $h_B \neq \emptyset$.

When $\alpha=\lambda$ is a limit ordinal, the proof is analogous.

$2)$: If $A\in \mc A_\alpha$ is a broom-extension of $B$, we observe that \err{$\D(A) = \textnormal{cl}_{\Tr}(B)$:}
\begin{align}
\D(A) & =  & \{ t\in \seq | \ & \textnormal{cl}_{\Tr}(A) \text{ contains infinitely many incomparable} \nonumber \\
							&&& \text{ extensions of $t$ of different length} \} \nonumber \\
	& = & \{ t\in \seq | \ & A \text{ contains infinitely many extensions of $t$} \}
	\label{equation: rewriting derivative of A} \\
	& = & \{ t\in \seq | \ & \text{the set }  \{ s\ext f^s_n \ext \nu^s_n | \ s\in B, n\in \omega \}
		\text{ contains} \nonumber \\
	&	&	& \text{infinitely many extensions of $t$} \} \nonumber \\	
	& = & \{ t \in \seq | \ & t \sqsubset s \text{ holds for some } s\in B \}
	= \textnormal{cl}_{\Tr}(B) , \nonumber 
\end{align}
The result then follows from $1)$, because the derivatives of $\textnormal{cl}_{\Tr}(B)$ are the same as those of $B$, and $B$ is an element of $\mc B_\alpha$.
\end{proof}


Finally, we describe a particular family of broom sets, used in \cite{talagrand1985choquet}, which has some additional properties.

Let $\varphi_n : \left( \mc B_{\omega_1} \right)^n \rightarrow \mc S$, for $n\in\omega$, be certain functions (to be specified later), such that for each $(B_n)_{n\in\omega}$, $\left( \varphi_n ( (B_i)_{i<n} \right))_{n\in \omega}$ is a forking sequence.
The collection $\mc B^T$ of \emph{finite Talagrand's brooms} is defined as the closure of $\mc B_0 = \{ \, \{ \emptyset \} \, \}$ in $\mc B_{\omega_1}$ with respect to the operations $B \mapsto h \ext B$, for $h \in \seq$, and
\begin{equation}\label{equation: construction of B T}
(B_n)_{n\in \omega} \mapsto \bigcup_{n\in\omega} f_n \ext B_n, \ \text{ where } f_n = \varphi_n ( (B_i)_{i<n} ) .
\end{equation}
For $\alpha \leq \omega_1$, we set $\mc B_\alpha^T := \mc B^T \cap \mc B_\alpha$.
Clearly, finite Talagrand's brooms satisfy an analogy of \eqref{equation: broom hierarchy}.

The family $\mc A^T$ of \emph{infinite Talagrand's brooms} consists of broom-extensions of elements of $\mc B^T$, where for each $B$, only certain combinations (to be specified later) of $(f^s_n)_n$, $s\in B$, and $\{ \nu^s_n | \ s\in B, n\in \omega \}$ are allowed.
For $\alpha \leq \omega_1$, we set $\mc A_\alpha^T := \mc A^T \cap \mc A_\alpha$.

The precise form of the as-of-yet-unspecified parameters above can be found in \cite{talagrand1985choquet}. However, our only concern is that the following lemma holds. (Recall that a family $\mc A$ is almost-disjoint if the intersection of any two distinct elements of $\mc A$ is finite.)

\begin{lemma}[The key properties of $\mc A^T$ and $\mc B^T$] \label{lemma: properties of Talagrands brooms}
The functions $\varphi_n$ and the ``allowed combinations'' above can be chosen in such a way that the following two properties hold:
\begin{enumerate}[(i)]
\item $\mc A^T$ is almost-disjoint; \label{case: A T is AD}
\item Let $B\in \mc B^T$ and let $L_s$, $s\in B$, be some sets. If each $L_s \cap \mc N(s)$ is $\tau_p$-dense in $\mc N(s)$, then there is some $\mc A^T$-extension $A$ of $B$ s.t. each $A\cap L_s$ is infinite. \label{case: A T is rich}
\end{enumerate}
\end{lemma}

\begin{proof}
$(i)$ follows from \cite[Lemma 4]{talagrand1985choquet}. $(ii)$ follows from \cite[Lemma 3]{talagrand1985choquet}.
\end{proof}

A third important property of Talagrand's broom sets is the following:

\begin{remark}
For Talagrand's brooms, the conclusion of Lemma \ref{lemma: rank of broom sets} is optimal.
\end{remark}

\noindent By ``conclusion being optimal'' we mean that for every $B\in \mc B_\alpha^T \setminus \mc B^T_{<\alpha}$ the derivative $\D^{\alpha'}(B)$ is non-empty, and similarly $\D^{\alpha'}(\D(A))$ is non-empty for every $A\in \mc A_\alpha^T \setminus \mc A^T_{<\alpha}$).
The ``importance'' of this property only appears implicitly -- if this property \err{didn't hold}, Talagrand's brooms would not be useful for our purposes. From this reason, we do not include the proof here, as we will not need the result in the following text. However, it can be proven analogously to Proposition 4.11~(ii) from \cite{kovarik2018brooms}.

\subsection{Broom spaces} \label{section: broom space properties}

Next, we introduce the class of broom spaces and state the main result of Section \ref{section: complexity of brooms}.
A useful related concept is that of a space with a single non-isolated point:

\begin{definition}[Space with a single non-isolated point] \label{definition: 1 non isolated point}
Let $\Gamma$ be a set, $\mc A\subset \mc P(\Gamma)$ and $\infty$ a point not in $\Gamma$. We define the \emph{space with a single non-isolated point} corresponding to $\Gamma$ and $\mc A$ as $\left(\Gamma\cup\{\infty\},\tau(\mc A)\right)$, where $\tau (A)$ is the topology in which
\begin{itemize}
\item each $\gamma\in \Gamma$ is isolated,
\item the neighborhood subbasis of $\infty$ consists of all sets of the form
\begin{align*}
& \{ \infty\} \cup \Gamma \setminus \{\gamma\}	 \ \ \ \text{ for } \gamma \in \Gamma \text{ and} \\
& \{ \infty\} \cup \Gamma \setminus A 			 \ \ \ \ \ \, \text{ for } A \in \mc A .
\end{align*}
\end{itemize}
\end{definition}

Broom spaces are a particular type of spaces with a single non-isolated point:

\begin{definition}[Broom spaces] \label{definition: AD topology}
A \emph{broom space} $T_{\mc A}$ corresponding to a non-empty family $\mc A \subset \mc A_{\omega_1}$ of infinite broom sets is defined as
\[ T_{\mc A} := \left( \baire \cup \{\infty\}, \tau(\mc A) \right) ,\]
that is, as the space with a single non-isolated point corresponding to $\Gamma=\baire$ and $\mc A$.
If the family $\mc A$ is almost-disjoint, the corresponding $T_{\mc A}$ is said to be an \emph{AD broom space}.

We define \emph{Talagrand's broom spaces} as
\begin{align*}
\mathbf{T}_\alpha & :=  T_{\mc A^T_{<\alpha}} && \text{for } \alpha \in [1,\omega_1] \text{ and } \\
\mathbf{S}_\alpha & :=  T_{\mc A^T_{<\alpha} \setminus \mc A^T_{<(\alpha-1)}}
	= T_{\mc A^T_{\alpha-1} \setminus \mc A^T_{<(\alpha-1)}}
	&& \text{for non-limit } \alpha \in [1,\omega_1] .
\end{align*}
\emph{The Talagrand's broom space} $\mathbf{T}$ is defined as $\mathbf{T}:=\mathbf{T}_{\omega_1}$ (emphasis on `\emph{the}'). 
\end{definition}

Formally, the main result of Section \ref{section: complexity of brooms} is the following theorem. (But as we mentioned earlier, we consider the methods to be the more interesting part.)
\errata{The proof is presented in Sections~\ref{section: cT} and \ref{section: dT}.}

\begin{theorem}[Complexity of Talagrand's broom spaces]\label{theorem: complexity of talagrands brooms}
Talagrand's broom spaces satisfy
\begin{align*}
\alpha \in [2,\omega_1] \implies &
\textnormal{Compl}\left(\mathbf{T}_\alpha\right) = \left[2,\alpha\right] \\
\alpha \in [2,\omega_1] \text{ is non-limit}  \implies &
\textnormal{Compl}\left(\mathbf{S}_\alpha\right) = \left[2,\alpha\right] .
\end{align*}
\end{theorem}

For every topological space \err{$Y$ which contains $\baire$ as a \emph{subset} (no topological assumptions)} and $y\in Y$, we set
\begin{align} \label{equation: S of y}
S(y) := \left\{ s\in\seq | \ \overline{\mc N(s)}^{Y} \ni y \right\}.
\end{align}
\noindent The following computation appears several times throughout the paper, so we formulate it separately.

\begin{lemma}[Sufficient condition for being $\Fa$]\label{lemma: S(y) and complexity of T}
Let $T$ be a broom space, $Y\supset T$ a topological space, and $\alpha\in [2,\omega_1]$. Suppose that for each $y\in Y\setminus T$, $S(y)$ is either finite or it can be covered by finitely many sets $\D(A)$, $A\in \mc A_{<\alpha}$. Then $T\in \mc F_\alpha(Y)$.
\end{lemma}

\noindent The proof of this lemma is the only place where Section \ref{section: complexity of brooms} relies on the concept of regular $\Fa$-representations and the results of Section \ref{section: regular representations}. Reader willing to treat this lemma as a ``black box'' is invited to ignore Section \ref{section: regular representations} altogether, and skip ahead to Proposition \ref{proposition: basic broom complexities}.

In the terminology of Section \ref{section: regular representations}, Lemma \ref{lemma: S(y) and complexity of T} actually shows that $\mc N$ is a regular $\Fa$-representation of $T$ in $Y$, in the sense of Definition \ref{definition: regular representations}. This  immediately follows from the remark just above Corollary \ref{corollary: sufficient condition for Fa}.

To prove Lemma \ref{lemma: S(y) and complexity of T} for a broom space $T$ contained in a topological space $Y$, we first need to have a Suslin scheme $\mc C$ in $T$ which satisfies $\mc A(\overline{\mc C}^Y\!\!) = T$.
It turns out that the basic open sets $\mc N(s) = \{ \sigma\in\baire | \ \sigma \sqsupset s \}$ from the product topology of $\baire$ provide a canonical solution:

\begin{lemma}[$\mc N$ as universally useful Suslin scheme] \label{lemma: A of N}
For every broom space $T$, the Suslin scheme $\mc N = (\mc N(s) )_{s\in\seq}$ on $\baire\subset T$ satisfies $\mc A(\overline{\mc N}^Y\!\!) = T$ for every topological space $Y\supset T$.
\end{lemma}

\begin{proof}
This holds, for example, by \cite[Lemma 5.10]{kovarik2018brooms}.
\end{proof}

Note that $S(\cdot) = S_{\mc N}(\cdot)$, so $S(\cdot)$ is a special case of the general definition in \eqref{equation: S C of y}. By definition of broom topology, the closure of any uncountable subset of $\baire$ contains $\infty$.
It follows that $\overline{\mc N(s)}^{Y} = \overline{\mc N(s)\cup \{\infty\}}^{Y}$ holds for every $s\in \seq$. In other words, we could equally well work with the non-disjoint Suslin scheme $(\mc N(s) \cup \{\infty\})_s$ on $T$.

\begin{proof}[Proof of Lemma \ref{lemma: S(y) and complexity of T}]
Let $T\subset Y$ and $\alpha$ be as in the lemma.
We shall verify that every $y\in Y \setminus T$ satisfies the assumptions of Corollary \ref{corollary: sufficient condition for Fa}, obtaining $T\in \Fa(Y)$ as a result.
In particular, we shall show that $\D^{\alpha'} (S(y))$ is empty when $\alpha$ is even, resp. that it can only contain the empty sequence when $\alpha$ is odd.

Let $y\in Y\setminus T$ \errata{and let $\mc A'\subset \mc A_{<\alpha}$ be a finite collection s.t. $\bigcup_{\mc A'} \D(A) \supset S(y)$.}
When $S(y)$ is finite, we have \errata{$\D^{\alpha'}(S(y)) = \emptyset$}, independently of the parity of $\alpha$ (because $\alpha\geq 2$ implies $\alpha' \geq 1$ and $\D^{\alpha'}(\cdot) \subset \D(\cdot)$, and $\D$-derivative of a finite set is empty).

Suppose that $\alpha$ is odd. We then have
\[ \D^{\alpha'}(S(y)) \subset \D^{\alpha'} ( \bigcup_{\mc A'} \D(A) )
 = \bigcup_{\mc A'} \D^{\alpha'} ( \D(A) ) 
 \overset{\alpha-1}{\underset{\leq \alpha}=}  \bigcup_{\mc A'} \D^{(\alpha-1)'} ( \D(A) ).\]
Since $\mc A' \subset \mc A_{<\alpha} = \mc A_{\alpha-1}$, it follows from $2)$\,$(ii)$ in Lemma \ref{lemma: rank of broom sets} that each set $\D^{(\alpha-1)'} ( \D(A) )$ is either empty or only contains the empty sequence.
This gives $\D^{\alpha'}(S(y)) \subset \{ \emptyset \}$.

Suppose that $\alpha$ is even. Let $A\in \mc A'$. Since $\mc A' \subset \mc A_{<\alpha}$, we have $A \in \mc A_{\beta}$ for some $\beta<\alpha$.
We have $\beta' < \alpha'$ (because $\alpha$ is even), which implies $\D^{\alpha'}(\cdot) \subset \D (\D^{\beta'} (\cdot))$.
Since $\D^{\beta'} (\D (A))$ is finite (by $2)$\,$(i)$ from Lemma \ref{lemma: rank of broom sets}), $\D(\D^{\beta'} (\D (A)))$ is empty, and thus $D^{\alpha'} ( \D (A) )$ is empty as well.
As in the case of odd $\alpha$, this gives $\D^{\alpha'}(S(y)) = \emptyset$.
\end{proof}

\noindent Actually, the proof of Lemma \ref{lemma: S(y) and complexity of T} did not use any special properties of $T$ and $S(\cdot)$, so the lemma holds when $T$ is replaced by an abstract topological space $X$ and $S(\cdot)$ is replaced by $S_{\mc C}(\cdot)$ for some complete Suslin scheme $\mc C$ on $X$.

\bigskip
We are now in a position to obtain the following (upper) bounds on the complexity of broom spaces:

\begin{proposition}[Basic broom space complexity results]\label{proposition: basic broom complexities}
Any broom space $T_{\mc A}$ satisfies:
\begin{enumerate}[(i)]
\item \cite{talagrand1985choquet} $T_{\mc A}$ is $\fsd$ in $\beta T_{\mc A}$, but it is not $\sigma$-compact;
\item If $\mc A \subset \mc A_{<\alpha}$ holds for some $\alpha\in[2,\omega_1]$, then $T_{\mc A}$ is an absolute $\Fa$ space.
\end{enumerate}
\end{proposition}

(In the language of Section \ref{section: regular representations}, (ii) even proves that $\mc N$ is a universal regular representation of $T_{\mc A}$, in the sense of Definition \ref{definition: regular representations}. This follows from the remark below Lemma \ref{lemma: S(y) and complexity of T}.)

\begin{proof}
(i) is proven in \cite[p.\,197]{talagrand1985choquet}, so it remains to prove (ii).
For even $\alpha$, the result \errata{follows from the proof of} \cite[Prop.\,5.13]{kovarik2018brooms} (but the proof below also applies).

Suppose we have $\mc A \subset \mc A_{<\alpha}$ for $\alpha\in[2,\omega_1]$ and let $Y\supset T_{\mc A}$ be a topological space.
Let $y\in Y \setminus T_{\mc A}$.
By Lemma \ref{lemma: S(y) and complexity of T}, it suffices to prove that $S(y)$ can be covered by finitely many sets $\D(A)$, $A\in \mc A \subset \mc A_{<\alpha}$.

Since $y\neq \infty$, there must be some open neighborhood $V$ of $y$ in $Y$ which satisfies $\infty \notin\overline{V}^Y$.
By definition of topology $\tau (\mc A)$, $\overline{V}^Y \cap T_{\mc A}$ can be covered by a union of some finite $\mc A' \subset \mc A$ and a finite set $F \subset \baire$. It follows that $S(y) = \bigcup_{\mc A'} \D(A)$ holds for some finite $\mc A' \subset \mc A$:
\begin{align} \label{equation: S of y is D of A}
& s \in S(y)
 \overset{\eqref{equation: S of y}}\iff y \in \overline{\mc N(s)}^{Y} \nonumber
 \iff y \in \overline{\mc N(s)}^{Y} \cap \overline{V}^Y
 \overset{\text{def.}}{\underset{\text{of }\tau(\mc A)}\iff} \\
& \overset{\text{def.}}{\underset{\text{of }\tau(\mc A)}\iff} 
 y \in \overline{\mc N(s)\cap \left(\bigcup \mc A' \cup F \right)}^{Y}
\overset{y\notin T}\implies \mc N(s)\cap \left(\bigcup \mc A' \cup F \right) \textrm{ is infinite} \\
& \iff \mc N(s)\cap A \textrm{ is inf. for some } A \in \mc A' \nonumber \\
& \iff \text{some } A\in \mc A' \textrm{ contains infinitely many incomparable extensions of } s \nonumber \\
& \overset{\eqref{equation: rewriting derivative of A}}{\iff} s \in \D(A) \text{ for some } A\in \mc A'
\iff s \in \bigcup_{\mc A'} \D(A). \nonumber
\end{align}
\end{proof}

\subsection{Talagrand's lemma}\label{section: Talagrand lemma}

In this subsection, we consider the following general problem:
We are given a $\mc K$-analytic space $X$, and a space $Y$ which contains it. We ``know'' $X$, but not its complexity in $Y$ -- we would like to obtain a lower bound on $\Compl{X}{Y}$.
To solve the problem, we take an arbitrary ``unkown'' set $Z\subset Y$, about which we however \emph{do} know that it belongs to $\Fa(Y)$. The mission shall be successful if we show that any such $Z$ which contains a big enough part of $X$ must also contain a part of $Y\setminus X$.
(Of course, the values of $\alpha$ for which this can work depend on the topology of $X$ and $Y$, so some more assumptions will be needed. In Sections \ref{section: cT} and \ref{section: dT}, we consider some suitable combinations of $X$ and $Y$.)

The method we adopt closely mimics the tools used in \cite{talagrand1985choquet} (namely, Lemma 1 and Lemma 3), and refines their conclusions. The author of the present article has found it hard to grasp the meaning of these tools. From this reason, Notations \ref{notation: witness} and \ref{notation: B corresponding to Z} introduce some auxiliary notions, which -- we hope -- will make the statements easier to parse.

The setting we consider is as follows:
\begin{itemize}
\item $X$ is a $\mc K$-analytic space of cardinality at least continuum.
	\begin{itemize}
	\item We assume that $X$ contains $\baire$ as a \emph{subset}.
	\item By default, $\baire$ is equipped with the subspace topology inherited from $X$. 	This subspace topology need not have any relation with the product topology $\tau_p$ 		of $\baire$; whenever $\tau_p$ is used, it shall be explicitly mentioned.
	\end{itemize}
\item $Y$ is a topological space containing $X$.
\item $Z$ is a subset of $Y$.
\end{itemize}

To get a quicker grasp of the notions, we can just imagine that
\[ \baire \subset X = Z \in \Fa(Y) \ \text{ for some } \alpha<\omega_1 ,\]
and we are aiming to arrive at a contradiction somewhere down the road.

The first auxiliary notion is that of ``witnessing'':\footnote{The author is aware of the fact that the word ``witness'' from Notation~\ref{notation: witness} is neither very original, nor particularly illuminating. Any ideas for a more fitting terminology would be appreciated. The same remark applies to the notion of ``correspondence'' from Notation~\ref{notation: B corresponding to Z}.}

\begin{notation}[Witnessing]\label{notation: witness}
$W\subset Y$ is a \emph{witness of $Z$ in $Y$} if there exists an indexing set $I$ and a family $\{ L_i | \ i\in I \}$ of closed subsets of $Y$ satisfying
\begin{enumerate}[(i)]
\item $\bigcap_{i\in I} L_i \subset Z$;
\item $\left( \forall i \in I \right) : L_i \cap W \text{ is infinite}$.
\end{enumerate}
\end{notation}

Of course, any infinite $Z$ which is closed in $Y$ is its own witness, since we can just set $I := \{i_0\}$ and $L_{i_0} := Z$.
But when the complexity of $Z$ is higher, we might need a much larger collection.
In general, there might even be no witnesses at all -- such as when $Z$ is discrete and $Y$ is its one-point compactification.

The key property of this notion is the following observation (used in \cite{talagrand1985choquet}). It shows that sets with certain properties cannot be witnesses, because they force the existence of points outside of $Z$.

\begin{lemma}[No discrete witnesses in certain spaces]\label{lemma: closed discrete witnesses}
Let $W$ be a closed discrete subset of $Z$. If $Y$ is s.t. $\overline{W}^Y$ is the one-point compactification of $W$, then $W$ cannot be a witness for $Z$ in $Y$.
\end{lemma}

For any discrete space $D$, we denote its one-point compactification as $\alpha D =: D\cup \{x_D\}$. We have
\begin{equation} \label{equation: closure in alpha D}
F\subset D \text{ is infinite } \implies \overline{F}^{\alpha D} = F \cup \{x_D\} .
\end{equation}

\begin{proof}[Proof of Lemma \ref{lemma: closed discrete witnesses}]
Suppose that $W$ is closed and discrete in $Z$ and $\overline{W}^Y$ is the one-point compactification of $W$ (and hence $x_W \notin Z$).
For contradiction, assume that there exists a family $\{ L_i \, | \ i\in I \}$ as in Notation~\ref{notation: witness}. Since each $L_i \cap W$ is infinite by $(ii)$, we have
\[ x_W \overset{\eqref{equation: closure in alpha D}}{\in} \overline{L_i \cap W}^Y
\subset \overline{L_i}^Y = L_i .\]
By $(i)$, we have $x_W \in \bigcap_I L_i \subset Z$ -- a contradiction.
\end{proof}

Note that whenever $W\subset Y$ is a witness of $Z$ in $Y$, then so is any $\widetilde W$ with $W\subset \widetilde W \subset Y$. However, it is more practical to have small witnesses, as it gives us more control over which elements appear in $\overline{L_i \cap W}^Y$.

Regarding positive results, we have the following refinement of \cite[Lemma 3]{talagrand1985choquet}:

\begin{lemma}[Existence of broom witnesses]\label{lemma: Talagrand lemma 3}
If $X\in \Fa(Y)$, then $X$ has an $\mc A^T_\alpha$-witness in $Y$.
\end{lemma}

As we will see later (in Lemma \ref{lemma: c T - lower bound} and Lemma \ref{lemma: d T - lower bound}), combining Lemma \ref{lemma: Talagrand lemma 3} and Lemma \ref{lemma: closed discrete witnesses} yields a lower bound on $\Compl{X}{Y}$ for certain $X$ and $Y$ (more precisely, for those $X\subset Y$ where each $A\in \mc A^T_\alpha$ is closed discrete in $X$, and satisfies $\overline{A}^Y = \alpha A$).
Before giving the proof of Lemma \ref{lemma: Talagrand lemma 3}, we first need some technical results.

Setting $W=\baire$ and $I=B\subset \seq$ in Notation \ref{notation: witness}, and strengthening the condition (ii), we obtain the following notion of ``correspondence''.
Its purpose is to allow the construction of $\mc A^T$-witnesses via Lemma \ref{lemma: properties of Talagrands brooms}\,\eqref{case: A T is rich}.

\begin{notation}[Correspondence]\label{notation: B corresponding to Z}
A set $B\subset \seq$ \emph{corresponds to $Z$} (in $Y$) if there exists a family $\{ L_s | \ s\in B \}$ of closed subsets of $Y$ satisfying
\begin{enumerate}[(i)]
\item $\bigcap_{s\in B} L_s \subset Z$;
\item $\left( \forall s \in B \right) : L_s \cap \mc N(s) \text{ is $\tau_p$-dense in } \mc N(s)$.
\end{enumerate}
\end{notation}

Talagrand's Lemma 1 \errata{(\cite{talagrand1985choquet})} says that when $X$ is $\mc F$-Borel in $Y$, there is a $\mc B^T$ broom set which corresponds to $X$. We refine this result to include the exact relation between complexity of $X$ and the ``rank'' of the corresponding broom set. The proof itself is identical to the one used in \cite{talagrand1985choquet} -- the non-trivial part was finding the ``right definitions'' for $\Fa$ and $\mc B_\alpha$ such that the correspondence holds.

\begin{lemma}[Exicence of corresponding brooms]\label{lemma: Talagrand lemma 1}
Let $\alpha<\omega_1$. If $X \in \Fa(Y)$, then there is some $B\in \mc B^T_\alpha$ which corresponds to $X$ in $Y$.
\end{lemma}

\begin{proof}
Recall that if $B\in \mc B^T_{\alpha+1}$ holds for even $\alpha$ and $h_B= \emptyset$, we actually have $B\in \mc B^T_{\alpha}$.
We shall prove the following stronger result. (Setting $Z:= X$ and $h:=\emptyset$ in Claim \ref{claim: Talagrand lemma 1} gives the conclusion of Lemma \ref{lemma: Talagrand lemma 1}.)

\begin{claim}\label{claim: Talagrand lemma 1}
Suppose that $Z\in\Fa(Y)$ and $Z\cap \mc N(h)$ is $\tau_p$-dense in $\mc N(h)$ for some $h\in \seq$.
\begin{enumerate}[(i)]
\item For odd $\alpha$, there is $B\in \mc B^T_\alpha$ with $h_B \sqsupset h$ which corresponds to $Z$ in $Y$.
\item For even $\alpha$, there is $B\in \mc B^T_{\alpha+1}$ with $h_B = h$ which corresponds to $Z$ in $Y$.
\end{enumerate}
\end{claim}

Let $Z$ and $h$ be as in the assumption of the claim. We shall prove the conclusion by transfinite induction.
For $\alpha=0$, $Z$ is closed, so we set $B:=\{h\}$ and $L_h := Z$.

Suppose that $\alpha$ is odd and the claim holds for $\alpha-1$. We have $Z=\bigcup_n Z_n$ for some $Z_n \in \mc F_{\alpha-1}(Y)$.
By the Baire theorem, some $Z_{n_0}$ is non-meager in $(\mc N(h),\tau_p)$. It follows that $Z_{n_0}$ is $\tau_p$ dense in $\mc N(h_0)$ for some $h_0 \sqsupset h$.
By the induction hypothesis, there is some $B\in \mc B^T_{(\alpha-1)+1} = \mc B^T_{\alpha}$ with $h_B = h_0 \sqsupset h$ which corresponds to $Z_{n_0}$. In particular, this $B$ also corresponds to $Z\supset Z_{n_0}$.

Suppose that $\alpha\in (0,\omega_1)$ is even and the claim holds for every $\beta < \alpha$.
We have $Z=\bigcap_n Z_n$, where $Z_n \in \mc F_{\alpha_n}(Y)$ for some odd $\alpha_n < \alpha$.
Suppose we have already constructed $(f_i)_{i<n}$, $(B'_i)_{i<n}$ and $(B_i)_{i<n}$ for some $n\in\omega$.
Let $f_n := \varphi_n((B'_i)_{i<n})$ (where $\varphi_n$ is the function from Lemma \ref{lemma: properties of Talagrands brooms}). By the induction hypothesis, there is some $B_n \in \mc B^T_{\alpha_n}$ which corresponds to $Z_n$ and satisfies $h_{B_n} \sqsupset h\ext f_n$.
As noted in Section \ref{section: broom sets}, $B_n$ can be rewritten as $h\ext f_n \ext B'_n$, where $B'_n \in \mc B^T_{\alpha_n} \subset \mc B^T_{<\alpha}$.

Once we have $f_n$ and $B'_n$ for every $n\in\omega$, we get a broom set $B' := \bigcup_n f_n\ext B'_n$. Since we have both $B' \in \mc B_\alpha$ and $B'\in \mc B^T$, $B$ belongs to $\mc B^T_\alpha$. It follows that $B:=h\ext B' \in \mc B^T_{\alpha+1}$.

It remains to prove that $B$ corresponds to $Z$. Since each $B_n$ corresponds to $Z_n$, there are some closed sets $L_s$, $s\in B_n$, such that $\bigcap_{B_n} L_s \subset Z_n$ (and $(ii)$ from Notation \ref{notation: B corresponding to Z} holds). Since $B=\bigcup_n h\ext f_n \ext B'_n = \bigcup_n B_n$, we have
\[ \bigcap_{s\in B} L_s = \bigcap_{n\in \omega} \bigcap_{s\in B_n} L_s \subset \bigcap_{n\in\omega} Z_n \subset Z ,\]
which shows that $B$ corresponds to $Z$.
\end{proof}

Lemma \ref{lemma: Talagrand lemma 3} is now a simple corollary of Lemma \ref{lemma: Talagrand lemma 1}:

\begin{proof}[Proof of Lemma \ref{lemma: Talagrand lemma 3}]
Let $B \in \mc B^T_\alpha$ be the set which corresponds to $X$ in $Y$ by Lemma \ref{lemma: Talagrand lemma 1}. Let $L_s$, $s\in B$, be the closed subsets of $Y$ as in Notation \ref{notation: B corresponding to Z}.
By Lemma \ref{lemma: properties of Talagrands brooms}, there is an $\mc A^T$-extension $A$ of $B$, such that each $L_s \cap A$ is infinite.
Setting $I := B$, we see that $A$ is a witness for $X$ in $Y$. Since we also have $A\in \mc A_\alpha$ (by definition of $\mc A_\alpha$), we have $A\in \mc A^T_\alpha$ and the proof is complete.
\end{proof}

We will also need the following technical version of Lemma \ref{lemma: Talagrand lemma 3}. As it contains no particularly novel ideas, we recommend skipping it on the first reading.

\begin{lemma}[Existence of broom witnesses -- technical version]\label{lemma: technical version of Talagrand lemma 3}
Suppose that $X\in \mc F_\eta (Y)$ holds for some $\eta < \omega_1$. Then for every $\delta \in [\eta,\omega_1)$, there is some $A\in \mc A^T_\delta \setminus A^T_{<\delta}$ and $h_0\in \seq$ s.t. $A \cap \mc N(h_0) \in \mc A^T_\eta$ is a witness for $X$ in $Y$.
\end{lemma}

\begin{proof}
Assume that $X\in \mc F_\eta (Y)$ and $\eta \leq \delta < \omega_1$.
First, construct a set $B\in \mc B^T_\delta \setminus \mc B^T_{<\delta}$ and $h_0 \in \seq$, such that the following set $B_0$ belongs to $\mc B^T_\eta$ and corresponds to $X$:
\[ B_0 := \{ s\in B| \ s\sqsupset h_0 \} .\]

Denote $h_{\textnormal{even}}:=\emptyset$ and choose an arbitrary $h_{\textnormal{odd}} \in \seq \setminus \{\emptyset\}$ (this notation is chosen merely so that it is simpler to explain how the construction differs depending on the parity of $\delta$).

Suppose first that $\delta$ is even.
Recall that $\varphi_n$, $n\in\omega$, are the functions from Lemma \ref{lemma: properties of Talagrands brooms}.
Set $f_0 := \varphi_0 (\emptyset)$, $h_0 := h_{\textnormal{even}}\ext f_0$.
By Claim \ref{claim: Talagrand lemma 1}, there is some $B_0 \in \mc B^T_\eta$ which corresponds to $Z$ and satisfies $h_{B_0} \sqsupset h_0$. Denote by $B'_0$ the $\mc B^T_\eta$-set satisfying $B_0 = h\ext f_0 \ext B'_0$.
We either have $B_0 \in \mc B^T_\delta \setminus \mc B^T_{<\delta}$ (in which case we set $B := B_0$ and the construction is complete), or $B_0, B'_0 \in \mc B^T_{<\delta}$.

Assume the second variant is true.
Let $(B'_n)_{n=1}^\infty$ be such that $\delta$ is the smallest ordinal s.t. each $B'_n$, is contained in $\mc B^T_{<\delta}$.
For $n\geq 1$, we set $f_n := \varphi_n ( (B'_i)_{i<n})$.
By \eqref{equation: construction of B T}, the following set $B$ belongs to $\mc B^T_{\delta}$:
\[ B := \bigcup_n f_n \ext B'_n = \bigcup_n h_{\textnormal{even}}\ext f_n \ext B'_n .\]
The choice of $(B'_n)_{n=1}^\infty$ ensures that $B \in \mc B^T_{\delta} \setminus \mc B^T_{<\delta}$.

To get the result for odd $\delta=\tilde \delta +1$, we just repeat the above process with $h_{\textnormal{odd}}$ in place of $h_{\textnormal{even}}$ and $\tilde \delta$ in place of $\delta$.

We now construct $A$ with the desired properties.
Let $L_s$, $s\in B_0$, be some sets which ensure that $B_0$ corresponds to $X$ and denote $L_s := \overline{X}^Y$ for $s\in B\setminus B_0$.
By Lemma~\ref{lemma: properties of Talagrands brooms}, there is some $\mc A^T$-extension $A$ of $B$ such that each $L_s \cap A$ is infinite.
By definition of $\mc A^T_{(\cdot)}$, we have $A\in \mc A^T_\delta \setminus A^T_{<\delta}$.
Moreover, $A \cap \mc N(h_0)$ is a broom-extension of $B_0\in \mc B^T_\eta$, which gives $A \cap \mc N(h_0) \in \mc A^T_\eta$.
Finally, $L_s \cap A = L_s \cap (A\cap \mc N(h_0))$ is infinite for each $s\in B_0$, and $\bigcap L_s \subset X$ holds even when the intersection is taken over $s\in B_0$. This proves that $A \cap \mc N(h_0)$ is a witness of $X$ in $Y$.
\end{proof}

In the~next part, we study broom spaces and their compactifications in a~more abstract setting. The purpose is to isolate the~few key properties which are required to obtain the~results we need, while ignoring all the~other details.

\subsection{Amalgamation spaces} \label{section: amalgamation spaces}

In Section~\ref{section: amalgamation spaces} (and only here), we retract our standing assumption that every topological space is Tychonoff.

In Section~\ref{section: zoom spaces}, we were able to take a~space $Z(cY,\mc X)$ and construct its compactification $Z(cY,c\mc X)$ by separately extending each $X_i$ into $cX_i$. Recall that this was easily doable, since the~spaces $X_i$ were pairwise disjoint and clopen in $Z(cY,X_i)$.
Our goal is to take an AD broom space $T_{\mc A}$, extend each $A\in \mc A$ separately into a~compactification $cA$, and thus obtain a~compactification of $T_{\mc A}$. However, the~family $\mc A$ is not disjoint, so we need a~generalization of the~approach from Section~\ref{section: zoom spaces}.

We will only need the~following properties of broom spaces:

\begin{lemma}[Example: Broom spaces]\label{lemma: broom spaces and A 1-4}
Any AD broom space $X=T_{\mc A}$ satisfies the~following four conditions:
\begin{itemize}
	\item[$(\mc A 1)$] $\mc A$ consists of clopen subsets of $X$.
	\item[$(\mc A 2)$] For distinct $A,A'\in \mc A$, the~intersection $A\cap A'$ is compact.
	\item[$(\mc A 3)$] $K := X \setminus \bigcup \mc A$ is compact.
	\item[$(\mc A 4)$] Whenever $\mc U$ is a~collection of open subsets of $X$ which covers $K$, there exists a~finite family $\mc A' \subset \mc A$, s.t. for every $A\in \mc A \setminus \mc A'$, we have $U \cup \bigcup \mc A' \supset A$ for some $U\in \mc U$.
\end{itemize}
\end{lemma}

To summarize the~properties, we can say that the~family $\mc A$ consists of clopen sets with small intersections. The third and fourth condition then ensure, in somewhat technical manner, that the~only parts where the~space $X$ is not compact are the~sets $A\in \mc A$.	

\begin{proof}
Let $T_{\mc A}$ be an AD broom space.
$(\mc A 1)$ holds by the~definition of topology $\tau(\mc A)$ on $T_{\mc A}$.
$(\mc A 2)$ because the~intersection of two distinct elements of $\mc A$ is, in fact, even finite.

$(\mc A 3)$: Often, $\mc A$ will cover the~whole space $T_{\mc A}$ except for $\infty$. In this case, $(\mc A3)$ holds trivially.
However even in the~non-trivial case where $K = T_{\mc A} \setminus \bigcup \mc A$ is infinite, the~subspace topology on $K = T_{\mc A} \setminus \bigcup \mc A$ coincides with the~topology of one-point compactification of $\baire \setminus \bigcup \mc A$ (by Definition~\ref{definition: 1 non isolated point}). This shows that $K$ is compact.

$(\mc A 4)$: Whenever $\mc U$ is an open cover of $K$, there is some $U \in \mc U$ which contains $\infty$. By definition of $\tau(\mc A)$, $U$ contains some basic open set $T_{\mc A} \setminus (\bigcup \mc A' \cup F)$, where $\mc A' \subset \mc A$ and $F \subset K$ are finite.
The fact that $T_{\mc A}$ satisfies $(\mc A4)$ for $\mc U$ is then witnessed by $\mc A'$ and $U$. Note that $U$ is a~universal witness, as we have $U \supset A \setminus \bigcup \mc A'$ for \emph{every} $A\in \mc A\setminus A'$.
\end{proof}

In the~remainder of this section, we will work with an abstract topological space $X$ and a~fixed family $\mc A \subset \mc P(X)$ s.t. the~conditions $(\mc A1)-(\mc A4)$ from Lemma~\ref{lemma: broom spaces and A 1-4} hold.

Let $\mc E=\left( E(A)\right)_{A\in \mc A}$ be a~collection of topological spaces such that for each $A\in \mc A$, $A$ is a~dense subset of $E(A)$.
Informally speaking, our goal is to find a~space whose ``local behavior'' is ``$E(A)$-like'', but the~``global properties'' are similar to those of $X$.

Without yet defining any topology on it, we set $\Amg{X}{\mc E} := X\cup\bigcup_{\mc A} E(A)$, assuming that each $A$ is ``extended into $E(A)$ separately'':
\begin{align} \label{equation: amalgamation set}
(\forall A \in \mc A) : 	& \ X\cap E(A) = A 	\nonumber \\
(\forall A, A' \in \mc A) : & \ A\neq A' \implies E(A)\cap E(A') = A \cap A'
\end{align}

To define the~topology on $\Amg{X}{\mc E}$, we first need the~following lemma.

\begin{lemma}[Largest open set with given trace] \label{lemma: W A U}
Let $P\subset Q$ be topological spaces and $G,G'$ open subsets of $P$.
Denote by $W^Q_P(G)$ the~largest open subset $W$ of $Q$ which satisfies $W \cap P = G \cap P$.
\begin{enumerate}[(i)]
\item $W^Q_P(G)$ is well defined. \label{case: W A U well defined}
\item $G\subset G' \implies W^Q_P(G) \subset W^Q_P(G')$. \label{case: W A U monotonicity}
\item $W^Q_P(G\cap G') = W^Q_P(G) \cap W^Q_P(G')$ .\label{case: W A U intersection}
\item If $P$ is dense in $Q$, we have $W^Q_P(G) \subset \textnormal{Int}_Q \, \overline{G}^Q$.
	\label{case: W A U is subset of ...}
\item For any compact $C\subset P$, we have $W^Q_P(P\setminus C) = Q \setminus C$.
	\label{case: W A U with compact complement}
\end{enumerate}
\end{lemma}

\begin{proof}
Let $P,Q,G$ and $G'$ be as in the~statement.

\eqref{case: W A U well defined}: Since $P$ is a~topological subspace of $Q$, there always exists \emph{some} open subset $W$ which satisfies $W \cap P = G$.
Consequently, we can define $W^Q_P(G)$ as
\begin{equation}\label{equation: W A U formula}
W^Q_P(G) := \bigcup \{ W \subset Q \text{ open } | \ W \cap A = G \} .
\end{equation}

\eqref{case: W A U monotonicity}:
Suppose that $G\subset G'$. The set $W' := W^Q_P(G) \cup W^Q_P(G')$ is open in $Q$ and satisfies
\[ W' \cap P = (W^Q_P(G) \cap P ) \cup (W^Q_P(G') \cap P)
= G \cup G' = G' .\]
Applying \eqref{equation: W A U formula} to $G'$, we have $W' \subset W^Q_P(G')$. It follows that $W^Q_P(G) \subset W^Q_P(G')$.

\eqref{case: W A U intersection}: 
``$\subset$'' follows from \eqref{case: W A U monotonicity}.
``$\supset$'' holds by \eqref{equation: W A U formula}, since $W^Q_P(G) \cap W^Q_P(G')$ is open in $Q$ and satisfies
\[ W^Q_P(G) \cap W^Q_P(G') \cap P = (W^Q_P(G) \cap P) \cap ( W^Q_P(G') \cap P) = G \cap G' .\]

\eqref{case: W A U is subset of ...}:
$P$ is dense in $Q$ and $W^Q_P(G) \subset Q$ is open.
Consequently, $W^Q_P(G) \cap P$ is dense in $W^Q_P(G)$ and we get
\[ W^Q_P(G) \subset \overline{W^Q_P(G) \cap P}^Q = \overline{G}^Q .\]
Since $W^Q_P(G)$ is open, the~conclusion follows.

\eqref{case: W A U with compact complement}: This is immediate, because $Q \setminus C$ is open in $Q$.
\end{proof}

For an open subset $U$ of $X$, we define
\[ V_U := U \cup \bigcup \{ W^{E(A)}_A(U\cap A) \, | \ A\in \mc A \} .\]
It follows from Lemma~\ref{lemma: W A U} that these sets satisfy
\begin{align}
U, U' \subset X \text{ are open in } X \implies \label{equation: V U cap}
& V_{U\cap U'} = V_U \cap V_{U'} .
\end{align}

\begin{definition}[Amalgamation space]\label{definition: amalgamation space}
The \emph{amalgamation of $X$ and $\mc E$} is defined as the~set $\Amg{X}{\mc E}$, equipped with the~topology whose basis\footnote{We do not claim that it is obvious that $\mc B$ is a~basis of topology. This is the~content of Lemma~\ref{lemma: Amg definition is correct}.} $\mc B$ consists of all sets of the~form
\begin{itemize}
\item $W\subset E(A)$, where $W$ is open in $E(A)$ and $A\in \mc A$;
\item $V_U$, where $U\subset X$ is open in $X$.
\end{itemize}
\end{definition}

In Lemma~\ref{lemma: Amg definition is correct}, we show that the~system $\mc B$ is closed under intersections, and therefore the~topology of $\Amg{X}{\mc E}$ is defined correctly.
We then follow with Lemma~\ref{lemma: basic properties of amalgamation spaces}, which captures the~basic and ``local'' properties of $\Amg{X}{\mc E}$. The ``global properties'' of amalgamation spaces are used implicitly in Proposition~\ref{proposition: compactifications of amalgamations}, where amalgamations are used to compactify $X$.

\begin{lemma}\label{lemma: Amg definition is correct}
The topology of $\Amg{X}{\mc E}$ is correctly defined.
\end{lemma}

\begin{proof}
Let $\Amg{X}{\mc E}$ be an amalgamation space and $\mc B$ the~system above, which we claim is a~basis of topology.
Since $\mc B$ obviously covers $\Amg{X}{\mc E}$, it remains to show that intersection of any two elements of $\mc B$ is again in $\mc B$.

This trivially holds when $W_0,W_1$ are two open subsets of the~same $E(A)$.
When $W$ is an open subset of $E(A)$ and $V_U$ corresponds to some open $U\subset X$, we have $V_U \cap E(A) = W^A_U$, so $W\cap V_U$ is again an open subset of $E(A)$.
When $V_U$ and $V_{U'}$ correspond to some opens subsets $U,U'$ of $X$, we have $V_U \cap V_{U'} = V_{U\cap U'} \in \mc B$ by \eqref{equation: V U cap}.

It remains to consider the~situation when $W$ is an open subset of some $E(A)$ and $W'$ is an open some $E(A')$, $A'\neq A$.
We need the~following claim:

\begin{claim}\label{claim: A cap A' is clopen}
For every distinct $A,A'\in \mc A$, $A'\cap E(A)$ is clopen in $E(A)$.
\end{claim}

\begin{proof}[Proof of the~claim]
Note that $A'\cap E(A) =A'\cap A$ is open in $X$, hence in $A$.
Further, we have
\[ A'\cap A \subset W_A^{E(A)} (A'\cap A)
\overset{L\ref{lemma: W A U}}{\underset{\eqref{case: W A U is subset of ...}}\subset}
\overline{A'\cap A}^{E(A)}
\overset{(\mc A2)} = A'\cap A .\]
Thus $A'\cap A$ is open in $E(A)$. It is also closed, as it is compact by $(\mc A2)$.
\end{proof}

Since $W'$ is open in $E(A)$, $W' \cap E(A)$ is open in $E(A') \cap E(A) = A' \cap A$. This set is, in turn, open in $E(A)$ (by the~claim). It follows that $W' \cap E(A)$ is open in $E(A')$, and thus $W\cap W' \in \mc B$.
\end{proof}

\begin{lemma}[Basic properties of amalgamations]\label{lemma: basic properties of amalgamation spaces}
The space $\Amg{X}{\mc E}$ has the~following properties:
\begin{enumerate}[(i)]
\item Each $E(A)$ is clopen in $\Amg{X}{\mc E}$. \label{case: E(A) is clopen in amalgamation}
\item The subspace topologies (inherited from $\Amg{X}{\mc E}$) on $X$ and $E(A)$ for $A\in \mc A$ coincide with the~original topologies of these spaces. \label{case: Amg and E(A) and X}
\item In particular, $X=\Amg{X}{\mc A}$. \label{case: X as Amg}
\item The space $\Amg{X}{\mc E}$ satisfies conditions $(\mc A1)-(\mc A4)$ for $\mc E$. 
	\label{case: Amg and E satisfies A 1 - 4}
\item The space $\Amg{X}{\mc E}$ is Hausdorff, provided that $X$ and each $E\in \mc E$ are Hausdorff \label{case: Amg is T 3}.
\end{enumerate}
\end{lemma}

\begin{proof}
\eqref{case: E(A) is clopen in amalgamation}: Let $A\in \mc A$. $E(A)$ is open in $\Amg{X}{\mc E}$ by definition (since it is open in itself).

To show that $\Amg{X}{\mc E} \setminus E(A)$ is open, we first prove $E(A) = V_A$.
For any $A'\neq A$, we have 
\begin{equation}\label{equation: W of A cap A'}
W^{E(A')}_{A'}(A\cap A') \overset{L\ref{lemma: W A U}}{\underset{\eqref{case: W A U is subset of ...}}\subset}
\overline{A\cap A'}^{E(A)} \overset{(\mc A 2)}{=} A\cap A' \subset A .
\end{equation}
Since $E(A) = W^{E(A)}_A(A)$, it follows that $E(A) = V_A$: 
\begin{align*}
E(A) = W^{E(A)}_A(A) \subset V_A
= A \cup W^{E(A)}_A (A) \cup \bigcup_{A'\neq A} W^{E(A')}_{A'} (A \cap A')
\overset{\eqref{equation: W of A cap A'}}\subset \\
\overset{\eqref{equation: W of A cap A'}}\subset
A \cup E(A) \cup \bigcup_{A'\neq A} A = E(A) .
\end{align*}

It remains to show that every point from the~set
\[ \Amg{X}{\mc E} \setminus E(A) = (X\setminus A) \cup \bigcup_{A'\neq A} E(A') \setminus A \]
is contained in some open set disjoint with $E(A)$.
For $x\in E(A') \setminus A$, the~set $E(A') \setminus E(A) = E(A')\setminus A$ is open in $E(A')$ (by Claim~\ref{claim: A cap A' is clopen}), and hence in $\Amg{X}{\mc E}$ as well.
For $x \in X \setminus A$, we have
\[ x \in X\setminus A \subset V_{X\setminus A}
\overset{\eqref{equation: V U cap}}{\subset} \Amg{X}{\mc E} \setminus V_A
 = \Amg{X}{\mc E} \setminus E(A) .\]

\eqref{case: Amg and E(A) and X}: Let $A\in \mc A$. For any basic open set $B$ in $\Amg{X}{\mc E}$, the~intersection $B\cap E(A)$ is open $E(A)$ (by definition of topology of $\Amg{X}{\mc E}$). Moreover, any $W\subset E(A)$ which is open in $E(A)$ is, by definition, also open in $\Amg{X}{\mc E}$. This shows that the~subspace topology of $E(A) \subset \Amg{X}{\mc E}$ coincides with the~original topology of $E(A)$.

We now show that the~subspace topology of $X \subset \Amg{X}{\mc E}$ coincides with the~original topology of $X$. Again, let $B$ be a~basic open subset of $\Amg{X}{\mc E}$.
If $B= V_U$ holds for some open $U\subset X$, we have $B\cap X = U$. When $B$ is an open subset of some $E(A)$, $B\cap X = B \cap A$ is open in $A$, and therefore also open in $X$ (because $A$ is open in $X$). This shows that the~original topology of $X$ is finer than the~subspace topology.

Conversely, for any open subset $U$ of $X$, we have $V_U\cap X=U$, which proves that the~subspace topology of $X$ is finer than the~original topology.

\eqref{case: X as Amg}: By (ii), $X$ is embedded in $\Amg{X}{\mc E}$ for any $\mc E$. Since $\Amg{X}{\mc A}$ adds no new points, this ``canonical'' embedding is a~homeomorphism.

\eqref{case: Amg and E satisfies A 1 - 4}:
$(\mc A1)$ for $\mc E$ (and $\Amg{X}{\mc E}$) is equivalent to \eqref{case: E(A) is clopen in amalgamation}.
By \eqref{equation: amalgamation set}, distinct sets $E(A), E(A') \in \mc E$ satisfy $E(A)\cap E(A') = A \cap A'$. Since $A \cap A'$ is compact by $(\mc A2)$, we get $(\mc A2)$ for $\mc E$ as well.
By \eqref{equation: amalgamation set}, we have
\[ \Amg{X}{\mc E} \setminus \bigcup \mc E = X \setminus \bigcup \mc A = K ,\]
which gives $(\mc A3)$ for $\mc E$.

To prove $(\mc A4)$ for $\mc E$, let $\mc V$ be a~collection of open subsets of $\Amg{X}{\mc E}$ covering $K$.
It suffices to work with a~suitable refinement -- in particular, we can assume that $\mc V$ consists of basic open sets. Denote
\begin{align*}
\mc V_K 		& := \mc V \cap \{ V_U |\, U\subset X \text{ open}\}
.
\end{align*}

By \eqref{equation: amalgamation set}, $\mc V_K$ is an open (in $\Amg{X}{\mc E}$) cover of $K$ and
\[ \mc U_K := \{ U\subset X | \ V_U \in \mc V_K \} \]
is an open (in $X$) cover of $K$. By $(\mc A4)$, there is a~finite family $\mc A' \subset \mc A$ s.t. for every $A\in \mc A \setminus A'$ we have $U_A \supset A \setminus \bigcup \mc A'$ for some $U_A \in \mc U$.
For any $A\in \mc A \setminus \mc A'$, the~set $C_A := A \cap \bigcup \mc A'$ is compact by $(\mc A2)$. Consequently, we have
\[ V_{U_A} \supset W^{E(A)}_A( U_A \cap A)
\overset{L\ref{lemma: W A U}}{\underset{\eqref{case: W A U monotonicity}}\supset}
W^{E(A)}_A(A \setminus C_A)
\overset{L\ref{lemma: W A U}}{\underset{\eqref{case: W A U with compact complement}}=}
E(A) \setminus C_A .\]
It follows that $V := V_{U_A}$ is an element of $\mc V$ satisfying $V \cup \bigcup \mc A' \supset E(A)$, which shows that $(\mc A4)$ holds for $\mc E$.

\eqref{case: Amg is T 3}:
First, we show that $\Amg{X}{\mc E}$ is Hausdorff.
Let $x,y\in\Amg{X}{\mc E}$ be distinct. It suffices to consider the~cases where
\begin{enumerate}[1)]
\item $x\in E(A)$ and $y\notin E(A)$ for some $A\in\mc A$,
\item $x,y \in E(A)$ for some $A\in\mc A$ and
\item $x,y \in K = X \setminus \bigcup \mc A$.
\end{enumerate}
In the~first case, the~open sets separating $x$ from $y$ are $\Amg{X}{\mc E}\setminus E(A)$ and $E(A)$ (by (i)).
In the~second case, we use the~fact that $E(A)$ is Hausdorff to find open $W,W'\subset E(A)$ which separate $x$ and $y$. By (i), $W$ and $W'$ are open in $\Amg{X}{\mc E}$ as well.

In the~last case, we use the~fact that $X$ is Hausdorff to get some disjoint open subsets $U, U'$ of $X$ for which $x \in U$ and $y\in U'$. By \eqref{equation: V U cap}, $V_U$ and $V_{U'}$ are disjoint as well.

To prove that $\Amg{X}{\mc E}$ is Tychonoff, it suffices to show that it is a~subspace of some compact space.
By \eqref{case: Amg and E satisfies A 1 - 4}, we can construct the~amalgamation $C := \Amg{\Amg{X}{\mc E}}{\beta \mc E}$.
By \eqref{case: X as Amg}, $\Amg{X}{\mc E}$ is a~subspace of $C$.
In Proposition~\ref{proposition: compactifications of amalgamations}, we show that the~amalgamation $C$ is compact (obviously, without relying on the~fact that $\Amg{X}{\mc E}$ is Tychonoff).
\end{proof}

\begin{proposition}[Compactifications of amalgamation spaces] \label{proposition: compactifications of amalgamations}
Let $X$ be a~Tychonoff space, $\mc A\subset \mc P(X)$ a~family satisfying $(\mc A1)$-$(\mc A4)$ and $cA$, $dA$ (Hausdorff) compactifications of $A$ for every $A\in\mc A$.
\begin{enumerate}[(i)]
\item For any regular topological space $Z$, a~function $f: \Amg{X}{\mc E} \rightarrow Z$ is continuous if and only if all the~restrictions $f|_X$ and $f|_{E(A)}$, $A\in\mc A$, are continuous. \label{case: continuous functions on amalgamation space}
\item $\Amg{X}{c\mc A}$ is a~compactification of $X$. In particular, $\Amg{X}{\mc A}$ is Tychonoff.
	\label{case: compactification of Amg}
\item $\Amg{X}{c\mc A} \preceq \Amg{C}{d\mc A} $ holds whenever $cA\preceq dA$ for each $A\in\mc A$.
	\label{case: smaller compactifications of Amg}
\item For every compactification $cX$ of $X$, we have $\Amg{X}{\overline{\mc A}^{cX}} \succeq cX$.\footnote{Analogously to Notation~\ref{notation: c A} we define $\overline{\mc A}^{cX} := \{ \overline{A}^{cX} | \ A\in\mc A \}$}
	\label{case: larger compactification of Amg type}
\item In particular, $\beta X = \Amg{X}{\beta \mc A}$.
	\label{case: beta compactification as Amg}
\end{enumerate}
\end{proposition}

\begin{proof}
\eqref{case: continuous functions on amalgamation space}: It remains to prove ``$\Leftarrow$''.
Let $f: \Amg{X}{\mc E} \rightarrow Z$ and suppose that all the~restrictions are continuous. We need to prove that $f$ is continuous at each point of $\Amg{X}{\mc E}$.

Let $x\in E(A)$ for some $A\in \mc A$.
Since $E(A)$ is clopen in $\Amg{X}{\mc E}$ (by \eqref{case: E(A) is clopen in amalgamation} of Lemma~\ref{lemma: basic properties of amalgamation spaces}) and $f|_{E(A)}$ is continuous (by the~assumption), $f$ is continuous at $x$.

Let $x\in K$. Let $G$ be an open neighborhood of $f(x)$ in $Z$ and let $H$ be an open neighborhood of $f(x)$ satisfying $\overline{H} \subset G$.
By continuity of $f|_X$, there is some open neighborhood $U$ of $x$ in $X$ which satisfies $f(U) \subset H$.
In particular, we have $f(U\cap A) \subset H$ for any $A\in \mc A$. Since $\overline{U \cap A} \subset E(A)$ and $f|_{E(A)}$ is continuous, we have
\[ f(\overline{U \cap A}) \subset \overline{ f( U \cap A )} \subset \overline{H} \subset G .\]
It follows that $f( U \cup \bigcup_{\mc A} \overline{U \cap A} ) \subset G$.
This proves that $V_U$ is an open neighborhood of $x$ which is mapped into $G$:
\[ x \in V_U \overset{\text{def.}}{=} U \cup \bigcup_{A \in \mc A} W^{E(A)}_A(U\cap A)
\overset{L\ref{lemma: W A U}}{\underset{\eqref{case: W A U is subset of ...}}\subset}
U \cup \bigcup_{A \in \mc A} \overline{U \cap A}
\ \ \& \ \ f( U \cup \bigcup_{\mc A} \overline{U \cap A} ) \subset G .\]

\eqref{case: compactification of Amg}:
By Lemma~\ref{lemma: basic properties of amalgamation spaces}, we already know that $\Amg{X}{c\mc A}$ is Hausdorff. Once we know that $\Amg{X}{c\mc A}$ is compact, we get that it is Tychonoff for free, which proves the~``in particular'' part.

Let $\mc V$ be an open cover. As in the~proof of Lemma~\ref{lemma: basic properties of amalgamation spaces} \eqref{case: Amg and E satisfies A 1 - 4}, we can assume that $\mc V$ consists of basic open sets, denoting
\begin{align*}
& \mc V_K := \mc V \cap \{ V_U |\, U\subset X \text{ open}\} \text{ and} \\
& \mc V_{cA} := \mc V \cap \{ W | \, W \cap cA \neq \emptyset \}, \ A\in\mc A .
\end{align*}

Clearly, we have $\bigcup \mc V_{cA} \supset cA$ for every $A\in \mc A$, so there exist some finite subfamilies $\mc V'_{cA}$ of $\mc V_{cA}$ satisfying $\bigcup \mc V'_{cA} \supset cA$.
Similarly we have $\bigcup \mc V_K \supset K$ and we denote by $\mc V'_K$ be some finite subfamily of $\mc V_K$ satisfying $\bigcup V'_K \supset K$.

By Lemma~\ref{lemma: basic properties of amalgamation spaces} \eqref{case: Amg and E satisfies A 1 - 4}, $(\mc A4)$ holds for $\Amg{X}{c\mc A}$ and $c\mc A$.
Applying $(\mc A4)$ yields a~finite family $c\mc A' \subset c \mc A$, such that every $cA \in c \mc A \setminus c \mc A'$ satisfies $V \cup \bigcup c\mc A' \supset cA$ for some $V \in \mc V'_K$.
In particular, $\mc V'_K$ covers the~whole space $\Amg{X}{c\mc A}$ except for $\bigcup c\mc A'$.

It follows that $\mc V' := \mc V'_K \cup \bigcup_{c\mc A'} \mc V'_{cA}$ is a~finite subcover of $\Amg{X}{c\mc A}$.

\eqref{case: smaller compactifications of Amg}: For $A\in \mc A$, denote by $q_A$ the~mapping witnessing that $cA\preceq dA$ and define $\varphi: \Amg{T}{d\mc A} \rightarrow \Amg{T}{c\mc A}$ as
\begin{equation*}
\varphi (x) :=
\begin{cases}
x & \text{  } x \in K , \\
q_A(x) & \text{ for } x\in dA, \ A\in \mc A.
\end{cases}
\end{equation*}
Clearly, $\varphi$ satisfies $\varphi|_X = \text{id}_X$.
By \eqref{case: continuous functions on amalgamation space}, the~mapping $\varphi$ is continuous.
This proves that $\varphi$ witnesses  $\Amg{T}{c\mc A} \preceq \Amg{T}{d\mc A}$.

\eqref{case: larger compactification of Amg type}: Let $X$ and $cX$ be as in the~statement. We denote by $i_A : \bar A \rightarrow cX $ the~identity mapping between $\bar{A} \subset \Amg{X}{\overline{\mc A}^{cX}}$ and $\bar{A} \subset cX$. We also denote as $i_X : X \rightarrow cX$ the~identity between $X \subset \Amg{X}{\overline{\mc A}^{cX}}$ and $X \subset cX$. By definition of topology on the~amalgamation space, $i_X$ and each $i_A$ is an embedding. We define a~mapping $\varphi : \Amg{X}{\overline{\mc A}^{cX}} \rightarrow cX$ as
\begin{equation*}
\varphi (x) :=
\begin{cases}
i_X(x) & \text{  } x \in X , \\
i_A(x) & \text{ for } x\in \bar A, \ A\in \mc A.
\end{cases}
\end{equation*}
The mapping $\varphi$ is well-defined and continuous (by \eqref{case: continuous functions on amalgamation space}). In particular, $\varphi$ witnesses that $cX \preceq \Amg{X}{\overline{\mc A}^{cX}}$.

\eqref{case: beta compactification as Amg}: This is an immediate consequence of \eqref{case: larger compactification of Amg type}.
\end{proof}



In the~following two subsections, we show two different ways of constructing compactifications of a~broom space $T$, and compute the~complexity of $T$ in each of them.

\subsection{Compactifications \texorpdfstring{$c_\gamma T$}{cT} and Broom Spaces \texorpdfstring{$\mathbf{T}_\alpha$}{T alpha}}
	\label{section: cT}

In one type of broom space compactifications, which we call $c_\gamma T$, the~closures of each $A\in \mc A$ are either the~smallest possible compactification $A$ (the Alexandroff one-point compactification), or the~largest possible one (the Čech-stone compactification):

\begin{notation}[Compactifications $c_\gamma T$]\label{notation:cT}
Let $\gamma \leq \omega_1$. For $A\in \mc A_{\omega_1}$, we denote
\[ c_\gamma A :=
\begin{cases}
	\alpha A, 	& \textrm{for } A\in \mc A_{<\gamma}\\
	\beta A, 	& \textrm{for } A\in \mc A_{\omega_1} \setminus \mc A_{<\gamma} ,
\end{cases} \]
where $A$ is endowed with the~discrete topology.
For an AD broom space $T = T_{\mc A}$, we set $c_\gamma T := \Amg{T}{c_\gamma \mc A}$.
\end{notation}

By Proposition~\ref{proposition: compactifications of amalgamations}, $c_\gamma T$ is a~compactification of $T$.
Since $\mc A_{<0} = \emptyset$, $c_0(\cdot)$ assigns to each $A$ its Čech-Stone compactification. It follows that $c_0 T = \beta T$:
\[ c_0 T \overset{\text{def.}}{\underset{\text{of }c_0 T}=} 
\Amg{T}{c_0 \mc A} \overset{\mc A_{<0}}{\underset{=\emptyset}=}
\Amg{T}{\beta \mc A} \overset{P\ref{proposition: compactifications of amalgamations}}{\underset{\eqref{case: beta compactification as Amg}}=}
\beta T .\]
On the~other hand, when $\gamma$ is such that the~whole $\mc A$ is contained in $\mc A_{<\gamma}$, every $c_\gamma A$ will be equal to $\alpha A$. This gives the~second identity in the~following observation:\footnote{$\Amg{T_{\mc A}}{\alpha \mc A}$ actually corresponds to the~compactification from in \cite{talagrand1985choquet}, used to prove the~existence of a~non-absolute $\fsd$ space.}
\[							\mc A \subset \mc A_{<\alpha} \ \& \ \gamma \geq \alpha \implies
c_\gamma T					\overset{\text{def.}}{\underset{\text{of }c_\gamma T}=}
\Amg{T}{c_\gamma \mc A}		\overset{\text{def.}}{\underset{\text{of }c_\gamma}=}
\Amg{T}{\alpha \mc A} .\]
By Proposition~\ref{proposition: compactifications of amalgamations}\,\eqref{case: smaller compactifications of Amg}, we obtain the~following chain of compactifications
\[ \beta T = c_0 T \succeq c_1 T \succeq \dots \succeq c_\gamma T
\succeq \dots \succeq c_\alpha T = c_{\alpha+1} T = \dots = c_{\omega_1} T , \]
which stabilizes at the~first ordinal $\alpha$ for which $\mc A \subset \mc A_{<\alpha}$.

The next lemma is the~only part which is specific to Talagrand's broom spaces:

\begin{lemma}[Lower bound on the~complexity in $c_\gamma T$] \label{lemma: c T - lower bound}
Let \errata{$\gamma \in [2,\omega_1]$}. If an AD broom space $T = T_\mc A$ satisfies $\mc A \supset \mc A^T_{<\gamma}$, then we have $T \notin \mc F_{<\gamma} (c_\gamma T)$.
\end{lemma}

\begin{proof}
Assuming for contradiction that $T \in \mc F_{<\gamma}(c_\gamma T)$, we can apply Lemma~\ref{lemma: Talagrand lemma 3} to $X=T$ and $Y=c_\gamma T$. It follows that there is an $\mc A^T_{<\gamma}$-witness $A$ for $T$ in $c_\gamma T$.
Since $A$ belongs to $\mc A^T_{<\gamma} \subset \mc A$, it is closed and discrete in $T$. By definition of $c_\gamma T$, we have $\overline{A}^{c_\gamma T} = c_\gamma A = \alpha A$ -- this contradicts Lemma~\ref{lemma: closed discrete witnesses}.
\end{proof}

In particular, this yields the~following result which we promised in Section~\ref{section:overview}:

\begin{corollary}[Spaces of \err{arbitrary} absolute complexity]\label{corollary: T alpha for odd}
For \errata{any} $\beta\in [2,\omega_1)$, there exists a~space $X^\beta_2$ satisfying
\[ \left\{ 2, \beta \right\} \subset \textnormal{Compl}\left( X^\beta_2 \right) \subset [2, \beta] .\]
\end{corollary}

\begin{proof}
Let $\beta\in [2,\omega_1)$ and set $X^\beta_2 := \mathbf{T}_\beta$.
Since $\mathbf{T}_\beta = T_{\mc A_{<\beta}}$ holds by definition, we can apply Proposition~\ref{proposition: basic broom complexities}\,$(ii)$ to get $\textnormal{Compl}(\mathbf{T}_\beta) \subset [0,\beta]$.
By Proposition~\ref{proposition: basic broom complexities}\,$(i)$, we have $2 \in \textnormal{Compl}(\mathbf{T}_\beta)$ and hence (by Proposition~\ref{proposition: basic attainable complexities}\,\eqref{case: compact and sigma compact})
\[ \left\{ 2 \right\} \subset \textnormal{Compl}\left( \mathbf{T}_\beta \right) \subset [2, \beta] .\]
By Lemma~\ref{lemma: c T - lower bound}, we have $\Compl{\mathbf{T}_\beta}{c_\beta \mathbf{T}_\beta} \geq \beta$, which implies that $\textnormal{Compl}(\mathbf{T}_\beta)$ contains $\beta$. This shows that $\textnormal{Compl}\left( \mathbf{T}_\beta \right) \subset [2, \beta]$ and concludes the~proof.
\end{proof}

The next lemma proves an upper estimate on the~complexity of $T$ in $c_\gamma T$.

\begin{lemma}[Upper bound on the~complexity in $c_\gamma T$] \label{lemma: c T - upper bound}
For any AD broom space $T$ and $\gamma\in [2,\omega_1]$, we have $T \in \mc F_\gamma (c_\gamma T)$.
\end{lemma}

\begin{proof}
Let $T=T_{\mc A}$ be an AD broom space and $\gamma\in [2,\omega_1]$.
We shall prove that $T$, $c_{\gamma} T$ and $\gamma$ satisfy the~assumptions of Lemma~\ref{lemma: S(y) and complexity of T} (which gives $T \in \mc F_\gamma (c_\gamma T)$).
To apply Lemma~\ref{lemma: S(y) and complexity of T}, we need to show that for every $x\in c_\gamma T \setminus T$, the~following set $S(x)$ is either finite or it can be covered by finitely sets $\D(A)$, $A\in \mc A_{<\gamma}$:
\begin{align*}
S(x) & = \left\{ s\in\seq | \ \overline{\mc N(s)}^{c_\gamma T} \ni x \right\} 
\end{align*}

First, we observe that $\overline{\mc N(s)}^{c_\gamma T} \cap \overline{A}^{c_\gamma T} = \overline{\mc N(s) \cap A}^{c_\gamma T}$ holds for any $s\in\seq$ and $A\in \mc A$.
The inclusion ``$\supset$'' is trivial.
For the~converse inclusion, let $x\in \overline{\mc N(s)}^{c_\gamma T} \cap \overline{A}^{c_\gamma T}$ and let $U$ be an open neighborhood of $x$ in $c_\gamma T$.
Since $\overline{A}^{c_\gamma T} = c_\gamma A$ is open in $c_\gamma T$ (Lemma~\ref{lemma: basic properties of amalgamation spaces}\,\eqref{case: E(A) is clopen in amalgamation}), $U\cap \overline{A}^{c_\gamma T}$ is an open neighborhood of $x$.
In particular, $x\in \overline{\mc N(s)}^{c_\gamma T}$ gives $U \cap \overline{A}^{c_\gamma T} \cap \mc N(s) \neq \emptyset$.
Since $\overline{A}^{c_\gamma T} \cap \mc N(s) = A \cap \mc N(s)$, it follows that $U$ intersects $A \cap \mc N(s)$.
This proves that $x$ belongs to $\overline{A \cap \mc N(s)}^{c_\gamma T}$.

Recall that $c_\gamma T \setminus T = \bigcup_{\mc A} \overline{A}^{c_\gamma T} \setminus A$. 
As the~first case, we shall assume that $x \in \overline{A}^{c_\gamma T} \setminus A$ for some $A \in \mc A \setminus \mc A_{<\gamma}$, and prove that the~set $S(x)$ is finite.
Suppose that such an $x$ satisfies $x\in \overline{\mc N(s)}^{c_\gamma T}$ and $x\in \overline{\mc N(t)}^{c_\gamma T}$ for two sequences $s$ and $t$.
By definition of $c_\gamma T$, we have $\overline{A}^{c_\gamma T} = c_\gamma A = \beta A$. It follows that
\begin{align}\label{equation: x and N of s}
x \in \ & \overline{\mc N(s)}^{c_\gamma T} \cap \overline{\mc N(t)}^{c_\gamma T} \cap \overline{A}^{c_\gamma T} = \overline{\mc N(s)\cap A}^{c_\gamma T} \cap \overline{\mc N(t) \cap A}^{c_\gamma T} = \nonumber \\
& = \overline{\mc N(s)\cap A}^{\beta A} \cap \overline{\mc N(t) \cap A}^{\beta A} .
\end{align}
Recall that for any normal topological space $X$ and closed $E,F\subset X$, we have $\overline{E}^{\beta X} \cap \overline{F}^{\beta X} = \overline{E \cap F}^{\beta X}$. In particular, this holds for $A$ (which is discrete).
Applying this to \eqref{equation: x and N of s} yields
\[ x \in \overline{\mc N(s) \cap \mc N(t) \cap A}^{\beta A} \subset
 \overline{\mc N(s) \cap \mc N(t) \cap A}^{c_\gamma T} \subset
 \overline{\mc N(s) \cap \mc N(t)}^{c_\gamma T} .\]
Since $x$ does not belong to $T \supset \mc N(s) \cap \mc N(t)$, the~set $\mc N(s) \cap \mc N(t)$ must in particular be non-empty.
The only way this might happen is when the~sequences $s$ and $t$ are comparable.
Consequently, $S(x)$ consists of a~single branch.
Since $x \notin \mc A(\overline{\mc N}^{c_\gamma T}\!)$ by Lemma~\ref{lemma: A of N}, this branch must necessarily be finite (by Lemma~\ref{lemma: IF trees}).

The remaining case is when $x$ belongs to $\overline{A}^{c_\gamma T} \setminus A$ for some $A \in \mc A \cap \mc A_{<\gamma}$, that is, when we have $x=x_A$.
We claim that such $x$ satisfies $S(x_A)=\D(A)$. Indeed, any $s\in\seq$ satisfies
\begin{align} \label{equation: S of x_A}
s \in \ & S(x_A) \iff x_A\in \overline{\mc N(s)}^{c_\gamma T} \iff x_A\in \overline{\mc N(s) \cap A}^{c_\gamma T} \iff \nonumber \\
 & \iff \mc N(s) \cap A \textrm{ is infinite} .
\end{align}
Clearly, \eqref{equation: S of x_A} is further equivalent to $\mc A$ containing infinitely many distinct extensions of $s$.
This happens precisely when $\cltr{A}$ contains infinitely many incomparable extensions of $s$.
Since each branch of $\cltr{A}$ is infinite, we can assume that each two of these extensions have different lengths. In other words, $s$ belongs to $S(x_A)$ precisely when it belongs to $\D(A)$.
\end{proof}

Finally, we apply all the~results to the~particular case of $T=\mathbf{T}_\alpha$:

\begin{proof}[Proof of the~``$\mathbf{T}_\alpha$'' part of Theorem~\ref{theorem: complexity of talagrands brooms}]
Let $\alpha \in [2,\omega_1]$. Since $\mathbf{T}_\alpha$ corresponds to the~family $\mc A = \mc A^T_{<\alpha} \subset \mc A_{<\alpha}$, we can apply Proposition~\ref{proposition: basic broom complexities} to get
\[ \textnormal{Compl}(\mathbf{T}_\alpha) \subset [2,\alpha] .\]
For any $\gamma \in [2,\alpha]$, $\mathbf{T}_\alpha$ is an $\mc F_\gamma$ subset of $c_\gamma \mathbf{T}_\alpha$ (by Lemma~\ref{lemma: c T - upper bound}), but it is \emph{not} its $\mc F_{<\gamma}$ subset (by Lemma~\ref{lemma: c T - lower bound}). It follows that $\Compl{\mathbf{T}_\alpha}{c_\gamma \mathbf{T}_\alpha} = \gamma$ and therefore $\gamma \in \textnormal{Compl}(\mathbf{T}_\alpha)$.
Since $\gamma \in [2,\alpha]$ was arbitrary, we get the~desired result:
\[ \textnormal{Compl}(\mathbf{T}_\alpha) = [2,\alpha] .\]
\end{proof}

\subsection{Compactifications \texorpdfstring{$d_\gamma T$}{dT} and Broom Spaces \texorpdfstring{$\mathbf{S}_\alpha$}{S alpha}} \label{section: dT}

The construction from Section~\ref{section: cT} is more involved than gluing together pre-existing examples, but nonetheless, it has somewhat similar flavor to the~approach from Section~\ref{section: topological sums}. We will now show that we can get the~same result by ``relying on the~same sets the~whole time''.

For $A\in \mc A_{\omega_1}$ and $h\in \seq$, we shall write
\[ A(h):=\left\{ \sigma \in A | \ \sigma \sqsupset h \right\} .\]

\noindent When $h=\emptyset$, we have $A(\emptyset) = A$ and $A(h)$ is ``as complicated as it can be'', in the~sense that its rank $\rank$ is the~same as the~rank of $A$. Conversely, when $h$ belongs to the~$B$ of which $A$ is a~broom-extension, then $A(h)$ is ``as simple as it can be'' -- it belongs to $\mc A_1$. The following lemma shows that the~intermediate possibilities are also possible, \err{and} defines the~corresponding set of ``$<\!\gamma$-handles'' of $A$.

\begin{lemma}[$H_{<\gamma}(\cdot)$, the~set $<\!\gamma$-handles]\label{lemma: H gamma}
Let $A$ be a~broom extension of $B$. For $\gamma\in [2,\omega_1]$, denote \emph{the set of $<\!\gamma$-handles of $A$} as
\[ H_{<\gamma}(A) := \{ h\in \cltr{B} | \ A(h)\in \mc A_{<\gamma} \ \& \text{ no other } s\sqsubset h \text{ satisfies } A(s) \in \mc A_{<\gamma} \}  .\]
\begin{enumerate}[(i)]
\item $A$ is the~disjoint union of sets $A(h)$, $h\in H_{<\gamma}(A)$.
\item If $A(h_0) \in \mc A_{<\gamma}$ holds for some $h_0\in \seq$, then we have $h\sqsubset h_0$ for some $h\in H_{<\gamma}(A)$.
\end{enumerate}
\end{lemma}

\begin{proof}
Let $A$, $B$ and $\gamma$ be as in the~statement.
First, we give a~more practical description of $H_{<\gamma}(A)$.
For each $s\in B$, $A(s)$ is a~broom extension of $\{s\} \in \mc B_1 \subset \mc B_{<\gamma}$. It follows that $A(s) \in \mc A_{<\gamma}$.
Moreover, we have
\[ \left( \forall u \sqsubset v \in \cltr{B} \right)
 \left(\forall \alpha<\omega_1\right): 
 A(u) \in \mc A_\alpha \implies A(v) \in \mc A_\alpha .\]
Indeed, this follows from the~definition of $\mc B_\alpha$.
Consequently, for each $s\in B$, there exists a~minimal $n\in\omega$ for which $A(s|n_s) \in \mc A_{<\gamma}$.
We claim that
\[ H_{<\gamma}(A) = \{ s|n_s \ | \ s\in B \} .\]
Indeed, each $s|n_s$ belongs to $H_{<\gamma}(A)$ by minimality of $n_s$. Conversely, any $h \in \cltr{B}$ satisfies $h \sqsubset s$ for some $s\in B$, so $h$ is either equal to $s|n_s$, or it does not belong to $H_{<\gamma}(A)$.

(i):
For each $\sigma\in A$, there is some $s\in B$ s.t. $s \sqsubset \sigma$.
It follows that $\sigma \in A(s|n_s)$, which proves that $A$ is covered by the~sets $A(h)$, $h\in H_{<\gamma}(A)$.

For distinct $s, t \in B$, the~initial segments $s|n_s $ and $t|n_t$ are either equal, or incomparable (by minimality of $n_s$ and $n_t$).
It follows that distinct $g,h\in H_{<\gamma}(A)$ are incomparable, which means that distinct $A(g)$ and $A(h)$ are disjoint.

(ii) follows from the~minimality of $n_s$ (and the~trivial fact that for each $h\in \cltr{B}$, there is some $s\in B$ with $s \sqsupset h$).
\end{proof}

By default, we equip each $A\in \mc A_{\omega_1}$ with a~discrete topology.
Using the~sets of $<\!\gamma$-handles, we define new compactifications of broom spaces. Unlike the~compactifications from $c_\gamma T$ from Section~\ref{section: cT}, these are no longer ``all or nothing'', but there are many intermediate steps between $\alpha A$ on one end and $\beta A$ on the~other:

\begin{notation}[Compactifications $d_\gamma T$]
For $\gamma \in [2,\omega_1]$ and $\mc A \in \mc A_{\omega_1}$, we define
\[ E_\gamma(A) := \bigoplus_{h\in H_{<\gamma}(A)} \alpha A(h)
 = \bigoplus_{h\in H_{<\gamma}(A)} \left( A(h) \cup \{ x_{A(h)} \} \right). \]
For an AD broom space $T$, we define a~compactification $d_\gamma T$ as
\[ d_\gamma T := \Amg{T}{d_\gamma \mc A} \text{, where } d_\gamma A := \beta E_\gamma(A) .\]
\end{notation}

By Lemma~\ref{lemma: H gamma}(ii), we have $E_\gamma(A) \supset A$, as well as the~fact that for distinct $g,h\in H_{<\gamma}(A)$, we have $x_{A(h)}\neq x_{A(g)}$.

We now show that the~compactifications $d_\gamma T$ have similar properties as $c_\gamma T$.

\begin{lemma}[Lower bound on the~complexity in $d_\gamma T$]\label{lemma: d T - lower bound}
Let $\gamma \in [2,\omega_1)$. If an AD system $\mc A\subset \mc A_{\omega_1}$ contains $\mc A^T_{<\alpha} \setminus \mc A^T_{<\alpha-1}$ for some non-limit $\alpha\geq \gamma$, then we have $T\notin \mc F_{<\gamma}(d_\gamma T)$.
\end{lemma}

\begin{proof}
Assuming for contradiction that $T \in \mc F_\eta (d_\gamma T)$ holds for some $\eta < \gamma$, we can apply Lemma~\ref{lemma: technical version of Talagrand lemma 3} (with `$\delta$'$=\alpha-1$).
It follows that there is some $A \in \mc A^T_{\alpha-1} \setminus \mc A^T_{<\alpha-1} \subset \mc A$ and $h_0\in \seq$, such that $A(h_0) \in \mc A^T_\eta \subset \mc A^T_{<\gamma}$ is a~witness for $T$ in $d_\gamma T$.

By Lemma~\ref{lemma: H gamma} (iii), there is some $h\in H_{<\gamma}(A)$ s.t. $h_0 \sqsupset h$.
Since $A$ belongs to $\mc A$, the~sets $A(h_0) \subset A(h) \subset A$ are all closed discrete in $T$. Moreover, we have $\overline{A(h)}^{d_\gamma T} = \alpha A(h)$ by definition of $d_\gamma T$. Because $A(h_0)$ is infinite, we get
\[ \overline{A(h_0)}^{d_\gamma T} = \overline{A(h_0)}^{\alpha A(h)} = A(h_0) \cup \{ x_{A(h)} \} .\]

We have found a~closed (in $T$) discrete witness for $T$ in $d_\gamma T$, whose closure in $d_\gamma T$ is homeomorphic to its one-point compactification -- a~contradiction with Lemma~\ref{lemma: closed discrete witnesses}.
\end{proof}

\begin{lemma}[Upper bound on the~complexity in $d_\gamma T$] \label{lemma: d T - upper bound}
For any AD broom space $T$ and $\gamma\in [2,\omega_1]$, $T \in \mc F_\gamma (d_\gamma T)$.
\end{lemma}

\begin{proof}
Let $x\in d_\gamma T \setminus T$.
When $x$ belongs to $\beta E_\gamma(A) \setminus  E_\gamma(A)$ for some $A\in \mc A$, we use the~exact same method as in Lemma~\ref{lemma: c T - upper bound} to prove that $S(x)$ is finite.

When $x$ belongs to $E_\gamma(A) \setminus A$ for some $A\in \mc A$, we have $x=x_{A(h)}$ for some $h\in H_{<\gamma } (A)$. By definition of $H_{<\gamma}(A)$, we have $A(h) \in \mc A_{<\gamma}$.
The approach from Lemma~\ref{lemma: c T - upper bound} yields $ S(x_{A(h)}) = \D(A\left(h\right))$.

We have verified the~assumptions of Lemma~\ref{lemma: S(y) and complexity of T}, which gives $T \in \mc F_\gamma (d_\gamma T)$.
\end{proof}

We have all the~ingredients necessary to finish the~proof of Theorem~\ref{theorem: complexity of talagrands brooms}:

\begin{proof}[Proof of the~``$\mathbf{S}_\alpha$'' part of Theorem~\ref{theorem: complexity of talagrands brooms}]
Let $\alpha \in [2,\omega_1)$ be a~non-limit ordinal. Since $\mathbf{S}_\alpha$ corresponds to the~family
\[ \mc A = \mc A^T_{<\alpha} \setminus \mc A^T_{<\alpha-1} \subset \mc A_{<\alpha} ,\]
we can apply Proposition~\ref{proposition: basic broom complexities} to get
\[ \textnormal{Compl}(\mathbf{S}_\alpha) \subset [2,\alpha] .\]
For any $\gamma \in [2,\alpha]$, $\mathbf{S}_\alpha$ satisfies $\Compl{\mathbf{S}_\alpha}{d_\gamma \mathbf{S}_\alpha}\leq\gamma $ (by Lemma~\ref{lemma: d T - upper bound}) and $\Compl{\mathbf{S}_\alpha}{d_\gamma \mathbf{S}_\alpha}\geq\gamma $ (by Lemma~\ref{lemma: d T - lower bound}). It follows that $\Compl{\mathbf{S}_\alpha}{d_\gamma \mathbf{S}_\alpha} = \gamma$ and therefore $\gamma \in \textnormal{Compl}(\mathbf{S}_\alpha)$.
Since $\gamma \in [2,\alpha]$ was arbitrary, we get
\[ \textnormal{Compl}(\mathbf{S}_\alpha) \supset [2,\alpha] ,\]
which concludes the~proof.
\end{proof}
%

\section*{Acknowledgment}
I would like to thank my supervisor, Ondřej Kalenda, for numerous very helpful
suggestions and fruitful consultations regarding this paper.
I am grateful to Adam Bartoš for discussions related to this paper.
This work was supported by the research grants  GAČR 17-00941S and  GA UK No. 915.


\putbib[refs]
\end{bibunit}

\chapwithtoc{List of papers included in the thesis}

\begin{enumerate}[{[1]}]
\item Kalenda, Ondřej, and Kovařík, Vojtěch. "Absolute $\fsd$ spaces." \emph{Topology and its Applications} 233 (2018): 44-51.
\item Kovařík, Vojtěch. "Absolute $\mathcal F $-Borel classes." \emph{Fundamenta Mathematicae}, \emph{published electronically}.
\item Kovařík, Vojtěch. "Complexities and Representations of $\mc F$-Borel Spaces" \emph{Submitted, available at arXiv:1804.08367}.
\end{enumerate}





\openright
\end{document}